\documentclass[leqno, openany, draft]{book}

\usepackage{makeidx}
\makeindex

\def\diam{\mathop{\rm diam}}

\def\dim{\mathop{\rm dim}}
\def\im{\mathop{\rm Im}}
\def\rad{\mathop{\rm rad}}
\def\re{\mathop{\rm Re}}
\def\supp{\mathop{\rm supp}}
\def\tr{\mathop{\rm tr}}

\newtheorem{theorem}{Theorem}
\newtheorem{lemma}[theorem]{Lemma}
\newtheorem{proposition}[theorem]{Proposition}
\newtheorem{definition}[theorem]{Definition}
\newtheorem{corollary}[theorem]{Corollary}
\newtheorem{problem}[theorem]{Problem}
\newtheorem{remark}[theorem]{Remark}
\newtheorem{examples}[theorem]{Examples}
\newtheorem{question}[theorem]{Question}

\newcommand{\begintheorem}{\addtocounter{equation}{1}\begin{theorem}}
\newcommand{\beginlemma}{\addtocounter{equation}{1}\begin{lemma}}
\newcommand{\beginproposition}{\addtocounter{equation}{1}\begin{proposition}}
\newcommand{\begindefinition}{\addtocounter{equation}{1}\begin{definition}}
\newcommand{\begincorollary}{\addtocounter{equation}{1}\begin{corollary}}
\newcommand{\beginproblem}{\addtocounter{equation}{1}\begin{problem}}
\newcommand{\beginremark}{\addtocounter{equation}{1}\begin{remark}}
\newcommand{\beginexamples}{\addtocounter{equation}{1}\begin{examples}}
\newcommand{\beginquestion}{\addtocounter{equation}{1}\begin{question}}

\renewcommand{\theequation}{\arabic{chapter}.\arabic{equation}}

\begin{document}

\frontmatter

\title{A Sampler in Analysis}

\author{Stephen Semmes \\
        Rice University}

\date{}

\maketitle

\chapter*{Preface}

\noindent
The mathematical area of analysis is often described as the study of
limits, continuity, and convergence, as in calculus.  It is also very
much concerned with estimates, whether or not a limit is involved.
Here we look at several topics involving norms on vector spaces and
linear mappings, with prerequisites along the lines of advanced
calculus and basic linear algebra.  The main idea is to explore
intermediate ranges of abstraction and sophistication, without getting
bogged down with too many technicalities.  Lebesgue integrals are not
required, but could easily be incorporated by readers familiar with
that theory.

        We begin with some inequalities related to convexity in the
first chapter, which can be applied to sums or integrals.  The next
three chapters focus on finite-dimensional vector spaces and linear
transformations between them.  Some properties of infinite sums are
described in Chapter \ref{sums, l^p spaces}, as well as a class of
infinite-dimensional spaces known as $\ell^p$ spaces.  The latter give
examples of Banach and Hilbert spaces, which are considered more
abstractly in Chapter \ref{banach, hilbert spaces}.  An important tool
for dealing with bounded linear operators on $\ell^p$ or $L^p$ spaces
is Marcel Riesz' convexity theorem, presented in Chapter
\ref{interpolation}.  As a further introduction to real-variable
methods in harmonic analysis, estimates for dyadic maximal and square
functions are discussed in Chapter \ref{dyadic analysis}.  A brief
review of some basic notions about metric spaces is included in
Appendix \ref{metric spaces}.

        Of course, there are numerous excellent texts on these and
related subjects, a selection of which can be found in the
bibliography.  Indeed, it is hoped that readers might pursue specific
topics more fully, according to their interests.  Here one might find
a few tricks of the trade, or simplified special cases, which
illustrate broader concepts.  I would like to dedicate this book to my
fellow students from Washington University in St Louis, and to the
faculty there, from whom we learned a great deal.

\tableofcontents

\mainmatter

\chapter{Preliminaries}
\label{preliminaries}

\section{Real and complex numbers}
\label{real, complex numbers}

	As usual, the real line is denoted ${\bf R}$, the complex
plane is denoted ${\bf C}$, and the set of integers is denoted ${\bf
Z}$.  If $A$ is a subset of ${\bf R}$ and $b$ is a real number such
that $a \le b$ for all $a \in A$, then $b$ is said to be an
\emph{upper bound}\index{upper bounds} for $A$.  A real number $c$ is
said to be the \emph{least upper bound}\index{least upper bound} or
\emph{supremum}\index{supremum} of $A$ if $c$ is an upper bound for
$A$ and $c \le b$ for every real number $b$ which is an upper bound
for $A$.  One version of the \emph{completeness}\index{completeness of
the real numbers} of the real numbers asserts that a nonempty subset
$A$ of ${\bf R}$ with an upper bound has a least upper bound.  It is
easy to see from the definition that the supremum $\sup A$ of $A$ is
unique when it exists.

	Similarly, if $A \subseteq {\bf R}$ and $y \in {\bf R}$
satisfy $y \le x$ for every $x \in A$, then $y$ is said to be a
\emph{lower bound}\index{lower bounds} of $A$.  If $z$ is a real
number such that $z$ is a lower bound for $A$ and $y \le z$ for every
real number $y$ which is a lower bound for $A$, then $z$ is said to be
a \emph{greatest lower bound}\index{greatest lower bound} or
\emph{infimum}\index{infimum} of $A$.  It follows from the
completeness of the real numbers that every nonempty subset $A$ of
${\bf R}$ with a lower bound has a greatest lower bound.  This can be
obtained as the supremum of the set of lower bounds for $A$, or as the
negative of the supremum of
\begin{equation}
        -A = \{-a : a \in A\}.
\end{equation}
 Again, it is easy to see directly from the definition that the
infimum $\inf A$ of $A$ is unique when it exists.

	It is sometimes convenient to use extended real numbers, which
are real numbers together with $+\infty$, $-\infty$, with standard
conventions concerning arithmetic operations and ordering.  More precisely,
\begin{equation}
        -\infty < x < +\infty
\end{equation}
and
\begin{equation}
        x + (+\infty) = (+\infty) + x = +\infty, \quad
         x + (-\infty) = (-\infty) + x = -\infty
\end{equation}
for every $x \in {\bf R}$.  If $x$ is a positive real number, then
\begin{equation}
        x \cdot (+\infty) = (+\infty) \cdot x = +\infty, \quad
         x \cdot (-\infty) = (-\infty) \cdot x = -\infty
\end{equation}
while the product of $\pm \infty$ with a negative real number changes
the sign.  The product of $\pm \infty$ with $\pm \infty$ is defined to
be $\pm \infty$, where the signs are multiplied in the usual way.  One
can also define $x / \pm \infty$ to be $0$ for every $x \in {\bf R}$,
but expressions such as $\infty - \infty$, $\infty / \infty$, and
$0/0$ are not defined.  If one allows extended real numbers, then
every nonempty set $A \subseteq {\bf R}$ has a supremum and an
infimum, where $\sup A = +\infty$ if $A$ does not have a finite upper
bound, and $\inf A = -\infty$ if $A$ does not have a finite lower
bound.  In situations where all of the quantities of interest are
nonnegative, it may be appropriate to interpret $1/0$ as being equal
to $+\infty$.
	
	If $a$ and $b$ are real numbers with $a < b$, then there are
four types of \emph{intervals}\index{intervals} in the real line with
endpoints $a$ and $b$, i.e., the open interval $(a, b)$, the
half-open, half-closed intervals $(a, b]$, $[a, b)$, and the closed
interval $[a, b]$.  These four types of intervals are defined as follows:
\begin{eqnarray*}
	(a,b) & = & \{x \in {\bf R} : a < x < b\};		\\  
	(a,b] & = & \{x \in {\bf R} : a < x \le b\};		\\  \ 
	[a,b) & = & \{x \in {\bf R} : a \le x < b\};		\\  \ 
	[a,b] & = & \{x \in {\bf R} : a \le x \le b\}.
\end{eqnarray*}
The \emph{length}\index{lengths of intervals} of each of these
intervals is defined to be $b-a$, and the length of an interval $I$
may be denoted $|I|$.

	We also consider $[a, b]$ to be defined when $a = b$, in which
event the interval consists of a single point and has length equal to
$0$.  For an interval which is open at the left endpoint $a$, we may
allow $a = -\infty$, and for an interval which is open at the right
endpoint $b$, we may allow $b = +\infty$.  Hence the real line may be
expressed as $(-\infty, +\infty)$.  In these cases, we say that the
interval is unbounded, while an interval with finite endpoints is said
to be bounded.

	If $x$ is a real number, then the absolute value of $x$ is
denoted $|x|$ and defined to be $x$ when $x \ge 0$ and to be $-x$ when
$x \le 0$.  Thus $|x|$ is always a nonnegative real number, $|x| = 0$
if and only if $x = 0$, and
\begin{equation}
	|x + y| \le |x| + |y|
\end{equation}
and
\begin{equation}
	|x \cdot y| = |x| \cdot |y|
\end{equation}
for every $x, y \in {\bf R}$.  These properties are not difficult to
verify.

	Suppose that $z = x + i y$ is a complex number, where $x, y
\in {\bf R}$.  One may refer to $x$, $y$ as the \emph{real}\index{real
part} and \emph{imaginary parts}\index{imaginary part} of $z$, denoted
$\re z$, $\im z$.  The \emph{complex conjugate}\index{complex
conjugation} of $z$ is denoted $\overline{z}$ and defined by
\begin{equation}
	\overline{z} = x - i y.
\end{equation}
It is easy to see that
\begin{equation}
\label{overline{z + w} = overline{z} + overline{w}}
	\overline{z + w} = \overline{z} + \overline{w}
\end{equation}
and
\begin{equation}
\label{overline{z cdot w} = overline{z} cdot overline{w}}
	\overline{z \cdot w} = \overline{z} \cdot \overline{w}
\end{equation}
for every $z, w \in {\bf C}$.  Note that
\begin{equation}
\label{z + overline{z} = 2 re z and z - overline{z} = 2 i im z}
 z + \overline{z} = 2 \re z \quad\hbox{and}\quad z - \overline{z} = 2 i \im z
\end{equation}
for every $z \in {\bf C}$, and that the complex conjugate of
$\overline{z}$ is equal to $z$.

        The \emph{modulus}\index{modulus of a complex number} of $z =
x + i y \in {\bf C}$, $x, y \in {\bf R}$, is denoted $|z|$ and defined
to be the nonnegative real number given by
\begin{equation}
\label{|z| = sqrt{x^2 + y^2}}
	|z| = \sqrt{x^2 + y^2}.
\end{equation}
Thus the modulus of $z$ is the same as the absolute value of $z$ when
$z \in {\bf R}$, and 
\begin{equation}
\label{|re z|, |im z| le |z|}
        |\re z|, |\im z| \le |z|
\end{equation}
for every $z \in {\bf C}$.  Of course, the modulus of $z$ is the same
as the modulus of the complex conjugate of $z$, and it is easy to see that
\begin{equation}
\label{|z|^2 = z cdot overline{z}}
	|z|^2 = z \cdot \overline{z}
\end{equation}
for every $z \in {\bf C}$.  This implies that
\begin{equation}
\label{|z cdot w| = |z| cdot |w|}
	|z \cdot w| = |z| \cdot |w|
\end{equation}
for every $z, w \in {\bf C}$, because of (\ref{overline{z cdot w} =
overline{z} cdot overline{w}}).  

        Similarly, we would like to check that
\begin{equation}
\label{|z + w| le |z| + |w|}
	|z + w| \le |z| + |w|
\end{equation}
for every $z, w \in {\bf C}$.  Using (\ref{|z|^2 = z cdot
overline{z}}) applied to $z + w$ and then (\ref{overline{z + w} =
overline{z} + overline{w}}), we get that
\begin{equation}
 |z + w|^2 = (z + w) \, (\overline{z} + \overline{w})
           = |z|^2 + z \, \overline{w} + w \, \overline{z} + |w|^2.
\end{equation}
We also have that
\begin{equation}
        z \, \overline{w} + w \, \overline{z}
         = z \, \overline{w} + \overline{(z \, \overline{w})} 
         = 2 \re z \, \overline w \le 2 \, |z \, \overline{w}|
         = 2 \, |z| \, |w|,
\end{equation}
and hence that
\begin{equation}
        |z + w|^2 \le |z|^2 + 2 \, |z| \, |w| + |w|^2 = (|z| + |w|)^2.
\end{equation}
This implies (\ref{|z + w| le |z| + |w|}), as desired.

	A sequence $\{z_n\}_{n=1}^\infty$ of complex numbers is said to
\emph{converge}\index{convergent sequences} to another complex number $z$
if for every $\epsilon > 0$ there is a positive integer $N$ such that
\begin{equation}
	|z_n - z| < \epsilon
\end{equation}
for every $n \ge N$.  One can check that the limit $z$ of the sequence
$\{z_n\}_{n=1}^\infty$ is unique when it exists, in which case we put
\begin{equation}
	\lim_{n \to \infty} z_n = z.
\end{equation}
If $\{w_n\}_{n=1}^\infty$, $\{z_n\}_{n=1}^\infty$ are two sequences of
complex numbers which converge to the complex numbers $w$, $z$,
respectively, then the sequences $\{w_n + z_n\}_{n=1}^\infty$, $\{w_n
\cdot z_n\}_{n=1}^\infty$ of sums and products converge to the sum $w
+ z$ and product $w \cdot z$ of the limits, respectively.  A sequence
$\{z_n\}_{n=1}^\infty$ of complex numbers converges to a complex
number $z$ if and only if the sequences of real and imaginary parts of
the $z_n$'s converge to the real and imaginary parts of $z$.

	Let $\{x_n\}_{n=1}^\infty$ be a sequence of real numbers which
is \emph{monotone increasing}, which is to say that $x_n \le x_{n+1}$
for each $n$.  One can check that $\{x_n\}_{n=1}^\infty$ converges if
and only if the set of $x_n$'s has an upper bound, in which case the
limit of the sequence is equal to the supremum of this set.  For any
sequence $\{x_n\}_{n=1}^\infty$ of real numbers,
\begin{equation}
	x_j \to +\infty \hbox{ as } j \to \infty
\end{equation}
if for each $L \ge 0$ there is a positive integer $N$ such that
\begin{equation}
	x_n \ge L
\end{equation}
for every $n \ge N$.  If $\{x_n\}_{n=1}^\infty$ is an unbounded
monotone increasing sequence of real numbers, then $x_n \to +\infty$
as $n \to \infty$.  Similar remarks apply to monotone decreasing
sequences of real numbers.

	Let $\{a_n\}_{n=1}^\infty$ be a sequence of real numbers.  For
each positive integer $k$, put
\begin{equation}
	A_k = \sup \{a_n : n \ge k\},
\end{equation}
which may be $+\infty$.  Thus $A_{k+1} \le A_k$ for every $k$.  The
\emph{upper limit} of $\{a_n\}_{n=1}^\infty$ is denoted $\limsup_{n
\to \infty} a_n$ and defined to be the infimum of the $A_k$'s, which
may be $\pm \infty$.  Similarly, if
\begin{equation}
	B_l = \inf \{a_n : n \ge l\},
\end{equation}
then $B_l \le B_{l+1}$ for every $l$, and the \emph{lower limit}
$\liminf_{n \to \infty} a_n$ of $\{a_n\}_{n=1}^\infty$ is defined to
be the supremum of the $B_l$'s.  By construction, $B_l \le A_k$ for
every $k$ and $l$, and hence
\begin{equation}
	\liminf_{n \to \infty} a_n \le \limsup_{n \to \infty} a_n.
\end{equation}
One can check that $a_n \to a$ as $n \to \infty$ if and only if
\begin{equation}
	\liminf_{n \to \infty} a_n = \limsup_{n \to \infty} a_n = a.
\end{equation}

	A sequence $\{z_n\}_{n=1}^\infty$ of complex numbers is said
to be a \emph{Cauchy sequence}\index{Cauchy sequences} if for every
$\epsilon > 0$ there is a positive integer $N$ such that
\begin{equation}
	|z_l - z_n| < \epsilon
\end{equation}
for each $l, n \ge N$.  It is easy to see that $\{z_n\}_{n =
1}^\infty$ is a Cauchy sequence if and only if the corresponding
sequences of real and imaginary parts of the $z_n$'s are Cauchy
sequences, and that convergent sequences are automatically Cauchy
sequences.  It is not difficult to show that the upper and lower
limits of a Cauchy sequence of real numbers are finite and equal, and
hence that every Cauchy sequence of real numbers converges.  It
follows that every Cauchy sequence of complex numbers converges too.

	An infinite series of complex numbers $\sum_{j=0}^\infty a_j$
is said to \emph{converge}\index{convergent series} if the
corresponding sequence of partial sums $\sum_{j=0}^n a_j$ converges,
in which case the sum of the series is defined to be the limit of the
sequence of partial sums.  If $\sum_{j=0}^\infty a_j$ converges, then
\begin{equation}
	\lim_{j \to \infty} a_j = 0.
\end{equation}
The partial sums of an infinite series whose terms are nonnegative
real numbers are monotone increasing, and therefore the series
converges if and only if the partial sums are bounded.  An infinite
series $\sum_{j=0}^\infty a_j$ of complex numbers is said to converge
\emph{absolutely}\index{absolutely convergent series} if
\begin{equation}
	\sum_{j=0}^\infty |a_j|
\end{equation}
converges.  One can check that the partial sums of an absolutely
convergent series form a Cauchy sequence, and therefore converge.

	If $A$ is a subset of a set $X$, then ${\bf 1}_A(x)$ denotes
the \emph{indicator function}\index{indicator functions} of $A$ on
$X$.  This is the function equal to $1$ when $x \in A$ and to $0$ when
$x \in X \backslash A$, and it is sometimes called the
\emph{characteristic function}\index{characteristic functions}
associated to $A$.  A function on the real line, or on an interval in
the real line, is called a \emph{step function}\index{step functions}
if it is a finite linear combination of indicator functions of
intervals.  Equivalently, this means that there is a finite partition
of the domain into intervals on which the function is constant.  In
this book, one is normally welcome to restrict one's attention to
functions on the real line that are step functions, at least in the
context of integrating functions on ${\bf R}$.  Step functions are
convenient because their integrals can be reduced immediately to
finite sums.  Results about other functions can often be derived from
those for step functions by approximation.

\section{Convex functions}
\label{convex functions}

	Let $I$ be an open interval in the real line, which may be
unbounded.  A real-valued function $\phi(x)$ on $I$ is said to be
\emph{convex}\index{convex functions} if
\begin{equation}
\label{convexity inequality for phi}
	\phi(\lambda \, x + (1-\lambda) \, y)
		\le \lambda \, \phi(x) + (1-\lambda) \, \phi(y)
\end{equation}
for every $x, y \in I$ and $\lambda \in [0,1]$.

	If $\phi(x)$ is an affine function, which is to say that
$\phi(x) = a \, x + b$ for some real numbers $a$ and $b$, then
$\phi(x)$ is a convex function on the whole real line, with equality
in (\ref{convexity inequality for phi}) for all $x$, $y$, and
$\lambda$.  Equivalently, both $\phi(x)$ and $-\phi(x)$ are convex,
which characterizes affine functions.  It is easy to see that $\phi(x)
= |x|$ is a convex function on the whole real line too.  If $\phi(x)$
is an arbitrary convex function on $I$, and if $c$ is a real number,
then the translation $\phi(x-c)$ of $\phi(x)$ is a convex function on
\begin{equation}
	I + c = \{x + c : x \in I\}.
\end{equation}
In particular, for each real number $c$, $|x-c|$ defines a convex function
on ${\bf R}$.  

\beginlemma
\label{convexity and increasing difference quotients}
A real-valued function $\phi(x)$ on $I$ is convex if and only if
\begin{equation}
\label{increasing difference quotients inequalities}
	\frac{\phi(t) - \phi(s)}{t-s} \le \frac{\phi(u) - \phi(s)}{u-s}
					\le \frac{\phi(u) - \phi(t)}{u-t}
\end{equation}
for every $s, t, u \in I$ with $s < t < u$.
\end{lemma}

	If $s$, $t$, $u$ are as in the lemma,
then
\begin{equation}
\label{t = frac{t-s}{u-s} u + frac{u-t}{u-s} s}
	t = \frac{t-s}{u-s} \, u + \frac{u-t}{u-s} \, s,
\end{equation}
where
\begin{equation}
	0 < \frac{t-s}{u-s} < 1
\end{equation}
and
\begin{equation}
	\frac{u-t}{u-s} = 1 - \frac{t-s}{u-s}.
\end{equation}
If $\phi(x)$ is convex, then
\begin{equation}
\label{convexity inequality rewritten}
	\phi(t) \le \frac{t-s}{u-s} \, \phi(u) + \frac{u-t}{u-s} \, \phi(s).
\end{equation}
One can rewrite this in two different ways to get (\ref{increasing
difference quotients inequalities}).  Conversely, one can work
backwards, and rewrite either of the inequalities in (\ref{increasing
difference quotients inequalities}) to get (\ref{convexity inequality
rewritten}), which gives (\ref{convexity inequality for phi}) when
$s$, $t$, and $u$ correspond to $x$, $y$, and $\lambda$ as in (\ref{t
= frac{t-s}{u-s} u + frac{u-t}{u-s} s}).

\beginlemma
\label{convexity and comparison with affine functions}
A function $\phi(x)$ on $I$ is convex if and only if for each $t \in
I$ there is a real-valued affine function $A(x)$ on ${\bf R}$ such
that $A(t) = \phi(t)$ and $A(x) \le \phi(x)$ for every $x \in I$.
\end{lemma}

	To see that this condition is sufficient for $\phi$ to be
convex, let $x$, $y$, and $\lambda$ be given in the usual way.  If $A$
is an affine function associated to
\begin{equation}
	t = \lambda \, x + (1 - \lambda) \, y
\end{equation}
as in the statement of the lemma, then
\begin{eqnarray}
\label{using A to get convexity inequality}
	\phi(\lambda \, x + (1 - \lambda) \, y)
		& = & A(\lambda \, x + (1 - \lambda) \, y)	\\
		& = &  \lambda \, A(x) + (1 - \lambda) \, A(y)	\nonumber \\
		& \le &  \lambda \, \phi(x) + (1 - \lambda) \, \phi(y).
								\nonumber
\end{eqnarray}

	Conversely, suppose that $\phi(x)$ is convex, and let $t \in I$
be given.  We would like to choose a real number $a$ so that
\begin{equation}
\label{A(x) = phi(t) + a (x-t)}
	A(x) = \phi(t) + a \, (x-t)
\end{equation}
satisfies $A(x) \le \phi(x)$ for all $x \in I$, which is equivalent to
\begin{equation}
\label{a (x-t) le phi(x) - phi(t)}
	a \, (x-t) \le \phi(x) - \phi(t)
\end{equation}
for $x \in I$.  This is trivial when $x = t$, and otherwise we can
rewrite (\ref{a (x-t) le phi(x) - phi(t)}) as
\begin{equation}
\label{a le frac{phi(x) - phi(t)}{x-t}}
	a \le \frac{\phi(x) - \phi(t)}{x-t}
\end{equation}
when $x > t$, and as
\begin{equation}
\label{frac{phi(t) - phi(x)}{t-x} le a}
	\frac{\phi(t) - \phi(x)}{t-x} \le a
\end{equation}
when $x < t$.  It follows from Lemma \ref{convexity and increasing
difference quotients} that
\begin{equation}
	\frac{\phi(t) - \phi(s)}{t - s} 
			\le \frac{\phi(u) - \phi(t)}{u - t}
\end{equation}
for every $s, u \in I$ such that $s < t < u$.  Hence
\begin{equation}
	D_l = \sup \bigg\{\frac{\phi(t) - \phi(s)}{t-s} 
						: s \in I, s < t \bigg\}
\end{equation}
and
\begin{equation}
	D_r = \inf \bigg\{\frac{\phi(u) - \phi(t)}{u-t}
						: u \in I, u > t \bigg\}
\end{equation}
are well-defined and satisfy
\begin{equation}
	D_l \le D_r.
\end{equation}
To get (\ref{a le frac{phi(x) - phi(t)}{x-t}}) and (\ref{frac{phi(t) -
phi(x)}{t-x} le a}), it suffices to choose $a \in {\bf R}$ such that
\begin{equation}
	D_l \le a \le D_r.
\end{equation}
This completes the proof of Lemma \ref{convexity and comparison with
affine functions}.

	A real-valued function $\phi(x)$ on $I$ is said to be
\emph{strictly convex}\index{strict convexity} if
\begin{equation}
\label{strict convexity inequality}
	\phi(\lambda \, x + (1-\lambda) \, y) <
		\lambda \, \phi(x) + (1-\lambda) \, \phi(y)
\end{equation}
for every $x, y \in I$ such that $x \ne y$ and each $\lambda \in
(0,1)$.

\beginlemma
\label{strict convexity and comparison with affine functions}
A real-valued function $\phi$ on $I$ is strictly convex if and only if
for every point $t \in I$ there is a real-valued affine function
$A(x)$ on ${\bf R}$ such that $A(t) = \phi(t)$ and $A(x) < \phi(x)$
for all $x \in I \backslash \{t\}$.
\end{lemma}

	This is the analogue of Lemma \ref{convexity and comparison
with affine functions} for strictly convex functions, which can be
obtained in practically the same manner as before.  For the existence
of $A$ when $\phi$ is strictly convex, one can start with $A$ as in
the previous lemma, and use strict convexity to show that $A(x) \ne
\phi(x)$ when $x \ne t$.

	The convexity of a real-valued function $\phi$ on $I$ can also
be characterized by the property that for each $x, y \in I$ with $x < y$,
\begin{equation}
	\phi(t) \le B(t) \hbox{ for every } t \in [x, y],
\end{equation}
where $B$ is the affine function on the real line which is equal to
$\phi$ at $x$ and $y$.  Strict convexity corresponds to
\begin{equation}
	\phi(t) < B(t) \hbox{ when } t \in (x, y).
\end{equation}
This is easy to check, just using the definitions.

	  Note that convex functions are automatically continuous.
This follows by trapping a convex function on both sides of a point
between affine functions with the same value at that point.

\beginlemma
\label{lemma about convexity conditions}
If $\phi$ is a continuous real-valued function on $I$, and if for each
$x$, $y$ in $I$ there is a $\lambda_{x,y} \in (0,1)$ such that
(\ref{convexity inequality for phi}) holds with $\lambda = \lambda_{x,
y}$, then $\phi$ is convex.
\end{lemma}

	This is often stated in the special case where $\lambda_{x, y}
= 1/2$ for every $x, y \in I$, in which event one can iterate the
condition and pass to a limit to get the desired inequality for
arbitrary $\lambda$.  Alternatively, for each $x, y \in I$ with $x <
y$, let $L(x, y)$ be the set of $\lambda \in [0, 1]$ such that
(\ref{convexity inequality for phi}) holds, which is a closed set when
$\phi$ is continuous, and which automatically contains $0$ and $1$.
If $L(x, y) \ne [0, 1]$, then one can get a contradiction under the
conditions of the lemma, by considering a maximal open interval in
$[0, 1] \backslash L(x, y)$, and showing that it has to contain an
element of $L(x, y)$.

\section{Some related inequalities}
\label{related inequalities}

	Suppose that $\phi$ is a convex function on an open interval
$I \subseteq {\bf R}$, as in the previous section.  If $K$ is an
interval in ${\bf R}$ of positive length and $f$ is an integrable
function on $K$ such that $f(x) \in I$ for all $x \in K$, then
\begin{equation}
        |K|^{-1} \int_K f(x) \, dx \in I,
\end{equation}
and
\begin{equation}
\label{phi(|K|^{-1} int_K f(x) dx) le |K|^{-1} int_K phi(f(x)) dx}
	\phi\Big(|K|^{-1} \int_K f(x) \, dx \Big)
		\le |K|^{-1} \int_J \phi(f(x)) \, dx.
\end{equation}
This is called \emph{Jensen's inequality}.\index{Jensen's inequality}

	Let us first consider the analogous statement for finite sums.
If $x_1, x_2, \ldots, x_n$ are elements of $I$ and $\lambda_1,
\lambda_2, \ldots, \lambda_n$ are nonnegative real numbers such that
$\sum_{i=1}^n \lambda_i = 1$, then
\begin{equation}
        \sum_{i=1}^n \lambda_i \, x_i \in I,
\end{equation}
and
\begin{equation}
\label{generalized convexity inequality}
	\phi\Big(\sum_{i=1}^n \lambda_i \, x_i \Big)
		\le \sum_{i=1}^n \lambda_i \, \phi(x_i).
\end{equation}
This is the same as (\ref{convexity inequality for phi}) when $n = 2$,
and one can apply (\ref{convexity inequality for phi}) repeatedly to
get the general case.  One can also use the characterization of
convexity in Lemma \ref{convexity and comparison with affine
functions}, as in (\ref{using A to get convexity inequality}).  If $f$
is a step function, then (\ref{phi(|K|^{-1} int_K f(x) dx) le |K|^{-1}
int_K phi(f(x)) dx}) follows directly from (\ref{generalized convexity
inequality}).  In general, one can reduce to the case of finite sums
through suitable approximations, or employ Lemma \ref{convexity and
comparison with affine functions} in the same way as for sums.

	It is well known that $\phi(t) = |t|^p$ is a convex function
on the real line when $p$ is a real number such that $p \ge 1$, and
moreover that $|t|^p$ is strictly convex when $p > 1$.  In particular,
\begin{equation}
\label{| |K|^{-1} int_K f(x) dx|^p le |K|^{-1} int_K f(x)^p dx}
 \biggl||K|^{-1} \int_K f(x) \, dx\biggr|^p \le |K|^{-1} \int_K |f(x)|^p \, dx
\end{equation}
for real-valued functions $f$ on an interval $K$ of positive length.

	Let $p$, $q$ be real numbers such that $p, q \ge 1$ and
\begin{equation}
\label{1/p + 1/q = 1}
	\frac{1}{p} + \frac{1}{q} = 1.
\end{equation}
In this event we say that $p$ and $q$ are \emph{conjugate
exponents}.\index{conjugate exponents} If $f, g$ are nonnegative
real-valued functions on an interval $K$, then \emph{H\"older's
inequality}\index{H\"older's inequality} states that
\begin{equation}
\label{Holder's inequality written out}
	\int_K f(x) \, g(x) \, dx
		\le \Big(\int_K f(y)^p \, dy\Big)^{1/p}
			\Big(\int_K g(z)^q \, dz\Big)^{1/q}.
\end{equation}
We can also allow $p$ or $q$ to be $1$ and the other to be $+\infty$,
which is consistent with (\ref{1/p + 1/q = 1}).  If $p = 1$ and $q =
+\infty$, then the substitute for (\ref{Holder's inequality written
out}) is
\begin{equation}
\label{Holder's inequality, p=1, q=infty}
	\int_K f(x) \, g(x) \, dx
		\le \Big(\int_K f(y) \, dy \Big)
			\Big(\sup_{z \in K} \, g(z) \Big).
\end{equation}

	Let us now prove (\ref{Holder's inequality written out}) when
$p, q > 1$, beginning with some initial reductions.  The inequality is
trivial if $f$ or $g$ is identically $0$, or zero ``almost
everywhere'', since the left side of (\ref{Holder's inequality written
out}) is then equal to $0$.  Thus we may suppose that
\begin{equation}
\label{(int_K f(y)^p dy)^{1/p}, (int_K g(z)^q dz)^{1/q}}
	\Big(\int_K f(y)^p \, dy\Big)^{1/p} \quad\hbox{and}\quad
		\Big(\int_K g(z)^q \, dz\Big)^{1/q}
\end{equation}
are nonzero.  We may suppose further that these expressions are both
equal to $1$, because the general case would follow by multiplying $f$
and $g$ by positive constants.  For any nonnegative real numbers $s$,
$t$,
\begin{equation}
\label{s t le frac{s^p}{p} + frac{t^q}{q}}
	s \, t \le \frac{s^p}{p} + \frac{t^q}{q}.
\end{equation}
This is a version of the geometric-arithmetic mean inequalities, which
can be treated as an exercise in calculus, or derived from the
convexity of the exponential function.  Note that the inequality is
strict when $s^p \ne t^q$.  Applying (\ref{s t le frac{s^p}{p} +
frac{t^q}{q}}) to $s = f(x)$ and $t = g(x)$ and then integrating in
$x$, we get that
\begin{equation}
	\int_K f(x) \, g(x) \, dx
		\le \frac{1}{p} \int_K f(x)^p \, dx
			+ \frac{1}{q} \int_K g(x)^q \, dx.
\end{equation}
This implies (\ref{Holder's inequality written out}) when the
integrals of $f^p$ and $g^q$ are equal to $1$, as desired.

	Similarly,
\begin{equation}
\label{Holder's inequality for sums}
\index{H\"older's inequality}
	\sum_{j=1}^n a_j \, b_j \le \Big(\sum_{k=1}^n a_k^p\Big)^{1/p}
				\Big(\sum_{l=1}^n b_l^q\Big)^{1/q}
\end{equation}
when $a_1, \ldots, a_n$, $b_1, \ldots, b_n$ are nonnegative real
numbers and $p, q \ge 1$ are conjugate exponents.  If $p = 1$ and $q =
\infty$, then this should be interpreted as
\begin{equation}
	\sum_{j=1}^n a_j b_j \le \Big(\sum_{k=1}^n a_k\Big) 
				\Big(\max \{b_l : 1 \le l \le n \}\Big).
\end{equation}

	Let $f$ and $g$ be nonnegative functions on an interval $K$
again, and let $p$ be a real number, $p \ge 1$.  \emph{Minkowski's
inequality}\index{Minkowski's inequality} states that
\begin{equation}
\label{Minkowski's inequality, written out}
	\quad  \Big(\int_K (f(x) + g(x))^p \, dx \Big)^{1/p}
	         \le \Big(\int_K f(x)^p \, dx \Big)^{1/p}
		    + \Big(\int_K g(x)^p \, dx \Big)^{1/p}.
\end{equation}
The analogue of (\ref{Minkowski's inequality, written out}) for $p =
+\infty$ is the elementary inequality
\begin{equation}
	\sup_{x \in K} \, (f(x) + g(x))
		\le \sup_{x \in K} \, f(x) + \sup_{x \in K} \, g(x).
\end{equation}

	Let us suppose that $1 < p < + \infty$, since
(\ref{Minkowski's inequality, written out}) is trivial when $p = 1$.
We begin with
\begin{eqnarray}
\lefteqn{\int_K (f(x) + g(x))^p \, dx}	\\
	&& = \int_K f(x) \, (f(x) + g(x))^{p-1} \, dx
		+ \int_K g(x) \, (f(x) + g(x))^{p-1} \, dx.	\nonumber
\end{eqnarray}
If $q > 1$ is the conjugate exponent of $p$, then H\"older's
inequality implies that
\begin{eqnarray}
\lefteqn{\int_K f(x) \, (f(x) + g(x))^{p-1} \, dx}	\\
	&& \le \Big(\int_K f(y)^p \, dy \Big)^{1/p}
	        \Big(\int_K (f(z) + g(z))^{q(p-1)} \, dz \Big)^{1/q}
						\nonumber  \\
	&& = \Big(\int_K f(y)^p \, dy \Big)^{1/p}
	        \Big(\int_K (f(z) + g(z))^p \, dz \Big)^{1 - 1/p}.
						\nonumber
\end{eqnarray}
There is an analogous estimate for $\int_K g(x) \, (f(x) + g(x))^{p-1}
\, dx$, which leads to
\begin{eqnarray}
\lefteqn{\qquad \int_K (f(x) + g(x))^p \, dx}		\\
	&& \le
	\bigg\{\Big(\int_K f(y)^p \, dy \Big)^{1/p}
			+ \Big(\int_J g(y)^p \, dy \Big)^{1/p}\bigg\}
	        \Big(\int_K (f(z) + g(z))^p \, dz \Big)^{1 - 1/p}.
					\nonumber
\end{eqnarray}
It is easy to derive (\ref{Minkowski's inequality, written out}) from
this.

	Minkowski's inequality\index{Minkowski's inequality} for
finite sums can be expressed as
\begin{equation}
\label{Minkowski's inequality for sums}
	\Big(\sum_{j=1}^n (a_j + b_j)^p \Big)^{1/p}
		\le \Big(\sum_{j=1}^n a_j^p \Big)^{1/p}
			+ \Big(\sum_{j=1}^n b_j^p \Big)^{1/p}
\end{equation}
when $1 \le p < \infty$, and
\begin{equation}
\label{max of a sum}
	\qquad \max \{ a_j + b_j : 1 \le j \le n \}
	 \le \max \{ a_j : 1 \le j \le n\} + \max \{ b_j : 1 \le j \le n\}
\end{equation}
when $p = \infty$, where $a_1, \ldots, a_n$, $b_1, \ldots, b_n$ are
nonnegative real numbers.  These inequalities can be shown in the same
way as for integrals.  As an alternate approach, fix $p$, $1 < p <
\infty$, since the $p = 1$ and $p = \infty$ cases are easy, and
suppose for the moment that
\begin{equation}
\label{(sum_j a_j^p)^{1/p} = (sum_j b_j^p)^{1/p} = 1}
	\Big(\sum_{j=1}^n a_j^p \Big)^{1/p}
			= \Big(\sum_{j=1}^n b_j^p \Big)^{1/p} = 1.
\end{equation}
If $t$ is a real number such that $0 \le t \le 1$, then
\begin{equation}
\label{(sum_j (t a_j + (1-t) b_j)^p)^{1/p} le 1}
	\Big(\sum_{j=1}^n (t \, a_j + (1-t) \, b_j)^p \Big)^{1/p} \le 1.
\end{equation}
To see this, rewrite (\ref{(sum_j a_j^p)^{1/p} = (sum_j b_j^p)^{1/p} = 1})
and (\ref{(sum_j (t a_j + (1-t) b_j)^p)^{1/p} le 1}) as
\begin{equation}
\label{sum_j a_j^p = sum_j b_j^p = 1}
	\sum_{j=1}^n a_j^p = \sum_j b_j^p = 1
\end{equation}
and
\begin{equation}
\label{sum_j (t a_j + (1-t) b_j)^p le 1}
	\sum_{j=1}^n (t \, a_j + (1-t) \, b_j)^p \le 1,
\end{equation}
respectively.  To go from (\ref{sum_j a_j^p = sum_j b_j^p = 1}) to
(\ref{sum_j (t a_j + (1-t) b_j)^p le 1}), it suffices to know that
\begin{equation}
	(t \, a_j + (1-t) \, b_j)^p \le t \, a_j^p + (1-t) \, b_j^p
\end{equation}
for each $j$, which follows from the convexity of the function
$\phi(x) = x^p$, $x \ge 0$.  Once one has (\ref{(sum_j (t a_j + (1-t)
b_j)^p)^{1/p} le 1}) under the assumption (\ref{(sum_j a_j^p)^{1/p} =
(sum_j b_j^p)^{1/p} = 1}), it is not difficult to derive
(\ref{Minkowski's inequality for sums}) in the general case.
Basically, the parameter $t$ compensates for $(\sum_{j=1}^n
a_j^p)^{1/p}$ and $(\sum_{j=1}^n b_j^p)^{1/p}$ not being equal.
	
	Fix a positive integer $n$, and suppose that $\{a_j\}_{j=1}^n$
is a finite sequence of nonnegative real numbers.  Let $p$ and $q$ be
positive real numbers, with $p < q$.  Clearly
\begin{equation}
\label{max_{1 le j le n} a_j le (sum_{j=1}^n a_j^p)^{1/p}}
	\max_{1 \le j \le n} \, a_j \le \Big(\sum_{j=1}^n a_j^p \Big)^{1/p}.
\end{equation}
Moreover,
\begin{equation}
\label{(sum_{j=1}^n a_j^q)^{1/q} le (sum_{j=1}^n a_j^p)^{1/p}}
	\Big(\sum_{j=1}^n a_j^q \Big)^{1/q}
				\le \Big(\sum_{j=1}^n a_j^p \Big)^{1/p},
\end{equation}
because
\begin{equation}
	\sum_{j=1}^n a_j^q 
		\le \Big(\max_{1 \le k \le n} a_k \Big)^{q-p}
			\, \Big(\sum_{l=1}^n a_l^p \Big)		
	  	\le  \Big(\sum_{r=1}^n a_r^p \Big)^{1 + (q-p)/p}
\end{equation}
and $1 + (q-p)/p = q/p$.

	In the other direction,
\begin{equation}
\label{(sum_{j=1}^n a_j^p)^{1/p} le n^{1/p} max_{1 le j le n} a_j}
	\Big(\sum_{j=1}^n a_j^p \Big)^{1/p} 
		\le n^{1/p} \, \max_{1 \le j \le n} a_j,
\end{equation}
and
\begin{equation}
\label{(sum_{j=1}^n a_j^p)^{1/p} le n^{(1/p) - (1/q)}(sum_{j=1}^n a_j^q)^{1/q}}
	\Big(\sum_{j=1}^n a_j^p \Big)^{1/p}
		\le n^{(1/p) - (1/q)} \, \Big(\sum_{j=1}^n a_j^q \Big)^{1/q}.
\end{equation}
The first inequality is trivial, and the second can be rewritten as
\begin{equation}
	\Big(\frac{1}{n} \sum_{j=1}^n a_j^p \Big)^{q/p}
		\le \frac{1}{n} \sum_{j=1}^n a_j^q,
\end{equation}
which is an instance of (\ref{generalized convexity inequality})
applied to $\phi(x) = x^{q/p}$.

	If $0 < p < 1$ and $u$, $v$ are nonnegative real numbers,
then
\begin{equation}
	(u + v)^p \le u^p + v^p.
\end{equation}
This is a special case of (\ref{(sum_{j=1}^n a_j^q)^{1/q} le
(sum_{j=1}^n a_j^p)^{1/p}}), with $q = 1$ and $n = 2$.  This 
leads to
\begin{equation}
\label{sum_{j=1}^n (b_j + c_j)^p le sum_{j=1}^n b_j^p + sum_j c_j^p}
	\sum_{j=1}^n (b_j + c_j)^p \le \sum_{j=1}^n b_j^p + \sum_j c_j^p,
\end{equation}
for nonnegative real numbers $b_1, \ldots, b_n$ and $c_1, \ldots, c_n$, and
\begin{equation}
\label{int_K (f(x) + g(x))^p dx le int_K f(x)^p dx + int_K g(x)^p dx }
	\int_K (f(x) + g(x))^p \, dx 
	  \le \int_K f(x)^p \, dx + \int_K g(x)^p \, dx 
\end{equation}
for nonnegative functions $f$, $g$ on an interval $K$.

	Suppose that $0 < p, q, r, < \infty$ and
\begin{equation}
	\frac{1}{r} = \frac{1}{p} + \frac{1}{q}.
\end{equation}
If $a_1, \ldots, a_n$, $b_1, \ldots, b_n$ are nonnegative real
numbers, then
\begin{equation}
	\Big(\sum_{j=1}^n (a_j \, b_j)^r \Big)^{1/r}
		\le \Big(\sum_{j=1}^n a_j^p \Big)^{1/p}
			\, \Big(\sum_{j=1}^n b_j^q \Big)^{1/q}.
\end{equation}
This follows from H\"older's inequality.  Similarly, for nonnegative
functions $f$, $g$ on an interval $K$,
\begin{equation}
	\Big(\int_K (f(x) \, g(x))^r \, dx \Big)^{1/r}
		\le \Big(\int_K f(x)^p \, dx \Big)^{1/p}
			\, \Big(\int_K g(x)^q \, dx \Big)^{1/q}.
\end{equation}
One can also allow for infinite exponents in the usual way.

\chapter{Norms on vector spaces}
\label{norms on vector spaces}

        In this book, all vector spaces use the real or complex
numbers as their underlying scalar field.  We may sometimes wish to
restrict ourselves to one or the other, but frequently both are fine.
Let us make the standing assumption that all vector spaces are
finite-dimensional in this and the next two chapters.

\section{Definitions and examples}
\label{definitions, examples}

	Let $V$ be a real or complex vector space.  By a
\emph{norm}\index{norms} on $V$ we mean a nonnegative real-valued
function $\| \cdot \|$ on $V$ such that $\|v\| = 0$ if and only if $v$
is the zero vector in $V$,
\begin{equation}
\label{||t v|| = |t| ||v||}
        \|t \, v\| = |t| \, \|v\|
\end{equation}
for every $v \in V$ and $t \in {\bf R}$ or ${\bf C}$, as appropriate,
and
\begin{equation}
\label{||v + w|| le ||v|| + ||w||}
        \|v + w\| \le \|v\| + \|w\|
\end{equation}
for every $v, w \in V$.  As a basic class of examples, let $V$ be
${\bf R}^n$ or ${\bf C}^n$, and consider
\begin{equation}
\label{def of ||v||_p}
	\|v\|_p = \Big(\sum_{j=1}^n |v_j|^p\Big)^{1/p}
\end{equation}
when $1 \le p < \infty$, and
\begin{equation}
\label{def of ||v||_{infty}}
	\|v\|_{\infty} = \max_{1 \le j \le n} \, |v_j|.
\end{equation}
The triangle inequality for these norms follows from (\ref{Minkowski's
inequality for sums}) and (\ref{max of a sum}).

        If
\begin{equation}
\label{B_1 = {v in V : ||v|| le 1}}
	B_1 = \{v \in V : \|v\| \le 1\}
\end{equation}
is the closed unit ball corresponding to a norm $\|v\|$ on $V$, then
it is easy to see that $B_1$ is a convex\index{convex sets} set in $V$.
This means that
\begin{equation}
	t \, v + (1-t) \, w \in B_1
\end{equation}
whenever $v, w \in B_1$ and $t$ is a real number such that $0 \le t
\le 1$, which follows from (\ref{||t v|| = |t| ||v||}) and (\ref{||v +
w|| le ||v|| + ||w||}).  Conversely, if $\|v\|$ is a nonnegative
real-valued function on $V$ such that $\|v\| = 0$ if and only if $v =
0$, $\|v\|$ satisfies (\ref{||t v|| = |t| ||v||}), and the unit ball
$B_1$ is convex, then one can show $\|v\|$ also satisfies (\ref{||v +
w|| le ||v|| + ||w||}), and hence that $\|v\|$ is a norm on $V$.  In
effect, this was mentioned already in Section \ref{related
inequalities}, as an alternate approach to Minkowski's inequality for
finite sums.	

	If $V$ is a vector space, and $\|\cdot\|$ is a norm on $V$,
then
\begin{equation}
\label{| ||v|| - ||w|| | le ||v - w||}
	\bigl|\|v\| - \|w\| \bigr| \le \|v - w\|
\end{equation}
for every $v, w \in V$.  This follows from
\begin{equation}
\label{||v|| le ||w|| + ||v - w||}
	\|v\| \le \|w\| + \|v-w\|,
\end{equation}
and the analogous inequality with the roles of $v$ and $w$
interchanged.  Suppose that $V = {\bf R}^n$ or ${\bf C}^n$ for some
positive integer $n$, which is not a real restriction, since every
real or complex vector space of positive finite dimension is
isomorphic to one of these.  Let $|x|$ denote the standard Euclidean
norm on ${\bf R}^n$ or ${\bf C}^n$, which is the same as the norm
$\|x\|_2$ in (\ref{def of ||v||_p}).  One can check that there is a
positive constant $C$ such that
\begin{equation}
\label{||v|| le C |v|}
	\|v\| \le C \, |v|
\end{equation}
for every $v \in V$, by expanding $v$ in the standard basis for $V =
{\bf R}^n$ or ${\bf C}^n$, and using the triangle inequality and
homogeneity of $\|v\|$.  This and (\ref{| ||v|| - ||w|| | le ||v -
w||}) imply that $\|v\|$ is a continuous real-valued function on $V$,
with respect to the standard Euclidean metric and topology.  Thus the
minimum $b > 0$ of $\|v\|$ among the vectors $v \in V$ with $|v| = 1$
is attained, by well-known results about continuity and compactness,
and $b > 0$.  It follows that
\begin{equation}
\label{b |v| le ||v||}
	b \, |v| \le \|v\|
\end{equation}
for every $v \in V$, because of the homogeneity property of the norms
$\|v\|$ and $|v|$.

\section{Dual spaces and norms}
\label{dual spaces, norms}

	Let $V$ be a vector space, real or complex.  By a \emph{linear
functional}\index{linear functionals} on $V$ we mean a linear mapping
from $V$ into the field of scalars, i.e., the real or complex numbers,
as appropriate.  The \emph{dual of $V$} is the space of linear
functionals on $V$, which is a vector space over the same field of
scalars as $V$, with respect to pointwise addition and scalar
multiplication.  The dual of $V$ is denoted $V^*$,\index{dual spaces}
and it is well known that $V^*$ is also finite-dimensional when $V$
is, with the same dimension as $V$.

	If $\|\cdot \|$ is a norm on $V$, then the corresponding
\emph{dual norm}\index{dual norms} $\|\cdot \|^*$ on $V^*$ is defined
as follows.  If $\lambda$ is a linear functional on $V$, then
\begin{equation}
\label{||lambda||^* = sup {|lambda(v)| : v in V, ||v|| le 1}}
	\|\lambda\|^* = \sup \{|\lambda(v)| : v \in V, \, \|v\| \le 1\}.
\end{equation}
Equivalently,
\begin{equation}
	|\lambda(v)| \le \|\lambda\|^* \, \|v\|
\end{equation}
for every $v \in V$, and $\|\lambda\|^*$ is the smallest nonnegative
real number with this property.  It is not difficult to verify that
$\| \cdot \|^*$ defines a norm on $V^*$.  In particular, the
finiteness of $\|\lambda\|^*$ can be derived from the remarks at the
end of the previous section.

	For example, suppose that $V = {\bf R}^n$ or ${\bf C}^n$ for
some positive integer $n$.  We can identify $V^*$ with ${\bf R}^n$ or
${\bf C}^n$, respectively, by associating to each $w$ in ${\bf R}^n$
or ${\bf C}^n$ the linear functional $\lambda_w$ on $V$ given by
\begin{equation}
	\lambda_w(v) = \sum_{j=1}^n w_j \, v_j.
\end{equation}
Let $1 \le p, q \le \infty$ be conjugate exponents, which is to say
that $1/p + 1/q = 1$, and let us check that $\|\lambda_w \|^* =
\|w\|_q$ is the dual norm for $\|v\| = \|v\|_p$.

	First, we have that
\begin{equation}
	|\lambda_w(v)| \le \|w\|_q \, \|v\|_p
\end{equation}
for all $w$ and $v$, by H\"older's inequality.  To show that
$\|\lambda_w\|^* = \|w\|_q$, we would like to check that for each $w$
there is a nonzero $v$ such that
\begin{equation}
\label{|lambda_w(v)| = ||w||_q ||v||_p}
	|\lambda_w(v)| = \|w\|_q \, \|v\|_p.
\end{equation}

	Let $w$ be given.  We may as well assume that $w \ne 0$, since
otherwise any $v$ would do.  Let us also assume for the moment that $p
> 1$, so that $q < \infty$.  Under these conditions, we can define $v$
by
\begin{equation}
	v_j = \overline{w_j} \, |w_j|^{q-2}
\end{equation}
when $w_j \ne 0$, and by $v_j = 0$ when $w_j = 0$.  Here
$\overline{w_j}$ is the complex conjugate of $w_j$, which is not
needed when we are working with real numbers instead of complex
numbers.  With this choice of $v$, we have that
\begin{equation}
	\lambda_w(v) = \sum_{j=1}^n |w_j|^q = \|w\|_q^q.
\end{equation}
It remains to check that
\begin{equation}
\label{||v||_p = ||w||_q^{q-1}}
	\|v\|_p = \|w\|_q^{q-1}.
\end{equation}

	If $q = 1$, then $p = \infty$, and (\ref{||v||_p =
||w||_q^{q-1}}) reduces to
\begin{equation}
\label{max_{1 le j le n} |v_j| = 1}
	\max_{1 \le j \le n} |v_j| = 1.
\end{equation}
In this case $|v_j| = 1$ for each $j$ such that $v_j \ne 0$, which
holds for at least one $j$ because $w \ne 0$.  Thus we get
(\ref{max_{1 le j le n} |v_j| = 1}).  If $q > 1$, then $|v_j| =
|w_j|^{q-1}$ for each $j$, and one can verify (\ref{||v||_p =
||w||_q^{q-1}}) using the identity $p (q-1) = q$.

	Finally, if $p = 1$, and hence $q = \infty$, then choose $l$,
$1 \le l \le n$, such that
\begin{equation}
	|w_l| = \max_{1 \le j \le n} |w_j| = \|w\|_{\infty}.
\end{equation}
Define $v$ by $v_l = \overline{w_l} / |w_l|^{-1}$ and $v_j = 0$ when
$j \ne l$.  This leads to
\begin{equation}
	\lambda_w(v) = |w_l| = \|w\|_{\infty}
\end{equation}
and $\|v\|_1 = 1$, as desired.

\section{Second duals}
\label{second duals}

	Let $V$ be a vector space, and $V^*$ its dual space.  The dual
of $V^*$ is denoted $V^{**}$.

	There is a canonical mapping from $V$ into $V^{**}$, defined
as follows.  Let $v \in V$ be given.  For each $\lambda \in V^*$, we
get a scalar by taking $\lambda(v)$.  The mapping $\lambda \mapsto
\lambda(v)$ is a linear functional on $V^*$, and hence an element of
$V^{**}$.  Since we can do this for every $v \in V$, we get a mapping
from $V$ into $V^{**}$.  One can check that this mapping is linear and
an isomorphism from $V$ onto $V^{**}$.  For instance, everything can
be expressed in terms of a basis for $V$.

	Now suppose that we have a norm $\|\cdot \|$ on $V$.  This
leads to a dual norm $\|\cdot \|^*$ on $V^*$, as in the previous
section, and a double dual norm $\|\cdot \|^{**}$ on $V^{**}$.  Using
the canonical isomorphism between $V$ and $V^{**}$ just described, we
can think of $\| \cdot \|^{**}$ as defining a norm on $V$.  
We would like to show that
\begin{equation}
\label{||v||^{**} = ||v|| for every v in V}
	\|v\|^{**} = \|v\|   \quad\hbox{for every } v \in V.
\end{equation}
Note that this holds for the $p$-norms $\| \cdot \|_p$ on ${\bf R}^n$
and ${\bf C}^n$, by the analysis of their duals in the preceding
section.

	Let $v \in V$ be given.  By definition of the dual norm
$\|\lambda\|^*$, we have that
\begin{equation}
	|\lambda(v)| \le \|\lambda\|^* \, \|v\|
\end{equation}
for every $\lambda \in V^*$, and hence that
\begin{equation}
        \|v\|^{**} \le \|v\|.
\end{equation}
It remains to show that the opposite inequality holds, which is
trivial when $v = 0$.  Thus it suffices to show that there is a
nonzero $\lambda_0 \in V^*$ such that
\begin{equation}
\label{lambda_0(v) = ||lambda_0||^* ||v||}
	\lambda_0(v) = \|\lambda_0\|^* \, \|v\|
\end{equation}
when $v \ne 0$.

\begintheorem
\label{extension theorem}
Let $V$ be a real or complex vector space, and let $\|\cdot \|$ be a
norm on $V$.  If $W$ is a linear subspace of $V$ and $\mu$ is a linear
functional on $W$ such that
\begin{equation}
	|\mu(w)| \le \|w\| \quad\hbox{for every } w \in W,
\end{equation}
then there is a linear functional $\widehat{\mu}$ on $V$ such that
$\widehat{\mu} = \mu$ on $W$ and $\|\widehat{\mu}\|^* \le 1$.
\end{theorem}

	The existence of a nonzero $\lambda_0 \in V^*$ satisfying
(\ref{lambda_0(v) = ||lambda_0||^* ||v||}) follows easily from this,
by first defining $\lambda_0$ on the span of $v$ so that $\lambda_0(v)
= \|v\|$, and then extending to a linear functional on $V$ with norm
$1$.

	To prove the theorem, let us begin by assuming that $V$ is a
real vector space.  Afterwards, we shall discuss the complex case.

	Let $W$ and $\mu$ be given as in the theorem, and let $\dim Z$
be the dimension of a linear subspace $Z$ of $V$.  For each integer
$j$ such that $\dim W \le j \le \dim V$, we would like to show that
there is a linear subspace $W_j$ of $V$ and a linear functional
$\mu_j$ on $W_j$ such that $W \subseteq W_j$, $\dim W_j = j$, $\mu_j =
\mu$ on $W$, and
\begin{equation}
\label{|mu_j(w)| le ||w|| for every w in W_j}
	|\mu_j(w)| \le \|w\| \quad\hbox{for every } w \in W_j.
\end{equation}
If we can do this with $j = \dim V$, then $W_j = V$, and this would
give a linear functional on $V$ with the required properties.

	Let us show that we can do this by induction.  For the base
case $j = \dim W$, we simply take $W_j = W$ and $\mu_j = \mu$.
Suppose that $\dim W \le j < \dim V$, and that $W_j$, $\mu_j$ are as
above.  We would like to choose $W_{j+1}$ and $\mu_{j+1}$ with the
analogous properties for $j+1$ instead of $j$.  To be more precise, we
shall choose them in such a way that $W_j \subseteq W_{j+1}$ and
$\mu_{j+1}$ is an extension of $\mu_j$ to $W_{j+1}$.

	Under these conditions, $W_j$ is a proper subspace of $V$, and
hence there is a $z \in V \backslash W_j$.  Fix any such $z$, and take
$W_{j+1}$ to be the span of $W_j$ and $z$.  Thus
\begin{equation}
	\dim W_{j + 1} = \dim W_j + 1 = j + 1.
\end{equation}

	Let $\alpha$ be a real number, to be chosen later in the
argument.  If we set $\mu_{j+1}(z)$ equal to $\alpha$, then
$\mu_{j+1}$ is determined on all of $W_{j+1}$ by linearity and the
condition that $\mu_{j+1}$ be an extension of $\mu_j$.  Specifically,
each $w \in W_{j + 1}$ can be expressed in a unique way as $x + t \,
z$ for some $x \in W_j$ and $t \in {\bf R}$, and
\begin{equation}
	\mu_{j + 1}(w) = \mu_j(x) + t \, \alpha.
\end{equation}

	It remains to choose $\alpha$ so that $\mu_{j+1}$ satisfies
the analogue of (\ref{|mu_j(w)| le ||w|| for every w in W_j}) for
$j+1$, which is to say that
\begin{equation}
	|\mu_{j+1}(w)| \le \|w\|  
		\quad\hbox{for every } w \in W_{j+1}.
\end{equation}
Equivalently, we would like to choose $\alpha$ so that
\begin{equation}
\label{|mu_j(x) + t alpha| le ||x + t z||, x in W_j, t in {bf R}}
	|\mu_j(x) + t \, \alpha| \le \|x + t \, z\|
		\quad\hbox{for every } x \in W_j
		\hbox{ and } t \in {\bf R}.
\end{equation}
It suffices to show that
\begin{equation}
\label{|mu_j(x) + alpha| le ||x + z|| for every x in W_j}
	|\mu_j(x) + \alpha| \le \|x + z\| 
			\quad\hbox{for every } x \in W_j,
\end{equation}
since the case where $t = 0$ in (\ref{|mu_j(x) + t alpha| le ||x + t
z||, x in W_j, t in {bf R}}) corresponds exactly to our induction
hypothesis (\ref{|mu_j(w)| le ||w|| for every w in W_j}), and one can
eliminate $t \ne 0$ using homogeneity.  Let us rewrite (\ref{|mu_j(x)
+ alpha| le ||x + z|| for every x in W_j}) as
\begin{equation}
\label{- mu_j(x) -||x + z|| le alpha le -mu_j(x) + ||x + z||, x in W_j}
  - \mu_j(x) -\|x + z\| \le \alpha \le -\mu_j(x) + \|x + z\|
		\quad\hbox{for every } x \in W_j.
\end{equation}

	It follows from (\ref{|mu_j(w)| le ||w|| for every w in W_j})
that
\begin{equation}
	\mu_j(x - y) \le \|x - y\|
		\quad\hbox{for every } x, y \in W_j.
\end{equation}
Using the triangle inequality, we get that
\begin{equation}
	\mu_j(x - y) \le \|x + z\| + \|y + z\|
		\quad\hbox{for every } x, y \in W_j,
\end{equation}
and hence
\begin{equation}
   -\mu_j(y) - \|y + z\| \le - \mu_j(x) + \|x + z\|
		\quad\hbox{for every } x, y \in W_j.
\end{equation}
If $A$ is the supremum of the left side of this inequality over $y
\in W_j$, and $B$ is the infimum of the right side of this inequality
over $x \in W_j$, then $A \le B$, and any $\alpha \in {\bf R}$ such
that $A \le \alpha \le B$ satisfies (\ref{- mu_j(x) -||x + z|| le
alpha le -mu_j(x) + ||x + z||, x in W_j}).  This finishes the
induction argument, and the proof of Theorem \ref{extension theorem}
when $V$ is a real vector space.

	Consider now the case of a complex vector space $V$.  The real
part of a linear functional on $V$ is also a linear functional on $V$
as a real vector space, i.e., forgetting about multiplication by $i$.
Conversely, if $\phi$ is a real-valued function on $V$ which is linear
with respect to vector addition and scalar multiplication by real
numbers, then there is a unique complex linear functional $\psi$ on
$V$ whose real part is $\phi$, given by
\begin{equation}
	\psi(v) = \phi(v) - i \, \phi(i \, v).
\end{equation}
For any complex number $\zeta$,
\begin{equation}
	|\zeta| = \sup \{\re (a \, \zeta) : a \in {\bf C}, |a| \le 1\}.
\end{equation}
If $V$ is equipped with a norm $\|v\|$, then the norm of a complex
linear functional $\lambda$ on $V$ can be expressed as
\begin{equation}
	\|\lambda\|^* = \sup \{\re (a \, \lambda(v)) :
		v \in V, a \in {\bf C}, \|v\| \le 1, |a| \le 1\}.
\end{equation}
By linearity, $\lambda(a \, v) = a \, \lambda(v)$, which implies that
\begin{equation}
	\|\lambda\|^* = \sup \{\re \lambda(v) : v \in V, \|v\| \le 1\}.
\end{equation}
Thus the norm of a complex linear functional on $V$ is the same as the
norm of its real part, as a linear functional on the real version of
$V$.  To prove the extension theorem in the complex case, one can
apply the extension theorem in the real case to the real part of the
given complex linear functional on a complex linear subspace, and then
complexify the real extension to get a complex linear extension with
the same estimate for the norm.

\section{Linear transformations}
\label{linear transformations}

	Let $V_1$ and $V_2$ be vector spaces, both real or both
complex, equipped with norms $\|\cdot \|_1$ and $\|\cdot \|_2$,
respectively.  Here the subscripts are merely labels to distinguish
these norms, rather than referring to the $p$-norms described in
Section \ref{definitions, examples}.  The corresponding \emph{operator
norm}\index{operator norms} $\|T\|_{op}$ of a linear transformation
$T$ from $V_1$ into $V_2$ is defined by
\begin{equation}
\label{defn of ||T||_{op}}
	\|T\|_{op} = \sup \{\|T(v)\|_2 : v \in V_1, \, \|v\|_1 \le 1\}.
\end{equation}
Equivalently,
\begin{equation}
\label{bound from ||T||_{op}}
	\|T(v)\|_2 \le \|T\|_{op} \, \|v\|_1
		\quad\hbox{for every } v \in V_1,
\end{equation}
and $\|T\|_{op}$ is the smallest nonnegative real number with this
property.  The finiteness of $\|T\|_{op}$ is easy to check using the
remarks at the end of Section \ref{definitions, examples}.

	The space $\mathcal{L}(V_1, V_2)$ of all linear
transformations from $V_1$ into $V_2$ is a vector space in a natural
way, using pointwise addition and scalar multiplication of linear
transformations, with the same scalar field as for $V_1$ and $V_2$.
One can verify that the operator norm $\| \cdot \|_{op}$ on
$\mathcal{L}(V_1, V_2)$ is a norm on this vector space.  Note that the
dual $V^*$ of a vector space $V$ is the same as $\mathcal{L}(V, {\bf
R})$ or $\mathcal{L}(V, {\bf C})$, as appropriate, and the dual norm
on $V^*$ associated to a norm on $V$ is the same as the operator norm
with respect to the standard norm on ${\bf R}$ or ${\bf C}$.

	Suppose that $V_3$ is another vector space, with the same
field of scalars as $V_1$ and $V_2$, and equipped with a norm $\|\cdot
\|_3$.  If $T_1 : V_1 \to V_2$ and $T_2 : V_2 \to V_3$ are linear
mappings, then the composition $T_2 \circ T_1$ is the linear mapping
from $V_1$ to $V_3$ given by
\begin{equation}
	(T_2 \circ T_1)(v) = T_2(T_1(v)).
\end{equation}
It is easy to see that the operator norm of $T_2 \circ T_1$ is less
than or equal to the product of the operator norms of $T_1$ and $T_2$
with respect to the given norms on $V_1$, $V_2$, and $V_3$.

	Let $V_1$ and $V_2$ be vector spaces, both real or both
complex, and let $T$ be a linear transformation from $V_1$ into $V_2$.
There is a canonical dual linear transformation $T^* : V_2^* \to
V_1^*$ corresponding to $T$, defined by
\begin{equation}
	T^*(\mu) = \mu \circ T  \quad\hbox{for every } \mu \in V_2^*.
\end{equation}
In other words, if $\mu$ is a linear functional on $V_2$, then $\mu
\circ T$ is a linear functional on $V_1$, and $T^*(\mu)$ is this
linear functional.  If $R, T : V_1 \to V_2$ are linear mappings and
$a$, $b$ are scalars, then
\begin{equation}
	(a \, R + b \, T)^* = a \, R^* + b \, T^*.
\end{equation}
If $V_3$ is another vector space with the same field of scalars as
$V_1$ and $V_2$, and if $T_1 : V_1 \to V_2$, $T_2 : V_2 \to V_3$ are
linear mappings, then
\begin{equation}
	(T_2 \circ T_1)^* = T_1^* \circ T_2^*.
\end{equation}

	If $T$ is a linear mapping from $V_1$ to $V_2$, then we can
pass to the second duals to get a linear transformation $T^{**} :
V_1^{**} \to V_2^{**}$.  As in Section \ref{second duals}, there are
canonical isomorphisms between $V_1$ and $V_1^{**}$, and between $V_2$
and $V_2^{**}$, which allow one to identify $T^{**}$ with a linear
mapping from $V_1$ to $V_2$.  It is easy to see that this mapping is
the same as $T$.

	The identity transformation $I = I_V$ on a vector space $V$ is
the mapping that takes each $v \in V$ to itself, and the dual of $I_V$
is equal to the identity mapping $I_{V^*}$ on $V^*$.  A one-to-one
linear transformation $T$ from $V_1$ onto $V_2$ is said to be
\emph{invertible}\index{invertible linear mappings}, which implies
that there is a linear transformation $T^{-1} : V_2 \to V_1$ such that
\begin{equation}
 T^{-1} \circ T = I_{V_1} \quad\hbox{and}\quad T \circ T^{-1} = I_{V_2}.
\end{equation}
One can check that $T : V_1 \to V_2$ is invertible if and only if $T^*
: V_2^* \to V_1^*$ is invertible, in which event
\begin{equation}
	(T^{-1})^* = (T^*)^{-1}.
\end{equation}
If $T_1 : V_1 \to V_2$ and $T_2 : V_2 \to V_3$ are both invertible, then
their composition $T_2 \circ T_1 : V_1 \to V_3$ is invertible too, with
\begin{equation}
        (T_2 \circ T_1)^{-1} = T_1^{-1} \circ T_2^{-1}.
\end{equation}

	Let $\|\cdot \|_1$ and $\|\cdot \|_2$ be norms on $V_1$ and
$V_2$ again, and let $\|\cdot \|_1^*$ and $\|\cdot \|_2^*$ be the
corresponding dual norms on $V_1^*$ and $V_2^*$.  We also have the
associated operator norm $\|\cdot \|_{op}$ on $\mathcal{L}(V_1, V_2)$,
and the operator norm $\|\cdot \|_{op*}$ on $\mathcal{L}(V_2^*,
V_1^*)$ determined by the dual norms on $V_1^*$ and $V_2^*$.  It is
easy to see that
\begin{equation}
\label{||T^*||_{op*} le ||T||_{op}}
        \|T^*\|_{op*} \le \|T\|_{op}
\end{equation}
for each linear mapping $T : V_1 \to V_2$, directly from the
definitions.  Using linear functionals as in (\ref{lambda_0(v) =
||lambda_0||^* ||v||}), one can show that the opposite inequality
holds, so that
\begin{equation}
\label{||T||_{op} = ||T^*||_{op*}}
	\|T\|_{op} = \|T^*\|_{op*}.
\end{equation}
Alternatively, one can get the opposite inequality by applying
(\ref{||T^*||_{op*} le ||T||_{op}}) to $T^*$ instead of $T$ and
identifying $T^{**}$ with $T$.

\section{Some special cases}
\label{special cases}

	Let $V$ be ${\bf R}^n$ or ${\bf C}^n$ for some positive
integer $n$, and let $\|\cdot\|_p$ be the norm on $V$ described in
Section \ref{definitions, examples} for some $p$, $1 \le p \le
\infty$.  If $T$ is a linear mapping from $V$ into $V$, then we can
express $T$ in terms of an $n \times n$ matrix $(a_{j,k})$ of real or
complex numbers, as appropriate, through the formula
\begin{equation}
\label{(T(v))_j = sum_{k=1}^n a_{j,k} v_k}
	(T(v))_j = \sum_{k=1}^n a_{j,k} \, v_k.
\end{equation}
Here $v_k$ is the $k$th component of $v \in V$, $(T(v))_j$ is the
$j$th component of $T(v)$, and conversely any $n \times n$ matrix
$(a_{j,k})$ of real or complex numbers determines such a linear
transformation $T$.  Let us write $\|T\|_{op, pp}$ for the operator
norm of $T$ with respect to the norm $\|\cdot\|_p$, used on $V$ both
as the domain and range of $T$.  These operator norms can be given
explicitly when $p = 1, \infty$, by
\begin{equation}
\label{||T||_{op, 11} = ....}
	\|T\|_{op, 11} = \max_{1 \le k \le n} \ \sum_{j=1}^n |a_{j,k}|
\end{equation}
and
\begin{equation}
\label{||T||_{op, infty infty} = ....}
   \|T\|_{op, \infty \infty} = \max_{1 \le j \le n} \ \sum_{k=1}^n |a_{j,k}|.
\end{equation}

	To see this, let $e_1, \ldots, e_n$ denote the standard basis
vectors of $V$, so that the $k$th coordinate of $e_k$ is equal to $1$
and the other coordinates are equal to $0$.  The right side of
(\ref{||T||_{op, 11} = ....}) is the same as
\begin{equation}
	\max_{1 \le k \le n} \|T(e_k)\|_1.
\end{equation}
This is obviously less than or equal to $\|T\|_{op, 11}$, by
definition of the operator norm.  The opposite inequality can be
derived by expressing any $v \in V$ as a linear combination of the
$e_l$'s and estimating $\|T(v)\|_1$ in terms of the $\|T(e_l)\|_1$'s.
Similarly, for $p = \infty$, we use the fact that
\begin{equation}
	\|T(w)\|_{\infty} = \max_{1 \le j \le n} |(T(w))_j|
\end{equation}
for every $w \in V$, by the definition of the $\|\cdot\|_{\infty}$ norm.
Clearly
\begin{equation}
\label{|(T(w))_j| le sum_{k=1}^n |a_{j,k}|}
	|(T(w))_j| \le \sum_{k=1}^n |a_{j,k}|
\end{equation}
when $w \in V$ and $\|w\|_{\infty} \le 1$, so that
\begin{equation}
 \|T\|_{op, \infty \infty} \le \max_{1 \le j \le n} \sum_{k = 1}^n |a_{j, k}|.
\end{equation}
To get the opposite inequality, one can observe that for each $j$
there is a $w \in V$ such that $\|w\|_\infty = 1$ and equality holds
in (\ref{|(T(w))_j| le sum_{k=1}^n |a_{j,k}|}).

	If $(a_{j, k})$ happens to be a diagonal matrix, so that
$a_{j, k} = 0$ when $j \ne k$, then $\|T\|_{op, pp}$ is equal to the
maximum of $|a_{j, j}|$, $1 \le j \le n$, for every $p$.  Otherwise,
it may not be so easy to compute $\|T\|_{op, pp}$ when $1 < p <
\infty$.  A famous theorem of Schur states that
\begin{equation}
\label{Schur's inequality}
	\|T\|_{op, pp} 
   \le \|T\|_{op, 11}^{1/p} \, \|T\|_{op, \infty \infty}^{1 - 1/p}.
\end{equation}
To show this, fix $p \in (1, \infty)$, and observe that
\begin{eqnarray}
	|(T(v))_j|^p & \le & 
 \Big(\sum_{k=1}^n |a_{j,k}| \Big)^{p-1} \sum_{l=1}^n |a_{j,l}| \, |v_l|^p \\
 & \le & \|T\|_{op, \infty \infty}^{p-1} \sum_{l=1}^n |a_{j,l}| \, |v_l|^p.
								\nonumber
\end{eqnarray}
for each $j = 1, \ldots, n$ and $v \in V$.  This uses H\"older's
inequality or simply the convexity of $\phi(r) = |r|^p$ for the first
inequality, and (\ref{||T||_{op, infty infty} = ....}) for the second.
Therefore
\begin{eqnarray}
	\sum_{j=1}^n |(T(v))_j|^p 
   & \le & \|T\|_{op, \infty \infty}^{p-1} 
		\sum_{j=1}^n \sum_{l=1}^n |a_{j,l}| \, |v_l|^p	\\
   & \le & \|T\|_{op, \infty \infty}^{p-1} 
		\, \|T\|_{op, 11} \sum_{l=1}^n |v_l|^p.		\nonumber
\end{eqnarray}
Schur's theorem follows by taking the $p$th root of both sides of this
inequality.

\section{Inner product spaces}
\label{inner product spaces}

	Let $V$ be a real or complex vector space.  An \emph{inner
product}\index{inner products} on $V$ is a scalar-valued function
$\langle \cdot, \cdot \rangle$ on $V \times V$ with the following
properties: (1) for each $w \in V$, $v \mapsto \langle v, w \rangle$
is a linear functional on $V$; (2) if $V$ is a real vector space, then
\begin{equation}
\label{langle w, v rangle = langle v, w rangle for every v, w in V}
	\langle w, v \rangle = \langle v, w \rangle 
				  \quad\hbox{for every } v, w \in V,
\end{equation}
and if $V$ is a complex vector space, then
\begin{equation}
\label{langle w, v rangle = overline{langle v, w rangle} for every v, w in V}
	\langle w, v \rangle = \overline{\langle v, w \rangle}
				  \quad\hbox{for every } v, w \in V;
\end{equation}
(3) the inner product is positive definite, in the sense that $\langle
v, v \rangle$ is a positive real number for every $v \in V$ such that
$v \ne 0$.  Note that $\langle v, w \rangle = 0$ whenever $v = 0$ or
$w = 0$, and that $\langle v, v \rangle \in {\bf R}$ for every $v \in
V$ even when $V$ is complex.

	A vector space with an inner product is called an \emph{inner
product space}.  If $(V, \langle \cdot, \cdot \rangle)$ is an inner
product space, then we put
\begin{equation}
	\|v\| = \langle v, v \rangle^{1/2}
\end{equation}
for every $v \in V$.  The \emph{Cauchy--Schwarz
inequality}\index{Cauchy--Schwarz inequality} says that
\begin{equation}
\label{Cauchy--Schwarz inequality in inner product spaces}
	|\langle v, w \rangle| \le \|v\| \, \|w\|
\end{equation}
for every $v, w \in V$, and it can be proved using the fact that
\begin{equation}
	\langle v + a \, w, v + a \, w \rangle \ge 0
\end{equation}
for all scalars $a$.  One can show that $\|\cdot \|$ satisfies the
triangle inequality, and is therefore a norm on $V$, by expanding
$\|v+w\|^2$ as a sum of inner products and applying the
Cauchy--Schwarz inequality.  For each positive integer $n$, the
standard inner products on ${\bf R}^n$ and ${\bf C}^n$ are given by
\begin{equation}
	\langle v, w \rangle = \sum_{j=1}^n v_j \, w_j
\end{equation}
on ${\bf R}^n$ and
\begin{equation}
	\langle v, w \rangle = \sum_{j=1}^n v_j \, \overline{w_j}
\end{equation}
on ${\bf C}^n$, and the associated norms are the standard Euclidean
norms on ${\bf R}^n$, ${\bf C}^n$.

	Let $(V, \langle \cdot, \cdot \rangle)$ be a real or complex
inner product space.  A pair of vectors $v, w \in V$ are said to be
\emph{orthogonal}\index{orthogonal vectors} if
\begin{equation}
\label{langle v, w rangle = 0}
	\langle v, w \rangle = 0,
\end{equation}
which may be expressed symbolically by $v \perp w$.
This condition is symmetric in $v$ and $w$, and implies that
\begin{equation}
\label{||v + w||^2 = ||v||^2 + ||w||^2}
	\|v + w\|^2 = \|v\|^2 + \|w\|^2.
\end{equation}
A collection $v_1, \ldots, v_n$ of vectors in $V$ is said to be
\emph{orthonormal}\index{orthonormal vectors} if $v_j \perp v_l$ when
$j \ne l$ and $\|v_j\| = 1$ for each $j$.  In this case, if $c_1,
\ldots, c_n$ are scalars and
\begin{equation}
	w = c_1 \, v_1 + \cdots + c_n \, v_n,
\end{equation}
then
\begin{equation}
        c_j = \langle w, v_j \rangle
\end{equation}
for each $j$, and
\begin{equation}
	\|w\|^2 = \sum_{j=1}^n |c_j|^2.
\end{equation}
An \emph{orthonormal basis}\index{orthonormal bases} for $V$ is an
orthonormal collection of vectors in $V$ whose linear span is equal to
$V$.  For example, the standard bases in ${\bf R}^n$ and ${\bf C}^n$
are orthonormal with respect to the standard inner products.

	Suppose that $v_1, \ldots, v_n$ are orthonormal vectors in
$V$, and define a linear transformation $P : V \to V$ by
\begin{equation}
\label{P(w) = sum_{j=1}^n langle w, v_j rangle v_j}
	P(w) = \sum_{j=1}^n \langle w, v_j \rangle \, v_j.
\end{equation}
Observe that
\begin{equation}
        P(w) = w
\end{equation}
when $w \in V$ is a linear combination of $v_1, \ldots, v_n$.
If $w$ is any vector in $V$, then
\begin{equation}
	\langle P(w), v_j \rangle = \langle w, v_j \rangle
\end{equation}
for each $j$, and hence $(w - P(w)) \perp v_j$ for each $j$.
Thus $(w - P(w)) \perp w$, and
\begin{equation}
\label{||w||^2 = ||P(w)||^2 + ||w - P(w)||^2 = ...}
	\|w\|^2 = \|P(w)\|^2 + \|w - P(w)\|^2
	         = \sum_{j=1}^n |\langle w, v_j \rangle|^2 + \|w - P(w)\|^2.
\end{equation}

	Suppose that $w$ is an element of $V$ which is not in the span
of $v_1, \ldots, v_n$.  This implies that $w - P(w) \ne 0$, and
\begin{equation}
	u = \frac{w - P(w)}{\|w - P(w)\|}
\end{equation}
satisfies $\|u\| = 1$ and $u \perp v_j$ for each $j$.  It follows that
$v_1, \ldots, v_n, u$ is an orthonormal collection of vectors in $V$
whose linear span is the same as the span of $v_1, \ldots, v_n, w$.
By repeating the process, we can extend $v_1, \ldots, v_n$ to an
orthonormal basis of $V$.  In particular, every finite-dimensional
inner product space has an orthonormal basis.

	The \emph{orthogonal complement}\index{orthogonal complements}
$W^\perp$ of a linear subspace $W$ of $V$ is defined by
\begin{equation}
	W^\perp = \{v \in V : \langle v, w \rangle = 0
				\ \hbox{ for every } w \in W\},
\end{equation}
and is also a linear subspace of $V$.  Note that
\begin{equation}
	W \cap W^\perp = \{0\},
\end{equation}
since $v \perp v$ if and only if $v = 0$.  Let $v_1, \ldots, v_n$ be
an orthonormal basis for $W$, and let $P : V \to V$ be defined as in
(\ref{P(w) = sum_{j=1}^n langle w, v_j rangle v_j}).  Thus
\begin{equation}
\label{P(v) in W and v - P(v) in W^perp}
	P(v) \in W \quad\hbox{and}\quad v - P(v) \in W^\perp
\end{equation}
for every $v \in V$.  If $v \in V$ and $x, y \in W$ satisfy $v - x, v
- y \in W^\perp$, then $x - y \in W \cap W^\perp$, and hence $x - y =
0$.  Therefore $P(v)$ is uniquely determined by (\ref{P(v) in W and v
- P(v) in W^perp}), and does not depend on the choice of orthonormal
basis $v_1, \ldots, v_n$ for $W$.  This linear transformation is
called the \emph{orthogonal projection}\index{orthogonal projections}
of $V$ onto $W$, and may be denoted $P_W$.

	Note that
\begin{equation}
	\lambda_w(v) = \langle v, w \rangle
\end{equation}
defines a linear functional on $V$ for each $w \in V$, and that
\begin{equation}
	|\lambda_w(v)| \le \|v\| \, \|w\|
\end{equation}
for every $v \in V$, by the Cauchy--Schwarz inequality.  Thus the dual
norm of $\lambda_w$ corresponding to the norm $\|\cdot \|$ on $V$ is
less than or equal to $\|w\|$.  In fact, the dual norm of $\lambda_w$
is equal to $\|w\|$, because $\lambda_w(w) = \|w\|^2$.

        Conversely, every linear functional $\lambda$ on $V$ can be
represented as $\lambda_w$ for some $w \in V$.  To see this, let $v_1,
\ldots, v_n$ be an orthonormal basis of $V$.  If
\begin{equation}
        w = \sum_{j = 1}^n \lambda(v_j) \, v_j,
\end{equation}
then
\begin{equation}
        \langle v_j, w \rangle = \lambda(v_j)
\end{equation}
for each $j$, and hence $\lambda(v) = \lambda_w(v)$ for every $v \in V$.
It is easy to reverse this argument to show that $w$ is uniquely
determined by $\lambda$.

\beginremark
\label{polarization, parallelogram law}
{\rm If $(V, \langle \cdot, \cdot \rangle)$ is a real or complex inner
product space, and if $\|\cdot\|$ is the norm associated to the inner
product, then there is a simple formula for the inner product in terms
of the norm, through \emph{polarization}.  Specifically, the
polarization identities\index{polarization identities} are
\begin{equation}
\label{4 langle v, w rangle, real case}
	4 \, \langle v, w \rangle = \|v + w\|^2 - \|v - w\|^2
\end{equation}
in the real case, and
\begin{equation}
\label{4 langle v, w rangle, complex case}
	4 \, \langle v, w \rangle 
		= \|v + w\|^2 - \|v - w\|^2 
			+ i \, \|v + i \, w\|^2 - i \, \|v - i \, w\|^2
\end{equation}
in the complex case.  The norm also satisfies the \emph{parallelogram
law}\index{parallelogram law}
\begin{equation}
\label{parallelogram law}
	\|v + w\|^2 + \|v - w\|^2 = 2 \, (\|v\|^2 + \|w\|^2)
		\quad\hbox{for every } v, w \in V.
\end{equation}
Conversely, if $V$ is a vector space and $\|\cdot\|$ is a norm on $V$
which satisfies the parallelogram law, then there is an inner product
on $V$ for which $\|\cdot\|$ is the associated norm.  This is a
well-known fact, which can be established using (\ref{4 langle v, w
rangle, real case}) or (\ref{4 langle v, w rangle, complex case}), as
appropriate, to define $\langle v, w \rangle$, and using the
parallelogram law to show that this is an inner product.}
\end{remark}

\section{Some more special cases}
\label{more special cases}

        Let $V_1$, $V_2$ be vector spaces, both real or both complex,
equipped with norms $\|\cdot \|_{V_1}$, $\|\cdot \|_{V_2}$, respectively.
Suppose first that $V_1$ is ${\bf R}^n$ or ${\bf C}^n$ for some positive
integer $n$, and that $\|\cdot \|_{V_1}$ is the norm $\|\cdot \|_1$ from
Section \ref{definitions, examples}.  Let $e_1, \ldots, e_n$ be the
standard basis vectors in $V_1$, so that the $k$th coordinate of $e_k$
is equal to $1$ for each $k$, and the rest of the coordinates are equal
to $0$.  If $T$ is any linear mapping from $V_1$ into $V_2$, then
\begin{equation}
\label{||T||_{op} = max_{1 le k le n} ||T(e_k)||_{V_2}}
        \|T\|_{op} = \max_{1 \le k \le n} \|T(e_k)\|_{V_2}.
\end{equation}
This reduces to (\ref{||T||_{op, 11} = ....}) when $V_2 = V_1$ with
the norm $\|\cdot \|_1$ from Section \ref{definitions, examples}, and
essentially the same argument works for any norm on any $V_2$.  As
before,
\begin{equation}
\label{||T(e_k)||_{V_2} le ||T||_{op}}
        \|T(e_k)\|_{V_2} \le \|T\|_{op}
\end{equation}
for each $k$ by definition of the operator norm, since $e_k$ has norm $1$ in
$V_1$ for each $k$, which implies that $\|T\|_{op}$ is less than or equal to
the right side of (\ref{||T||_{op} = max_{1 le k le n} ||T(e_k)||_{V_2}}).
To get the opposite inequality, one can express any $v \in V_1$ as
$\sum_{k = 1}^n v_k \, e_k$, where $v_1, \ldots, v_n$ are the coordinates
of $v$ in $V_1 = {\bf R}^n$ or ${\bf C}^n$, and observe that
\begin{equation}
        \|T(v)\|_{V_2} \le \sum_{k = 1}^n |v_k| \, \|T(e_k)\|_{V_2}
          \le \Big(\max_{1 \le k \le n} \|T(e_k)\|_{V_2}\Big) \, \|v\|_1.
\end{equation}

        Now let $V_1$ be any vector space with any norm $\|\cdot
\|_{V_1}$, and let $V_2$ be ${\bf R}^n$ or ${\bf C}^n$ with the norm
$\|\cdot \|_\infty$ from Section \ref{definitions, examples}.  Let $T$
be a linear mapping from $V_1$ into $V_2$ again, and let
$\lambda_j(v)$ be the $j$th component of $T(v)$ in ${\bf R}^n$ or
${\bf C}^n$, as appropriate, for $j = 1, \ldots, n$.  Thus $\lambda_j$
is a linear functional on $V_1$, with a dual norm
$\|\lambda_j\|_{V_1}^*$ with respect to the norm $\|\cdot \|_{V_1}$ on
$V_1$ for each $j$.  In this case, it is easy to see that
\begin{equation}
\label{||T||_{op} = max_{1 le j le n} ||lambda_j||_{V_1}^*}
        \|T\|_{op} = \max_{1 \le j \le n} \|\lambda_j\|_{V_1}^*,
\end{equation}
using the definitions of the dual norm, the operator norm, and the
norm $\|\cdot \|_\infty$ on $V_2$.  If $V_1 = V_2$ equipped with the
norm $\|\cdot \|_\infty$ from Section \ref{definitions, examples},
then (\ref{||T||_{op} = max_{1 le j le n} ||lambda_j||_{V_1}^*})
reduces to (\ref{||T||_{op, infty infty} = ....}), because of the
standard identification of the dual of the norm $\|\cdot \|_\infty$ on
${\bf R}^n$ or ${\bf C}^n$ with the norm $\|\cdot \|_1$ from Section
\ref{definitions, examples}, as in Section \ref{dual spaces, norms}.
Note that this case is dual to the previous one, in the sense that it
can be applied to the dual of a linear mapping as in the previous
paragraph.  Similarly, the remarks in the previous paragraph can be
applied to the dual of a linear mapping as in this paragraph.

\section{Quotient spaces}
\label{quotients}

	Let $V$ be a real or complex vector space, and let $W$ be a
linear subspace of $V$.  The quotient $V / W$ of $V$ by $W$ is defined
by identifying $v, v' \in V$ when $v - v' \in W$.  More precisely, one
can define an equivalence relation $\sim$ on $V$ by saying that $v
\sim v'$ when $v - v' \in W$, and the elements of $V / W$ correspond
to equivalence classes in $V$ determined by $\sim$.  By standard
arguments, $V / W$ is a vector space in a natural way, and there is a
canonical quotient mapping $q$ from $V$ onto $V / W$ that sends each
$v \in V$ to the equivalence class that contains it, and which is a
linear mapping from $V$ onto $V / W$ whose kernel is $W$.

        If $V$ is equipped with a norm $\|\cdot\|$, then there is a
natural quotient norm $\|\cdot\|_Q$ on $V / W$ defined by
\begin{equation}
	\|q(v)\|_Q = \inf \{\|v + w\| : w \in W\}.
\end{equation}
It is not too difficult to show that this does determine a norm on $V
/ W$.  More precisely, to check that $\|q(v)\|_Q > 0$ when $v \in V
\backslash W$ and hence $q(v) \ne 0$ in $V / W$, one can use the
remarks at the end of Section \ref{definitions, examples}, and the
fact that linear subspaces of ${\bf R}^n$ and ${\bf C}^n$ are closed
with respect to the standard topology on those spaces.  Note that the
operator norm of $q$ is less than or equal to $1$ with respect to the
given norm on $V$ and the corresponding quotient norm on $V / W$, at
that it is equal to $1$ when $W \ne V$.

	The dual $(V / W)^*$ of $V / W$ can be identified with a
subspace of $V^*$ in a natural way.  Of course, every linear
functional on $V / W$ determines a linear functional on $V$, by
composition with the quotient mapping $q$.  The linear functionals on
$V$ that occur in this way are exactly those that are equal to $0$ on
$W$.  If $\lambda$ is a linear functional on $V$ that is equal to $0$
on $W$ and $k$ is a nonnegative real number, then the statements
\begin{equation}
	|\lambda(v)| \le k \, \|v\|
		\quad\hbox{for every } v \in V
\end{equation}
and
\begin{equation}
	|\lambda(v)| \le k \, \inf \{\|v + w\| : w \in W\}
		\quad\hbox{for every } v \in V
\end{equation}
are equivalent to each other.  This means that the dual norm of
$\lambda$ as a linear functional on $V$ with respect to $\|\cdot \|$
is the same as the dual norm of the linear functional on $V / W$ that
corresponds to $\lambda$ under the quotient mapping with respect to
the quotient norm on $V / W$.

        If $v$ is any element of $V$, then there is a $w_0 \in W$ such that
\begin{equation}
\label{||v + w_0|| le ||v + w||}
        \|v + w_0\| \le \|v + w\|
\end{equation}
for every $w \in W$, so that the infimum in the definition of the
quotient norm is attained.  To see this, we may as well suppose that
$V = {\bf R}^n$ or ${\bf C}^n$ for some positive integer $n$, since
every real or complex vector space of positive finite dimension is
isomorphic to one of these.  As in Section \ref{definitions,
examples}, the norm $\|\cdot \|$ defines a continuous function on
${\bf R}^n$ or ${\bf C}^n$, as appropriate, with respect to the
standard Euclidean metric and topology.  It is also well known that
linear subspaces of ${\bf R}^n$ and ${\bf C}^n$ are closed sets with
respect to the standard Euclidean metric and topology.  Remember too
that $\|\cdot \|$ is bounded from below by a positive constant
multiple of the standard Euclidean metric on ${\bf R}^n$ or ${\bf
C}^n$, as appropriate, as in (\ref{b |v| le ||v||}).  Using this, it
suffices to consider a bounded subset of $W$ when minimizing $\|v +
w\|$ over $w \in W$.  This permits the existence of the minimum to be
derived from well-known results about minimizing continuous functions
on compact sets, because closed and bounded subsets of ${\bf R}^n$ and
${\bf C}^n$ are compact.

\section{Projections}
\label{projections}

        Let $V$ be a real or complex vector space, and let $U$ and $W$
be linear subspaces of $V$.  Suppose that $U \cap W = \{0\}$, and that
every $v \in V$ can be expressed as
\begin{equation}
\label{v = u + w}
        v = u + w
\end{equation}
for some $u \in U$ and $w \in W$.  If $u' \in U$ and $w' \in W$ also
satisfy $v = u' + w'$, then
\begin{equation}
        u - u' = w' - w,
\end{equation}
and this implies that $u = u'$ and $w = w'$, because $u - u' \in U$,
$w - w' \in W$, and $U \cap W = \{0\}$.  In this case, $U$ and $W$ are
said to be \emph{complementary}\index{complementary linear subspaces}
in $V$.

        Consider the mapping $P : V \to V$ defined by
\begin{equation}
        P(v) = u
\end{equation}
for each $v \in V$, where $u \in U$ is as in (\ref{v = u + w}).  It is
easy to see that $P$ is a linear mapping of $V$ onto $U$ with kernel
equal to $W$, and that
\begin{equation}
        P(u) = u
\end{equation}
for every $u \in U$.  Conversely, suppose that $P$ is a linear mapping
from $V$ onto a linear subspace $U$ of $V$ such that the restriction
of $P$ to $U$ is equal to the identity mapping on $U$.  If $W$ is the
kernel of $P$, then
\begin{equation}
        v - P(v) \in W
\end{equation}
for every $v \in V$, and it follows that $U$ and $W$ are complementary
in $V$.

        A linear mapping $P : V \to V$ is said to be a
\emph{projection}\index{projections} if
\begin{equation}
        P \circ P = P.
\end{equation}
This is equivalent to saying that the restriction of $P$ to $U = P(V)$
is the identity mapping on $U$, as in the previous paragraph, so that
$U$ is complementary to the kernel $W$ of $P$.  If $P$ is a projection
on $V$, then it is easy to see that $I - P$ is also a projection on
$V$, where $I$ denotes the identity mapping on $V$, because
\begin{equation}
\label{(I - P) circ (I - P) = I - P - P + P circ P = I - P}
        (I - P) \circ (I - P) = I - P - P + P \circ P = I - P.
\end{equation}
More precisely, $I - P$ maps $V$ onto the kernel $W$ of $P$, and the
kernel of $I - P$ is $U = P(V)$.  Of course, orthogonal projections
onto linear subspaces of inner product spaces are projections in this
sense.

        Let $U$ and $W$ be linear subspaces of $V$ again, and let $q$
be the canonical quotient mapping from $V$ onto $V / W$, as in the
previous section.  It is easy to see that $U$ and $W$ are
complementary in $V$ if and only if the restriction of $q$ to $U$ is a
one-to-one mapping from $U$ onto $V / W$.  Suppose that this is the
case, and let $P$ be the corresponding projection of $V$ onto $U$ with
kernel $W$.

        Let $\|\cdot \|$ be a norm on $V$, let $\|\cdot \|_{op}$ be
the corresponding operator norm for linear mappings on $V$, and let
$\|\cdot \|_Q$ be the corresponding quotient norm on $V / W$.
If $u \in U$ and $w \in W$, then $P(u + w) = u$, and hence
\begin{equation}
\label{||u|| = ||P(u + w)|| le ||P||_{op} ||u + w||}
	\|u\| = \|P(u + w)\| \le \|P\|_{op} \, \|u + w\|.
\end{equation}
This implies that
\begin{equation}
	\|P\|_{op}^{-1} \, \|u\| \le \|q(u)\|_Q \le \|u\|
\end{equation}
for every $u \in U$.

        Observe that
\begin{equation}
\label{||P||_{op} = ||P circ P||_{op} le ||P||_{op}^2}
        \|P\|_{op} = \|P \circ P\|_{op} \le \|P\|_{op}^2.
\end{equation}
Thus
\begin{equation}
        \|P\|_{op} \ge 1
\end{equation}
when $P \ne 0$.  If $\|P\|_{op} = 1$, then we get that
\begin{equation}
        \|q(u)\|_Q = \|u\|
\end{equation}
for every $u \in U$.  Orthogonal projections onto nontrival subspaces
of inner product spaces have operator norm equal to $1$, for instance.

        Suppose that $V = {\bf R}^n$ or ${\bf C}^n$ for some positive
integer $n$, and let $I$ be a subset of the set $\{1, \ldots, n\}$ of
positive integers less than or equal to $n$.  Let $U_I$ be the linear
subspace of $V$ consisting of vectors $u$ such that $u_j = 0$ when $j
\not\in I$, and let $W_I$ be the complementary subspace consisting of
vectors $w$ such that $w_j = 0$ when $j \in I$.  The associated
projection $P_I$ of $V$ onto $U_I$ with kernel $W_I$ sends $v \in V$
to the vector whose $j$th coordinate is equal to $v_j$ when $j \in I$,
and to $0$ when $j \not\in I$.  If $V$ is equipped with a norm
$\|\cdot \|_p$ as in Section \ref{definitions, examples} for some $p$,
$1 \le p \le \infty$, then the corresponding operator norm of $P_I$ is
equal to $1$ when $I \ne \emptyset$, so that $P_I \ne 0$.

        Let $V$ be any real or complex vector space with a norm
$\|\cdot \|$ again, and let $u_1$ be an element of $V$ such that
$\|u_1\| = 1$.  As in Section \ref{second duals}, there is a linear
functional $\lambda$ on $V$ such that
\begin{equation}
        \lambda(u_1) = 1
\end{equation}
and the dual norm of $\lambda$ with respect to the given norm $\|\cdot
\|$ on $V$ is equal to $1$.  Under these conditions, one can check that
\begin{equation}
\label{P(v) = lambda(v) u_1}
        P(v) = \lambda(v) \, u_1
\end{equation}
is a projection of $V$ onto the $1$-dimensional linear subspace $U$ of
$V$ spanned by $u_1$ with operator norm equal to $1$.

\section{Extensions and liftings}
\label{extensions, liftings}

	Let $V_1$ and $V_2$ be vector spaces, both real or both
complex, and equipped with norms.  If $U_1$ is a linear subspace of
$V_1$, and $T$ is a linear mapping from $U_1$ into $V_2$, then it is
easy to see that there is an extension $\widehat{T}$ of $T$ to a
linear mapping from $V_1$ into $V_2$.  The operator norm of
$\widehat{T}$ is automatically greater than or equal to the operator
norm of $T$, and one would like to choose $\widehat{T}$ so that its
operator norm is as small as possible.  If $P_1$ is a projection from
$V_1$ onto $U_1$, then the composition $T \circ P_1$ is an extension
of $T$ to $V_1$ whose operator norm is less than or equal to the
product of the operator norm of $T$ on $U_1$ and the operator norm of
$P_1$ on $V_1$.  Conversely, if $V_2 = U_1$ and $T$ is the identity
mapping on $U_1$, then an extension of $T$ to a linear mapping from
$V_1$ into $V_2$ is the same as a projection from $V_1$ onto $U_1$.

	If $V_2$ is $1$-dimensional, then this extension problem is
equivalent to the one for linear functionals discussed in Section
\ref{second duals}.  Similarly, suppose that $V_2$ is ${\bf R}^n$ or
${\bf C}^n$ for some positive integer $n$, equipped with the norm
$\|\cdot\|_{\infty}$ defined in Section \ref{definitions, examples}.
In this case, a linear mapping $T$ from another vector space into
$V_2$ is equivalent to $n$ linear functionals on the vector space, and
the operator norm of a $T$ is equal to the maximum of the dual norms
of the corresponding $n$ linear functionals, as in the second part of
Section \ref{more special cases}.  This permits the extension problem
for $T$ to be reduced to its counterpart for linear functionals again.

        Now let $W_2$ be a linear subspace of $V_2$, and let $L$ be a
linear mapping from $V_1$ into $V_2 / W_2$.  It is easy to see that
there is a linear mapping $\widetilde{L}$ from $V_1$ into $V_2$ whose
composition with the canonical quotient mapping $q_2$ from $V_2$ onto
$V_2 / W_2$ is equal to $L$, and one would like to choose
$\widetilde{L}$ so that its operator norm is as small as possible.
This problem is dual to the extension problem discussed in the
previous paragraphs.  Of course, the operator norm of $\widetilde{L}$
is greater than or equal to the operator norm of $L$, with respect to
the quotient norm on $V_2 / W_2$ that corresponds to the given norm on
$V_2$.  One way to approach this problem is to use a linear subspace
$U_2$ of $V_2$ which is complementary to $W_2$, so that the
restriction of $q_2$ to $U_2$ is a one-to-one linear mapping of $U_2$
onto $V_2 / W_2$.  In this case, one can get a lifting $\widetilde{L}$
of $L$ to $V_2$ by composing $L$ with the inverse of the restriction
of $q_2$ to $U_2$.  Conversely, if $V_1 = V_2 / W_2$ and $L$ is the
identity mapping on $V_2 / W_2$, then a lifting of $L$ to a linear
mapping $\widetilde{L}$ from $V_2 / W_2$ into $V_2$ whose composition
with $q_2$ is the identity mapping on $V_2 / W_2$ would map $V_2 /
W_2$ onto a linear subspace $U_2$ of $V_2$ which is complementary to $W_2$.

        Suppose that $V_1$ has dimension $1$, and let $v_1$ be an
element of $V_1$ with norm $1$.  Let $L$ be a linear mapping from
$V_1$ into $V_2 / W_2$, and let $v_2$ be an element of $V_2$ such that
$q_2(v_2) = L(v_1)$ and the norm of $v_2$ in $V_2$ is equal to the
quotient norm of $L(v_1)$ in $V_2 / W_2$.  The existence of $v_2$
follows from the discussion of minimization at the end of Section
\ref{quotients}.  If $\widetilde{L}$ is the linear mapping from $V_1$
into $V_2$ that sends $v_1$ to $v_2$, then $\widetilde{L}$ is a
lifting of $L$ with the same operator norm as $L$.  Similarly, if
$V_1$ is ${\bf R}^n$ or ${\bf C}^n$ with the norm $\|\cdot \|_1$ from
Section \ref{definitions, examples}, then one can get a lifting
$\widetilde{L}$ of $L$ to a linear mapping from $V_1$ into $V_2$ with
the same operator norm as $L$ by lifting the $n$ vectors $L(e_j)$ in
$V_2 / W_2$ to $V_2$ for each of the standard basis vectors $e_j$ in
${\bf R}^n$ or ${\bf C}^n$, since the operator norm of a linear
mapping on $V_1$ may be computed as in the first part of Section
\ref{more special cases}.

\section{Minimizing distances}
\label{minimizing distances}

        Let $V$ be a real or complex vector space with a norm $\|\cdot
\|$, and let $W$ be a linear subspace of $V$.  If $v$ is any element
of $V$, then there is a $w_1 \in W$ such that
\begin{equation}
\label{||v - w_1|| le ||v - w||}
        \|v - w_1\| \le \|v - w\|
\end{equation}
for every $w \in W$.  This is equivalent to the minimization problem
discussed at the end of Section \ref{quotients}, with $w$ replaced by
$-w$ and $w_1 = - w_0$

        Suppose for the moment that there is an inner product $\langle
\cdot, \cdot \rangle$ on $V$ for which $\|\cdot \|$ is the
corresponding norm, and let $P_W(v)$ be the orthogonal projection of
$v$ onto $W$, as in Section \ref{inner product spaces}.  Thus $P_W(v)
\in W$ and $v - P_W(v) \in W^\perp$, as in (\ref{P(v) in W and v -
P(v) in W^perp}).  If $w$ is any element of $W$, then it follows that
$P_W(v) - w \in W$, and hence
\begin{equation}
\label{||v - w||^2 = ||v - P_W(v)||^2 + ||P_W(v) - w||^2}
        \|v - w\|^2 = \|v - P_W(v)\|^2 + \|P_W(v) - w\|^2,
\end{equation}
because $(v - P_W(v)) \perp (P_W(v) - w)$.  This implies that $P_W(v)$
minimizes the distance to $v$ among elements of $W$ in this case, and
that $P_W(v)$ is the only element of $W$ with this property.

        Let $\|\cdot \|$ be any norm on $V$ again, and let $\|\cdot
\|_{op}$ be the corresponding operator norm for linear mappings on
$V$.  If $P_1$ is a projection of $V$ onto $W$, then the kernel of
$I - P_1$ is equal to $W$, and hence
\begin{eqnarray}
        \|v - P_1(v)\| = \|(I - P_1)(v)\| & = & \|(I - P_1)(v - w)\|  \\
                           & \le & \|I - P_1\|_{op} \, \|v - w\| \nonumber
\end{eqnarray}
for every $w \in W$.  If $\|I - P_1\|_{op} = 1$, then it follows that
\begin{equation}
\label{||v - P_1(v)|| le ||v - w||}
        \|v - P_1(v)\| \le \|v - w\|
\end{equation}
for every $w \in W$, so that $w_1 = P_1(v) \in W$ satisfies (\ref{||v
- w_1|| le ||v - w||}).

        Suppose for the sake of convenience now that $V = {\bf R}^n$
or ${\bf C}^n$ for some positive integer $n$, which is not a real
restriction, since every real or complex vector space of positive
finite dimension is isomorphic to one of these.  Let $E$ be a nonempty
subset of $V$ which is closed with respect to the standard Euclidean
metric and topology on $V$.  If $v \in V$, then there is a $w_1 \in E$
which minimizes the distance to $v$ with respect to the norm $\|\cdot
\|$ on $V$, in the sense that (\ref{||v - w_1|| le ||v - w||}) holds
for every $w \in E$.  This follows from the same type of argument
using continuity and compactness as before.  More precisely, although
$E$ may not be bounded, and hence not compact, it suffices to consider
a bounded subset of $E$ for this minimization problem.

        Let $V$ be any real or complex vector space again, and let
$B_1$ be the closed unit ball associated to the norm $\|\cdot \|$ on
$V$, as in Section \ref{definitions, examples}.  Let us say that $B_1$
is \emph{strictly convex}\index{strict convexity} if 
\begin{equation}
	\|t \, v + (1 - t) \, w\| < 1.
\end{equation}
for every $v, w \in V$ with $\|v\| = \|w\| = 1$ and $v \ne w$ and
every real number $t$ with $0 < t < 1$.  The unit ball in any inner
product space is strictly convex, as one can show by analyzing the
case of equality in the proof of the triangle inequality.  If $V =
{\bf R}^n$ or ${\bf C}^n$ with the norm $\|\cdot \|_p$ as in Section
\ref{definitions, examples} for some $p$, $1 < p < \infty$, then one
can check that the unit ball is strictly convex, using the strict
convexity of the function $|r|^p$.  If $n \ge 2$ and $p = 1$ or
$\infty$, then it is easy to see that the unit ball is not strictly
convex.

        Let $E$ be a nonempty convex set in $V$, and let $v$ be an
element of $V$.  Suppose that $w_1, w_2 \in E$ both minimize the
distance to $v$ with respect to $\|\cdot \|$ in $V$, in the sense that
\begin{equation}
\label{||v - w_1|| = ||v - w_2|| le ||v - w||}
        \|v - w_1\| = \|v - w_2\| \le \|v - w\|
\end{equation}
for every $w \in E$.  If $0 < t < 1$, then $w = t \, w_1 + (1 - t) \,
w_2 \in E$, because $E$ is convex.  However, if $B_1$ is strictly
convex and $w_1 \ne w_2$, then the norm of
\begin{equation}
        v - w = t \, (v - w_1) + (1 - t) \, (v - w_2),
\end{equation}
is strictly less than the common value of the norms of $v - w_1$ and
$v - w_2$, contradicting (\ref{||v - w_1|| = ||v - w_2|| le ||v - w||}).
This shows that $w_1 = w_2$ under these conditions when $B_1$ is strictly
convex.

        Note that we could simply take $t = 1/2$ in the preceding
argument.  If the norm on $V$ is associated to an inner product, then
the strict convexity property of the unit ball with $t = 1/2$ follows
from the parallelogram law (\ref{parallelogram law}).

\chapter{Structure of linear operators}
\label{structure of linear operators}

        In this chapter, we continue to restrict our attention to
finite-dimensional vector spaces.

\section{The spectrum and spectral radius}
\label{spectrum, spectral radius}

	Let $V$ be a complex vector space with positive dimension, and
let $T$ be a linear operator from $V$ into $V$.  The
\emph{spectrum}\index{spectrum} of $T$ is the set of complex numbers
$\alpha$ such that $\alpha$ is an \emph{eigenvalue}\index{eigenvalues}
of $T$, which is to say that there is a $v \in V$ such that $v \ne 0$ and
\begin{equation}
	T(v) = \alpha v.
\end{equation}
In this case, $v$ is said to be an
\emph{eigenvector}\index{eigenvectors} of $T$ with eigenvalue
$\alpha$.  If $\alpha$ is an eigenvalue of $T$, then
\begin{equation}
\label{{v in V : T(v) = alpha v}}
	\{v \in V : T(v) = \alpha \, v\}
\end{equation}
is a nontrivial linear subspace of $V$, called the
\emph{eigenspace}\index{eigenspaces} of $T$ associated to $\alpha$.

        It is well known that a one-to-one linear mapping $R : V \to
V$ automatically maps $V$ onto itself, and hence is invertible,
because $R(V)$ is a linear subspace of $V$ with the same dimension as
$V$.  If $R$ is not invertible on $V$, then it follows that $R$ is not
one-to-one, so that the kernel of $R$ is nontrivial.  By definition,
$\alpha \in {\bf C}$ is an eigenvalue of $T$ when the kernel of $T -
\alpha \, I$ is nontrivial, where $I$ denotes the identity
transformation on $V$.  Equivalently, $\alpha \in {\bf C}$ is not in
the spectrum of $T$ when $T - \alpha \, I$ is an invertible linear
operator on $V$.

        A famous theorem states that every linear operator $T$ on $V$
has at least one eigenvalue.  To see this, note that $\alpha \in {\bf
C}$ lies in the spectrum of $T$ exactly when the determinant of $T -
\alpha \, I$ is $0$.  The determinant of $T - \alpha \, I$ is a
polynomial in $\alpha$, whose degree is equal to the dimension of $V$.
By the ``Fundamental Theorem of Algebra'', $\det(T - \alpha \, I)$ has
at least one root, as desired.

        This argument also shows that the number of distinct
eigenvalues of $T$ is less than or equal to the dimension of $V$,
since a polynomial of degree $n$ has at most $n$ roots.  The
\emph{spectral radius}\index{spectral radius} $\rad(T)$ of $T$ is
defined to be the maximum of $|\alpha|$, where $\alpha \in {\bf C}$ is
an eigenvalue of $T$.  Thus $T - \alpha \, I$ is invertible when
$|\alpha| > \rad(T)$, and $\rad(T)$ is the largest nonnegative real
number with this property.

	Let $\|\cdot \|$ be a norm on $V$, and let $\|\cdot \|_{op}$
be the corresponding operator norm for linear transformations on $V$,
as in Section \ref{linear transformations}.  If $\alpha \in {\bf C}$
is an eigenvalue of $T$, and $v \in V$ is a nonzero eigenvector
corresponding to $\alpha$, then
\begin{equation}
	|\alpha| \, \|v\| \le \|T\|_{op} \, \|v\|,
\end{equation}
and hence $|\alpha| \le \|T\|_{op}$.  It follows that $\rad(T) \le
\|T\|_{op}$.

        Now let $n$ be a positive integer, and let us check that
\begin{equation}
        \rad(T^n) = \rad(T)^n.
\end{equation}
If $\alpha$ is an eigenvalue of $T$, then $\alpha^n$ is
obviously an eigenvalue of $T^n$ for each $n$, and hence
\begin{equation}
        \rad(T)^n \le \rad(T^n).
\end{equation}
To get the opposite inequality, suppose that $\beta$ is an eigenvalue
of $T^n$, and let us show that $\alpha$ is an eigenvalue of $T$ for
some complex number $\alpha$ such that $\alpha^n = \beta$.  Let
$\alpha_1, \ldots, \alpha_n$ be the $n$th roots of $\beta$, so that
\begin{equation}
	z^n - \beta = (z - \alpha_1) \cdots (z - \alpha_n).
\end{equation}
This implies that
\begin{equation}
\label{T^n - beta I = (T - alpha_1 I) cdots (T - alpha_n I)}
	T^n - \beta \, I = (T - \alpha_1 \, I) \cdots (T - \alpha_n \, I),
\end{equation}
where the product of linear operators on $V$ is defined by their
composition.  If $T - \alpha_j \, I$ is invertible on $V$ for each $j
= 1, \ldots, n$, then it follows that $T^n - \beta \, I$ is also
invertible on $V$, because the composition of invertible operators is
invertible.  If $\beta$ is an eigenvalue of $T$, then $T^n - \beta \,
I$ is not invertible, and hence $T - \alpha_j \, I$ is not invertible
for some $j$.  This says exactly that $\alpha_j$ is an eigenvalue of
$T$ for some $j$, as desired.

\section{Adjoints}
\label{adjoints}

	In this section, both real and complex vector spaces are
allowed.  Let $(V, \langle \cdot, \cdot \rangle)$ be an inner product
space.  If $T$ is a linear operator on $V$, then there is a unique
linear operator $T^*$ on $V$ such that
\begin{equation}
\label{characterization of T^*}\index{adjoint operators}
	\langle T(v), w \rangle = \langle v, T^*(w) \rangle
\end{equation}
for every $v, w \in V$, called the \emph{adjoint} of $T$.  More
precisely, for each $w \in V$,
\begin{equation}
\label{mu_w(v) = langle T(v), w rangle}
        \mu_w(v) = \langle T(v), w \rangle
\end{equation}
defines a linear functional on $V$.  This implies that there is a
unique element $T^*(w)$ of $V$ such that
\begin{equation}
\label{mu_w(v) = langle v, T^*(w) rangle}
        \mu_w(v) = \langle v, T^*(w) \rangle
\end{equation}
for every $v \in V$, as in Section \ref{inner product spaces}.  One
can check that $T^*$ is linear on $V$, using the fact that $T^*(w)$ is
uniquely determined by (\ref{mu_w(v) = langle v, T^*(w) rangle}).
Alternatively, if we fix an orthonormal basis for $V$, then we can
express $T$ in terms of a matrix relative to this basis.  Remember
that the \emph{transpose}\index{transpose of a matrix} of a matrix
$(a_{j,l})$ is the matrix $(b_{j, l})$ given by $b_{j, l} = a_{l, j}$.
If $V$ is a real vector space, then $T^*$ is the linear transformation
on $V$ that corresponds to the transpose of the matrix for $T$ with
respect to the same basis.  In the complex case, the entries of the
matrix for $T^*$ are the complex conjugates of the entries of the
transpose of the matrix for $T$.  It is easy to see that $T^*$ is
uniquely determined by (\ref{characterization of T^*}), so that
different orthonormal bases for $V$ lead to the same linear
transformation $T^*$ when one computes $T^*$ in terms of matrices.

        Although we are using the same notation here for the adjoint
as we did in Section \ref{linear transformations} for dual linear
mappings, we should be careful about some of the differences.  In the
real case, we can identify $V$ with its dual space $V^*$, since every
linear functional on $V$ can be represented as
\begin{equation}
        \lambda_w(v) = \langle v, w \rangle
\end{equation}
for some $w \in V$.  In this case, it is easy to see that the adjoint
of $T$ corresponds exactly to the dual linear transformation defined
previously.  However, in the complex case, the mapping from $w \in V$
to the linear functional $\lambda_w \in V^*$ is not quite linear, but
rather conjugate-linear, in the sense that multiplication of $w$ by a
complex number $a$ corresponds to multiplying $\lambda_w$ by the
complex conjugate $\overline{a}$ of $a$.  Thus the adjoint of $T$ is
not quite the same as the dual linear transformation defined earlier
in the complex case, which is also reflected in the linearity
properties of the mapping from $T$ to $T^*$ discussed next.

	Note that $I^* = I$, where $I$ is the identity transformation
on $V$.  If $S$ and $T$ are linear transformations on $V$ and $a$, $b$
are scalars, then
\begin{equation}
	(a \, S + b \, T)^* = a \, S^* + b \, T^*
\end{equation}
when $V$ is a real vector space, and
\begin{equation}
	(a \, S + b \, T)^* = 
			\overline{a} \, S^* + \overline{b} \, T^*
\end{equation}
when $V$ is a complex vector space.  Also, $(T^*)^* = T$. and
\begin{equation}
\label{(S T)^* = T^* S^*}
	(S \, T)^* = T^* \, S^*.
\end{equation}
If $T$ is invertible, then $T^*$ is invertible, and
\begin{equation}
	(T^*)^{-1} = (T^{-1})^*.
\end{equation}
It follows that $T - \lambda \, I$ is invertible if and only if $T^* -
\lambda \, I$ is invertible for each $\lambda \in {\bf R}$ in the real
case, and similarly that $T - \lambda \, I$ is invertible if and only
if $T^* - \overline{\lambda} \, I$ is invertible for each $\lambda \in
{\bf C}$ in the complex case.
	
	Let $\|\cdot \|$ be the norm on $V$ associated to the inner
product $\langle \cdot, \cdot \rangle$, and let $\|\cdot \|_{op}$ be
the corresponding operator norm for linear transformations on $V$ with
respect to $\|\cdot \|$.  Let us check that
\begin{equation}
\label{formula for ||T||_{op} in terms of the inner product}
	\|T\|_{op} = \sup \{ |\langle T(v), w \rangle| : 
			    v, w \in V, \ \|v\|, \|w\| \le 1\}
\end{equation}
for any linear transformation $T$ on $V$.  The right side of
(\ref{formula for ||T||_{op} in terms of the inner product}) is
clearly less than or equal to the operator norm of $T$, because of the
Cauchy--Schwarz inequality.  To get the opposite inequality, one can
take $w = T(v)/\|T(v)\|$ in the right side of (\ref{formula for
||T||_{op} in terms of the inner product}) when $T(v) \ne 0$.
It follows that
\begin{equation}
\label{||T^*||_{op} = ||T||_{op}}
	\|T^*\|_{op} = \|T\|_{op},
\end{equation}
because the right side of (\ref{formula for ||T||_{op} in terms of the
inner product}) is equal to the analogous quantity for $T^*$.

\section{Self-adjoint linear operators}
\label{self-adjoint linear operators}

        Let $(V, \langle \cdot, \cdot \rangle)$ be a real or complex
inner product space, as in the previous section.  A linear
transformation $A$ on $V$ is said to be
\emph{self-adjoint}\index{self-adjoint linear operators} if
\begin{equation}
\label{A^* = A}
        A^* = A.
\end{equation}
This is equivalent to the condition that
\begin{equation}
	\langle A(v), w \rangle = \langle v, A(w) \rangle
\end{equation}
for every $v, w \in V$.  Thus the identity operator $I$ on $V$ is
self-adjoint.

        Let $W$ be a linear subspace of $V$, and let $P_W$ be the
orthogonal projection of $V$ onto $W$, as in Section \ref{inner
product spaces}.  Remember that $P_W(v)$ is characterized by the
conditions $P_W(v) \in W$ and $v - P_W(v) \in W^\perp$, for each $v
\in V$.  This implies that
\begin{equation}
\label{langle P_W(v), w rangle = ... = langle v, P_W(w) rangle}
	\langle P_W(v), w \rangle = \langle P_W(v), P_W(w) \rangle
		= \langle v, P_W(w) \rangle
\end{equation}
for every $v, w \in V$, and hence that $P_W$ is self-adjoint on $V$.

        If $A$ and $B$ are self-adjoint linear operators on $V$, then
their sum $A + B$ is self-adjoint.  Similarly, if $A$ is a
self-adjoint linear operator on $V$ and $t$ is a real number, then $t
\, A$ is self-adjoint as well.  Note that it is important to take $t
\in {\bf R}$ here, even when $V$ is a complex vector space.

	If $A$ is a self-adjoint linear operator on $V$ and $V$ is
complex, then it is easy to see that
\begin{equation}
\label{langle A(v), v rangle in {bf R}}
        \langle A(v), v \rangle \in {\bf R}
\end{equation}
for every $v \in V$.  Using this, one can check that the eigenvalues
of $A$ are also real numbers.

	Suppose that $A$ is a self-adjoint linear operator on a real
or complex inner product space $(V, \langle \cdot, \cdot \rangle)$,
and that $v \in V$ is an eigenvector of $A$ with eigenvalue $\lambda$.
If $y \in V$ and $y \perp v$, then
\begin{equation}
\label{langle v, A(y) rangle = ... = lambda langle v, y rangle = 0}
	\langle v, A(y) \rangle = \langle A(v), y \rangle
		= \lambda \, \langle v, y \rangle = 0,
\end{equation}
so that $A(y) \perp v$.  If $w \in V$ is an eigenvector of $A$ with
eigenvalue $\mu \ne \lambda$, then
\begin{eqnarray}
\label{lambda langle v, w rangle = ... = mu langle v, w rangle}
	\lambda \, \langle v, w \rangle & = &
	\langle \lambda \, v, w \rangle = \langle A(v), w \rangle \\
  & = & \langle v, A(w) \rangle = \langle v, \mu \, w \rangle 
	= \mu \, \langle v, w \rangle,		\nonumber
\end{eqnarray}
which implies that $\langle v, w \rangle = 0$.

        If the dimension $n$ of $V$ is positive, then one can use an
orthonormal basis of $V$ to show that $V$ is isomorphic to ${\bf R}^n$
or ${\bf C}^n$ with its standard inner product, and we may as well
take $V = {\bf R}^n$ or ${\bf C}^n$ for the moment.  By the usual
considerations of continuity and compactness, $\langle A(v), v
\rangle$ attains its maximum and minimum on the unit sphere
\begin{equation}
\label{{v in V : ||v|| = 1}}
        \{v \in V : \|v\| = 1\}.
\end{equation}
It is well known that the critical points of $\langle A(v), v \rangle$
on the unit sphere are exactly the eigenvectors of $A$ with norm $1$,
and hence that the maximum and minimum of $\langle A(v), v \rangle$ on
the unit sphere are attained at eigenvectors of $A$.  In particular,
$A$ has a nonzero eigenvector $v$, and $A$ maps
\begin{equation}
        W = \{w \in V : \langle w, v \rangle = 0\}
\end{equation}
to itself, as in the previous paragraph.  By repeating the process,
one can show that there is an orthonormal basis of $V$ consisting of
eigenvectors of $A$.

	A linear operator $T$ on a complex inner product space $V$ is
said to be \emph{normal}\index{normal operators} if $T$ and $T^*$
commute, which is to say that
\begin{equation}
	T \circ T^* = T^* \circ T.
\end{equation}
If $T$ can be diagonalized in an orthonormal basis, then $T^*$ is
diagonalized by the same basis, and $T$ is normal.  Conversely, one
can show that a normal operator $T$ on $V$ can be diagonalized in an
orthonormal basis, as follows.  Any linear operator $T$ on $V$ can be
expressed as $A + i B$, where
\begin{equation}
\label{A = frac{T + T^*}{2} and B = frac{T - T^*}{2 i}}
	A = \frac{T + T^*}{2} \quad\hbox{and}\quad B = \frac{T - T^*}{2 i}
\end{equation}
are self-adjoint.  Thus $A$ and $B$ can each be diagonalized in an
orthonormal basis of $V$, and one would like to show that they can
both be diagonalized by the same orthonormal basis when $T$ is normal,
which implies that $A$ and $B$ commute.  If $A$ and $B$ are commuting
linear transformations on any vector space, then it is easy to see
that the eigenspaces of $A$ are invariant under $B$.  To diagonalize
$T$ in an orthonormal basis, one can first use a diagonalization of
$A$ to decompose $V$ into an orthogonal sum of eigenspaces of $A$, and
then diagonalize the restriction of $B$ to each of the eigenspaces of
$A$.

	A linear operator $T$ on $V$ is said to be an \emph{orthogonal
transformation}\index{orthogonal transformations} when $V$ is real, or
a \emph{unitary transformation}\index{unitary transformations} when
$V$ is complex, if
\begin{equation}
\label{langle T(v), T(w) rangle = langle v, w rangle}
 	\langle T(v), T(w) \rangle = \langle v, w \rangle
\end{equation}
for every $v, w \in V$.  This implies that
\begin{equation}
\label{||T(v)|| = ||v||}
	\|T(v)\| = \|v\|
\end{equation}
for every $v \in V$, and the converse holds because of polarization,
as in Remark \ref{polarization, parallelogram law}.  This condition
obviously implies that the kernel of $T$ is trivial, and hence that
$T$ is invertible, because $V$ is supposed to be finite-dimensional.
More precisely, it is easy to see that $T$ is orthogonal or unitary,
as appropriate, if and only if $T$ is invertible and
\begin{equation}
        T^{-1} = T^*.
\end{equation}
In particular, unitary operators are normal, because $T$ automatically
commutes with $T^{-1}$.

\section{Anti-self-adjoint operators}
\label{anti-self-adjoint operators}

        Let $(V, \langle \cdot, \cdot \rangle)$ be a real or complex
inner product space, as before.  A linear operator $R$ on $V$ is said to
be \emph{anti-self-adjoint}\index{anti-self-adjoint linear operators} if
\begin{equation}
\label{R^* = -R}
        R^* = -R,
\end{equation}
which is equivalent to asking that
\begin{equation}
\label{langle R(v), w rangle = - langle v, R(w) rangle}
        \langle R(v), w \rangle = - \langle v, R(w) \rangle
\end{equation}
for every $v, w \in V$.  If $V$ is a complex vector space, then $R$ is
anti-self-adjoint if and only if $R = i B$ for some self-adjoint
linear operator $B$ on $V$.  Note that the sum of two
anti-self-adjoint linear operators on $V$ is also anti-self-adjoint,
as is the product of an anti-self-adjoint linear operator and a real
number.

        If $V$ is a real inner product space and $R$ is an
anti-self-adjoint linear operator on $V$, then
\begin{equation}
\label{langle R(v), v rangle = ... = - langle R(v), v rangle}
        \langle R(v), v \rangle = - \langle v, R(v) \rangle
                                = - \langle R(v), v \rangle
\end{equation}
for every $v \in V$, using the symmetry of the inner product on $V$ in
the second step.  This implies that
\begin{equation}
\label{langle R(v), v rangle = 0}
        \langle R(v), v \rangle = 0
\end{equation}
for every $v \in V$, and hence that any eigenvalue of $R$ must be
equal to $0$, if there is one.  The analogous argument in the complex
case would only give that
\begin{equation}
        \langle R(v), v \rangle
\end{equation}
is purely imaginary for each $v \in V$ when $R$ is anti-self-adjoint,
and hence that the eigenvalues of $R$ are purely imaginary, which also
follow from the representation of $R$ as $i B$ for some self-adjoint
linear operator $B$ on $V$.  As in the case of self-adjoint operators,
if $v, w \in V$ satisfy $R(v) = 0$ and $v \perp w$, then it is easy to
see that $v \perp R(w)$ too.

        If $R$ is an anti-self-adjoint linear operator on $V$, then
\begin{equation}
\label{langle R^2(v), w rangle = ... = langle v, R^2(w) rangle}
        \langle R^2(v), w \rangle = - \langle R(v), R(w) \rangle
                                  = \langle v, R^2(w) \rangle
\end{equation}
for every $v, w \in W$.  This shows that $R^2 = R \circ R$ is a
self-adjoint linear operator on $V$, and in particular that $R^2$ can
be diagonalized in an orthonormal basis in $V$, as in the previous
section.  If we take $v = w$ in (\ref{langle R^2(v), w rangle = ... =
langle v, R^2(w) rangle}), then we get that
\begin{equation}
\label{langle R^2(v), v rangle = - langle R(v), R(v) rangle = - ||R(v)||^2}
        \langle R^2(v), v \rangle = - \langle R(v), R(v) \rangle = - \|R(v)\|^2
\end{equation}
for every $v \in V$.  Of course, $R(v) = 0$ implies that $R^2(v) =
R(R(v)) = 0$ trivially, and (\ref{langle R^2(v), v rangle = - langle
R(v), R(v) rangle = - ||R(v)||^2}) shows that $R^2(v) = 0$ implies
that $R(v) = 0$ in this case.

        If $T$ is any linear operator on $V$, then $T$ can be expressed as
\begin{equation}
\label{T = A + R}
        T = A + R,
\end{equation}
where $A = (T + T^*)/2$ is self-adjoint, and $R = (T - T^*)/2$ is
anti-self-adjoint.  If $T$ commutes with $T^*$, then $A$ commutes with
$R$, and hence $A$ commutes with $R^2$.  As in the previous section,
one can show that there is an orthonormal basis of $V$ in which $A$
and $R^2$ are simultaneously diagonalized under these conditions.  One
also gets that the eigenspaces of $A$ are invariant under $R$, because
$A$ commutes with $R$.  In particular, these remarks can be applied to
the case of an orthogonal linear transformation $T$ on a real inner
product space $V$.

\section{The $C^*$-identity}
\label{section on C^*-identity}

	Any linear operator $T$ on a real or complex inner product
space $V$ satisfies
\begin{equation}
	\|T^* \, T\|_{op} \le \|T^*\|_{op} \, \|T\|_{op} = \|T\|_{op}^2,
\end{equation}
using (\ref{||T^*||_{op} = ||T||_{op}}) in the second step.
Conversely,
\begin{equation}
	\|T(v)\|^2 = \langle T(v), T(v) \rangle
		= \langle T^*(T(v)), v \rangle
\end{equation}
for every $v \in V$, and hence
\begin{equation}
	\|T(v)\|^2 \le \|(T^* \, T)(v)\| \, \|v\|
			\le \|T^* \, T\|_{op} \, \|v\|^2.
\end{equation}
Thus $\|T\|_{op}^2 \le \|T^* \, T\|_{op}$, which implies that
\begin{equation}
\label{C^*-identity}
	\|T^* \, T\|_{op} = \|T\|_{op}^2.
\end{equation}
This is known as the
\emph{$C^*$-identity}\index{C^*-identity@$C^*$-identity}.

	Let $(V_1, \langle \cdot, \cdot \rangle_1)$, $(V_2, \langle
\cdot, \cdot \rangle_2)$ be inner product spaces which are both real
or both complex, and with norms $\|\cdot\|_1$, $\|\cdot\|_2$
associated to their inner products, respectively.  Let $\|T\|_{op,ab}$
be the operator norm of a linear mapping $T : V_a \to V_b$ using
$\|\cdot\|_a$ on the domain and $\|\cdot\|_b$ on the range, where $a,
b = 1, 2$.  This can be characterized in terms of inner products by
\begin{equation}
\label{||T||_{op,ab} in terms of inner products}
	\|T\|_{op,ab}  =  \sup \{ |\langle T(v), w \rangle_b| : 
			     v \in V_a, \, w \in V_b, 
				\ \|v\|_a, \|w\|_b \le 1\},
\end{equation}
as in (\ref{formula for ||T||_{op} in terms of the inner product}).

        If $T$ is a linear mapping from $V_1$ into $V_2$, then there
is a unique linear mapping $T^* : V_2 \to V_1$ such that
\begin{equation}
	\langle T(v), w \rangle_2 = \langle v, T^*(w) \rangle_1
\end{equation}
for every $v \in V_1$ and $w \in V_2$, again called the \emph{adjoint}
of $T$.\index{adjoint operators}  As before, 
\begin{equation}
        \mu_w(v) = \langle T(v), w \rangle_2
\end{equation}
defines a linear functional on $V_1$ for each $w \in V_2$, which can
be represented as
\begin{equation}
        \mu_w(v) = \langle v, T^*(w) \rangle_1
\end{equation}
for a unique element $T^*(w)$ of $V_1$, and one can check that $T^* :
V_2 \to V_1$ is linear.  Otherwise, one can get $T^*$ using
orthonormal bases for $V_1$ and $V_2$, with respect to which the
matrix for $T^*$ is equal to the transpose or the complex conjugate of
the transpose of the corresponding matrix for $T$, depending on
whether $V_1$, $V_2$ are real or complex vector spaces.  As in Section
\ref{adjoints}, the adjoint of $T : V_1 \to V_2$ is very similar to
the dual linear mapping discussed in Section \ref{linear
transformations}, but there are some differences, especially when
$V_1$ and $V_2$ are complex.

	If $T : V_1 \to V_2$ is a linear mapping and if $a$ is a real
or complex number, as appropriate, then
\begin{equation}
\label{(a T)^* = a T^*}
        (a \, T)^* = a \, T^*
\end{equation}
in the real case, and
\begin{equation}
\label{(a T)^* = overline{a} T^*}
        (a \, T)^* = \overline{a} \, T^*
\end{equation}
in the complex case.  If $S, T : V_1 \to V_2$ are linear mappings, then
\begin{equation}
	(S + T)^* = S^* + T^*.
\end{equation}
It is easy to see that $(T^*)^* = T$ for every $T : V_1 \to V_2$.  If
$V_1$, $V_2$, and $V_3$ are inner product spaces, all real or all
complex, and if $T_1 : V_1 \to V_2$ and $T_2 : V_2 \to V_3$ are linear
mappings, then
\begin{equation}
	(T_2 \circ T_1)^* = T_1^* \circ T_2^*
\end{equation}
as linear mappings from $V_3$ into $V_1$.  A linear mapping $T : V_1
\to V_2$ is invertible if and only if $T^* : V_2 \to V_1$ is
invertible, in which case
\begin{equation}
\label{(T^{-1})^* = (T^*)^{-1}}
        (T^{-1})^* = (T^*)^{-1}.
\end{equation}
Using (\ref{||T||_{op,ab} in terms of inner products}), we get that
\begin{equation}
\label{||T^*||_{op,21} = ||T||_{op,12}}
	\|T^*\|_{op,21} = \|T\|_{op,12},
\end{equation}
for any linear mapping $T : V_1 \to V_2$, as in
(\ref{||T^*||_{op} = ||T||_{op}}).  The
$C^*$-identities\index{C^*-identity@$C^*$-identity}
\begin{equation}
\label{||T^* T||_{op,11} = ||T T^*||_{op,22} = ||T||_{op,12}^2}
        \|T^* \, T\|_{op,11} = \|T \, T^*\|_{op,22} = \|T\|_{op,12}^2
\end{equation}
can be verified in this setting as well.

\section{The trace norm}
\label{trace norm}

	Let $V_1$ and $V_2$ be vector spaces, both real or both
complex, equipped with norms $\|\cdot\|_1$ and $\|\cdot\|_2$,
respectively.  Let $V_1^*$, $V_2^*$ and $\|\cdot\|_1^*$,
$\|\cdot\|_2^*$ be the corresponding dual spaces and norms, as in
Section \ref{dual spaces, norms}.

	Any linear mapping $T : V_1 \to V_2$ can be expressed as
\begin{equation}
\label{T(v) = sum_{j=1}^N lambda_j(v) w_j}
	T(v) = \sum_{j=1}^N \lambda_j(v) \, w_j,
\end{equation}
where $N$ is a positive integer, $\lambda_1, \ldots, \lambda_N \in
V_1^*$, and $w_1, \ldots, w_N \in V_2$.  The \emph{trace
norm}\index{trace norm} of $T$ relative to $\|\cdot\|_1$ and
$\|\cdot\|_2$ is defined to be the infimum of
\begin{equation}
\label{sum_{j=1}^N ||lambda_j||_1^* ||w_j||_2}
	\sum_{j=1}^N \|\lambda_j\|_1^* \, \|w_j\|_2
\end{equation}
over all such representations of $T$, and is denoted $\|T\|_{tr}$, or
$\|T\|_{tr, 12}$ to indicate the role of the norms $\|\cdot\|_1$,
$\|\cdot\|_2$.

	If $\lambda \in V_1$, $w \in V_2$, and $A(v) = \lambda(v) \,
w$, then
\begin{equation}
\label{||A||_{op, 12} = ||lambda||_1^* ||w||_2, A = lambda w}
	\|A\|_{op, 12} = \|\lambda\|_1^* \, \|w\|_2,
\end{equation}
where $\|A\|_{op, 12}$ is the operator norm of $A$ with respect to the
norms $\|\cdot \|_1$ and $\|\cdot \|_2$.  Thus
\begin{equation}
	\|T\|_{op, 12} \le \sum_{j=1}^N \|\lambda_j\|_1^* \, \|w_j\|_2
\end{equation}
for each representation (\ref{T(v) = sum_{j=1}^N lambda_j(v) w_j}) of
$T$, and therefore
\begin{equation}
\label{||T||_{op, 12} le ||T||_{tr, 12}}
	\|T\|_{op, 12} \le \|T\|_{tr, 12}.
\end{equation}
Using this, one can check that the trace norm is a norm on the vector
space of linear mappings from $V_1$ to $V_2$.  If $A(v) = \lambda(v)
\, w$, where $\lambda \in V_1^*$ and $w \in V_2$, then
\begin{equation}
	\|\lambda\|_1^* \, \|w\|_2 = \|A\|_{op, 12} \le \|A\|_{tr, 12}
		\le \|\lambda\|_1^* \, \|w\|_2,
\end{equation}
and hence the operator and trace norms of $A$ are the same.

	Suppose that $V_3$ is another vector space which is real or
complex depending on whether $V_1$, $V_2$ are real or complex, and
that $\|\cdot\|_3$ is a norm on $V_3$.  If $T_1 : V_1 \to V_2$, $T_2 :
V_2 \to V_3$ are linear mappings, then
\begin{equation}
	\|T_2 \circ T_1\|_{tr, 13} \le \|T_1\|_{op, 12} \, \|T_2\|_{tr, 23}
\end{equation}
and
\begin{equation}
	\|T_2 \circ T_1\|_{tr, 13} \le \|T_1\|_{tr, 12} \, \|T_2\|_{op, 23}.
\end{equation}
Here $\|\cdot\|_{op, ab}$ and $\|\cdot\|_{tr, ab}$ are the operator
and trace norms for linear mappings from $V_a$ to $V_b$, $a, b = 1, 2,
3$.  This follows by converting representations of the form (\ref{T(v)
= sum_{j=1}^N lambda_j(v) w_j}) for $T_1$ or $T_2$ into similar
representations for $T_1 \circ T_2$.

	Let us briefly review the notion of the trace of a linear
mapping.  Fix a vector space $V$, and suppose that $A$ is a linear
mapping from $V$ to itself.  Let $v_1, \ldots, v_n$ be a basis for
$V$, so that every element of $V$ can be expressed as a linear
combination of the $v_j$'s in exactly one way.  With respect to this
basis, $A$ can be described by an $n \times n$ matrix $(a_{j,k})$ of
real or complex numbers, as appropriate, through the formula
\begin{equation}
	A(v_k) = \sum_{j = 1}^n a_{j,k} \, v_j.
\end{equation}
The \emph{trace}\index{trace} of $A$ is denoted $\tr A$ and defined by
\begin{equation}
	\tr A = \sum_{j = 1}^n a_{j,j}.
\end{equation}
Clearly $\tr A$ is linear in $A$, and one can check that
\begin{equation}
	\tr \, (A \circ B) = \tr \, (B \circ A)
\end{equation}
for any linear transformations $A$ and $B$ on $V$.  In particular,
\begin{equation}
	\tr \, (T \circ A \circ T^{-1}) = \tr A
\end{equation}
for every invertible linear transformation $T$ on $V$.  This implies
that the trace does not depend on the choice of basis for $V$.

	Suppose further that $V$ is equipped with a norm $\|\cdot\|$,
and let $\|\cdot\|_{tr}$ be the corresponding trace norm for operators
on $V$.  If $A$ is any linear transformation on $V$, then
\begin{equation}
\label{|tr A| le ||A||_{tr}}
	|\tr A| \le \|A\|_{tr}.
\end{equation}
To prove this, it suffices to show that 
\begin{equation}
\label{|tr A| le sum_{l = 1}^N ||lambda_l||^* ||w_l||}
	|\tr A| \le \sum_{l = 1}^N \|\lambda_l\|^* \, \|w_l\|
\end{equation}
whenever $\lambda_1, \ldots, \lambda_N \in V^*$, $w_1, \ldots, w_N \in
V$, and $A(v) = \sum_{l = 1}^N \lambda_l(v) \, w_l$.  By linearity, it
is enough to check that
\begin{equation}
\label{|tr A| le ||lambda||^* ||w||}
	|\tr A| \le \|\lambda\|^* \, \|w\|
\end{equation}
when $\lambda \in V^*$, $w \in V$, and $A(v) = \lambda(v) \, w$.  In
this case, $\tr A = \lambda(w)$, and $|\lambda(w)| \le \|\lambda\|^*
\, \|w\|$ by definition of the dual norm.

	Let us return to the setting of two vector spaces $V_1$,
$V_2$, with norms $\|\cdot\|_1$, $\|\cdot\|_2$.  If $T_1 : V_1 \to
V_2$ and $T_2 : V_2 \to V_1$ are linear mappings, then
\begin{equation}
	{\tr}_{V_1} (T_2 \circ T_1) = {\tr}_{V_2} (T_1 \circ T_2).
\end{equation}
Here the trace on the left applies to linear operators on $V_1$, and
the trace on the right applies to linear operators on $V_2$, as
indicated by the notation.  By (\ref{|tr A| le ||A||_{tr}}),
\begin{equation}
\label{|{tr}_{V_1} (T_2 circ T_1)| = |{tr}_{V_2} (T_1 circ T_2)|}
	|{\tr}_{V_1} (T_2 \circ T_1)| = |{\tr}_{V_2} (T_1 \circ T_2)|
\end{equation}
is less than or equal to the trace norm of the composition of $T_1$
and $T_2$ in either order.  Hence it is less than or equal to the
product of the trace norm of $T_1$ and the operator norm of $T_2$, or
the operator norm of $T_1$ times the trace norm of $T_2$.

	Remember that $\mathcal{L}(V_1, V_2)$ denotes the vector space
of linear mappings from $V_1$ into $V_2$.  If $R$ is a linear mapping
from $V_2$ into $V_1$, then
\begin{equation}
\label{T mapsto tr_{V_1} (R circ T)}
	T \mapsto {\tr}_{V_1} (R \circ T)
\end{equation}
is a linear functional $\mathcal{L}(V_1, V_2)$.  This defines a linear
isomorphism from $\mathcal{L}(V_2, V_1)$ onto $\mathcal{L}(V_1,
V_2)^*$.

	Fix a linear mapping $R : V_2 \to V_1$, and let us check that
the dual norm of (\ref{T mapsto tr_{V_1} (R circ T)}) with respect to
the trace norm on $\mathcal{L}(V_1, V_2)$ is equal to $\|R\|_{op,
21}$.  We have already seen that
\begin{equation}
\label{|{tr}_{V_1} (R circ T)| le ||R||_{op, 21} ||T||_{tr, 12}}
	|{\tr}_{V_1} (R \circ T)| \le \|R\|_{op, 21} \, \|T\|_{tr, 12},
\end{equation}
which says exactly that the aforementioned dual norm of (\ref{T mapsto
tr_{V_1} (R circ T)}) is less than or equal to $\|R\|_{op, 21}$.  To
establish the opposite inequality, let $\lambda \in V_1^*$ and $w \in
V_2$ be given, and put $T_0(v) = \lambda(v) \, w$.  Thus $T_0$ is a
linear mapping from $V_1$ to $V_2$,
\begin{equation}
	(R \circ T_0)(v) = \lambda(v) \, R(w),
\end{equation}
and
\begin{equation}
	{\tr}_{V_1} (R \circ T_0) = \lambda(R(w)).
\end{equation}
By definition, $|{\tr}_{V_1} (R \circ T_0)|$ is less than or equal to the
dual norm of (\ref{T mapsto tr_{V_1} (R circ T)}) times the trace norm of
$T_0$.  The trace norm of $T_0$ is equal to $\|\lambda\|_1^* \,
\|w\|_2$, and hence $|\lambda(R(w))|$ is less than or equal to the
dual norm of (\ref{T mapsto tr_{V_1} (R circ T)}) times $\|\lambda\|_1^* \,
\|w\|_2$.  Since $\lambda \in V_1^*$ and $w \in V_2$ are arbitrary,
this implies that $\|R\|_{op, 21}$ is less than or equal to the dual norm
of (\ref{T mapsto tr_{V_1} (R circ T)}), so that the two are the same.

	Similarly, if $T$ is a linear mapping from $V_1$ into $V_2$, then
\begin{equation}
\label{R mapsto {tr}_{V_1} (R circ T)}
	R \mapsto {\tr}_{V_1} (R \circ T)
\end{equation}
is a linear functional on $\mathcal{L}(V_2, V_1)$, and this defines a linear
isomorphism from $\mathcal{L}(V_1, V_2)$ onto $\mathcal{L}(V_2, V_1)^*$.
It follows from the previous discussion that the dual norm of (\ref{R
mapsto {tr}_{V_1} (R circ T)}) with respect to the operator norm on
$\mathcal{L}(V_2, V_1)$ is equal to $\|T\|_{tr, 12}$, by the results
in Section \ref{second duals}.  In other words, we just saw that the
operator norm on $\mathcal{L}(V_2, V_1)$ corresponds to the dual of the
trace norm on $\mathcal{L}(V_1, V_2)$, and this implies that the dual
of the operator norm on $\mathcal{L}(V_2, V_1)$ corresponds to the trace
norm on $\mathcal{L}(V_1, V_2)$, as in Section \ref{second duals}.

\section{The Hilbert--Schmidt norm}
\label{hilbert--schmidt norm}

	Let $(V_1, \langle \cdot, \cdot \rangle_1)$, $(V_2, \langle
\cdot, \cdot \rangle_2)$ be inner product spaces, both real or both
complex, with norms $\|\cdot \|_1$, $\|\cdot \|_2$ associated to their
inner products, as usual.  If $\{a_j\}_{j=1}^p$, $\{b_k\}_{k=1}^q$ are
orthonormal bases for $V_1$, $V_2$, respectively, then
\begin{equation}
	v = \sum_{j=1}^p \langle v, a_j \rangle_1 \, a_j, \qquad
	  w = \sum_{k=1}^q \langle w, b_k \rangle_2 \, b_k
\end{equation}
for every $v \in V_1$, $w \in V_2$, and we can express a linear
mapping $T : V_1 \to V_2$ as
\begin{equation}
	T(v) = \sum_{j=1}^p \sum_{k=1}^q \langle v, a_j \rangle_1 \,
		\langle T(a_j), b_k \rangle_2 \, b_k.
\end{equation}
Because $\{a_j\}_{j=1}^p$ and $\{b_k\}_{k=1}^q$ are orthonormal bases
for $V_1$ and $V_2$, we get that
\begin{eqnarray}
	\sum_{j=1}^p \|T(a_j)\|_2^2 & = &
	\sum_{j=1}^p \sum_{k=1}^q |\langle T(a_j), b_k \rangle_2|^2	\\
	& = & \sum_{j=1}^p \sum_{k=1}^q |\langle a_j, T^*(b_k) \rangle_1|^2
		= \sum_{k=1}^q \|T^*(b_k)\|_1^2.		\nonumber
\end{eqnarray}
The \emph{Hilbert--Schmidt norm}\index{Hilbert--Schmidt norm} of $T$
is denoted $\|T\|_{HS}$ and defined to be the square root of the
common value of these sums.  It follows that the Hilbert--Schmidt
norms of $T$ and $T^*$ are the same, and do not depend on the
particular choices of orthonormal bases for $V_1$ and $V_2$.

	If we express $T(v)$ as
\begin{equation}
	T(v) = \sum_{j=1}^p \langle v, a_j \rangle_1 \, T(a_j),
\end{equation}
and $(T^* \circ T)(v)$ as
\begin{equation}
	(T^* \circ T)(v) = \sum_{j=1}^p \sum_{l=1}^p
 \langle v, a_j \rangle_1 \, \langle T(a_j), T(a_l) \rangle_1 \, a_l,
\end{equation}
then we see that
\begin{equation}
	{\tr}_{V_1} (T^* \circ T) = \sum_{j=1}^p \|T(a_j)\|_1^2 = \|T\|_{HS}^2.
\end{equation}
Here the left side is the trace of $T^* \circ T$ as an operator on
$V_1$, and similarly the trace of $T \circ T^*$ on $V_2$ is equal to
$\|T\|_{HS}^2$.  One can check that
\begin{equation}
\label{langle A, B rangle_{mathcal{L}(V_1, V_2)} = ...}
	\langle A, B \rangle_{\mathcal{L}(V_1, V_2)}
		= {\tr}_{V_1} (B^* \circ A) = {\tr}_{V_2} (A \circ B^*)
\end{equation}
defines an inner product on the vector space $\mathcal{L}(V_1, V_2)$
of linear transformations from $V_1$ into $V_2$, so that the
Hilbert--Schmidt norm is exactly the norm on $\mathcal{L}(V_1, V_2)$
associated to this inner product.  More precisely, if $R$ is a linear
transformation on a real or complex inner product space $V$, then
\begin{equation}
        {\tr}_V R^* = {\tr}_V R
\end{equation}
in the real case, and
\begin{equation}
        {\tr}_V R^* = \overline{{\tr}_V R}
\end{equation}
in the complex case, as one can verify using the description of the
matrix of $R^*$ with respect to an orthonormal basis for $V$ mentioned
in Section \ref{adjoints}.  This implies that (\ref{langle A, B
rangle_{mathcal{L}(V_1, V_2)} = ...}) satisfies the symmetry property
required to be an inner product.

\section{Schmidt decompositions}
\label{schmidt decompositions}

        Let $(V_1, \langle \cdot, \cdot \rangle_1)$, $(V_2, \langle
\cdot, \cdot \rangle_2)$ be inner product spaces again, both real or
both complex, with norms $\|\cdot \|_1$, $\|\cdot \|_2$ associated to
their inner products.  A \emph{Schmidt decomposition}\index{Schmidt
decompositions} for a linear mapping $T : V_1 \to V_2$ is a
representation of $T$ as
\begin{equation}
\label{schmidt decomposition}
	T(v) = \sum_{j=1}^r \lambda_j \, \langle v, u_j \rangle_1 \, w_j,
\end{equation}
where $r$ is a positive integer, $u_1, \ldots, u_r$ and $w_1, \ldots,
w_r$ are orthonormal vectors in $V_1$ and $V_2$, respectively, and
$\lambda_1, \ldots, \lambda_r$ are scalars.  The existence of a
Schmidt decomposition uses the fact that $T^* \circ T$ is self-adjoint
on $V_1$, and hence can be diagonalized in an orthonormal basis.
It also uses the observation that
\begin{equation}
\label{T(y) perp T(z)}
        T(y) \perp T(z)
\end{equation}
in $V_2$ when $y, z \in V_1$, $y \perp z$, and $y$ is an eigenvector
for $T^* \circ T$.

        If $T$ has a Schmidt decomposition (\ref{schmidt
decomposition}), then it is easy to see that
\begin{equation}
	\|T\|_{op} = \max (|\lambda_1|, \ldots, |\lambda_r|),
\end{equation}
and
\begin{equation}
	\|T\|_{HS} = \Big(\sum_{j=1}^r |\lambda_j|^2 \Big)^{1/2}.
\end{equation}
Let us check that
\begin{equation}
\label{||T||_{tr} = sum_{j=1}^r |lambda_j|}
	\|T\|_{tr} = \sum_{j=1}^r |\lambda_j|.
\end{equation}
Clearly
\begin{equation}
\label{||T||_{tr} le sum_{j=1}^r |lambda_j|}
        \|T\|_{tr} \le \sum_{j=1}^r |\lambda_j|,
\end{equation}
by the definition of the trace norm, and
\begin{equation}
\label{sum_{j=1}^r |lambda_j| = sum_{j=1}^r |langle T(u_j), w_j rangle_2|}
	\sum_{j=1}^r |\lambda_j| =
		\sum_{j=1}^r |\langle T(u_j), w_j \rangle_2|.
\end{equation}

        Let $y_1, \ldots, y_k \in V_1$ and $z_1, \ldots, z_k \in V_2$
be arbitrary orthonormal collections of vectors, and let us check that
\begin{equation}
\label{sum_{l=1}^k |langle R(y_l), z_l rangle_2| le ||R||_{tr}}
	\sum_{l=1}^k |\langle R(y_l), z_l \rangle_2| \le \|R\|_{tr}
\end{equation}
for every linear mapping $R : V_1 \to V_2$.  If $R$ is of the form
$R(v) = \langle v, a \rangle_1 \, b$ for some $a \in V_1$ and $b \in
V_2$, then
\begin{equation}
	\sum_{l=1}^k |\langle R(y_l), z_l \rangle_2| = 
 \sum_{l=1}^k |\langle y_l, a \rangle_1| \, |\langle b, z_l \rangle_2|
		\le \|a\|_1 \, \|b\|_2,
\end{equation}
using the Cauchy--Schwarz inequality in the second step.  Thus the
left side of (\ref{sum_{l=1}^k |langle R(y_l), z_l rangle_2| le
||R||_{tr}}) is less than or equal to the operator norm of $R$ when
$R$ has rank $1$, which implies (\ref{sum_{l=1}^k |langle R(y_l), z_l
rangle_2| le ||R||_{tr}}) in general.  If $T$ has Schmidt
decomposition (\ref{schmidt decomposition}), then we can apply
(\ref{sum_{l=1}^k |langle R(y_l), z_l rangle_2| le ||R||_{tr}}) with
$R = T$, $y_l = u_l$, and $z_l = w_l$, to get that the right side of
(\ref{sum_{j=1}^r |lambda_j| = sum_{j=1}^r |langle T(u_j), w_j
rangle_2|}) is less than or equal to the trace norm of $T$, as
desired.

\section{$\mathcal{S}_p$ norms}
\label{mathcal{S}_p norms}

	Let $(V_1, \langle \cdot, \cdot \rangle_1)$, $(V_2, \langle
\cdot, \cdot \rangle_2)$ be inner product spaces again, both real or
both complex, with norms $\|\cdot \|_1$, $\|\cdot \|_2$ associated to
their inner products, and let $p$ be a real number, $1 < p < \infty$.
If $T$ is a linear mapping from $V_1$ to $V_2$ with Schmidt
decomposition (\ref{schmidt decomposition}), and $y_1, \ldots, y_k \in
V_1$, $z_1, \ldots, z_k \in V_2$ are orthonormal collections of
vectors, then
\begin{equation}
\label{... le (sum_{j = 1}^r |lambda_j|^p)^{1/p}}
	\Big(\sum_{h = 1}^k |\langle T(y_h), z_h \rangle_2|^p \Big)^{1/p}
 		\le \Big(\sum_{j = 1}^r |\lambda_j|^p \Big)^{1/p}.
\end{equation}
Equivalently,
\begin{equation}
\label{... le (sum_{j = 1}^r |lambda_j|^p)^{1/p}, 2}
	\Big(\sum_{h = 1}^k \biggl|\sum_{l = 1}^r \lambda_l \,
 \langle y_h, u_l \rangle_1 \, \langle w_l, z_h \rangle_2 \biggr|^p \Big)^{1/p}
		\le \Big(\sum_{j = 1}^r |\lambda_j|^p \Big)^{1/p}.
\end{equation}
To see this, observe that
\begin{equation}
\label{(sum_{l = 1}^r |langle y_h, u_l rangle_1|^2)^{1/2} le ||y_h||_1 = 1}
	\Big(\sum_{l = 1}^r |\langle y_h, u_l \rangle_1|^2 \Big)^{1/2} 
		\le \|y_h\|_1 = 1
\end{equation}
and
\begin{equation}
\label{(sum_{l = 1}^r |langle w_l, z_h rangle_2|^2)^{1/2} le ||z_h||_2 = 1}
	\Big(\sum_{l = 1}^r |\langle w_l, z_h \rangle_2|^2 \Big)^{1/2} 
		\le \|z_h\|_2 = 1
\end{equation}
for each $h$, and that
\begin{equation}
\label{(sum_{h = 1}^k |langle y_h, u_l rangle_1|^2)^{1/2} le ||u_l||_1 = 1}
	\Big(\sum_{h = 1}^k |\langle y_h, u_l \rangle_1|^2 \Big)^{1/2} 
		\le \|u_l\|_1 = 1,
\end{equation}
and
\begin{equation}
\label{(sum_{h = 1}^k |langle w_l, z_h rangle_2|^2)^{1/2} le ||w_l||_2 = 1}
	\Big(\sum_{h = 1}^k |\langle w_l, z_h \rangle_2|^2 \Big)^{1/2} 
		\le \|w_l\|_2 = 1
\end{equation}
for each $l$.  Hence
\begin{equation}
\label{sum_{l = 1}^r |langle y_h, u_l rangle_1 langle w_l, z_h rangle_2| le 1}
  \sum_{l = 1}^r |\langle y_h, u_l \rangle_1 \, \langle w_l, z_h \rangle_2|
		\le 1
\end{equation}
for each $h$, and
\begin{equation}
\label{sum_{h = 1}^k |langle y_h, u_l rangle_1 langle w_l, z_h rangle_2| le 1}
  \sum_{h = 1}^k |\langle y_h, u_l \rangle_1 \, \langle w_l, z_h \rangle_2|
		\le 1
\end{equation}
for each $l$, by the Cauchy--Schwarz inequality.  The desired estimate
(\ref{... le (sum_{j = 1}^r |lambda_j|^p)^{1/p}, 2}) can now be
derived from Schur's theorem in Section \ref{special cases}.

	The $\mathcal{S}_p$ norm\index{Sp norm@$\mathcal{S}_p$ norm}
of $T$ is denoted $\|T\|_{\mathcal{S}_p}$ and defined by
\begin{equation}
	\|T\|_{\mathcal{S}_p} = \Big(\sum_{j = 1}^r |\lambda_j|^p \Big)^{1/p}
\end{equation}
when $1 \le p < \infty$, and
\begin{equation}
	\|T\|_{\mathcal{S}_\infty} = \max (|\lambda_1|, \ldots, |\lambda_r|).
\end{equation}
This is equal to the trace norm of $T$ when $p = 1$, the
Hilbert--Schmidt norm of $T$ when $p = 2$, and the operator norm of
$T$ when $p = \infty$.  Equivalently,
\begin{equation}
\label{||T||_{S_p} = sup (sum_{h = 1}^k |langle T(y_h), z_h rangle_2|^p)^{1/p}}
	\|T\|_{\mathcal{S}_p} = 
	\sup \Big(\sum_{h = 1}^k |\langle T(y_h), z_h \rangle_2|^p \Big)^{1/p}
\end{equation}
when $1 \le p < \infty$, where the supremum is taken over all
orthonormal collections of vectors $y_1, \ldots, y_k$ and $z_1,
\ldots, z_k$ in $V_1$ and $V_2$, since (\ref{... le (sum_{j = 1}^r
|lambda_j|^p)^{1/p}}) shows that the supremum is attained by the
Schmidt decomposition.  Using (\ref{||T||_{S_p} = sup (sum_{h = 1}^k
|langle T(y_h), z_h rangle_2|^p)^{1/p}}), one can check that the
$\mathcal{S}_p$ norm satisfies the triangle inequality, and hence is a
norm on the vector space $\mathcal{L}(V_1, V_2)$ of linear mappings
from $V_1$ into $V_2$.

        Similarly, if $p \ge 2$ and $y_1, \ldots, y_k \in V_1$ are
orthonormal, then
\begin{equation}
\label{(sum_{h = 1}^k ||T(y_h)||_2^p )^{1/p} le ||T|_{mathcal{S}_p}}
	\Big(\sum_{h = 1}^k \|T(y_h)\|_2^p \Big)^{1/p}
		\le \|T\|_{\mathcal{S}_p}
\end{equation}
for every linear mapping $T : V_1 \to V_2$.  To see this, let
(\ref{schmidt decomposition}) be a Schmidt decomposition for $T$, and
observe that
\begin{equation}
	\|T(y_h)\|_2 = 
   \Big(\sum_{l = 1}^r |\lambda_l|^2 \, 
		|\langle y_h, u_l \rangle_1|^2 \Big)^{1/2}
\end{equation}
for each $h$.  If $q = p/2$ and $\mu_j = |\lambda_j|^2$, then
(\ref{(sum_{h = 1}^k ||T(y_h)||_2^p )^{1/p} le ||T|_{mathcal{S}_p}}) can
be re-expressed as
\begin{equation}
	\Big(\sum_{h = 1}^k \Big(\sum_{l = 1}^r \mu_l \, 
		|\langle y_h, u_l \rangle_1|^2 \Big)^{q} \Big)^{1/q}
		\le \Big(\sum_{j = 1}^r \mu_j^q \Big)^{1/q}.
\end{equation}
The orthonormality of $y_1, \ldots, y_k$ and $u_1, \ldots, u_r$ in
$V_1$ imply that (\ref{(sum_{l = 1}^r |langle y_h, u_l
rangle_1|^2)^{1/2} le ||y_h||_1 = 1}) holds for each $h$ and that
(\ref{(sum_{h = 1}^k |langle y_h, u_l rangle_1|^2)^{1/2} le ||u_l||_1
= 1}) holds for each $l$, as before.  The desired estimate again
follows from Schur's theorem in Section \ref{special cases}.

\section{Duality}
\label{duality}

        Let $(V_1, \langle \cdot, \cdot \rangle_1)$, $(V_2, \langle
\cdot, \cdot \rangle_2)$ be inner product spaces again, both real or
both complex, and with norms $\|\cdot \|_1$, $\|\cdot \|_2$ associated
to their inner products.  Let $T$ be a linear mapping from $V_1$ into
$V_2$, and let $R$ be a linear mapping from $V_2$ into $V_1$.  Suppose
that $T$ has Schmidt decomposition (\ref{schmidt decomposition}), so that
\begin{equation}
\label{(R circ T)(v) = sum_{j = 1}^r lambda_j langle v, u_j rangle_1 R(w_j)}
 (R \circ T)(v) = \sum_{j = 1}^r \lambda_j \langle v, u_j \rangle_1 \, R(w_j)
\end{equation}
for each $v \in V_1$.  Thus
\begin{equation}
 {\tr}_{V_1} (R \circ T) = \sum_{j = 1}^r \lambda_j
                                            \langle R(w_j), u_j \rangle_1.
\end{equation}
If $1 < p, q < \infty$ are conjugate exponents, then we get that
\begin{equation}
\label{|{tr}_{V_1} (R circ T)| le ...}
 |{\tr}_{V_1} (R \circ T)| \le \Big(\sum_{j = 1}^j |\lambda_j|^p \Big)^{1/p}
  \, \Big(\sum_{j = 1}^r |\langle R(w_j), u_j \rangle_1|^q \Big)^{1/q},    
\end{equation}
by H\"older's inequality.  This implies that
\begin{equation}
\label{|{tr}_{V_1} (R circ T)| le ||T||_{mathcal{S}_p} ||R||_{mathcal{S}_q}}
 |{\tr}_{V_1} (R \circ T)| \le \|T\|_{\mathcal{S}_p} \, \|R\|_{\mathcal{S}_q},
\end{equation}
using the definition of the $\mathcal{S}_p$ norm of $T$ in terms the
Schmidt decomposition, and the analogue of (\ref{||T||_{S_p} = sup
(sum_{h = 1}^k |langle T(y_h), z_h rangle_2|^p)^{1/p}}) for $R$ and
$q$.  This also works when $p = \infty$ or $q = \infty$, by the same
argument.

        The preceding inequality implies that
\begin{equation}
        T \mapsto {\tr}_{V_1} (R \circ T)
\end{equation}
has dual norm less than or equal to $\|R\|_{\mathcal{S}_q}$ with
respect to the norm $\|T\|_{\mathcal{S}_p}$ on $\mathcal{L}(V_1, V_2)$
when $1 \le p, q \le \infty$ are conjugate exponents.  To show that
the dual norm is equal to $\|R\|_{\mathcal{S}_q}$, it suffices to check that
\begin{equation}
\label{{tr}_{V_1} (R circ T) = ||T||_{mathcal{S}_p} ||R||_{mathcal{S}_q}}
 {\tr}_{V_1} (R \circ T) = \|T\|_{\mathcal{S}_p} \, \|R\|_{\mathcal{S}_q}
\end{equation}
for some $T \ne 0$.  Let us begin this time with a Schmidt decomposition
\begin{equation}
        R(w) = \sum_{j = 1}^r \mu_j \, \langle w, w_j \rangle_2 \, u_j
\end{equation}
for $R$, where $u_1, \ldots, u_r$ and $w_1, \ldots, w_r$ be
orthonormal vectors in $V_1$ and $V_2$, respectively, and $\mu_1,
\ldots, \mu_k$ are real or complex numbers, as appropriate.  Let us
also restrict our attention now to linear mappings $T$ from $V_1$ into
$V_2$ with Schmidt decomposition (\ref{schmidt decomposition}), using
the same orthonormal vectors $u_1, \ldots, u_r$ and $w_1, \ldots, w_r$
as for $R$.  In this case, the trace of $R \circ T$ reduces to
\begin{equation}
\label{{tr}_{V_1} (R circ T) = sum_{h = 1}^k mu_h lambda_h}
        {\tr}_{V_1} (R \circ T) = \sum_{h = 1}^k \mu_h \, \lambda_h,
\end{equation}
and for any $\mu_1, \ldots, \mu_k$ one can choose $\lambda_1, \ldots,
\lambda_k$, not all equal to $0$, such that (\ref{{tr}_{V_1} (R circ
T) = ||T||_{mathcal{S}_p} ||R||_{mathcal{S}_q}}) holds, as in Section
\ref{dual spaces, norms}.

\chapter{Seminorms and sublinear functions}
\label{seminorms, sublinear functions}

        As in the previous two chapters, we continue to restrict our
attention to finite-dimensional vector spaces in this chapter.

\section{Seminorms}
\label{seminorms}

        A \emph{seminorm}\index{seminorms} on a real or complex vector
space $V$ is a nonnegative real-valued function $N$ on $V$ such that
\begin{equation}
\label{N(t v) = |t| N(v)}
        N(t \, v) = |t| \, N(v)
\end{equation}
for every $v \in V$ and $t \in {\bf R}$ or ${\bf C}$, as appropriate,
and
\begin{equation}
\label{N(v + w) le N(v) + N(w)}
        N(v + w) \le N(v) + N(w)
\end{equation}
for every $v, w \in V$.  Note that $N(0) = 0$, as one can see by
applying (\ref{N(t v) = |t| N(v)}) with $t = 0$.  Thus a seminorm $N$
on $V$ is a norm when $N(v) > 0$ for every $v \in V$ with $v \ne 0$.
If $\lambda$ is a linear functional on $V$, then
\begin{equation}
\label{N_lambda(v) = |lambda(v)|}
        N_\lambda(v) = |\lambda(v)|
\end{equation}
is a seminorm on $V$, and the sum and maximum of finitely many
seminorms on $V$ are also seminorms on $V$.

        If $N$ is a seminorm on $V$ and $v, w \in V$, then
\begin{equation}
        N(v) - N(w) \le N(v - w)
\end{equation}
and
\begin{equation}
        N(w) - N(v) \le N(w - v) = N(v - w),
\end{equation}
by the triangle inequality.  Hence
\begin{equation}
\label{|N(v) - N(w)| le N(v - w)}
        |N(v) - N(w)| \le N(v - w).
\end{equation}
Suppose for the moment that $V = {\bf R}^n$ or ${\bf C}^n$ for some
positive integer $n$, and let $|v|$ be the standard Euclidean norm on $V$.
It is easy to see that there is a nonnegative real number $C$ such that
\begin{equation}
        N(v) \le C \, |v|
\end{equation}
for every $v \in V$, by expressing $v$ as a linear combination of the
standard basis vectors for $V$ and using (\ref{N(t v) = |t| N(v)}) and
(\ref{N(v + w) le N(v) + N(w)}).  This together with (\ref{|N(v) -
N(w)| le N(v - w)}) implies that $N$ is a continuous function on $V$
with respect to the standard Euclidean metric and topology.

        Suppose that $V$ is a real vector space, $N$ is a seminorm on
$V$, and $\lambda$ is a linear functional on $V$ such that
\begin{equation}
\label{lambda(v) le N(v)}
        \lambda(v) \le N(v)
\end{equation}
for every $v \in V$.  This implies that
\begin{equation}
        -\lambda(v) = \lambda(-v) \le N(-v) = N(v)
\end{equation}
for every $v \in V$, and hence that
\begin{equation}
\label{|lambda(v)| le N(v)}
        |\lambda(v)| \le N(v).
\end{equation}
Similarly, if $V$ is a complex vector space, $N$ is a seminorm on $V$,
and $\lambda$ is a linear functional on $V$ such that
\begin{equation}
\label{re lambda(v) le N(v)}
        \re \lambda(v) \le N(v)
\end{equation}
for every $v \in V$, then we get that
\begin{equation}
        \re t \, \lambda(v) = \re \lambda(t \, v) \le N(t \, v) = N(v)
\end{equation}
for every $v \in V$ and $t \in {\bf C}$ with $|t| = 1$.  This implies
again that (\ref{|lambda(v)| le N(v)}) holds for every $v \in V$.

\section{Sublinear functions}
\label{sublinear functions}

	A \emph{sublinear function}\index{sublinear functions} on a
real or complex vector space $V$ is a real-valued function $p(v)$ on
$V$ such that
\begin{equation}
\label{p(t v) = t p(v)}
	p(t \, v) = t \, p(v)
\end{equation}
for every $v \in V$ and nonnegative real number $t$, and
\begin{equation}
\label{p(v + w) le p(v) + p(w)}
	p(v + w) \le p(v) + p(w)
\end{equation}
for every $v, w \in V$.  In particular, $p(0) = 0$, as one can see by
applying (\ref{p(t v) = t p(v)}) with $t = 0$.  Note that seminorms
are sublinear functions, but sublinear functions are not required to
be nonnegative.  Linear functionals on real vector spaces are
sublinear functions, as are the real parts of linear functionals on
complex vector spaces.  The sum and maximum of finitely many sublinear
functions are sublinear functions, which includes the maximum of a
sublinear function and $0$, to get a nonnegative sublinear function.

        Let $p$ be a sublinear function on $V$, and observe that
\begin{equation}
        0 = p(0) \le p(v) + p(-v)
\end{equation}
for every $v \in V$.  If
\begin{equation}
        p(-v) = p(v)
\end{equation}
for every $v \in V$, then it follows that
\begin{equation}
\label{p(v) ge 0}
        p(v) \ge 0
\end{equation}
for every $v \in V$, and that $p$ is a seminorm on $V$ when $V$ is a
real vector space.  Similarly, if $V$ is a complex vector space, and
\begin{equation}
        p(t \, v) = p(v)
\end{equation}
for every $t \in {\bf C}$ with $|t| = 1$, then $p$ is a seminorm on $V$.

        If $p$ is a sublinear function on $V$, then
\begin{equation}
        p(v) - p(w) \le p(v - w)
\end{equation}
and
\begin{equation}
        p(w) - p(v) \le p(w - v)
\end{equation}
for every $v, w \in V$, and hence
\begin{equation}
        |p(v) - p(w)| \le \max(p(v - w), p(w - v))
\end{equation}
for every $v, w \in V$.  Note that
\begin{equation}
\label{N(v) = max (p(v), p(-v))}
        N(v) = \max (p(v), p(-v))
\end{equation}
is a seminorm on $V$ when $V$ is a real vector space, by the remarks
in the previous paragraphs.  Of course, a complex vector space may
also be considered to be a real vector space, by forgetting about
mulitplication by $i$.  If $V = {\bf R}^n$ or ${\bf C}^n$ for some
positive integer $n$, then $N(v)$ is bounded by a constant multiple of
the standard Euclidean norm of $v$, as in the previous section.  This
implies that $p(v)$ is continuous with respect to the standard
Euclidean metric and topology on $V$, as before.

\section{Another extension theorem}
\label{another extension theorem}

\begintheorem
\label{extension theorem, 2}
Let $V$ be a real vector space, and let $p$ be a sublinear function on
$V$.  If $W$ is a linear subspace of $V$ and $\mu$ is a linear
functional on $W$ such that
\begin{equation}
\label{mu(w) le p(w) for every w in W}
	\mu(w) \le p(w)  \quad\hbox{for every } w \in W,
\end{equation}
then there is a linear functional $\widehat{\mu}$ on $V$ which is
equal to $\mu$ on $W$ and satisfies
\begin{equation}
\label{widehat{mu}(v) le p(v) for every v in V}
	\widehat{\mu}(v) \le p(v) \quad\hbox{for every } v \in V.
\end{equation}
\end{theorem}

        This is the analogue of Theorem \ref{extension theorem} in
Section \ref{second duals} for sublinear functions instead of norms.
Note that the statement and proof of Theorem \ref{extension theorem}
already work in exactly the same way for seminorms instead of norms.
The case of sublinear functions is essentially the same, except for a
few simple changes following the differences in the statements of the
theorems.

        As in (\ref{|mu_j(x) + t alpha| le ||x + t z||, x in W_j, t in
{bf R}}), we would like to show that there is an $\alpha \in {\bf R}$
such that
\begin{equation}
\label{mu_j(x) + t alpha le p(x + t z), x in W_j, t in R}
	\mu_j(x) + t \, \alpha \le p(x + t z)
		\quad\hbox{for every } x \in W_j
		\hbox{ and } t \in {\bf R}.
\end{equation}
This is equivalent to
\begin{equation}
\label{mu_j(x) + alpha le p(x + z), mu_j(x) - alpha le p(x - z), x in W_j}
	\mu_j(x) + \alpha \le p(x + z), \ \mu_j(x) - \alpha \le p(x - z)
		\quad\hbox{for every } x \in W_j,
\end{equation}
because one can convert (\ref{mu_j(x) + t alpha le p(x + t z), x in
W_j, t in R}) into (\ref{mu_j(x) + alpha le p(x + z), mu_j(x) - alpha
le p(x - z), x in W_j}) when $t \ne 0$ using homogeneity, and the $t =
0$ case of (\ref{mu_j(x) + t alpha le p(x + t z), x in W_j, t in R})
follows from the induction hypothesis for $\mu_j$ and $W_j$.  Let us
rewrite (\ref{mu_j(x) + alpha le p(x + z), mu_j(x) - alpha le p(x -
z), x in W_j}) as
\begin{equation}
\label{mu_j(x) - p(x - z) le alpha le p(x + z) - mu_j(x), x in W_j}
	\mu_j(x) - p(x - z) \le \alpha \le p(x + z) - \mu_j(x)
		\hbox{ for every } x \in W_j.
\end{equation}

	To show that there is an $\alpha \in {\bf R}$ that satisfies
(\ref{mu_j(x) - p(x - z) le alpha le p(x + z) - mu_j(x), x in W_j}),
it suffices to verify
\begin{equation}
	\mu_j(x) - p(x - z) \le p(y + z) - \mu_j(y)
		\quad\hbox{for every } x, y \in W_j,
\end{equation}
which reduces to
\begin{equation}
	\mu_j(x + y) \le p(x - z) + p(y + z) 
		\quad\hbox{for every } x, y \in W_j.
\end{equation}
The subadditivity property of $p(v)$ implies that this condition holds
if
\begin{equation}
	\mu_j(x + y) \le p(x + y) 
		\quad\hbox{for every } x, y \in W_j,
\end{equation}
which is the same as
\begin{equation}
	\mu_j(w) \le p(w) \quad\hbox{for every } w \in W_j.
\end{equation}
This holds by induction hypothesis, which completes the proof of
Theorem \ref{extension theorem, 2}.

\section{Minkowski functionals}
\label{minkowski functionals}

        Let $V$ be a real or complex vector space, and let $A$ be a
subset of $V$ such that $0 \in A$.  Suppose also that $A$ has the
absorbing property that for each $v \in V$ there is a positive real
number $t$ such that
\begin{equation}
        t \, v \in A.
\end{equation}
If $V = {\bf R}^n$ or ${\bf C}^n$ for some positive integer $n$, and
if $0$ is an element of the interior of $A$ with respect to the
standard Euclidean metric and topology on $V$, then $A$ obviously has
this property.

        Under these conditions, the \emph{Minkowski functional}\index{Minkowski
functionals} on $V$ associated to $A$ is defined by
\begin{equation}
\label{N_A(v) = inf {r > 0 : r^{-1} v in A}}
        N_A(v) = \inf \{r > 0 : r^{-1} \, v \in A\}
\end{equation}
for each $v \in V$.  Note that $N(v) \ge 0$ for every $v \in V$, $N(0)
= 0$, and
\begin{equation}
\label{N_A(t v) = t N(v)}
        N_A(t \, v) = t \, N(v)
\end{equation}
for every nonnegative real number $t$.  By construction,
\begin{equation}
        N_A(v) \le 1
\end{equation}
for every $v \in A$.  If $V = {\bf R}^n$ or ${\bf C}^n$ for some
positive integer $n$, and if $A$ is an open set in $V$ with respect to
the standard Euclidean metric and topology, then
\begin{equation}
        N_A(v) < 1
\end{equation}
for every $v \in A$.  This is because $r^{-1} \, v \in A$ for all $r$
sufficiently close to $1$ when $v \in A$ and $A$ is an open set.

        Let $-A$ be the set of vectors in $V$ of the form $-v$ with $v
\in A$, and let $t \, A$ be the set of vectors in $V$ of the form $t
\, v$ with $v \in A$ for each real or complex number $t$, as
appropriate.  Thus
\begin{equation}
        N_A(v) = \inf \{r > 0 : v \in r \, A\}
\end{equation}
for each $v \in V$.  If $-A = A$, then
\begin{equation}
        N_A(-v) = N_A(v)
\end{equation}
for every $v \in V$.  Similarly, if $V$ is a complex vector space, and
if $t \, A = A$ for every complex number $t$ with $|t| = 1$, then
\begin{equation}
\label{N_A(t v) = N_A(v)}
        N_A(t \, v) = N_A(v)
\end{equation}
for every $v \in V$ and $t \in {\bf C}$ with $|t| = 1$.

        Let us suppose from now on in this section that $A$ is
star-like about $0$, which means that $A$ contains every line segment
in $V$ between $0$ and any other element of $A$.  Equivalently,
$A$ is star-like about $0$ if
\begin{equation}
\label{t A subseteq A}
        t \, A \subseteq A
\end{equation}
for every nonnegative real number $t$ with $t \le 1$.  If $v \in V$
satisfies $N_A(v) < 1$, then there is a positive real number $r < 1$
such that $v \in r \, A$, and it follows that $v \in A$.  If $V = {\bf
R}^n$ or ${\bf C}^n$ for some $n$ and $A$ is an open set in $V$ with
respect to the standard metric and topology, then we get that
\begin{equation}
\label{A = {v in V : N_A(v) < 1}}
        A = \{v \in V : N_A(v) < 1\}.
\end{equation}

        Similarly, if $A$ is a closed set in $V = {\bf R}^n$ or ${\bf C}^n$,
then
\begin{equation}
        A = \{v \in V : N_A(v) \le 1\}.
\end{equation}
To see this, it remains to check that $v \in A$ when $v \in V$
satisfies $N_A(v) = 1$.  In this case, $r^{-1} \, v \in A$ for some
positive real numbers $r$ that are arbitrarily close to $1$, which
implies that $v \in A$ when $A$ is closed.

        Suppose now that $A$ is a convex set in $V$.  Note that this
implies that $A$ is star-like about $0$, because $0 \in A$.  We would
like to show that
\begin{equation}
\label{N_A(v + w) le N_A(v) + N_A(w)}
        N_A(v + w) \le N_A(v) + N_A(w)
\end{equation}
for every $v, w \in V$.

        Let $v, w \in V$ be given, and let $r_v$, $r_w$ be positive real
numbers such that
\begin{equation}
        r_v > N_A(v), \ r_w > N_A(w).
\end{equation}
This implies that $r_v^{-1} \, v$, $r_w^{-1} \, w$ are elements of
$A$, because of the definition of $N_A$ and the fact that $A$ is
star-like about $0$.  If $A$ is convex, then it follows that
\begin{equation}
 (r_v + r_w)^{-1} \, (v + w) = \frac{r_v}{r_v + r_w} \, (r_v^{-1} \, v)
                                + \frac{r_w}{r_v + r_w} \, (r_w^{-1} \, w)
\end{equation}
is an element of $A$ too, so that
\begin{equation}
        N_A(v + w) < r_v + r_w.
\end{equation}
This implies (\ref{N_A(v + w) le N_A(v) + N_A(w)}), since $r_v$, $r_w$
can be arbitrarily close to $N_A(v)$, $N_A(w)$.

\section{Convex cones}
\label{cones}

	Let $V$ be a vector space over the real numbers.  A nonempty
set $E \subseteq V$ is said to be a \emph{cone}\index{cones} if
\begin{equation}
        t \, v \in E
\end{equation}
for every $v \in E$ and nonnegative real number $t$, or equivalently,
\begin{equation}
        t \, E \subseteq E
\end{equation}
for every $t \ge 0$.  In particular, $0 \in E$, since we can take $t = 0$,
and $E \ne \emptyset$.  Note that every linear subspace of $V$ is a cone.

        We shall be especially interested in \emph{convex
cones}\index{convex cones}, which are cones that are also convex sets.
Equivalently, a nonempty set $E \subseteq V$ is a convex cone if
\begin{equation}
        r \, v + t \, w \in E
\end{equation}
for every $v, w \in E$ and nonnegative real numbers $r$, $t$.  Thus
linear subspaces of $V$ are also convex cones.  If $V = {\bf R}^n$,
then it is easy to see that
\begin{equation}
\label{{v in {bf R}^n : v_j ge 0 for j = 1, ldots, n}}
        \{v \in {\bf R}^n : v_j \ge 0 \hbox{ for } j = 1, \ldots, n\}
\end{equation}
is a convex cone in $V$.

        If $E$ is any nonempty subset of $V$, then we can get a cone
$C(E)$ in $V$ from $E$, by taking
\begin{equation}
        C(E) = \bigcup_{t \ge 0} t \, E = \{t \, v : v \in E, \, t \ge 0\}.
\end{equation}
Of course, $C(E) = E$ when $E$ is already a cone in $V$.  One can
check that $C(E)$ is also convex in $V$ when $E$ is convex.

        If $p$ is a sublinear function on $V$, then it is easy to see
that
\begin{equation}
\label{C_p = {v in V : p(v) le 0}}
        C_p = \{v \in V : p(v) \le 0\}
\end{equation}
is a convex cone in $V$.  This cone has another important property,
which is that it is a closed set in a suitable sense.  To make this
precise, it is convenient to suppose that $V = {\bf R}^n$ for some
positive integer $n$.  As in Section \ref{sublinear functions},
sublinear functions on ${\bf R}^n$ are continuous with respect to the
standard Euclidean metric and topology, which implies that $C_p$ is a
closed set.

        Let $V$ be any real vector space again, and let $\|\cdot \|$
be a norm on $V$.  If $E$ is a nonempty subset of $V$, then put
\begin{equation}
        p_E(v) = \inf \{\|v - w\| : w \in E\}
\end{equation}
for each $v \in V$.  Thus $p_E(v) = 0$ for every $v \in E$, and
conversely $p_E(v) = 0$ implies that $v \in E$ when $E$ is a closed
set in a suitable sense.  In particular, this works when $V = {\bf
R}^n$ and $E$ is a closed set with respect to the standard Euclidean
metric and topology, because of the remarks at the end of Section
\ref{definitions, examples}, about the relationship between an
arbitrary norm on ${\bf R}^n$ and the standard Euclidean norm.  Of
course, if $V = {\bf R}^n$, then one can simply take $\|\cdot \|$
to be the standard Euclidean norm on ${\bf R}^n$.

        If $E$ is a convex cone in $V$, then one can check that $p_E$
is a sublinear function on $V$.  It follows that every closed convex
cone in ${\bf R}^n$ can be expressed as in (\ref{C_p = {v in V : p(v)
le 0}}) for some sublinear function $p$ on ${\bf R}^n$.

\section{Dual cones}
\label{dual cones}

        Let $V$ be a real vector space again, and let $E$ be a
nonempty subset of $V$.  Consider the set $E' \subseteq V^*$
consisting of all linear functionals $\lambda$ on $V$ such that
\begin{equation}
\label{lambda(v) ge 0}
        \lambda(v) \ge 0
\end{equation}
for every $v \in E$.  It is easy to see that this is a convex cone in
$V^*$, known as the \emph{dual cone}\index{dual cones} associated to
$E$.  One can also check that
\begin{equation}
        C(E)' = E',
\end{equation}
so that we may as well restrict our attention to convex cones $E$ in $V$.

        Suppose for the moment that $V = {\bf R}^n$ for some positive
integer $n$, and let us identify $V^*$ with ${\bf R}^n$ as in Section
\ref{dual spaces, norms}.  Observe that $E'$ is the same as the dual
cone associated to the closure $\overline{E}$ of $E$ with respect to
the standard Euclidean metric and topology.  Thus we may as well
restrict our attention to closed subsets of ${\bf R}^n$, and hence to
closed convex cones.  Similarly, the dual cone $E'$ of any set $E
\subseteq {\bf R}^n$ automatically corresponds to a closed subset of
${\bf R}^n$, and therefore to a closed convex cone.

        If $V$ is any real vector space and $E$ is a nonempty subset
of $V$, then let $E''$ be the dual cone associated to $E' \subseteq
V^*$, which is a convex cone in the second dual $V^{**}$ of $V$.  As
in Section \ref{second duals}, $V^{**}$ can be identified with $V$ in
a natural way.  Using this identification, one can check that
\begin{equation}
\label{E subseteq E''}
        E \subseteq E''.
\end{equation}

        Suppose for convenience that $V = {\bf R}^n$ for some $n$
again.  If $E$ is a closed convex cone in $V$, then
\begin{equation}
\label{E = E''}
        E = E''.
\end{equation}
To see this, it suffices to show that $E'' \subseteq E$, since we
already know (\ref{E subseteq E''}).  If $v$ is any element of ${\bf
R}^n$ that is not in $E$, then we would like to show that $v$ is also
not in $E''$.  Equivalently, we would like to show that there is a
linear functional $\lambda$ on $V$ such that $\lambda \ge 0$ on $E$
and $\lambda(v) < 0$.

        It is easy to find such a linear functional $\lambda$
initially on the linear span $W$ of $v$ in $V$, and we would like to
extend $\lambda$ to all of $V$ using Theorem \ref{extension theorem,
2} in Section \ref{another extension theorem}.  More precisely, let
$p$ be a sublinear function on $V$ such that $E$ is the set where $p
\le 0$, which exists when $E$ is a closed convex cone in $V = {\bf
R}^n$, as in the previous section.  It is convenient to ask also that
$p \ge 0$ everywhere on $V$, so that $E$ is the set where $p = 0$.
This holds by construction when $p = p_E$ as before, and otherwise one
can simply replace $p$ with $\max (p, 0)$.

        Let $\mu$ be the linear functional defined on the linear span
$W$ of $V$ by
\begin{equation}
\label{mu(t v) = t p(v)}
        \mu(t \, v) = t \, p(v)
\end{equation}
for each $t \in {\bf R}$.  Thus
\begin{equation}
        \mu(t \, v) = p(t \, v)
\end{equation}
when $t \ge 0$, and
\begin{equation}
        \mu(t \, v) = t \, p(v) \le 0 \le p(t \, v)
\end{equation}
when $t \le 0$, so that $\mu \le p$ on $W$.  Theorem \ref{extension
theorem, 2} implies that there is an extension $\widehat{\mu}$ of
$\mu$ to a linear functional on $V$ such that $\widehat{\mu} \le p$ on
all of $V$.  If $\lambda = - \widehat{\mu}$, then $\lambda(v) = - p(v)
< 0$ and $\lambda \ge - p$ on all of $V$, which implies that $\lambda
\ge 0$ on $E$, as desired.

        As an example, consider the case where $E$ is the set of $v
\in {\bf R}^n$ such that $v_j \ge 0$ for $j = 1, \ldots, n$, as in
(\ref{{v in {bf R}^n : v_j ge 0 for j = 1, ldots, n}}).  In this case,
one can check that $E' = E$, using the standard identification of
${\bf R}^n$ with its own dual space.  In particular, $E'' = E$.

\section{Nonnegative self-adjoint operators}
\label{nonnegative self-adjoint operators}

	Let $(V, \langle \cdot, \cdot \rangle)$ be a real or complex
inner product space.  A self-adjoint linear transformation $A$ on $V$
is said to be \emph{nonnegative}\index{nonnegative self-adjoint
operators} if
\begin{equation}
\label{langle A(v), v rangle ge 0}
	\langle A(v), v \rangle \ge 0
\end{equation}
for every $v \in V$.  The identity transformation obviously has this
property, for instance.  If $W$ is a linear subspace of $V$, then we
have seen that the orthogonal projection $P_W$ of $V$ onto $W$ is
self-adjoint, as in Section \ref{self-adjoint linear operators}.  In
this case, we also have that
\begin{equation}
 \langle P_W(v), v \rangle = \langle P_W(v), P_W(v) \rangle = \|P_W(v)\|^2
\end{equation}
for every $v \in V$, as in (\ref{langle P_W(v), w rangle = ... =
langle v, P_W(w) rangle}), and hence that $P_W$ is nonnegative.

        If $a \in V$, then it is easy to see that
\begin{equation}
\label{A_a(v) = langle v, a rangle a}
	A_a(v) = \langle v, a \rangle \, a
\end{equation}
defines a nonnegative self-adjoint linear operator on $V$, which is the
same as the orthogonal projection onto the span of $a$ when $\|a\| = 1$.
Of course, every self-adjoint linear operator $A$ on $V$ is a linear
combination of rank-one operators like this, as a consequence of
diagonalization.  If the eigenvalues of $A$ are nonnegative, then it
follows that $A$ is nonnegative, because $A$ can be expressed as a
linear combination of rank-one operators like this with nonnegative
coefficients, by diagonalization.  Conversely, one can check that the
eigenvalues of a nonnegative self-adjoint linear operator $A$ are
nonnegative, by applying the nonnegativity condition to the
eigenvectors of $A$.

	If $T$ is any linear operator on $V$, then $T^* \, T$ is
self-adjoint and nonnegative, because
\begin{equation}
	(T^* \, T)^* = T^* \, (T^*)^* = T^* \, T
\end{equation}
and
\begin{equation}
	\langle (T^* \, T)(v), v \rangle = \langle T(v), T(v) \rangle
		= \|T(v)\|^2 \ge 0
\end{equation}
for each $v \in V$.  In particular, if $B$ is a self-adjoint linear
operator on $V$, then $B^*$ is a nonnegative self-adjoint linear
operator on $V$.  Conversely, one can use diagonalizations to show
that every nonnegative self-adjoint linear operator $A$ on $V$ can be
expressed as $B^2$ for some nonnegative self-adjoint linear operator
$B$ on $V$.  Note that $-B^2$ is a nonnegative self-adjoint linear
operator on $V$ when $B$ is anti-self-adjoint.

	If $A$ is a nonnegative self-adjoint linear transformation on
$V$, then the trace of $A$ is equal to the sum of the eigenvalues of
$A$, with their appropriate multiplicity, because of diagonalization.
More precisely,
\begin{equation}
	\tr A = \|A\|_{\mathcal{S}_1} = \|A\|_{tr},
\end{equation}
as in Sections \ref{schmidt decompositions} and \ref{mathcal{S}_p norms}.
Similarly, let us check that
\begin{equation}
\label{tr A B ge 0}
	\tr A \, B \ge 0
\end{equation}
for every pair of nonnegative self-adjoint linear operators $A$, $B$
on $V$.  If $A = A_a$ is as in (\ref{A_a(v) = langle v, a rangle a}), then
\begin{equation}
	\tr A \, B = \langle B(a), a \rangle \ge 0.
\end{equation}
Otherwise, $A$ can be expressed as a sum of rank $1$ operators of this
type, as before, and (\ref{tr A B ge 0}) follows.

        Conversely, if $B$ is a self-adjoint linear operator on $V$
such that (\ref{tr A B ge 0}) holds for every nonnegative self-adjoint
linear operator $A$ on $V$, then $B$ is also nonnegative.  This
follows by applying this condition to $A = A_a$ as in (\ref{A_a(v) =
langle v, a rangle a}), as in the previous paragraph.

        Let $\mathcal{L}(V)$ be the vector space of all linear
operators on $V$, and let $\mathcal{L}_{sa}(V)$ be the collection of
self-adjoint linear operators on $V$.  Thus $\mathcal{L}_{sa}(V)$ is a
linear subspace of $\mathcal{L}(V)$ when $V$ is a real vector space.
If $V$ is a complex vector space, then $\mathcal{L}(V)$ is a complex
vector space, but $\mathcal{L}_{sa}(V)$ is a real vector space, which
may be considered as a real-linear subspace of $\mathcal{L}(V)$.

        As in Section \ref{hilbert--schmidt norm},
\begin{equation}
\label{langle A, B rangle_{mathcal{L}(V)} = tr A B^*}
        \langle A, B \rangle_{\mathcal{L}(V)} = \tr A \, B^*
\end{equation}
defines an inner product on $\mathcal{L}(V)$, for which the
corresponding norm is the Hilbert--Schmidt norm.  This reduces to
\begin{equation}
        \langle A, B \rangle_{\mathcal{L}_{sa}(V)} = \tr A \, B
\end{equation}
when $A$ and $B$ are self-adjoint.  Using this inner product, we can
identify $\mathcal{L}_{sa}(V)$ with its own dual space in the usual way.

        The set of nonnegative self-adjoint linear operators on $V$
forms a convex cone in $\mathcal{L}_{sa}(V)$, and it is also a closed
set in $\mathcal{L}_{sa}(V)$ in a suitable sense.  This cone is equal
to its own dual cone, when we identify the dual of
$\mathcal{L}_{sa}(V)$ with itself using the inner product in the
previous paragraph.  This follows from the earlier remarks about
traces of products of self-adjoint linear operators on $V$.

\chapter{Sums and $\ell^p$ spaces}
\label{sums, l^p spaces}

\section{Nonnegative real numbers}
\label{nonnegative real numbers}

        Let $E$ be a nonempty set, and let $f$ be a nonnegative
real-valued function on $E$.  If $A$ is a nonempty subset of $E$ with
only finitely many elements, then the sum
\begin{equation}
\label{sum_{x in A} f(x)}
        \sum_{x \in A} f(x)
\end{equation}
can be defined in the usual way.  The sum
\begin{equation}
\label{sum_{x in E} f(x)}
        \sum_{x \in E} f(x)
\end{equation}
is then defined as the supremum of the finite subsums (\ref{sum_{x in
A} f(x)}), over all finite nonempty subsets $A$ of $E$.  More
precisely, the sum (\ref{sum_{x in E} f(x)}) is considered to be
$+\infty$ as an extended real number when there is no finite upper
bound for the finite subsums (\ref{sum_{x in A} f(x)}).

        Of course, if $E$ is itself a finite set, then this definition
of the sum (\ref{sum_{x in E} f(x)}) reduces to the usual one.  If $E$
is the set ${\bf Z}_+$ of positive integers, then the usual definition
of the sum of an infinite series is equivalent to
\begin{equation}
\label{sum_{j = 1}^infty f(j) = sup_{n ge 1} sum_{j = 1}^n f(j)}
        \sum_{j = 1}^\infty f(j) = \sup_{n \ge 1} \sum_{j = 1}^n f(j),
\end{equation}
which is again interpreted as being $+\infty$ when there is no finite
upper bound for the partial sums.  In this case, this definition of
the infinite sum is equivalent to the previous one.  More precisely,
this definition of the sum is less than or equal to the previous one,
because the partial sums $\sum_{j = 1}^n f(j)$ are subsums of the form
(\ref{sum_{x in A} f(x)}) for each $n$.  Similarly, the previous definition
of the sum is less than or equal to this one, because every finite set 
$A \subseteq E = {\bf Z}_+$ is contained in a set of the form $\{1, \ldots, 
n\}$ for some $n$, so that the corresponding subsum (\ref{sum_{x in A} f(x)})
is less than or equal to the partial sum $\sum_{j = 1}^n f(j)$.

        Suppose that the finite subsums (\ref{sum_{x in A} f(x)}) are
bounded by a nonnegative real number $C$, so that the sum (\ref{sum_{x
in E} f(x)}) is also less than or equal to $C$.  In this case, it is
easy to see that
\begin{equation}
\label{E(f, epsilon) = {x in E : f(x) ge epsilon}}
        E(f, \epsilon) = \{x \in E : f(x) \ge \epsilon\}
\end{equation}
has at most $C / \epsilon$ elements for each $\epsilon > 0$.
Otherwise, if $f(x) \ge \epsilon$ for more than $C / \epsilon$
elements $x$ of $E$, then there would be a finite set $A \subseteq E$
such that (\ref{sum_{x in A} f(x)}) is larger than $C$.  
In particular, $E(f, \epsilon)$ has only finitely many elements for
each $\epsilon > 0$.  This implies that
\begin{equation}
\label{{x in E : f(x) > 0} = bigcup_{n = 1}^infty E(f, 1/n)}
        \{x \in E : f(x) > 0\} = \bigcup_{n = 1}^\infty E(f, 1/n)
\end{equation}
has only finitely or countably many elements, so that (\ref{sum_{x in
E} f(x)}) can be reduced to a finite sum or an ordinary infinite series.

        If $f$, $g$ are nonnegative real-valued functions on $E$, then
one can check that
\begin{equation}
\label{sum_{x in E} (f(x) + g(x)) = sum_{x in E} f(x) + sum_{x in E} g(x)}
 \sum_{x \in E} (f(x) + g(x)) = \sum_{x \in E} f(x) + \sum_{x \in E} g(x),
\end{equation}
by reducing to the case of finite sums.  This includes the case where
some of the sums are infinite, with the usual conventions for sums of
extended real numbers.  Similarly, if $a$ is a nonnegative real
number, then
\begin{equation}
\label{sum_{x in E} a f(x) = a sum_{x in E} f(x)}
        \sum_{x \in E} a \, f(x) = a \, \sum_{x \in E} f(x),
\end{equation}
with the convention that $a \cdot (+\infty) = +\infty$ when $a > 0$.
In this context, it is also appropriate to make the convention that $0
\cdot (+\infty) = 0$, since the left side of (\ref{sum_{x in E} a f(x)
= a sum_{x in E} f(x)}) is automatically equal to $0$ when $a = 0$.

\section{Summable functions}
\label{summable functions}

        A real or complex-valued function $f$ on a nonempty set $E$
is said to be \emph{summable}\index{summable functions} on $E$ if
\begin{equation}
\label{sum_{x in E} |f(x)|}
        \sum_{x \in E} |f(x)|
\end{equation}
is finite.  If $f$ and $g$ are summable functions on $E$, then it is
easy to see that $f + g$ is summable, since
\begin{equation}
\label{|f(x) + g(x)| le |f(x)| + |g(x)|}
        |f(x) + g(x)| \le |f(x)| + |g(x)|
\end{equation}
for each $x \in E$.  Similarly, if $f$ is a summable function on $E$
and $a$ is a real or complex number, as appropriate, then $a \, f$ is
a summable function on $E$ too.  Thus the real or complex-valued
summable functions on $E$ form a vector space over the real or complex
numbers, as appropriate.

        Let $f$ be a real or complex-valued summable function on $E$,
and let $\epsilon > 0$ be given.  By definition of (\ref{sum_{x in E}
|f(x)|}), there is a finite set $A_\epsilon \subseteq E$ such that
\begin{equation}
        \sum_{x \in E} |f(x)| \le \sum_{x \in A_\epsilon} |f(x)| + \epsilon.
\end{equation}
Equivalently,
\begin{equation}
        \sum_{x \in A} |f(x)| \le \sum_{x \in A_\epsilon} |f(x)| + \epsilon
\end{equation}
for every finite set $A \subseteq E$.  If $B$ is a finite subset of $E$
which is disjoint from $A_\epsilon$, then we get that
\begin{equation}
\label{sum_{x in B} |f(x)| le epsilon}
        \sum_{x \in B} |f(x)| \le \epsilon,
\end{equation}
by applying the previous inequality to $A = A_\epsilon \cup B$.  This
will be used frequently to approximate various sums by finite sums.

        We would like to define the sum $\sum_{x \in E} f(x)$ of a
real or complex-valued summable function $f$ on $E$.  One way to do
this is to express $f$ as a linear combination of nonnegative
real-valued summable functions on $E$, and then use the previous
definition of the sum for nonnegative real-valued functions on $E$.
Another way to do this is to use the fact that the set of $x \in E$
such that $f(x) \ne 0$ has only finitely or countably many elements,
as in the previous section, and then treat the sum as a finite sum or
an infinite series, by enumerating the elements of this set.  The
summability of $f$ on $E$ implies that such an infinite series would
be absolutely convergent, and hence convergent.  The value of the sum
would also not depend on the choice of the enumeration of the set
where $f \ne 0$, because absolutely-convergent infinite series are
invariant under rearrangements.

        Whichever way one uses to define the sum, a key point is that
it should be approximable by finite sums in a natural way.  More
precisely, for each $\epsilon > 0$, there should be a finite set
$A_\epsilon \subseteq E$ such that
\begin{equation}
\label{|sum_{x in E} f(x) - sum_{x in A} f(x)| le epsilon}
 \biggl|\sum_{x \in E} f(x) - \sum_{x \in A} f(x)\biggr| \le \epsilon
\end{equation}
for every finite set $A \subseteq E$ such that $A_\epsilon \subseteq
A$.  One can check that both of the approaches to defining the sum
described in the previous paragraph have this approximation property.
Practically any reasonable way of defining the sum as a limit of
finite subsums should also have this property, because of the
approximation property (\ref{sum_{x in B} |f(x)| le epsilon})
mentioned earlier.

        At any rate, it is easy to see that any definition of the sum
that satisfies (\ref{|sum_{x in E} f(x) - sum_{x in A} f(x)| le
epsilon}) is uniquely determined by this property.  In particular,
this can be helpful for showing that such a definition does not depend
on any auxiliary choices used to define the sum.  In effect,
(\ref{|sum_{x in E} f(x) - sum_{x in A} f(x)| le epsilon})
characterizes the sum as a somewhat fancy limit of the finite subsums
$\sum_{x \in A} f(x)$.  This can be made precise by treating these
subsums as a net, which is indexed by the collection of all finite
subsets $A$ of $E$.  This also uses the natural partial ordering on
the collection of all finite subsets of $E$ by inclusion.

        If $f$ is any real or complex-valued summable function on $E$, then
\begin{equation}
\label{|sum_{x in E} f(x)| le sum_{x in E} |f(x)|}
        \biggl|\sum_{x \in E} f(x)\biggr| \le \sum_{x \in E} |f(x)|.
\end{equation}
This follows easily by approximating the sum by finite subsums, and
using the triangle inequality.  If $f$ and $g$ are summable functions
on $E$, and $a$ and $b$ are real or complex numbers, as appropriate,
then one can also check that
\begin{equation}
        \sum_{x \in E} (a \, f(x) + b \, g(x))
                      = a \, \sum_{x \in E} f(x) + b \, \sum_{x \in E} g(x).
\end{equation}
As usual, this can be shown by approximating the sums by finite
subsums.  Using these two properties of the sum, we get that
\begin{equation}
\label{|sum_{x in E} f(x) - sum_{x in E} g(x)| le sum_{x in E} |f(x) - g(x)|}
        \biggl|\sum_{x \in E} f(x) - \sum_{x \in E} g(x)\biggr|
 = \biggl|\sum_{x \in E} (f(x) - g(x))\biggr| \le \sum_{x \in E} |f(x) - g(x)|
\end{equation}
for any pair of summable functions $f$, $g$ on $E$.

        If $f$ is any real or complex-valued function on $E$, then the
\emph{support}\index{support} of $f$ is the set $\supp f$ of $x \in E$
such that $f(x) \ne 0$.  If the support of $f$ has only finitely many
elements, then the sum $\sum_{x \in E} f(x)$ reduces to an ordinary
finite sum.  The sum $\sum_{x \in E} f(x)$ of a summable function $f$
on $E$ can be characterized as the unique linear functional on the
vector space of all summable functions on $E$ that satisfies
(\ref{|sum_{x in E} f(x)| le sum_{x in E} |f(x)|}) and reduces to the
ordinary finite sum when $f$ has finite support.  This follows by
approximating an arbitrary summable function on $E$ by functions with
finite support using (\ref{sum_{x in B} |f(x)| le epsilon}), and then
using (\ref{|sum_{x in E} f(x) - sum_{x in E} g(x)| le sum_{x in E}
|f(x) - g(x)|}) to analyze the corresponding sums.

\section{Convergence theorems}
\label{convergence theorems}

	Let $\{f_j\}_{j=1}^\infty$ be a sequence of real or
complex-valued summable functions on a nonempty set $E$, and suppose that
\begin{equation}
\label{lim_{j to infty} f_j(x) = f(x)}
	\lim_{j \to \infty} f_j(x) = f(x)
\end{equation}
for every $x \in E$.  Of course,
\begin{equation}
\label{lim_{j to infty} sum_{x in E} f_j(x) = sum_{x in E} f(x)}
	\lim_{j \to \infty} \sum_{x \in E} f_j(x) = \sum_{x \in E} f(x)
\end{equation}
when $E$ is finite, but the situation is more complicated when $E$ is
infinite.  For example, if $E$ is the set of positive integers,
$f_j(j) = 1$ for each $j$, and $f_j(x) = 0$ when $x \ne j$, then
\begin{equation}
        \sum_{x \in E} f_j(x) = 1
\end{equation}
for each $j$, but $\{f_j(x)\}_{j = 1}^\infty$ converges to $0$ for
each $x \in E$.  Alternatively, if we take $f_j(x) = 1/j$ when $x \le
j$ and $f_j(x) = 0$ when $x > j$, then we get the same conclusions,
but with $\{f_j\}_{j = 1}^\infty$ converging to $0$ uniformly on $E$.

        However, there are some positive results, as follows.  Suppose
that $\{f_j(x)\}_{j = 1}^\infty$ is a monotone increasing sequence of
nonnegative real valued functions on $E$, so that
\begin{equation}
        f_j(x) \le f_{j + 1}(x)
\end{equation}
for every $x \in E$ and $j \ge 1$.  Put
\begin{equation}
        f(x) = \sup_{j \ge 1} f_j(x)
\end{equation}
for each $x \in E$, which may be $+\infty$.  Thus
\begin{equation}
        f_j(x) \to \sup_{l \ge 1} f_l(x) \quad\hbox{as } j \to \infty
\end{equation}
for every $x \in E$, with the usual interpretation when $f(x) =
+\infty$.  Under these conditions, the \emph{monotone convergence
theorem}\index{monotone convergence theorem} implies that
\begin{equation}
\label{sum_{x in E} f_j(x) to sum_{x in E} f(x) as j to infty}
 \sum_{x \in E} f_j(x) \to \sum_{x \in E} f(x) \quad\hbox{as } j \to \infty,
\end{equation}
with suitable interpretations when any of the quantities involved are infinite.

        Of course,
\begin{equation}
 \sum_{x \in E} f_j(x) \le \sum_{x \in E} f_{j + 1}(x) \le \sum_{x \in E} f(x)
\end{equation}
for each $j$, by monotonicity.  If $\sum_{x \in E} f_l(x) = +\infty$
for some $l$, then
\begin{equation}
        \sum_{x \in E} f_j(x) = \sum_{x \in E} f(x) = +\infty
\end{equation}
for each $j \ge l$, and (\ref{sum_{x in E} f_j(x) to sum_{x in E} f(x)
as j to infty}) is trivial.  Similarly, if $f(x_0) = +\infty$ for some
$x_0 \in E$, then $\sum_{x \in E} f(x)$ is interpreted as being equal
to $+\infty$, and
\begin{equation}
        \sum_{x \in E} f_j(x) \ge f_j(x_0)
\end{equation}
also tends to $+\infty$ as $j \to \infty$.  Suppose then that each
$f_j$ is summable, and that $f(x)$ is finite for every $x \in E$.  In
this case, one can get (\ref{sum_{x in E} f_j(x) to sum_{x in E} f(x)
as j to infty}) using the fact that
\begin{equation}
	\sum_{x \in A} f(x) = \lim_{j \to \infty} \sum_{x \in A} f_j(x)
			\le \lim_{j \to \infty} \sum_{x \in E} f_j(x)
\end{equation}
for every finite set $A \subseteq E$, since $\sum_{x \in E} f(x)$ is
equal to the supremum of $\sum_{x \in A} f(x)$ over all finite subsets
$A$ of $E$.

        Now let $\{f_j\}_{j = 1}^\infty$ be any sequence of
nonnegative real-valued functions on $E$, and put
\begin{equation}
        f(x) = \liminf_{j \to \infty} f_j(x).
\end{equation}
Under these conditions, \emph{Fatou's lemma}\index{Fatou's lemma}
states that
\begin{equation}
\label{sum_{x in E} f(x) le liminf_{j to infty} sum_{x in E} f_j(x)}
	\sum_{x \in E} f(x) \le \liminf_{j \to \infty} \sum_{x \in E} f_j(x),
\end{equation}
again with suitable interpretations when any of the quantities are infinite.
If $A$ is a finite subset of $E$, then
\begin{equation}
        \sum_{x \in A} f(x) \le \liminf_{j \to \infty} \sum_{x \in A} f_j(x)
\end{equation}
is a well known property of the lower limit.  This implies that
\begin{equation}
\label{sum_{x in A} f(x) le liminf_{j to infty} sum_{x in E} f_j(x)}
        \sum_{x \in A} f(x) \le \liminf_{j \to \infty} \sum_{x \in E} f_j(x),
\end{equation}
because the sum of $f_j(x)$ over $x \in A$ is less than or equal to
the sum of $f_j(x)$ over $x \in E$ when $A \subseteq E$.  In order to
get (\ref{sum_{x in E} f(x) le liminf_{j to infty} sum_{x in E}
f_j(x)}), it suffices to take the supremum of the left side of
(\ref{sum_{x in A} f(x) le liminf_{j to infty} sum_{x in E} f_j(x)})
over all finite subsets $A$ of $E$.

        Suppose now that $\{f_j\}_{j = 1}^\infty$ is a sequence of
real or complex valued functions on $E$ that converges pointwise to
another function $f$ on $E$.  Suppose also that $h$ is a nonnegative
real-valued function on $E$ which is summable, and that
\begin{equation}
\label{|f_j(x)| le h(x)}
        |f_j(x)| \le h(x)
\end{equation}
for every $x \in E$ and $j \ge 1$.  This implies that
\begin{equation}
\label{|f(x)| le h(x)}
        |f(x)| \le h(x)
\end{equation}
for every $x \in E$, because $\{f_j(x)\}_{j = 1}^\infty$ converges to
$f(x)$ for every $x \in E$ by hypothesis.  Note that $f_j$ is a
summable function on $E$ for each $j$, and that $f$ is a summable
function on $E$ too, since $h$ is summable.

        Under these conditions, the \emph{dominated convergence
theorem}\index{dominated convergence theorem} implies that
\begin{equation}
        \lim_{j \to \infty} \sum_{x \in E} f_j(x) = \sum_{x \in E} f(x).
\end{equation}
To prove this, it suffices to show that
\begin{equation}
	\lim_{j \to \infty} \sum_{x \in E} |f_j(x) - f(x)| = 0.
\end{equation}
Let $\epsilon > 0$ be given, and let us show that there is a $L \ge 1$
such that
\begin{equation}
\label{sum_{x in E} |f(x) - f_j(x)| < epsilon}
	\sum_{x \in E} |f(x) - f_j(x)| < \epsilon
\end{equation}
for every $j \ge L$.

        Because $h$ is summable on $E$, there is a finite set
$E_\epsilon \subseteq E$ such that
\begin{equation}
        \sum_{x \in E \backslash E_\epsilon} |h(x)| \le \frac{\epsilon}{3},
\end{equation}
as in (\ref{sum_{x in B} |f(x)| le epsilon}).  Thus
\begin{equation}
\label{sum_{x in E backslash E_epsilon} |f(x) - f_j(x)| le frac{2 epsilon}{3}}
	\sum_{x \in E \backslash E_\epsilon} |f(x) - f_j(x)|
 \le \sum_{x \in E \backslash E_\epsilon} 2 \, |h(x)| \le \frac{2 \epsilon}{3}
\end{equation}
for each $j$.  We also have that
\begin{equation}
\label{lim_{j to infty} sum_{x in E_epsilon} |f(x) - f_j(x)| = 0}
        \lim_{j \to \infty} \sum_{x \in E_\epsilon} |f(x) - f_j(x)| = 0,
\end{equation}
because $E_\epsilon$ is finite and $\{f_j(x)\}_{j = 1}^\infty$
converges to $f(x)$ for each $x \in E$.  This implies that there is an
$L \ge 1$ such that
\begin{equation}
\label{sum_{x in E_epsilon} |f(x) - f_j(x)| < frac{epsilon}{3}}
	\sum_{x \in E_\epsilon} |f(x) - f_j(x)| < \frac{\epsilon}{3}
\end{equation}
when $j \ge L$.  Combining (\ref{sum_{x in E backslash E_epsilon}
|f(x) - f_j(x)| le frac{2 epsilon}{3}}) and (\ref{sum_{x in E_epsilon}
|f(x) - f_j(x)| < frac{epsilon}{3}}), we get that (\ref{sum_{x in E}
|f(x) - f_j(x)| < epsilon}) holds when $j \ge L$, as desired.

\section{Double sums}
\label{double sums}

        Let $E_1$ and $E_2$ be nonempty sets, and let $E = E_1 \times
E_2$ be their Cartesian product.  If $f(x, y)$ is a nonnegative
real-valued function on $E$, then we would like to check that
\begin{equation}
\label{various sums}
\sum_{(x,y) \in E} f(x,y) = \sum_{x \in E_1} \Big(\sum_{y \in E_2} f(x,y)\Big)
	   = \sum_{y \in E_2} \Big(\sum_{x \in E_1} f(x,y)\Big).
\end{equation}
More precisely, it may be that $\sum_{y \in E_2} f(x', y) = +\infty$
for some $x' \in E_1$, or that $\sum_{x \in E_1} f(x, y') = +\infty$
for some $y' \in E_2$, in which case the corresponding iterated sum in
(\ref{various sums}) is considered to be $+\infty$ as well.

        Let $A$ be a finite subset of $E$, and $A_1$, $A_2$ be finite
subsets of $E_1$, $E_2$, respectively, such that
\begin{equation}
	A \subseteq A_1 \times A_2.
\end{equation}
Clearly
\begin{equation}
\label{sum_{(x, y) in A} f(x, y) le ...}
        \sum_{(x, y) \in A} f(x, y) 
              \le \sum_{x \in A_1} \Big(\sum_{y \in A_2} f(x, y)\Big)
		= \sum_{y \in A_2} \Big(\sum_{x \in A_1} f(x, y)\Big).
\end{equation}
This implies that $\sum_{(x, y) \in A} f(x, y)$ is less than or equal
to each of
\begin{equation}
	\sum_{x \in E_1} \Big(\sum_{y \in E_2} f(x, y)\Big)
		\quad\hbox{and}\quad
	\sum_{y \in E_2} \Big(\sum_{x \in E_1} f(x, y)\Big).
\end{equation}
It follows that $\sum_{(x, y) \in E} f(x, y)$ is less than or equal to
the same iterated sums, because $A$ is an arbitrary finite subset of $E$.

	In the other direction, if $A_1$ and $A_2$ are arbitary finite
subsets of $E_1$ and $E_2$, respectively, then
\begin{equation}
	\sum_{x \in A_1} \Big(\sum_{y \in A_2} f(x, y)\Big)
		= \sum_{y \in A_2} \Big(\sum_{x \in A_1} f(x, y)\Big)
		= \sum_{(x, y) \in A_1 \times A_2} f(x, y),
\end{equation}
and therefore
\begin{equation}
	\sum_{x \in A_1} \Big(\sum_{y \in A_2} f(x, y)\Big)
		= \sum_{y \in A_2} \Big(\sum_{x \in A_1} f(x, y)\Big)
		\le \sum_{(x, y) \in E} f(x, y).
\end{equation}
Using this, it is easy to see that
\begin{equation}
	\sum_{x \in A_1} \Big(\sum_{y \in E_2} f(x, y)\Big)
		\le \sum_{(x, y) \in E} f(x, y)
\end{equation}
and
\begin{equation}
	\sum_{y \in A_2} \Big(\sum_{x \in E_1} f(x, y)\Big)
		\le \sum_{(x, y) \in E} f(x, y).
\end{equation}
This implies that the iterated sums are less than or equal to the sum
of $f(x, y)$ over $E$, as desired, by taking the suprema over $A_1$
and $A_2$.

        Now let $f(x, y)$ be a real or complex-valued function on $E =
E_1 \times E_2$.  If any of the sums
\begin{equation}
\label{sums of |f(x, y)|}
	\sum_{(x, y) \in E} |f(x, y)|, 
	  \quad \sum_{x \in E_1} \Big(\sum_{y \in E_2} |f(x, y)|\Big),
	  \quad \sum_{y \in E_2} \Big(\sum_{x \in E_1} |f(x, y)|\Big)
\end{equation}
are finite, then they are all finite, and equal to each other, by the
previous discussion.  Suppose that this is the case, so that $f(x, y)$
is a summable function on $E$.  We also get that
\begin{equation}
	\sum_{y \in D_2} |f(x, y)| < +\infty
\end{equation}
for every $x \in D_1$, and
\begin{equation}
	\sum_{x \in D_1} |f(x, y)| < +\infty
\end{equation}
for every $y \in D_2$, so that $f(x, y)$ is a summable function on
$D_2$ for each $x \in D_1$, and $f(x, y)$ is a summable function on
$D_1$ for each $y \in D_2$.  Put
\begin{equation}
\label{f_1(x) = sum_{y in E_2} f(x, y)}
	f_1(x) = \sum_{y \in E_2} f(x, y)
\end{equation}
for each $x \in E_1$, and
\begin{equation}
\label{f_2(y) = sum_{x in E_1} f(x, y)}
	f_2(y) = \sum_{x \in E_1} f(x, y)
\end{equation}
for each $y \in E_2$.  Thus
\begin{equation}
	|f_1(x)| \le \sum_{y \in E_2} |f(x, y)|
\end{equation}
for every $x \in E_1$, and
\begin{equation}
	|f_2(y)| \le \sum_{x \in E_1} |f(x, y)|
\end{equation}
for every $y \in E_2$.  This implies that $f_1(x)$ is a summable
function on $E_1$, and that $f_2(y)$ is a summable function on $E_2$,
because of the finiteness of the iterated sums in (\ref{sums of |f(x,
y)|}).  Under these conditions, we also have that
\begin{equation}
\label{sums of f(x, y)}
        \sum_{(x, y) \in E} f(x, y)
              = \sum_{x \in E_1} f_1(x)
               = \sum_{y \in E_2} f_2(y).
\end{equation}
One way to show this is to express $f(x, y)$ as a linear combination
of nonnegative summable functions on $E$, and apply (\ref{various
sums}) in that case.  Alternatively, one can approximate $f(x, y)$ by
functions with finite support on $E$, as in (\ref{sum_{x in B} |f(x)|
le epsilon}).  The equality of the iterated and double sums is clear
for functions with finite support, and one can use the previous
results for nonnegative functions to estimate the errors in the
approximations.

\section{$\ell^p$ Spaces}
\label{ell^p spaces}

	Let $E$ be a nonempty set, and let $p$ be a positive real
number.  A real or complex-valued function $f$ on $E$ is said to be
\emph{$p$-summable}\index{p-summability@$p$-summability} if $|f(x)|^p$
is a summable function on $E$.  The spaces of real and complex-valued
$p$-summable functions on $D$ are denoted $\ell^p(E, {\bf R})$ and
$\ell^p(E, {\bf C})$, respectively, although we may also use
$\ell^p(E)$ to include both cases at the same time.  If $f$ and $g$
are $p$-summable functions on $E$, then it is easy to see that $f + g$
is also $p$-summable, using the observation that
\begin{equation}
        |f(x) + g(x)| \le |f(x)| + |g(x)| \le 2 \, \max(|f(x)|, |g(x)|)
\end{equation}
for every $x \in E$, and hence
\begin{equation}
        |f(x) + g(x)|^p \le 2^p \, \max(|f(x)|^p, |g(x)|^p)
                         \le 2^p \, (|f(x)|^p + |g(x)|^p).
\end{equation}
Similarly, $a \, f(x)$ is $p$-summable on $E$ when $f(x)$ is
$p$-summable on $E$ and $a$ is a real or complex number, as
appropriate, so that $\ell^p(E, {\bf R})$ and $\ell^p(E, {\bf C})$ are
vector spaces with respect to pointwise addition and scalar
multiplication of functions on $E$.

	If $f$ is a real or complex-valued $p$-summable function on
$E$, then we put
\begin{equation}
	\|f\|_p = \Big(\sum_{x \in E} |f(x)|^p \Big)^{1/p}.
\end{equation}
Thus
\begin{equation}
	\|a \, f\|_p = |a| \, \|f\|_p
\end{equation}
for every real or complex number $a$, as appropriate.  In this
context, Minkowski's inequality\index{Minkowski's inequality} states that
\begin{equation}
\label{||f_1 + f_2||_p le ||f_1||_p + ||f_2||_p}
	\|f_1 + f_2\|_p \le \|f_1\|_p + \|f_2\|_p
\end{equation}
when $p \ge 1$ and $f_1$, $f_2$ are $p$-summable functions on $E$.
This follows from the version (\ref{Minkowski's inequality for sums})
of Minkowski's inequality for finite sums, although analogous
arguments could be used more directly in this case.  If $0 < p \le 1$,
then
\begin{equation}
\label{||f_1 + f_2||_p^p le ||f_1||_p^p + ||f_2||_p^p}
	\|f_1 + f_2\|_p^p \le \|f_1\|_p^p + \|f_2\|_p^p,
\end{equation}
as in (\ref{sum_{j=1}^n (b_j + c_j)^p le sum_{j=1}^n b_j^p + sum_j c_j^p}).

	Let $\ell^\infty(E, {\bf R})$, $\ell^\infty(E, {\bf C})$ be
the spaces of bounded real and complex-valued functions $f$ on $E$,
respectively, which is to say that the values of $f$ on $E$ are
contained in a bounded subset of ${\bf R}$ or ${\bf C}$.  Of course,
the sum and product of two bounded functions is bounded.  If $f$ is a
bounded function on $E$, then put
\begin{equation}
	\|f\|_\infty = \sup \{|f(x)| : x \in E\}.
\end{equation}
Clearly
\begin{equation}
	\|a \, f\|_\infty = |a| \, \|f\|_\infty
\end{equation}
for every real or complex number $a$, as appropriate.  One can also check that
\begin{equation}
	\|f_1 + f_2\|_\infty \le \|f_1\|_\infty + \|f_2\|_\infty
\end{equation}
and
\begin{equation}
	\|f_1 \, f_2\|_\infty \le \|f_1\|_\infty \, \|f_2\|_\infty
\end{equation}
for all bounded real or complex-valued functions $f_1$, $f_2$ on $E$.

        If $f$ and $g$ are real or complex-valued $2$-summable functions
on $E$, then their product $f \, g$ is a summable function on $E$, because
\begin{equation}
\label{|f(x)| |g(x)| le max(|f(x)|^2, |g(x)|^2) le |f(x)|^2 + |g(x)|^2}
        |f(x)| \, |g(x)| \le \max(|f(x)|^2, |g(x)|^2) \le |f(x)|^2 + |g(x)|^2
\end{equation}
for every $x \in E$.  Alternatively, one can use the well-known fact that
\begin{equation}
        2 \, a \, b \le a^2 + b^2
\end{equation}
for all nonnegative real numbers $a$ and $b$ to get a better estimate.  Put
\begin{equation}
\label{langle f, g rangle = sum_{x in E} f(x) g(x)}
        \langle f, g \rangle = \sum_{x \in E} f(x) \, g(x)
\end{equation}
in the case of real-valued functions on $E$, and
\begin{equation}
\label{langle f, g rangle = sum_{x in E} f(x) overline{g(x)}}
        \langle f, g \rangle = \sum_{x \in E} f(x) \, \overline{g(x)}
\end{equation}
in the complex case.  It is easy to see that (\ref{langle f, g rangle
= sum_{x in E} f(x) g(x)}) defines an inner product on $\ell^2(E, {\bf
R})$, that (\ref{langle f, g rangle = sum_{x in E} f(x) overline{g(x)}})
defines an inner product on $\ell^2(E, {\bf C})$, and that the corresponding
norms are the same as the $\ell^2$ norm $\|f\|_2$ defined earlier.

\section{Additional properties}
\label{additional properties}

        If $f$ is a $p$-summable function on a nonempty set $E$ for
some $p > 0$, then it is easy to see that $f$ is bounded on $E$, and that
\begin{equation}
\label{||f||_infty le ||f||_p}
	\|f\|_\infty \le \|f\|_p,
\end{equation}
as in (\ref{max_{1 le j le n} a_j le (sum_{j=1}^n a_j^p)^{1/p}}).
Similarly, if $p$, $q$ are positive real numbers with $p \le q$, and
$f$ is a $p$-summable function on $E$, then $f$ is $q$-summable, and
\begin{equation}
\label{||f||_q le ||f||_p}
	\|f\|_q \le \|f\|_p,
\end{equation}
as in (\ref{(sum_{j=1}^n a_j^q)^{1/q} le (sum_{j=1}^n a_j^p)^{1/p}}).
        
	A real or complex-valued function $f$ on $E$ is said to
``vanish at infinity'' if for each $\epsilon > 0$ there is a finite
set $A_\epsilon \subseteq E$ such that
\begin{equation}
	|f(x)| < \epsilon
\end{equation}
for every $x \in E \backslash A_\epsilon$.  Of course, any function on
$E$ with finite support satisfies this condition.  Let $c_0(E, {\bf
R})$, $c_0(E, {\bf C})$ denote the spaces of real and complex-valued
functions on $E$, respectively, with this property.  As before, we may
also use the notation $c_0(E)$ to include both cases at the same time.
Note that these are linear subspaces of the corresponding spaces of
bounded functions on $E$.

        If $f$ is $p$-summable on $E$ for some $p > 0$, then $f$
vanishes at infinity on $E$.  Equivalently, if there is an $\epsilon >
0$ such that $|f(x)| \ge \epsilon$ for infinitely many $x \in E$, then
$f$ is not $p$-summable for any $p > 0$.  Note that a bounded function
$f$ on $E$ vanishes at infinity if and only if for each $\epsilon > 0$
there is a function $f_\epsilon$ with finite support on $E$ such that
\begin{equation}
\label{||f - f_epsilon||_infty < epsilon}
        \|f - f_\epsilon\|_\infty < \epsilon.
\end{equation}
If $f$ is a $p$-summable function on $E$, then one can check that for
each $\epsilon > 0$ there is a function $f_\epsilon$ with finite
support on $E$ such that
\begin{equation}
\label{||f - f_epsilon||_p < epsilon}
        \|f - f_\epsilon\|_p < \epsilon.
\end{equation}
In particular, this implies (\ref{||f - f_epsilon||_infty < epsilon})
in this case, because of (\ref{||f||_infty le ||f||_p}).

        Suppose now that $\{f_j\}_{j = 1}^\infty$ is a sequence of
real or complex-valued functions on $E$ that converges pointwise to a
function $f$ on $E$.  If there are positive real numbers $p$, $C$ such
that $f_j$ is $p$-summable for each $j$, with
\begin{equation}
        \|f_j\|_p \le C
\end{equation}
for each $j$, then $f$ is $p$-summable too, and
\begin{equation}
        \|f\|_p \le C.
\end{equation}
This can be derived from Fatou's lemma, as in Section \ref{convergence
theorems}.

        Similarly, if $f_j$ is a bounded function on $E$ for each $j$, with
\begin{equation}
        \|f_j\|_\infty \le C
\end{equation}
for each $j$, then $f$ is bounded on $E$ as well, and
\begin{equation}
        \|f\|_\infty \le C.
\end{equation}
This is easy to verify, directly from the definitions.  However, it is
easy to give examples where $f_j$ vanishes at infinity on $E$ for each
$j$, but $f$ does not vanish at infinity on $E$.  If $f_j$ vanishes at
infinity on $E$ for each $j$, and $\{f_j\}_{j = 1}^\infty$ converges
to $f$ uniformly on $E$, then $f$ vanishes at infinity on $E$, by
standard arguments.  Of course, $\{f_j\}_{j = 1}^\infty$ converges to
$f$ uniformly on $E$ if and only if $\{f_j\}_{j = 1}^\infty$ converges
to $f$ with respect to the $\ell^\infty$ norm, in the sense that
\begin{equation}
\label{lim_{j to infty} ||f_j - f||_infty = 0}
        \lim_{j \to \infty} \|f_j - f\|_\infty = 0.
\end{equation}

        Let $\{f_j\}_{j = 1}^\infty$ be a sequence of functions on $E$
in $\ell^p(E, {\bf R})$ or $\ell^p(E, {\bf C})$ for some $p$, $0 < p
\le \infty$.  As usual, we say that $\{f_j\}_{j = 1}^\infty$ is a
\emph{Cauchy sequence}\index{Cauchy sequences} in $\ell^p(E, {\bf R})$
or $\ell^p(E, {\bf C})$, as appropriate, if for each $\epsilon > 0$
there is an $L(\epsilon) \ge 1$ such that
\begin{equation}
\label{||f_j - f_l||_p < epsilon}
	\|f_j - f_l\|_p < \epsilon
\end{equation}
for every $j, l \ge L(\epsilon)$.  This is equivalent to saying that
$\{f_j\}_{j = 1}^\infty$ is a Cauchy sequence with respect to the metric
\begin{equation}
        d_p(g, h) = \|g - h\|_p
\end{equation}
on $\ell^p(E, {\bf R})$ or $\ell^p(E, {\bf C})$ when $p \ge 1$, and
with respect to the metric
\begin{equation}
        d_p(g, h) = \|g - h\|_p^p
\end{equation}
when $0 < p < 1$.  We would like to show that $\{f_j\}_{j = 1}^\infty$
converges to some function $f$ on $E$ in $\ell^p(E, {\bf R})$ or
$\ell^p(E, {\bf C})$, as appropriate, in the sense that
\begin{equation}
\label{lim_{j to infty} ||f_j - f||_p = 0}
        \lim_{j \to \infty} \|f_j - f\|_p = 0,
\end{equation}
and thereby conclude that $\ell^p(E, {\bf R})$ and $\ell^p(E, {\bf
C})$ are complete as metric spaces.

        To do this, observe first that $\{f_j(x)\}_{j=1}^\infty$ is a
Cauchy sequence of real or complex numbers, as appropriate, for every
$x \in E$, because
\begin{equation}
        |f_j(x) - f_l(x)| \le \|f_j - f_l\|_p
\end{equation}
for every $x \in E$, $j, l \ge 1$, and $0 < p \le \infty$.  Thus
$\{f_j(x)\}_{j = 1}^\infty$ converges to a real or complex number
$f(x)$, as appropriate, for each $x \in E$, since every Cauchy
sequence of real and complex numbers converges.  Note that
\begin{equation}
        \|f_j\|_p \le \|f_{L(1)}\|_p + 1
\end{equation}
for every $j \ge L(1)$ when $p \ge 1$, and that
\begin{equation}
        \|f_j\|_p^p \le \|f_{L(1)}\|_p^p + 1
\end{equation}
for every $j \ge L(1)$ when $0 < p \le 1$, by applying the Cauchy
condition (\ref{||f_j - f_l||_p < epsilon}) with $\epsilon = 1$ and $l
= L(1)$, and using the corresponding form of the triangle inequality.
This implies that $f \in \ell^p(E, {\bf R})$ or $\ell^p(E, {\bf C})$,
as appropriate, by the earlier remarks about pointwise convergent
sequences of functions on $E$.  We also get that
\begin{equation}
\label{||f_j - f||_p le epsilon}
        \|f_j - f\|_p \le \epsilon
\end{equation}
for every $j \ge L(\epsilon)$, by applying the earlier remarks to $f_j
- f_l$ as a sequence in $l$ for each $j$ and using the Cauchy
condition (\ref{||f_j - f_l||_p < epsilon}) again, so that
(\ref{lim_{j to infty} ||f_j - f||_p = 0}) holds, as desired.

\section{Bounded linear functionals}
\label{bounded linear functionals}

	Let $p$, $q$ be real numbers such that $1 < p, q < \infty$ and
\begin{equation}
	\frac{1}{p} + \frac{1}{q} = 1,
\end{equation}
so that $p$, $q$ are conjugate exponents.  If $f$, $g$ are real or
complex-valued functions on a nonempty set $E$ which are $p$, $q$-summable,
respectively, then their product $f \, g$ is a summable function on $E$, and
\begin{equation}
\label{sum_{x in E} |f(x)| |g(x)| le ||f||_p ||g||_q}
	\sum_{x \in E} |f(x)| \, |g(x)| \le \|f\|_p \, \|g\|_q.
\end{equation}
This is H\"older's inequality\index{H\"older's inequality} in the
present context, which follows easily from the version (\ref{Holder's
inequality for sums}) for finite sums.  We can also allow $p = 1$, $q
= \infty$ or $p = \infty$, $q = 1$, using bounded functions on $E$
when the corresponding exponent is infinite.

	A \emph{bounded linear functional}\index{bounded linear
functionals} on $\ell^p(E, {\bf R})$ or $\ell^p(E, {\bf C})$ is a
linear mapping $\lambda$ from this space into the real or complex
numbers, as appropriate, for which there is a nonnegative real number
$C$ such that
\begin{equation}
\label{|lambda(f)| le C ||f||_p}
	|\lambda(f)| \le C \, \|f\|_p
\end{equation}
for every $f \in \ell^p(E, {\bf R})$ or $\ell^p(E, {\bf C})$.  This
definition makes sense for every $p$ in the range $0 < p \le \infty$,
but let us focus first on the case where $p \ge 1$, and consider $p <
1$ afterwards.

        If $q$ is the conjugate exponent associated to $p \ge 1$, and
$g$ is a real or complex-valued function on $E$ in $\ell^q(E, {\bf
R})$ or $\ell^q(E, {\bf C})$, then
\begin{equation}
\label{lambda_g(f) = sum_{x in E} f(x) g(x)}
	\lambda_g(f) = \sum_{x \in E} f(x) \, g(x)
\end{equation}
defines a bounded linear functional on $\ell^p(E, {\bf R})$ of
$\ell^p(E, {\bf C})$, as appropriate, that satisfies (\ref{|lambda(f)|
le C ||f||_p}) with $C = \|g\|_q$, by H\"older's inequality.  One can
also check that $\|g\|_q$ is the smallest value of $C$ for which
(\ref{|lambda(f)| le C ||f||_p}) holds, which is analogous to the case
of finite sums discussed in Section \ref{dual spaces, norms}.

	Conversely, let $\lambda$ be a bounded linear functional on
$\ell^p(E, {\bf R})$ or $\ell^p(E, {\bf C})$, $1 \le p < \infty$, that
satisfies (\ref{|lambda(f)| le C ||f||_p}).  If $z \in E$, then let
$\delta_z$ be the function on $E$ defined by $\delta_z(y) = 1$ when $y
= z$, and $\delta_z(y) = 0$ otherwise.  Put
\begin{equation}
\label{g(z) = lambda(delta_z)}
	g(z) = \lambda(\delta_z)
\end{equation}
for every $z \in E$, so that
\begin{equation}
\label{lambda(f) = sum_{x in E} f(x) g(x)}
	\lambda(f) = \sum_{x \in E} f(x) \, g(x)
\end{equation}
when $f$ has finite support on $E$.  If $A$ is a finite subset of $E$,
then one can show that
\begin{equation}
\label{(sum_{x in A} |g(x)|^q)^{1/q} le C}
	\Big(\sum_{x \in A} |g(x)|^q \Big)^{1/q} \le C
\end{equation}
when $q < \infty$, and
\begin{equation}
\label{max_{x in A} |g(x)| le C}
	\max_{x \in A} |g(x)| \le C
\end{equation}
when $q = \infty$, using suitable choices of functions $f$ supported
on $A$.  This is also very similar to the discussion in Section
\ref{dual spaces, norms}.

        It follows that $g$ is $q$-summable when $q < \infty$, and
that $g$ is bounded when $q = \infty$, with
\begin{equation}
\label{||g||_q le C}
	\|g\|_q \le C
\end{equation}
in both cases.  Thus we can define $\lambda_g$ as a bounded linear
functional on $\ell^p(E, {\bf R})$ or $\ell^p(E, {\bf C})$, as
appropriate, as in (\ref{lambda_g(f) = sum_{x in E} f(x) g(x)}), and
\begin{equation}
\label{lambda(f) = lambda_g(f)}
        \lambda(f) = \lambda_g(f)
\end{equation}
when $f$ has finite support on $E$, as in (\ref{lambda(f) = sum_{x in
E} f(x) g(x)}).  To show that this holds for every $p$-summable
function $f$ on $E$, one can approximate $f$ by functions with finite
support on $E$, as in (\ref{||f - f_epsilon||_p < epsilon}) in the
previous section.  This also uses the fact that both $\lambda$ and
$\lambda_g$ are bounded linear functionals on $\ell^p(E, {\bf R})$ or
$\ell^p(E, {\bf C})$, as appropriate.

        If $p = \infty$, then it is better to consider bounded linear
functionals on $c_0(E, {\bf R})$, $c_0(E, {\bf C})$ instead of
$\ell^\infty(E, {\bf R})$, $\ell^\infty(E, {\bf C})$.  As before, a
bounded linear functional on $c_0(D, {\bf R})$ or $c_0(D, {\bf C})$ is
a linear mapping from this space to the real or complex numbers, as
appropriate, for which there is a nonnegative real number $C$ such that
\begin{equation}
\label{|lambda(f)| le C ||f||_infty}
	|\lambda(f)| \le C \, \|f\|_\infty
\end{equation}
for every function $f$ on $E$ that vanishes at infinity.  More
precisely, this is a bounded linear functional on $c_0(E, {\bf R})$ or
$c_0(E, {\bf C})$ with respect to the $\ell^\infty$ norm, which is the
natural norm in this case.

        If $g$ is a summable function on $D$, then we have seen that
(\ref{lambda_g(f) = sum_{x in E} f(x) g(x)}) defines a bounded linear
functional $\lambda_g$ on $\ell^\infty(E, {\bf R})$ or $\ell^\infty(E,
{\bf C})$, as appropriate, and that $\lambda_g$ satisfies
(\ref{|lambda(f)| le C ||f||_p}) with $p = \infty$ and $C = \|g\|_1$.
Hence the restriction of $\lambda_g$ to $c_0(E, {\bf R})$ or $c_0(E,
{\bf C})$, as appropriate, is a bounded linear functional that
satisfies (\ref{|lambda(f)| le C ||f||_infty}) with $C = \|g\|_1$.
One can check that $\|g\|_1$ is still the smallest value of $C$ for
which (\ref{|lambda(f)| le C ||f||_infty}) holds, even when we
restrict our attention to functions $f$ that vanish at infinity on
$E$, instead of considering all bounded functions on $E$.

        Conversely, suppose that $\lambda$ is a bounded linear
functional on $c_0(D, {\bf R})$ or $c_0(D, {\bf C})$ that satisfies
(\ref{|lambda(f)| le C ||f||_infty}).  Let $g$ be the function on $E$
defined by (\ref{g(z) = lambda(delta_z)}), as before, so that
(\ref{lambda(f) = sum_{x in E} f(x) g(x)}) holds when $f$ has finite
support on $E$.  If $A$ is a finite subset of $E$, then one can show
that
\begin{equation}
\label{sum_{x in A} |g(x)| le C}
	\sum_{x \in A} |g(x)| \le C,
\end{equation}
using suitable choices of functions $f$ supported on $A$.
This implies that $g$ is a summable function on $D$, with
\begin{equation}
\label{||g||_1 le C}
	\|g\|_1 \le C.
\end{equation}
To show that $\lambda(f) = \lambda_g(f)$ for every function $f$ that
vanishes at infinity on $E$, one can approximate such a function $f$
by functions with finite support on $E$ with respect to the
$\ell^\infty$ norm, as in (\ref{||f - f_epsilon||_infty < epsilon}) in
the previous section.

        If $g$ is a bounded real or complex-valued function on $E$,
then we have seen that (\ref{lambda_g(f) = sum_{x in E} f(x) g(x)})
defines a bounded linear functional $\lambda_g$ on $\ell^1(E, {\bf R})$
or $\ell^1(E, {\bf C})$, as appropriate, and that $\lambda_g$ satisfies
(\ref{|lambda(f)| le C ||f||_p}) with $p = 1$ and $C = \|g\|_\infty$.
If $0 < p < 1$, then we have also seen that every $p$-summable function
$f$ on $E$ is summable and satisfies
\begin{equation}
\label{||f||_1 le ||f||_p}
        \|f\|_1 \le \|f\|_p,
\end{equation}
as in (\ref{||f||_q le ||f||_p}) in the previous section, with $q = 1$.
It follows that the restriction of $\lambda_g$ to $\ell^p(E, {\bf R})$
or $\ell^p(E, {\bf C})$, as appropriate, is also a bounded linear functional
that satisfies (\ref{|lambda(f)| le C ||f||_p}) with $C = \|g\|_\infty$.
It is easy to see that $\|g\|_\infty$ is still the smallest value of
$C$ for which (\ref{|lambda(f)| le C ||f||_p}) holds, because
\begin{equation}
\label{||delta_z||p = 1}
        \|\delta_z\|p = 1
\end{equation}
for every $z \in E$.

        Conversely, suppose that $\lambda$ is a bounded linear
functional on $\ell^p(E, {\bf R})$ or $\ell^p(E, {\bf C})$ that
satisfies (\ref{|lambda(f)| le C ||f||_p}), where $0 < p < 1$.  As
usual, we can define a function $g$ on $E$ by (\ref{g(z) =
lambda(delta_z)}), so that (\ref{lambda(f) = sum_{x in E} f(x) g(x)})
holds when $f$ has finite support on $E$.  It is easy to see that $g$
is a bounded function on $E$, with
\begin{equation}
\label{||g||_infty le 1}
        \|g\|_\infty \le 1,
\end{equation}
by applying these conditions to $f = \delta_z$ for each $z \in E$.
One can then check that $\lambda(f) = \lambda_g(f)$ for every
$p$-summable function $f$ on $E$, by approximating $f$ by functions
with finite support on $E$, as in (\ref{||f - f_epsilon||_p <
epsilon}) in the previous section.

\section{Another convergence theorem}
\label{another convergence theorem}

        Let $E$ be a nonempty set, and let $1 \le p \le \infty$ be
given.  Also let $\{f_j\}_{j = 1}^\infty$ be a sequence of real or
complex-valued functions on $E$ that converges pointwise to another
function $f$ on $E$.  Suppose that $f_j$ is $p$-summable for each $j$
when $p < \infty$, and that $f_j$ is bounded for each $j$ when $p =
\infty$, with
\begin{equation}
        \|f_j\|_p \le C
\end{equation}
for some nonnegative real number $C$ and for every $j$ in both cases.
As in Section \ref{additional properties}, this implies that $f$ is
$p$-summable when $p < \infty$, and that $f$ is bounded when $p =
\infty$, with $\|f\|_p \le C$ in both cases.

        Let $1 \le q \le \infty$ be the exponent conjugate to $p$, so
that $1/p + 1/q = 1$, and let $g$ be a real or complex-valued function
on $E$ which is $q$-summable when $q < \infty$, and which vanishes at
infinity on $E$ when $q = \infty$.  In particular, $g$ is bounded when
$q = \infty$.  As in the previous section, $f_j \, g$ is summable on
$E$ for each $j$, as is $f \, g$.  Under these conditions, we would
like to show that
\begin{equation}
\label{lim_{j to infty} sum_{x in E} f_j(x) g(x) = sum_{x in E} f(x) g(x)}
        \lim_{j \to \infty} \sum_{x \in E} f_j(x) \, g(x)
                                   = \sum_{x \in E} f(x) \, g(x).
\end{equation}
The proof is similar to that of the dominated convergence theorem in
Section \ref{convergence theorems}, and in fact one can derive this
from the dominated convergence theorem when $p = \infty$.

        Equivalently, we would like to show that
\begin{equation}
\label{lim_{j to infty} sum_{x in E} (f_j(x) - f(x)) g(x) = 0}
        \lim_{j \to \infty} \sum_{x \in E} (f_j(x) - f(x)) \, g(x) = 0.
\end{equation}
Let $\epsilon > 0$ be given, and let $A$ be a finite subset of $E$ such that
\begin{equation}
\label{(sum_{x in E backslash A} |g(x)|^q)^{1/q} < epsilon}
        \Big(\sum_{x \in E \backslash A} |g(x)|^q\Big)^{1/q} < \epsilon
\end{equation}
when $q < \infty$, and $|g(x)| < \epsilon$ for every $x \in E
\backslash A$ when $q = \infty$.  Using this, it is easy to see that
\begin{equation}
\label{|sum_{x in E backslash A} (f_j(x) - f(x)) g(x)| le ... le 2 C epsilon}
        \biggl|\sum_{x \in E \backslash A} (f_j(x) - f(x)) \, g(x)\biggr|
                \le \epsilon \, \|f_j - f\|_p \le 2 \, C \, \epsilon,
\end{equation}
for each $j$, by H\"older's inequality.  We also have that
\begin{equation}
        \biggl|\sum_{x \in A} (f_j(x) - f(x)) \, g(x)\biggr| < \epsilon
\end{equation}
for all sufficiently large $j$, because $\{f_j\}_{j = 1}^\infty$
converges to $f$ pointwise on $E$, and because $A$ has only finitely
many elements.  The desired conclusion (\ref{lim_{j to infty} sum_{x
in E} (f_j(x) - f(x)) g(x) = 0}) follows by combining these two
statements.

        Note that (\ref{lim_{j to infty} sum_{x in E} (f_j(x) - f(x))
g(x) = 0}) follows directly from H\"older's inequality if we ask that
\begin{equation}
        \lim_{j \to \infty} \|f_j - f\|_p = 0.
\end{equation}
In this case, it would also have been sufficient to ask that $g$ be
bounded on $E$ when $q = \infty$.

\chapter{Banach and Hilbert spaces}
\label{banach, hilbert spaces}

\section{Basic concepts}
\label{basic concepts}

	Let $V$ be a vector space over the real or complex numbers,
and let $\|v\|$ be a norm on $V$.  In this chapter, $V$ is allowed to
be infinite-dimensional, but the definition of a norm on $V$ is the
same as in the finite-dimensional case.

        As usual, a sequence $\{v_j\}_{j=1}^\infty$ of elements of $V$
is said to \emph{converge}\index{convergent sequences} to $v \in V$ if
for each $\epsilon > 0$ there is an $L \ge 1$ such that
\begin{equation}
	\|v_j - v\| < \epsilon
\end{equation}
for every $j \ge L$.  In this case, we put
\begin{equation}
        \lim_{j \to \infty} v_j = v,
\end{equation}
and call $v$ the limit of the sequence $\{v_j\}_{j = 1}^\infty$.  It
is easy to see that the limit of a convergent sequence is unique when
it exists.

        Similarly, a sequence $\{v_j\}_{j=1}^\infty$ in $V$ is a
\emph{Cauchy sequence}\index{Cauchy sequences} if for each $\epsilon >
0$ there is an $L \ge 1$ such that
\begin{equation}
	\|v_j - v_l\| < \epsilon
\end{equation}
for every $j, l \ge L$.  Note that every convergent sequence is
automatically a Cauchy sequence, by a simple argument.

        If every Cauchy sequence of elements of $V$ converges to an
element of $V$, then we say that $V$ is \emph{complete}.  This is equivalent
to the compleness of $V$ as a metric space, with respect to the metric
\begin{equation}
\label{d(v, w) = ||v - w||}
        d(v, w) = \|v - w\|
\end{equation}
corresponding to the norm $\|v\|$ on $V$.  If $V$ is complete in this
sense, then we say that $V$ is a \emph{Banach space}.\index{Banach
spaces} If $\langle v, w \rangle$ is an inner product on $V$, and if
$V$ is complete with respect to the associated norm
\begin{equation}
        \|v\| = \langle v, v \rangle^{1/2},
\end{equation}
then we say that $V$ is a \emph{Hilbert space}.\index{Hilbert spaces}

        Of course, the real and complex numbers are complete with
respect to their standard norms, as in Section \ref{real, complex
numbers}.  If $V = {\bf R}^n$ or ${\bf C}^n$ for some positive integer
$n$, and $\|v\|$ is one of the norms $\|v\|_p$ defined in Section
\ref{definitions, examples}, $1 \le p \le \infty$, then it is easy to
see that $V$ is complete, by reducing to the $n = 1$ case.  If $V =
{\bf R}^n$ or ${\bf C}^n$ equipped with any norm $\|v\|$, then one can
also check that $V$ is complete with respect to $\|v\|$, by reducing
to the case of the standard norm on $V$ using the remarks at the end
of Section \ref{definitions, examples}.  This implies that any
finite-dimensional vector space $V$ over the real or complex numbers
is complete with respect to any norm $\|v\|$ on $V$, because $V$ is
isomorphic to ${\bf R}^n$ or ${\bf C}^n$ for some positive integer
$n$, or $V = \{0\}$.  If $E$ is a nonempty set and $V = \ell^p(E, {\bf
R})$ or $\ell^p(E, {\bf C})$ for some $1 \le p \le \infty$, then $V$
is complete with respect to the corresponding norm $\|f\|_p$,
as in Section \ref{additional properties}.

        Let $V$ be a real or complex vector space with a norm $\|v\|$.
A subset $W$ of $V$ is said to be a \emph{closed set}\index{closed
sets} in $V$ if for every sequence $\{w_j\}_{j = 1}^\infty$ of
elements of $W$ that converges to an element $w$ of $V$, we have that
$w \in W$.  This is equivalent to other standard definitions of closed
sets in metric spaces, in terms of a set containing all of its limit
points, or being a closed set when the complement is an open set.  If
$V$ is complete with respect to $\|v\|$, and $W$ is a closed linear
subspace of $V$, then it follows that $W$ is complete with respect to
the restriction of the norm $\|v\|$ to $v \in W$.  This is because a
Cauchy sequence in $W$ is also a Cauchy sequence in $V$ in this
situation, which converges to an element of $V$ when $V$ is complete,
and the limit is in $W$ when $W$ is a closed set in $V$.

        In particular, if $E$ is a nonempty set, then $c_0(E, {\bf
R})$ and $c_0(E, {\bf C})$ are closed linear subspaces of
$\ell^\infty(E, {\bf R})$ and $\ell^\infty(E, {\bf C})$, as in Section
\ref{additional properties}.  Hence $c_0(E, {\bf R})$ and $c_0(E, {\bf
R})$ are complete with respect to the $\ell^\infty$ norm, as in the
preceding paragraph, because $\ell^\infty(E, {\bf R})$ and
$\ell^\infty(E, {\bf C})$ are complete.

        Let $V$ be the vector space of continuous real or
complex-valued functions on the closed unit interval $[0, 1]$, with
respect to pointwise addition and scalar multiplication.  Remember
that continuous functions on $[0, 1]$ are automatically bounded,
because $[0, 1]$ is compact.  Thus
\begin{equation}
\label{||f|| = sup_{0 le x le 1} |f(x)|}
        \|f\| = \sup_{0 \le x \le 1} |f(x)|
\end{equation}
defines a norm on $V$, known as the \emph{supremum
norm}.\index{supremum norm} It is well known that $V$ is complete with
respect to this norm, because of the fact that the limit of a
uniformly-convergent sequence of continuous functions is also
continuous.  If $1 \le p < \infty$, then
\begin{equation}
\label{||f||_p = (int_0^1 |f(x)|^p dx)^{1/p}}
        \|f\|_p = \Big(\int_0^1 |f(x)|^p \, dx\Big)^{1/p}
\end{equation}
also defines a norm on $V$, because of Minkowski's inequality in
Section \ref{related inequalities}.  However, it is well known that
$V$ is not complete with respect to this norm for any $p < \infty$.
To get a complete space, one can use Lebesgue integrals.  Similarly,
the vector space of bounded continuous real or complex-valued
functions on any topological space is complete with respect to the
corresponding supremum norm.  There are also $L^p$ spaces associated
to any measure space, which are Banach spaces when $p \ge 1$, and
Hilbert spaces when $p = 2$.  Note that the $\ell^p$ spaces discussed
in the previous chapter may be considered as $L^p$ spaces with respect
to counting measure on a set $E$.

\section{Sequences and series}
\label{sequences, series}

        Let $V$ be a real or complex vector space with a norm $\|v\|$.
If $\{v_j\}_{j=1}^\infty$, $\{w_j\}_{j=1}^\infty$ are sequences in $V$
that converge to $v, w \in V$, respectively, then
\begin{equation}
	\lim_{j \to \infty} (v_j + w_j) = v + w.
\end{equation}
This can be shown in essentially the same way as for sequences of real
or complex numbers.  Similarly, if $\{t_j\}_{j=1}^\infty$ is a
sequence of real or complex numbers, as appropriate, that converges to
the real or complex number $t$, and if $\{v_j\}_{j=1}^\infty$ is a
sequence of vectors in $V$ that converges to $v \in V$, then
\begin{equation}
	\lim_{j \to \infty} t_j \, v_j = t \, v.
\end{equation}
Equivalently, this means that addition and scalar multiplication are
continuous on $V$ with respect to the metric associated to the norm.

        As in (\ref{| ||v|| - ||w|| | le ||v - w||}) in Section
\ref{definitions, examples}, one can check that
\begin{equation}
\label{| ||v|| - ||w|| | le ||v - w||, 2}
        \bigl| \|v\| - \|w\| \bigr| \le \|v - w\|
\end{equation}
for every $v, w \in V$, using the triangle inequality.  If
$\{v_j\}_{j=1}^\infty$ is a sequence in $V$ that converges to $v \in
V$, then it follows that
\begin{equation}
	\lim_{j \to \infty} \|v_j\| = \|v\|,
\end{equation}
as a sequence of real numbers.  This is the same as saying that
$\|v\|$ is a continuous real-valued function on $V$ with respect to
the metric associated to the norm.

        Let $\sum_{j = 1}^\infty a_j$ be an infinite series whose
terms $a_j$ are elements of $V$.  As usual, we say that $\sum_{j =
1}^\infty a_j$ \emph{converges}\index{convergent series} in $V$ if the
corresponding sequence of partial sums $\sum_{j = 1}^n a_j$ converges
in $V$, in which case we put
\begin{equation}
        \sum_{j = 1}^\infty a_j = \lim_{n \to \infty} \sum_{j = 1}^n a_j.
\end{equation}
If $\sum_{j = 1}^\infty a_j$ and $\sum_{j = 1}^\infty b_j$ are
convergent series with terms in $V$, then it is easy to see that
$\sum_{j = 1}^\infty (a_j + b_j)$ also converges, and that
\begin{equation}
        \sum_{j = 1}^\infty (a_j + b_j)
                 = \sum_{j = 1}^\infty a_j + \sum_{j = 1}^\infty b_j,
\end{equation}
because of the corresponding fact about sums of convergent sequences
mentioned earlier.  Similarly, if $\sum_{j = 1}^\infty a_j$ is a
convergent series with terms in $V$, and if $t$ is a real or complex
number, as appropriate, then $\sum_{j = 1}^\infty t \, a_j$ also
converges, and
\begin{equation}
        \sum_{j = 1}^\infty t \, a_j = t \, \sum_{j = 1}^\infty a_j.
\end{equation}

        An infinite series $\sum_{j = 1}^\infty a_j$ with terms in $V$ is
said to converge \emph{absolutely}\index{absolutely convergent series} if
\begin{equation}
\label{sum_{j = 1}^infty ||a_j||}
        \sum_{j = 1}^\infty \|a_j\|
\end{equation}
converges as an infinite series of nonnegative real numbers.  In this
case, one can show that the partial sums of $\sum_{j = 1}^\infty a_j$
form a Cauchy sequence, because
\begin{equation}
        \biggl\|\sum_{j = l}^n a_j\biggr\| \le \sum_{j = l}^n \|a_j\|
\end{equation}
for every $n \ge l \ge 1$.  If $V$ is complete, then it follows that
$\sum_{j = 1}^\infty a_j$ converges in $V$.  We also get that
\begin{equation}
\label{||sum_{j = 1}^infty a_j|| le sum_{j = 1}^infty ||a_j||}
 \biggl\|\sum_{j = 1}^\infty a_j\biggr\| \le \sum_{j = 1}^\infty \|a_j\|.
\end{equation}

        Conversely, suppose that $\{v_j\}_{j = 1}^\infty$ is a Cauchy
sequence of elements of $V$.  It is easy to see that there is a
subsequence $\{v_{j_l}\}_{l = 1}^\infty$ of $\{v_j\}_{j = 1}^\infty$
such that
\begin{equation}
        \|v_{j_l} - v_{j_{l + 1}}\| < 2^{-l}
\end{equation}
for each $l$, so that
\begin{equation}
\label{sum_{l = 1}^infty (v_{j_l} - v_{j_{l + 1}})}
        \sum_{l = 1}^\infty (v_{j_l} - v_{j_{l + 1}})
\end{equation}
converges absolutely.  Of course,
\begin{equation}
        \sum_{l = 1}^n (v_{j_l} - v_{j_{l + 1}}) = v_{j_1} - v_{j_{n + 1}}
\end{equation}
for each $n$, which implies that the series (\ref{sum_{l = 1}^infty
(v_{j_l} - v_{j_{l + 1}})}) converges in $V$ if and only if
$\{v_{j_l}\}_{l = 1}^\infty$ converges in $V$.  If $\{v_{j_l}\}_{l =
1}^\infty$ converges in $V$, then one can check that $\{v_j\}_{j =
1}^\infty$ also converges to the same element of $V$, because
$\{v_j\}_{j = 1}^\infty$ is a Cauchy sequence.  If every absolutely
convergent series in $V$ converges, then it follows that every Cauchy
sequence in $V$ converges, which is to say that $V$ is complete.

        Suppose now that the norm $\|v\|$ on $V$ is associated to an
inner product $\langle v, w \rangle$ on $V$ in the usual way.
Let $\sum_{j = 1}^\infty a_j$ be an infinite series whose terms are
pairwise-orthogonal vectors in $V$,\index{orthogonal vectors} in the sense that
\begin{equation}
\label{langle a_j, a_k rangle = 0}
        \langle a_j, a_k \rangle = 0
\end{equation}
when $j \ne k$.  In this case,
\begin{equation}
        \biggl\|\sum_{j = 1}^n a_j\biggr\|^2 = \sum_{j = 1}^n \|a_j\|^2
\end{equation}
for each $n$.  If $\sum_{j = 1}^\infty a_j$ converges in $V$, then
$\sum_{j = 1}^\infty \|a_j\|^2$ converges in ${\bf R}$, and
\begin{equation}
 \biggl\|\sum_{j = 1}^\infty a_j \biggr\|^2 = \sum_{j = 1}^\infty \|a_j\|^2.
\end{equation}

        Of course, the orthogonality condition (\ref{langle a_j, a_k
rangle = 0}) implies that
\begin{equation}
\label{||sum_{j = l}^n a_j||^2 = sum_{j = l}^n ||a_j||^2}
        \biggl\|\sum_{j = l}^n a_j\biggr\|^2 = \sum_{j = l}^n \|a_j\|^2
\end{equation}
for every $n \ge l \ge 1$.  If $\sum_{j = 1}^\infty \|a_j\|^2$
converges, then it follows that the partial sums of $\sum_{j =
1}^\infty a_j$ form a Cauchy sequence in $V$.  If $V$ is complete,
then $\sum_{j = 1}^\infty a_j$ converges in $V$ under these
conditions.

\section{Minimizing distances}
\label{minimizing distances, 2}

        Let $V$ be a real or complex vector space with an inner
product $\langle v, w \rangle$, and let $\|v\|$ be the corresponding
norm on $V$.  Also let $E$ be a nonempty subset of $V$, let $v$ be an
element of $V$, and put
\begin{equation}
        \rho = \inf \{\|v - w\| : w \in E\}.
\end{equation}
Thus for each positive integer $j$ there is a $w_j \in E$ such that
\begin{equation}
\label{||v - w_j|| < rho + frac{1}{j}}
        \|v - w_j\| < \rho + \frac{1}{j}.
\end{equation}
Note that
\begin{equation}
        \biggl\|\frac{x + y}{2}\biggr\|^2 + \biggl\|\frac{x - y}{2}\biggr\|^2
          = \frac{1}{2} \, (\|x\|^2 + \|y\|^2)
\end{equation}
for every $x, y \in V$, by applying the parallelogram law
(\ref{parallelogram law}) to $x/2$, $y/2$.  If we take $x = v - w_j$
and $y = v - w_l$, then we get that
\begin{equation}
 \biggl\|v - \Big(\frac{w_j + w_l}{2}\Big)\biggr\|^2
                                           + \frac{1}{4} \, \|w_j - w_l\|^2
           = \frac{1}{2} \, (\|v - w_j\|^2 + \|v - w_l\|^2).
\end{equation}
Combining this with (\ref{||v - w_j|| < rho + frac{1}{j}}) gives
\begin{equation}
 \quad  \biggl\|v - \Big(\frac{w_j + w_l}{2}\Big)\biggr\|^2
                                            + \frac{1}{4} \, \|w_j - w_l\|^2
           < \rho^2 + \rho \, \Big(\frac{1}{j} + \frac{1}{l}\Big)
                     + \frac{1}{2} \, \Big(\frac{1}{j^2} + \frac{1}{l^2}\Big).
\end{equation}

        Suppose now that $E$ is a convex set in $V$.  This implies that
\begin{equation}
        \frac{w_j + w_l}{2} \in E
\end{equation}
for each $j$, $l$, and hence that
\begin{equation}
        \biggl\|v - \frac{w_j + w_l}{2}\biggr\| \ge \rho.
\end{equation}
In this case, we get that
\begin{equation}
\label{||w_j - w_l||^2 / 4 < rho (1/j + 1/l) + 1/2 (1/j^2 + 1/l^2)}
 \frac{1}{4} \, \|w_j - w_l\|^2 < \rho \, \Big(\frac{1}{j} + \frac{1}{l}\Big)
                     + \frac{1}{2} \, \Big(\frac{1}{j^2} + \frac{1}{l^2}\Big)
\end{equation}
for each $j$, $l$.  Thus $\{w_j\}_{j = 1}^\infty$ is a Cauchy sequence
in $V$ under these conditions.

        If $V$ is complete, then it follows that $\{w_j\}_{j =
1}^\infty$ converges to an element $w$ of $V$.  If $E$ is a closed set
in $V$, then $w \in E$.  Moreover,
\begin{equation}
        \|v - w\| = \rho,
\end{equation}
so that $w$ minimizes the distance to $v$ among elements of $E$.

        This argument also works in a large class of Banach spaces.
More precisely, a norm $\|v\|$ on $V$ is said to be \emph{uniformly
convex}\index{uniform convexity} if for each $\epsilon > 0$
there is a $\delta(\epsilon) > 0$ such that $\delta(\epsilon) < 1$ and
\begin{equation}
\label{||frac{u + z}{2}|| le 1 - delta(epsilon)}
        \biggl\|\frac{u + z}{2}\biggr\| \le 1 - \delta(\epsilon)
\end{equation}
for every $u, z \in V$ such that $\|u\| = \|z\| = 1$ and $\|u - z\| >
\epsilon$.  Equivalently, this means that
\begin{equation}
\label{||u - z|| le epsilon}
        \|u - z\| \le \epsilon
\end{equation}
for every $u, z \in V$ such that $\|u\| = \|z\| = 1$ and $\|(u +
z)/2\| > 1 - \delta(\epsilon)$.  If $\|v\|$ is associated to an inner
product on $V$, then it is easy to see that $\|v\|$ is uniformly
convex, using the parallelogram law.  It is well known that the $L^p$
norm is uniformly convex when $1 < p < \infty$, as a consequence of
famous inequalities of Clarkson.  If $\|v\|$ is uniformly convex, then
one can modify the previous arguments to show that the minimum of the
distance from a point $v \in V$ to a nonempty closed convex set $E
\subseteq V$ is attained when $V$ is complete.  As before, the main
step is to show that a minimizing sequence $\{w_j\}_{j = 1}^\infty$
is a Cauchy sequence when $\|v\|$ is uniformly convex.

\section{Orthogonal projections}
\label{orthogonal projections}

        Let $V$ be a real or complex vector space with an inner
product $\langle v, w \rangle$, and let $\|v\|$ be the corresponding
norm on $V$.  Suppose that $V$ is complete, so that $V$ is a Hilbert
space, and let $W$ be a closed linear subspace of $V$.  If $v$ is any
element of $V$, then there is a $w \in W$ whose distance to $v$ is
minimal among elements of $W$, as in the previous section.
Equivalently,
\begin{equation}
        \|v - w\| \le \|v - w + z\|
\end{equation}
for every $z \in W$, because $W$ is a linear subspace of $V$.
Using the inner product, we get that
\begin{eqnarray}
        \|v - w\|^2 & \le & \|v - w + z\|^2                          \\
  & = & \|v - w\|^2 + \langle v - w, z \rangle + \langle z, v - w \rangle
                                                + \|z\|^2 \nonumber \\
  & = & \|v - w\|^2 + 2 \, \re \langle v - w, z \rangle + \|z\|^2. \nonumber
\end{eqnarray}
More precisely, it is not necessary to take the real part of $\langle
v - w, z \rangle$ in the last step when $V$ is a real vector space,
but this is needed when $V$ is complex.  Of course, this inequality
reduces to
\begin{equation}
\label{0 le 2 re langle v - w, z rangle + ||z||^2}
        0 \le 2 \, \re \langle v - w, z \rangle + \|z\|^2,
\end{equation}
by subtracting $\|v - w\|^2$ from both sides.

        Let $t$ be a real number, and put
\begin{equation}
        f(t) = 2 \, \re \langle v - w, t \, z \rangle + \|t \, z\|^2
              = 2 \, t \, \re \langle v - w, z \rangle + t^2 \, \|z\|^2.
\end{equation}
If $z \in W$, then $t \, z \in W$, and hence $f(t) \ge 0$, by (\ref{0
le 2 re langle v - w, z rangle + ||z||^2}).  Thus the minimum of
$f(t)$ is attained at $t = 0$, which implies that the derivative of
$f(t)$ at $t = 0$ is also equal to $0$.  This shows that
\begin{equation}
\label{re langle v - w, z rangle = 0}
        \re \langle v - w, z \rangle = 0
\end{equation}
for every $z \in W$, which is the same as saying that
\begin{equation}
\label{langle v - w, z rangle = 0}
        \langle v - w, z \rangle = 0
\end{equation}
for every $z \in W$ in the real case.  In the complex case, one can
get (\ref{langle v - w, z rangle = 0}) by applying (\ref{re langle v -
w, z rangle = 0}) to $z$ and to $i \, z$.

        Conversely, (\ref{langle v - w, z rangle = 0}) implies that
\begin{equation}
\label{||v - w + z||^2 = ||v - w||^2 + ||z||^2}
        \|v - w + z\|^2 = \|v - w\|^2 + \|z\|^2
\end{equation}
for every $z \in W$, and hence that $w$ minimizes the distance to $v$
among elements of $W$.  As in Section \ref{inner product spaces}, $w
\in W$ is uniquely determined by the condition that (\ref{langle v -
w, z rangle = 0}) holds for every $z \in W$.  Put $P_W(v) = w$, which
is the \emph{orthogonal projection}\index{orthogonal projections} of
$v$ onto $W$.  Note that
\begin{equation}
\label{||v||^2 = ||v - P_W(v)||^2 + ||P_W(v)||^2}
        \|v\|^2 = \|v - P_W(v)\|^2 + \|P_W(v)\|^2,
\end{equation}
which is the same as (\ref{||v - w + z||^2 = ||v - w||^2 + ||z||^2})
with $z = w$.  In particular,
\begin{equation}
        \|P_W(v)\| \le \|v\|
\end{equation}
for every $v \in V$.

        If $v_1, v_2 \in V$, then $P_W(v_1) + P_W(v_2) \in W$ and
\begin{equation}
        (v_1 + v_2) - (P_W(v_1) + P_W(v_2))
            = (v_1 - P_W(v_1)) + (v_2 - P_W(v_2))
\end{equation}
is orthogonal to every element of $W$, by the corresponding properties
of $P_W(v_1)$ and $P_W(v_2)$.  This implies that
\begin{equation}
\label{P_W(v_1 + v_2) = P_W(v_1) + P_W(v_2)}
        P_W(v_1 + v_2) = P_W(v_1) + P_W(v_2),
\end{equation}
because $P_W(v_1 + v_2)$ is characterized by these conditions, as in
the preceding paragraph.  Similarly,
\begin{equation}
        P_W(t \, v) = t \, P_W(v)
\end{equation}
for every $v \in V$ and $t \in {\bf R}$ or ${\bf C}$, as appropriate,
because $t \, P_W(v) \in W$ and
\begin{equation}
        t \, v - t \, P_W(v) = t \, (v - P_W(v))
\end{equation}
is orthogonal to every element of $W$, by the corresponding properties
of $P_W(v)$.  Thus $P_W(v)$ is a linear mapping from $V$ into $W$.  Of
course, $P_W(v) = v$ when $v \in W$.

\section{Orthonormal sequences}
\label{orthonormal sequences}

        Let $V$ be a real or complex vector space with an inner product
$\langle v, w \rangle$ again, and let $\|v\|$ be the corresponding norm
on $V$.  Suppose that $e_1, e_2, e_3, \ldots$ is an infinite sequence
of orthonormal vectors\index{orthonormal vectors} in $V$, so that
\begin{equation}
\label{langle e_j, e_k rangle = 0}
        \langle e_j, e_k \rangle = 0
\end{equation}
when $j \ne k$, and $\|e_j\| = 1$ for each $j$.  Note that any
sequence of vectors in $V$ can be modified to get an orthonormal
sequence with the same linear span, using the Gram--Schmidt process.
This was already used in Section \ref{inner product spaces} to show
that every finite-dimensional inner product space has an orthonormal
basis.

        Put
\begin{equation}
        P_n(v) = \sum_{j = 1}^n \langle v, e_j \rangle \, e_j
\end{equation}
for each $v \in V$ and positive integer $n$.  This is the same as the
orthogonal projection of $V$ onto the linear subspace $W_n$ spanned by
$e_1, \ldots, e_n$ in $V$, as in Section \ref{inner product spaces}.
Remember that
\begin{equation}
\label{||v||^2 = ||P_n(v)||^2 + ||v - P_n(v)||^2 = ...}
        \|v\|^2 = \|P_n(v)\|^2 + \|v - P_n(v)\|^2
          = \sum_{j = 1}^n |\langle v, e_j \rangle|^2 + \|v - P_n(v)\|^2,
\end{equation}
as in (\ref{||w||^2 = ||P(w)||^2 + ||w - P(w)||^2 = ...}).  In
particular,
\begin{equation}
\label{sum_{j = 1}^n |langle v, e_j rangle|^2 le ||v||^2}
        \sum_{j = 1}^n |\langle v, e_j \rangle|^2 \le \|v\|^2
\end{equation}
for every $v \in V$ and $n \ge 1$.  Remember also that $P_n(v) \in
W_n$ minimizes the distance to $v$ among elements of $W_n$, as in
Section \ref{minimizing distances}.

        Note that $\bigcup_{n = 1}^\infty W_n$ is a linear subspace of
$V$, because $W_n$ is a linear subspace of $V$ for each $n$, and
because $W_n \subseteq W_{n + 1}$ for each $n$, by construction.
Let $W$ be the closure of $\bigcup_{n = 1}^\infty W_n$ in $V$, which
is the set of $v \in V$ with the property that for each $\epsilon > 0$
there is a $w \in \bigcup_{n = 1}^\infty W_n$ such that
\begin{equation}
\label{||v - w|| < epsilon}
        \|v - w\| < \epsilon.
\end{equation}
Thus $\bigcup_{n = 1}^\infty W_n \subseteq W$ automatically, and one
can check that $W$ is a closed linear subspace of $V$.

        Equivalently, $v \in V$ is an element of $W$ if and only if
\begin{equation}
\label{lim_{n to infty} ||v - P_n(v)|| = 0}
        \lim_{n \to \infty} \|v - P_n(v)\| = 0.
\end{equation}
More precisely, if $v$ satisfies (\ref{lim_{n to infty} ||v - P_n(v)||
= 0}), then it is easy to see that $v \in W$, because $P_n(v) \in W_n$
for each $n$.  Conversely, suppose that $v \in W$, and let $\epsilon >
0$ be given.  By definition of $W$, there is a positive integer $k$
and a $w \in W_k$ such that (\ref{||v - w|| < epsilon}) holds.  This
implies that
\begin{equation}
\label{||v - P_n(v)|| le ||v - w|| < epsilon}
        \|v - P_n(v)\| \le \|v - w\| < \epsilon
\end{equation}
for every $n \ge k$, as desired, because $w_k \in W_n$ for each $n \ge
k$, and because $P_n(v)$ minimizes the distance to $v$ among elements
of $W_n$, as before.

        We would like to put
\begin{equation}
\label{P(v) = sum_j langle v, e_j rangle e_j = lim_{n to infty} P_n(v)}
        P(v) = \sum_{j = 1}^\infty \langle v, e_j \rangle \, e_j
              = \lim_{n \to \infty} P_n(v)
\end{equation}
for every $v \in V$, but we need to be careful about the existence of
the limit.  If $v \in W$, then (\ref{lim_{n to infty} ||v - P_n(v)|| =
0}) implies that $\{P_n(v)\}_{n = 1}^\infty$ converges to $v$, so that
the definition of $P(v)$ makes sense and $P(v) = v$.  Otherwise, if
$v$ is any element of $V$, then
\begin{equation}
        \sum_{j = 1}^\infty |\langle v, e_j \rangle|^2
\end{equation}
converges and is less than or equal to $\|v\|^2$, because of
(\ref{sum_{j = 1}^n |langle v, e_j rangle|^2 le ||v||^2}).  If $V$ is
complete, then it follows that the series in (\ref{P(v) = sum_j langle
v, e_j rangle e_j = lim_{n to infty} P_n(v)}) converges, as in Section
\ref{sequences, series}.  In this case, it is easy to see that $P(v)
\in W$ for every $v \in V$, because $P_n(v) \in W_n$ for each $n$.
One can also check that $P$ is a linear mapping from $V$ into $W$
under these conditions.  By construction, we also have that
\begin{equation}
        \|P(v)\|^2 = \sum_{j = 1}^\infty |\langle v, e_j \rangle|^2 \le \|v\|^2
\end{equation}
for every $v \in V$.

        Let us suppose from now on in this section that $V$ is
complete, and thus a Hilbert space.  Observe that
\begin{equation}
\label{langle P(v), e_l rangle = ... = langle v, e_l rangle}
        \langle P(v), e_l \rangle
 = \lim_{n \to \infty} \langle P_n(v), e_l \rangle = \langle v, e_l \rangle
\end{equation}
for every $v \in V$ and $l \ge 1$.  This uses the fact that
\begin{equation}
        \langle P_n(v), e_l \rangle = \langle v, e_l \rangle
\end{equation}
for each $n \ge l$, because of the orthonormality of the $e_j$'s.
This also implicitly uses the Cauchy--Schwarz inequality, in order to
take the limit outside of the inner product, which is basically the
same as the continuity of the inner product with respect to the
associated norm.  It follows that $v - P(v)$ is orthogonal to $e_l$
for each $l$, which implies that $v - P(v)$ is orthogonal to every
element of $\bigcup_{n = 1}^\infty W_n$, because of the linearity
properties of the inner product.  Using continuity of the inner
product again, we get that $v - P(v)$ is orthogonal to every element
of the closure $W$ of $\bigcup_{n = 1}^\infty W_n$.  As in Section
\ref{inner product spaces}, $P(v)$ is uniquely determined by the
conditions that $P(v) \in W$ and $v - P(v)$ is orthogonal to every
element of $W$, and hence $P$ is the same as the orthogonal projection
$P_W$ of $V$ onto $W$, as in the previous section.\index{orthogonal
projections}

        If $\{a_j\}_{j = 1}^\infty$ is any sequence of real or complex
numbers, as appropriate, such that $\sum_{j = 1}^\infty |a_j|^2$
converges in ${\bf R}$, then the same arguments show that
\begin{equation}
\label{sum_{j = 1}^infty a_j e_j}
        \sum_{j = 1}^\infty a_j \, e_j
\end{equation}
converges in $V$ to an element of $W$, and that
\begin{equation}
\label{||sum_{j = 1}^infty a_j e_j||^2 = sum_{j = 1}^infty |a_j|^2}
\biggl\|\sum_{j = 1}^\infty a_j \, e_j\biggr\|^2 = \sum_{j = 1}^\infty |a_j|^2.
\end{equation}
If $\{b_j\}_{j = 1}^\infty$ is another sequence of real or complex
numbers, as appropriate, such that $\sum_{j = 1}^\infty |b_j|^2$
converges, then one can check that
\begin{equation}
\label{langle sum_{j = 1}^infty a_j e_j, sum_{k = 1}^infty b_k e_k rangle, 1}
        \Big\langle \sum_{j = 1}^\infty a_j \, e_j,
                             \sum_{k = 1}^\infty b_k \, e_k \Big\rangle
                                            = \sum_{j = 1}^\infty a_j \, b_j
\end{equation}
in the real case, and
\begin{equation}
\label{langle sum_{j = 1}^infty a_j e_j, sum_{k = 1}^infty b_k e_k rangle, 2}
        \Big\langle \sum_{j = 1}^\infty a_j \, e_j,
                             \sum_{k = 1}^\infty b_k \, e_k \Big\rangle
                                  = \sum_{j = 1}^\infty a_j \, \overline{b_j}
\end{equation}
in the complex case, using the orthonormality of the $e_j$'s and the
continuity properties of the inner product, as before.  Note that the
infinite series on the right sides of (\ref{langle sum_{j = 1}^infty
a_j e_j, sum_{k = 1}^infty b_k e_k rangle, 1}) and (\ref{langle sum_{j
= 1}^infty a_j e_j, sum_{k = 1}^infty b_k e_k rangle, 2}) are absolutely
convergent under these conditions, as in Section \ref{ell^p spaces}.
If $V = W$, then the $e_j$'s are said to form an \emph{orthonormal
basis}\index{orthonormal bases} of $V$.  In this case, we get a natural
isomorphism between $V$ and $\ell^2({\bf Z}_+, {\bf R})$ or
$\ell^2({\bf Z}_+, {\bf C})$, as appropriate, associated to this
orthonormal basis for $V$.

\section{Bounded linear functionals}
\label{bounded linear functionals, 2}

        Let $V$ be a real or complex vector space with a norm $\|v\|$.
A linear functional $\lambda$ on $V$ is said to be \emph{bounded}\index{bounded
linear functionals} if there is a nonnegative real number $C$ such that
\begin{equation}
\label{|lambda(v)| le C ||v||}
        |\lambda(v)| \le C \, \|v\|
\end{equation}
for every $v \in V$.  If $V$ has finite dimension, then every linear
functional on $V$ is bounded, as in Section \ref{dual spaces, norms}.
If $\lambda$ is a bounded linear functional on $V$, then it is easy to
see that $\lambda$ is continuous on $V$ with respect to the metric
associated to the norm.  In particular, if $\{v_j\}_{j = 1}^\infty$ is
a sequence of vectors in $V$ that converges to another vector $v \in
V$, as in Section \ref{basic concepts}, then it is easy to see that
\begin{equation}
        \lim_{j \to \infty} \lambda(v_j) = \lambda(v)
\end{equation}
in ${\bf R}$ or ${\bf C}$, as appropriate.  Conversely, one can check
that a linear functional $\lambda$ on $V$ is bounded when it is
continuous at $0$.  Note that continuity of a linear functional on $V$
at $0$ implies continuity at every point in $V$, by linearity.

        Suppose for the moment that $V$ is equipped with an inner
product $\langle v, w \rangle$, and that $\|v\|$ is the norm
associated to this inner product.  If $w \in V$, then
\begin{equation}
        \lambda_w(v) = \langle v, w \rangle
\end{equation}
defines a bounded linear functional on $V$, since
\begin{equation}
        |\lambda_w(v)| \le \|v\| \, \|w\|
\end{equation}
for every $v \in V$, by the Cauchy--Schwarz inequality.  Conversely,
if $V$ is complete, and if $\lambda$ is a bounded linear functional on
$V$, then there is a unique $w \in V$ such that $\lambda(v) =
\lambda_w(v)$ for every $v \in V$.  The uniqueness of $w$ is a simple
exercise that does not use the completeness of $V$, and so we proceed
now to the proof of the existence of $w$.  This is trivial when
$\lambda(v) = 0$ for every $v \in V$, and hence we suppose that
$\lambda(v_0) \ne 0$ for some $v_0 \in V$.

        Let $Z$ be the kernel of $\lambda$, which is to say that
\begin{equation}
        Z = \{v \in V : \lambda(v) = 0\}.
\end{equation}
It is easy to see that $Z$ is a closed linear subspace of $V$, because
of the continuity of $\lambda$ that follows from boundedness.  Thus
the orthogonal projection $P_Z$ of $V$ onto $Z$ may be defined as in
Section \ref{orthogonal projections}.  Consider
\begin{equation}
        w_0 = v_0 - P_Z(v_0).
\end{equation}
Note that $w_0 \ne 0$, since $v_0 \not\in Z$ by hypothesis, and that
$w_0$ is orthogonal to every element of $Z$, as in Section
\ref{orthogonal projections}.  In particular,
\begin{equation}
  \langle v_0, w_0 \rangle = \langle v_0 - P_Z(v_0), w_0 \rangle 
                           = \langle w_0, w_0 \rangle = \|w_0\|^2 > 0.
\end{equation}
Put $w = \lambda(v_0) \, \|w_0\|^{-2} \, w_0$ in the real case, and $w
= \overline{\lambda(v_0)} \, \|w_0\|^{-2} \, w_0$ in the complex case.
By construction, $\lambda_w(v_0) = \lambda(v_0)$, and $\lambda_w(z) =
0$ for every $z \in Z$.  This implies that $\lambda_w(v) = \lambda(v)$
for every $v \in V$, as desired, because $V$ is spanned by $v_0$ and
$Z$ in this situation.

        Let $V$ be a real or complex vector space with an arbitrary
norm $\|v\|$ again.  Suppose that $W$ is a linear subspace of $V$, and
that $\lambda$ is a linear functional on $W$ that satisfies
(\ref{|lambda(v)| le C ||v||}) for some $C \ge 0$ and every $v \in V$.
Under these conditions, the Hahn--Banach theorem states that there is
an extension of $\lambda$ to a linear functional on $V$ that satisfies
(\ref{|lambda(v)| le C ||v||}) for every $v \in V$, with the same
constant $C$.  If $V$ is finite-dimensional, then this is the same in
essence as Theorem \ref{extension theorem} in Section \ref{second
duals}.  Otherwise, there is an argument using the axiom of choice,
with the previous construction as an important part of the proof.  In
some situations, one can use a sequence of extensions as before to
extend $\lambda$ to a dense linear subspace of $V$, and then extend
$\lambda$ to all of $V$ using continuity.  At any rate, an important
consequence of the Hahn--Banach theorem is that for each $v \in V$
with $v \ne 0$ there is a bounded linear functional $\lambda$ on $V$
such that $\lambda(v) \ne 0$.  More precisely, one can first define
$\lambda$ on the $1$-dimensional linear subspace of $V$ spanned by
$v$, and then use the Hahn--Banach theorem to extend $\lambda$ to a
bounded linear functional on all of $V$.

\section{Dual spaces}
\label{dual spaces, 2}

        Let $V$ be a real or complex vector space with a norm $\|v\|$
again, and let $V^*$ be the space of all bounded linear functionals on
$V$.  This is also a vector space over the real or complex numbers in
a natural way, because the sum of two bounded linear functionals on
$V$ is also bounded, as is the product of a bounded linear functional
on $V$ by a scalar.  If $\lambda$ is a bounded linear functional on
$V$, then the dual norm $\|\lambda\|^*$ of $\lambda$ is defined by
\begin{equation}
        \|\lambda\|^* = \sup \{|\lambda(v)| : v \in V, \, \|v\| \le 1\}.
\end{equation}
This is the same as (\ref{||lambda||^* = sup {|lambda(v)| : v in V,
    ||v|| le 1}}) in Section \ref{dual spaces, norms}, except that now
we need to ask that $\lambda$ be a bounded linear functional on $V$ to
ensure that the supremum is finite.  As before, $\lambda$ satisfies
(\ref{|lambda(v)| le C ||v||}) with $C = \|\lambda\|^*$, and this is
the smallest value of $C$ for which (\ref{|lambda(v)| le C ||v||})
holds.

        It is easy to see that $\|\lambda\|^*$ is a norm on $V^*$, as
in the finite-dimensional case.  Let us check that $V^*$ is
automatically complete with respect to the dual norm.  Let
$\{\lambda_j\}_{j = 1}^\infty$ be a sequence of bounded linear
functionals on $V$ which is a Cauchy sequence with respect to the dual
norm.  This means that for each $\epsilon > 0$ there is an
$L(\epsilon) \ge 1$ such that
\begin{equation}
\label{||lambda_j - lambda_l||^* < epsilon}
        \|\lambda_j - \lambda_l\|^* < \epsilon
\end{equation}
for every $j, l \ge L(\epsilon)$, and hence
\begin{equation}
\label{|lambda_j(v) - lambda_l(v)| le epsilon ||v||}
        |\lambda_j(v) - \lambda_l(v)| \le \epsilon \, \|v\|
\end{equation}
for every $v \in V$ and $j, l \ge L(\epsilon)$.  In particular,
$\{\lambda_j(v)\}_{j = 1}^\infty$ is a Cauchy sequence of real or
complex numbers, as appropriate, for each $v \in V$.  Thus
$\{\lambda_j(v)\}_{j = 1}^\infty$ converges to a real or complex
number $\lambda(v)$ for each $v \in V$, by the completeness of ${\bf
R}$, ${\bf C}$.  One can check that $\lambda$ defines a linear
functional on $V$, because $\lambda_j$ is linear on $V$ for each $j$.
We also have that
\begin{equation}
\label{|lambda_j(v) - lambda(v)| le epsilon ||v||}
        |\lambda_j(v) - \lambda(v)| \le \epsilon \, \|v\|
\end{equation}
for every $v \in V$ and $j \ge L(\epsilon)$, by taking the limit
as $l \to \infty$ in (\ref{|lambda_j(v) - lambda_l(v)| le epsilon ||v||}).
Applying this with $\epsilon = 1$ and $j = L(1)$, we get that
\begin{equation}
        |\lambda(v)| \le |\lambda_{L(1)}(v)| + \|v\|
                      \le (\|\lambda_{L(1)}\|^* + 1) \, \|v\|
\end{equation}
for every $v \in V$, so that $\lambda$ is a bounded linear functional
on $V$.  Using (\ref{|lambda_j(v) - lambda(v)| le epsilon ||v||}) again,
we get that $\{\lambda_j\}_{j = 1}^\infty$ converges to $\lambda$ with
respect to the dual norm, as desired.

        Let $V^{**}$ be the space of bounded linear functionals on
$V^*$, with respect to the dual norm $\|\lambda\|^*$ on $V^*$.  
If $v \in V$, then
\begin{equation}
\label{L_v(lambda) = lambda(v)}
        L_v(\lambda) = \lambda(v)
\end{equation}
defines a linear functional on $V^*$, which satisfies
\begin{equation}
\label{|L_v(lambda)| = |lambda(v)| le ||lambda||^* ||v||}
        |L_v(\lambda)| = |\lambda(v)| \le \|\lambda\|^* \, \|v\|
\end{equation}
for every $\lambda \in V^*$, by the definition of $\|\lambda\|^*$.
This implies that $L_v$ is a bounded linear functional on $V^*$.
More precisely, if $\|L\|^{**}$ is the dual norm of a bounded linear
functional $L$ on $V^*$ with respect to the dual norm $\|\lambda\|^*$
on $V^*$, then (\ref{|L_v(lambda)| = |lambda(v)| le ||lambda||^*
||v||}) implies that
\begin{equation}
\label{||L_v||^{**} le ||v||}
        \|L_v\|^{**} \le \|v\|
\end{equation}
for every $v \in V$.  Using the Hahn--Banach theorem, one can check
that
\begin{equation}
\label{||L_v||^{**} = ||v||}
        \|L_v\|^{**} = \|v\|
\end{equation}
for every $v \in V$.  The main point is to show that if $v \ne 0$,
then there is a $\lambda \in V^*$ such that $\|\lambda\|^* = 1$ and
$\lambda(v) = \|v\|$, so that equality holds in (\ref{|L_v(lambda)| =
|lambda(v)| le ||lambda||^* ||v||}).  As usual, one can start by
defining $\lambda$ on the $1$-dimensional subspace of $V$ spanned by
$v$, and then extend $\lambda$ to all of $V$ using the Hahn--Banach
theorem.

        A Banach space $V$ is said to be
\emph{reflexive}\index{reflexive Banach spaces} if every bounded
linear functional on $V^{**}$ is of the form $L_v$ for some $v \in V$.
Note that $V$ has to be complete for this to hold, since we already
know that $V^{**}$ is complete, because it is a dual space.  It is
easy to see that Hilbert spaces are reflexive, using the
characterization of their dual spaces in the previous section.  It is
also well known that $L^p$ spaces are reflexive when $1 < p < \infty$,
because the dual of $L^p$ can be identified with the corresponding
$L^q$ space, where $1 < q < \infty$ is conjugate to $p$ in the usual
sense that $1/p + 1/q = 1$.  In particular, $\ell^p$ spaces are
reflexive when $1 < p < \infty$, by the characterization of their dual
spaces in Section \ref{bounded linear functionals}.  We also saw in
Section \ref{bounded linear functionals} that the dual of $c_0(E)$ can
be identified with $\ell^1(E)$ for any nonempty set $E$, and that the
dual of $\ell^1(E)$ can be identified with $\ell^\infty(E)$.  If $E$
is an infinite set, then $c_0(E)$ is a proper linear subspace of
$\ell^\infty(E)$, and it follows that $c_0(E)$ is not reflexive.

\section{Bounded linear mappings}
\label{bounded linear mappings}

        Let $V_1$ and $V_2$ be vector spaces, both real or both
complex, and equipped with norms $\|\cdot \|_1$, $\|\cdot \|_2$,
respectively.  A linear mapping $T$ from $V_1$ into $V_2$ is said to
be \emph{bounded}\index{bounded linear mappings} if
\begin{equation}
\label{||T(v)||_2 le C ||v||_1}
        \|T(v)\|_2 \le C \, \|v\|_1
\end{equation}
for some $C \ge 0$ and every $v \in V_1$.  If $V_1$ has finite
dimension, then one can check that every linear mapping from $V_1$
into $V_2$ is bounded, using a basis for $V_1$ and the remarks at the
end of Section \ref{definitions, examples} to reduce to the case where
$V_1$ is ${\bf R}^n$ or ${\bf C}^n$ equipped with the standard norm.
If $V_2 = {\bf R}$ or ${\bf C}$, as appropriate, then a bounded linear
mapping from $V_1$ into $V_2$ is the same as a bounded linear
functional on $V_1$.  The boundedness of any linear mapping is
equivalent to suitable continuity conditions, as in the context of
linear functionals.

        Let $\mathcal{BL}(V_1, V_2)$ be the space of bounded linear
mappings from $V_1$ into $V_2$.  It is easy to see that this is a
vector space with respect to pointwise addition and scalar
multiplication.  If $T$ is a bounded linear mapping from $V_1$ into
$V_2$, then the operator norm of $T$ is defined by
\begin{equation}
\label{||T||_{op} = sup {||T(v)||_2 : v in V_1, ||v||_1 le 1}}
        \|T\|_{op} = \sup \{\|T(v)\|_2 : v \in V_1, \, \|v\|_1 \le 1\},
\end{equation}
as in (\ref{defn of ||T||_{op}}) in Section \ref{linear transformations}.
Equivalently, (\ref{||T(v)||_2 le C ||v||_1}) holds with $C = \|T\|_{op}$,
and this is the smallest value of $C$ for which (\ref{||T(v)||_2 le C ||v||_1})
holds.  One can check that (\ref{||T||_{op} = sup {||T(v)||_2 : v in
V_1, ||v||_1 le 1}}) defines a norm on $\mathcal{BL}(V_1, V_2)$.
If $V_2 = {\bf R}$ or ${\bf C}$, as appropriate, then the operator norm
reduces to the dual norm on $(V_1)^*$ defined in the previous section.
If $V_2$ is any vector space which is complete with respect to the norm
$\|\cdot \|_2$, then one can show that $\mathcal{BL}(V_1, V_2)$ is complete
with respect to the operator norm, in the same way as for dual spaces.

        Let $V_3$ be another vector space, which is real or complex
depending on whether $V_1$ and $V_2$ are real or complex, and let
$\|\cdot \|_3$ be a norm on $V_3$.  If $T_1$ is a bounded linear
mapping from $V_1$ into $V_2$, and $T_2$ is a bounded linear mapping
from $V_2$ into $V_3$, then it is easy to see that their composition
$T_2 \circ T_1$ is a bounded linear mapping from $V_1$ into $V_3$.
Moreover,
\begin{equation}
\label{||T_2 circ T_1||_{op, 13} le ||T_1||_{op, 12} ||T_2||_{op, 23}}
        \|T_2 \circ T_1\|_{op, 13} \le \|T_1\|_{op, 12} \, \|T_2\|_{op, 23},
\end{equation}
where $\|\cdot \|_{op, ab}$ is the operator norm for a linear mapping
from $V_a$ into $V_b$, with $a, b = 1, 2, 3$.

        A bounded linear mapping $T : V_1 \to V_2$ is said to be
\emph{invertible}\index{invertible linear mappings} if it is a
one-to-one linear mapping from $V_1$ onto $V_2$ whose inverse $T^{-1}$
is bounded as a linear mapping from $V_2$ into $V_1$.  Note that the
composition of invertible mappings is also invertible.  If $T$ is
invertible, then
\begin{equation}
\label{||T(v)||_2 ge c ||v||_1}
        \|T(v)\|_2 \ge c \, \|v\|_1
\end{equation}
for some $c > 0$ and every $v \in V_1$.  More precisely, this holds
with $c$ equal to the reciprocal of the operator norm of $T^{-1}$, by
applying the boundedness of $T^{-1}$ to $T^{-1}(T(v)) = v$.
Conversely, suppose that $T$ is a bounded linear mapping from $V_1$
into $V_2$ that satisfies (\ref{||T(v)||_2 ge c ||v||_1}).  In
particular, $v = 0$ when $T(v) = 0$, so that $T$ is one-to-one.  If
$T$ maps $V_1$ onto $V_2$, then (\ref{||T(v)||_2 ge c ||v||_1})
implies that $T^{-1}$ is bounded, with operator norm less than or
equal to $1/c$.

        If $V_1$ is complete and $T : V_1 \to V_2$ is a bounded linear
mapping that satisfies (\ref{||T(v)||_2 ge c ||v||_1}), then it is
easy to see that $T(V_1)$ is also complete.  This is because a sequence
$\{v_j\}_{j = 1}^\infty$ of elements of $V_1$ is a Cauchy sequence in $V_1$
if and only if $\{T(v_j)\}_{j = 1}^\infty$ is a Cauchy sequence in $V_2$,
and $\{v_j\}_{j = 1}^\infty$ converges to $v \in V_1$ if and only if 
$\{T(v_j)\}_{j = 1}^\infty$ converges to $T(v)$ in $V_2$.  In this case,
it follows that $T(V_1)$ is a closed linear subspace in $V_2$.  To see this,
let $\{v_j\}_{j = 1}^\infty$ be a sequence of elements of $V_1$ such that
$\{T(v_j)\}_{j = 1}^\infty$ converges to some $z \in V_2$, and let us check
that $z = T(v)$ for some $v \in V_1$.  Note that $\{T(v_j)\}_{j = 1}^\infty$
is a Cauchy sequence in $V_2$, since it converges in $V_2$.  As before,
this implies that $\{v_j\}_{j = 1}^\infty$ is a Cauchy sequence in $V_1$,
so that $\{v_j\}_{j = 1}^\infty$ converges to some $v \in V$, because $V$
is complete.  Thus $\{T(v_j)\}_{j = 1}^\infty$ converges to $T(v)$ in
$V_2$, because $T$ is bounded, and hence $z = T(v)$, as desired.

        Let $V$ be a real or complex vector space with a norm $\|v\|$,
and let $T$ be a bounded linear operator on $V$.  If $j$ is a positive
integer, then let $T^j$ be the composition of $j$ $T$'s, so that $T^1
= T$, $T^2 = T \circ T$, and so on.  It will be convenient to
interpret $T^j$ as being the identity operator $I$ on $V$ when $j =
0$.  Observe that
\begin{equation}
\label{(I - T) (sum_{j = 0}^n T^j) = (sum_{j = 0}^n T^j) (I - T) = I - T^{n+1}}
        (I - T) \, \Big(\sum_{j = 0}^n T^j\Big)
                              = \Big(\sum_{j = 0}^n T^j\Big) \, (I - T)
                                   = I - T^{n + 1}
\end{equation}
for each nonnegative integer $n$, as in the case of ordinary geometric
series of real and complex numbers.  Of course,
\begin{equation}
\label{||T^j||_{op} le ||T||_{op}^j}
        \|T^j\|_{op} \le \|T\|_{op}^j
\end{equation}
for each $j$, by (\ref{||T_2 circ T_1||_{op, 13} le ||T_1||_{op, 12}
||T_2||_{op, 23}}).  If $\|T\|_{op} < 1$, then we get that
\begin{equation}
\label{sum_{j = 0}^infty ||T^j||_{op} le ... = frac{1}{1 - ||T||_{op}}}
        \sum_{j = 0}^\infty \|T^j\|_{op} \le \sum_{j = 0}^\infty \|T\|_{op}^j
                                            = \frac{1}{1 - \|T\|_{op}},
\end{equation}
by the usual formula for the sum of an geometric series.

        This shows that the infinite series
\begin{equation}
\label{sum_{j = 0}^infty T^j}
        \sum_{j = 0}^\infty T^j
\end{equation}
converges absolutely in the vector space $\mathcal{BL}(V) = \mathcal{BL}(V,
V)$ of bounded linear operators on $V$ when $\|T\|_{op} < 1$.  If $V$
is complete, then we have seen that $\mathcal{BL}(V)$ is complete with
respect to the operator norm, and hence that (\ref{sum_{j = 0}^infty
T^j}) converges in $\mathcal{BL}(V)$.  We also get that
\begin{equation}
\label{(I - T) (sum_{j = 0}^infty T^j) = (sum_{j = 0}^infty T^j) (I - T) = I}
        (I - T) \, \Big(\sum_{j = 0}^\infty T^j\Big)
                  = \Big(\sum_{j = 0}^\infty T^j\Big) \, (I - T) = I,
\end{equation}
by taking the limit as $n \to \infty$ in (\ref{(I - T) (sum_{j = 0}^n
T^j) = (sum_{j = 0}^n T^j) (I - T) = I - T^{n+1}}), and using the fact
that $T^{n + 1} \to 0$ as $n \to \infty$ when $\|T\|_{op} < 1$.  Thus
$I - T$ is invertible on $V$ when $\|T\|_{op} < 1$ and $V$ is
complete.

        Let us continue to ask that $V$ be complete.  If $T$ is any
bounded linear operator on $V$, and $\lambda$ is a real or complex
number, as appropriate, such that $|\lambda| > \|T\|_{op}$, then
$\lambda \, I - T$ is invertible on $V$.  This follows from the
preceding argument applied to $\lambda^{-1} \, T$.  Similarly, if $R$
is a bounded linear operator on $V$ which is also invertible, and if
$T$ is a bounded linear operator on $V$ that satisfies
\begin{equation}
        \|R^{-1}\|_{op} \, \|T\|_{op} < 1,
\end{equation}
then
\begin{equation}
        R - T = R \, (I - R^{-1} \, T)
\end{equation}
is invertible on $V$.

        Suppose that $V$ is a complex Banach space, and let $T$ be a
bounded linear operator on $V$.  The \emph{spectrum}\index{spectrum}
of $T$ is the set of complex numbers $\lambda$ such that $\lambda \, I
- T$ is not invertible on $V$.  As usual, eigenvalues of $T$ are
elements of the spectrum, but the converse does not hold in infinite
dimensions.  Note that
\begin{equation}
\label{|lambda| le ||T||_{op}}
        |\lambda| \le \|T\|_{op},
\end{equation}
for every $\lambda \in {\bf C}$ in the spectrum of $T$, as in the
preceding paragraph.  If $\lambda \in {\bf C}$ is not in the spectrum
of $T$, so that $\lambda \, I - T$ is invertible on $V$, then $\mu \,
I - T$ is also invertible for every complex number $\mu$ sufficiently
close to $\lambda$, by the remarks in the previous paragraph.  This
implies that the spectrum of $T$ is a closed set in the complex plane.
A famous theorem states that the spectrum of $T$ is always nonempty.
The main idea in the proof is that otherwise $(\lambda \, I - T)^{-1}$
would be a holomorphic function of $\lambda$ on the complex plane that
tends to $0$ as $|\lambda| \to \infty$.

        Let $V_1$ and $V_2$ be Banach spaces, both real or both
complex, and with norms $\|\cdot \|_1$, $\|\cdot \|_2$, respectively.
Also let
\begin{equation}
        B_1 = \{v \in V_1 : \|v\|_1 \le 1\}
\end{equation}
be the closed unit ball in $V_1$.  A linear mapping $T$ from $V_1$
into $V_2$ is said to be \emph{compact}\index{compact linear mappings}
if the closure of $T(B_1)$ in $V_2$ a compact set.  This is equivalent
to asking that $T(B_1)$ be \emph{totally bounded}\index{totally
bounded sets} in $V_2$, which means that for each $\epsilon > 0$,
$T(B_1)$ can be covered by finitely many balls of radius $\epsilon$ in
$B_2$.  In particular, totally bounded sets are bounded, and hence
compact linear mappings are bounded.  It is easy to see that bounded
subsets of finite-dimensional spaces are totally bounded, so that
bounded linear mappings of finite rank are compact.  One can also
check that the composition of a bounded linear mapping with a compact
linear mapping is compact, where the compact operator is either first
or second in the composition.

        Let $\mathcal{CL}(V_1, V_2)$ be the space of compact linear
mappings from $V_1$ into $V_2$.  This is a linear subspace of the
vector space $\mathcal{BL}(V_1, V_2)$ of bounded linear mappings from
$V_1$ into $V_2$, which is closed with respect to the operator norm on
$\mathcal{BL}(V_1, V_2)$.  This means that if $\{T_j\}_{j = 1}^\infty$
is a sequence of compact linear mappings from $V_1$ into $V_2$ that
converges to a bounded linear mapping $T : V_1 \to V_2$ with respect
to the operator norm, then $T$ is compact too.  In particular, $T$ is
compact if it is the limit of a sequence of bounded linear mappings of
finite rank with respect to the operator norm.  In some situations,
including mappings between Hilbert spaces, one can show that every
compact linear mapping is the limit of a sequence of bounded linear
mappings of finite rank with respect to the operator norm.

        Let $T$ be a compact linear mapping from a Banach space $V$
into itself.  If $\lambda$ is a real or complex number, as
appropriate, such that $\lambda \ne 0$ and $\lambda \, I - T$ is not
invertible, then it can be shown that $\lambda$ is an eigenvalue of
$T$, and that the corresponding eigenspace is finite-dimensional.  It
can also be shown that for each $r > 0$, there are only finitely many
eigenvalues $\lambda$ with $|\lambda| \ge r$.

\section{Self-adjoint linear operators}
\label{self-adjoint linear operators, 2}

        Let $V$ be a real or complex vector space with an inner
product $\langle v, w \rangle$, and suppose that $V$ is complete with
respect to the corresponding norm $\|v\|$, so that $V$ is a Hilbert
space.  As in Section \ref{self-adjoint linear operators}, a bounded
linear operator $A$ on $V$ is said to be
\emph{self-adjoint}\index{self-adjoint linear operators} if
\begin{equation}
        \langle A(v), w \rangle = \langle v, A(w) \rangle
\end{equation}
for every $v, w \in V$.  As before, the identity operator $I$ on $V$
is self-adjoint, as is the orthogonal projection $P_W$ of $V$ onto a
closed linear subspace $W$ of $V$.  The sum of two bounded
self-adjoint linear operators on $V$ is also self-adjoint, and the
product of a bounded self-adjoint linear operator on $V$ and a real
number is self-adjoint too.

        Suppose for the moment that $V$ is a complex Hilbert space.
If $A$ is a bounded self-adjoint linear operator on $V$, then
\begin{equation}
        \langle A(v), v \rangle = \langle v, A(v) \rangle
                                = \overline{\langle A(v), v \rangle}
\end{equation}
for every $v \in V$, and hence
\begin{equation}
\label{langle A(v), v rangle in {bf R}, 2}
        \langle A(v), v \rangle \in {\bf R}
\end{equation}
for every $v \in V$.  Using this, it is easy to see that the
eigenvalues of $A$ are real numbers, as before.  Let us check that the
spectrum of $A$, as defined in the previous section, is also contained
in the real line under these conditions.  Equivalently, this means
that $\lambda \, I - A$ is invertible on $V$ for every complex number
$\lambda$ with nonzero imaginary part.

        Observe that
\begin{equation}
\label{im langle (lambda I - A)(v), v rangle = (im lambda) ||v||^2}
        \im \langle (\lambda \, I - A)(v), v \rangle = (\im \lambda) \, \|v\|^2
\end{equation}
for every $v \in V$, by (\ref{langle A(v), v rangle in {bf R}, 2}), so that
\begin{equation}
\label{|langle (lambda I - A)(v), v rangle| ge |im lambda| ||v||^2}
 |\langle (\lambda \, I - A)(v), v \rangle| \ge |\im \lambda| \, \|v\|^2
\end{equation}
for every $v \in V$.  The Cauchy--Schwarz inequality implies that
\begin{equation}
        |\langle (\lambda \, I - A)(v), v \rangle|
                \le \|(\lambda \, I - A)(v)\| \, \|v\|,
\end{equation}
from which we get that
\begin{equation}
\label{||(lambda I - A)(v)|| ge |im lambda| ||v||}
        \|(\lambda \, I - A)(v)\| \ge |\im \lambda| \, \|v\|
\end{equation}
for every $v \in V$.  This is the same type of condition as
(\ref{||T(v)||_2 ge c ||v||_1}) in the previous section, since $\im
\lambda \ne 0$, by hypothesis.  In order to show that $\lambda \, I -
A$ is invertible on $V$, it suffices to check that $\lambda \, I - A$
maps $V$ onto itself.

        Suppose for the sake of a contradiction that $W = (\lambda \,
I - A)(V)$ is a proper linear subspace of $V$.  Note that $W$ is a
closed linear subspace of $V$, because of (\ref{||(lambda I - A)(v)||
ge |im lambda| ||v||}) and the completeness of $V$, as in the previous
section.  Let $v$ be any element of $V \backslash W$, and put
\begin{equation}
\label{y = v - P_W(v)}
        y = v - P_W(v),
\end{equation}
where $P_W(v)$ is the orthogonal projection of $v$ onto $W$, as in
Section \ref{orthogonal projections}.  Thus $y \ne 0$, because $v
\not\in W$ and $P_W(v) \in W$, and $y$ is orthogonal to every element
of $W$.  The latter condition is the same as saying that
\begin{equation}
        \langle (\lambda \, I - A)(v), y \rangle = 0
\end{equation}
for every $v \in V$.  In particular, we can apply this to $v = y$,
to get that
\begin{equation}
        \langle (\lambda \, I - A)(y), y \rangle = 0.
\end{equation}
This implies that $y = 0$, by (\ref{|langle (lambda I - A)(v), v
rangle| ge |im lambda| ||v||^2}), contradicting the hypothesis that $y
\ne 0$.  It follows that $W = V$, so that $\lambda \, I - A$ is
invertible on $V$, as desired.

        Let $V$ be a real or complex Hilbert space again.  A bounded
self-adjoint linear operator $A$ on $V$ is said to be
\emph{nonnegative}\index{nonnegative self-adjoint operators} if
\begin{equation}
        \langle A(v), v \rangle \ge 0
\end{equation}
for every $v \in V$.  Suppose that $A$ satisfies the strict positivity
condition that
\begin{equation}
        \langle A(v), v \rangle \ge c \, \|v\|^2
\end{equation}
for some $c > 0$ and every $v \in V$, and let us check that $A$ is
invertible on $V$.  By the Cauchy-Schwarz inequality,
\begin{equation}
\label{langle A(v), v rangle le ||A(v)|| ||v||}
        \langle A(v), v \rangle \le \|A(v)\| \, \|v\|
\end{equation}
for every $v \in V$, and hence
\begin{equation}
\label{||A(v)|| ge c ||v||}
        \|A(v)\| \ge c \, \|v\|
\end{equation}
for every $v \in V$.  This is the same as (\ref{||T(v)||_2 ge c
||v||_1}) in this context, and it suffices to show that $A$ maps $V$
onto itself.

        As before, $W = A(V)$ is a closed linear subspace of $V$ under
these conditions.  If $W \ne V$, then there is a $y \in V$ such that
$y \ne 0$ and $y$ is orthogonal to every element of $W$.  Equivalently,
this means that
\begin{equation}
        \langle A(v), y \rangle = 0
\end{equation}
for every $v \in V$, and for $v = y$ in particular, so that $\langle
A(y), y \rangle = 0$.  This implies that $y = 0$, by the strict
positivity of $A$, contradicting the hypothesis that $y \ne 0$.  Thus
$A(V) = V$, and hence $A$ is invertible on $V$, as desired.

        Similarly, if $A$ is a bounded self-adjoint linear operator on
$V$ that satisfies (\ref{||A(v)|| ge c ||v||}) for some $c > 0$ and
every $v \in V$, then $A$ is invertible on $V$.  As before, $W = A(V)$
is a closed linear subspace of $V$ under these conditions, and we want
to show that $W = V$.  Otherwise, there is a $y \in V$ such that $y
\ne 0$ and $y$ is orthogonal to every element of $W$, so that
\begin{equation}
        \langle v, A(y) \rangle = \langle A(v), y \rangle = 0
\end{equation}
for every $v \in V$.  This implies that $A(y) = 0$, and hence that $y
= 0$, because of (\ref{||A(v)|| ge c ||v||}).  Thus $A(V) = V$, so
that $A$ is invertible on $V$, as desired.

        In analogy with the finite-dimensional case, it can be shown
that a compact self-adjoint linear operator $T$ on $V$ can be
diagonalized in an orthonormal basis for $V$.  Using this, one can
show that any compact linear mapping between Hilbert spaces has a
Schmidt decomposition\index{Schmidt decompositions} as in Section
\ref{schmidt decompositions}, but perhaps with infinite sequences of
orthonormal vectors, and coefficients $\lambda_j$ converging to $0$ as
$j \to \infty$.  If the $\lambda_j$'s are $p$-summable for some $p >
0$, then the operator is said to be in the $\mathcal{S}_p$
class.\index{Sp class@$\mathcal{S}_p$ class}

\chapter{Marcel Riesz' convexity theorem}
\label{riesz convexity theorem}

\noindent
Let $(a_{j, k})$ be an $n \times n$ matrix of complex numbers, and let
$A(x,y)$ be the bilinear form defined for $x, y \in {\bf C}^n$ by
\begin{equation}
\label{def of A(x,y)}
	A(x,y) = \sum_{j=1}^n \sum_{k=1}^n y_j \, a_{j, k} \, x_k.
\end{equation}
For $1 < p < \infty$, let $M_p$ be the quantity
\begin{equation}
\label{def of M_p, 1 < p < infty}
 \quad 	\sup \bigg\{ |A(x,y)| : x, y \in {\bf C}^n, 
		\, \Big(\sum_{k=1}^n |x_k|^p \Big)^{1/p} \le 1,
		\, \Big(\sum_{j=1}^n |y_j|^{p'} \Big)^{1/p'} \le 1 \bigg\},
\end{equation}
where $p'$ denotes the conjugate exponent of $p$, $1/p + 1/p' = 1$.
When $p = 1$, $p' = \infty$, put
\begin{equation}
\label{def of M_p, p = 1}
	M_1 = \sup \bigg\{ |A(x,y)| : x, y \in {\bf C}^n, 
		\, \sum_{k=1}^n |x_k| \le 1,
		\, \max_{1 \le j \le n} |y_j| \le 1 \bigg\},
\end{equation}
and when $p = \infty$, $p' = 1$, set
\begin{equation}
\label{def of M_p, p = infty}
	M_\infty = \sup \bigg\{ |A(x,y)| : x, y \in {\bf C}^n, 
		\, \max_{1 \le k \le n} |x_k| \le 1,
		\, \sum_{j=1}^n |y_j| \le 1 \bigg\}.
\end{equation}

\begintheorem 
\label{M. Riesz' convexity theorem}
As a function of $1/p \in [0,1]$, $\log M_p$ is convex.
\end{theorem}

	More precisely, if $1 \le p < q \le \infty$, $0 < t < 1$, $1 <
r < \infty$, and
\begin{equation}
\label{frac{1}{r} = frac{t}{p} + frac{1 - t}{q}}
        \frac{1}{r} = \frac{t}{p} + \frac{1 - t}{q},
\end{equation}
then
\begin{equation}
\label{M_r le M_p^t M_q^{1-t}}
	M_r \le M_p^t \, M_q^{1-t}.
\end{equation}
If $M_p = 0$ for some $p$, then $A \equiv 0$ and $M_p = 0$ for every
$p$, and hence we may as well assume that $A \not\equiv 0$ in the
arguments that follow.  The special case of $p = 1$, $q = \infty$
corresponds exactly to the theorem of Schur discussed in Section
\ref{special cases}.  Note that the analogous inequality holds when
the $a_{j, k}$'s are real numbers, and we use $x, y \in {\bf R}^n$ in
the definition of $M_p$, by the same proof.

	We can also describe $M_p$ as
\begin{equation}
\label{second version of M_p, p < infty}
	M_p = \sup \bigg\{ 
 \Big(\sum_{j=1}^n \, \biggl|\sum_{k=1}^n a_{j, k} \, x_k \biggr|^p\Big)^{1/p}
	   : x \in {\bf C}^n, \, \Big(\sum_{k=1}^n |x_k|^p \Big)^{1/p} \le 1
								\bigg\}
\end{equation}
when $1 \le p < \infty$, and
\begin{equation}
\label{second version of M_p, p = infty}
	M_\infty = \sup \bigg\{ 
 	\max_{1 \le j \le n} \, \biggl|\sum_{k=1}^n a_{j, k} \, x_k \biggr|
	   : x \in {\bf C}^n, \, \max_{1 \le k \le n} |x_k| \le 1 \bigg\}.
\end{equation}
This definition of $M_p$ is greater than or equal to the previous one
by H\"older's inequality, and for each $x \in {\bf C}^n$, there is a
$y \in {\bf C}^n$ for which equality holds and $\big(\sum_{j=1}
|y_j|^{p'} \big)^{1/p'}$ or $\max_{1 \le j \le n} |y_j|$ is equal to
$1$, according to whether $p' < \infty$ or $p' = \infty$.
Equivalently, $M_p$ is the operator norm of the linear transformation
on ${\bf C}^n$ associated to the matrix $(a_{j,k})$ with respect to
the $p$-norm $\|\cdot \|_p$ defined in Section \ref{definitions, examples}.
Similarly,
\begin{equation}
\label{third version of M_p, p > 1}
 \quad 	M_p = \sup \bigg\{ 
\Big(\sum_{k=1}^n \, 
	\biggl|\sum_{j=1}^n y_j \, a_{j, k} \biggr|^{p'}\Big)^{1/p'}
	   : y \in {\bf C}^n, \, \Big(\sum_{j=1}^n |y_j|^p \Big)^{1/p'} \le 1
								\bigg\}
\end{equation}
when $1 < p \le \infty$, and
\begin{equation}
	M_1 = \sup \bigg\{ 
 	\max_{1 \le k \le n} \, \biggl|\sum_{j=1}^n y_j \, a_{j, k} \biggr|
	   : y \in {\bf C}^n, \, \max_{1 \le j \le n} |y_j| \le 1 \bigg\},
\end{equation}
which says that $M_p$ is equal to the operator norm of the dual linear
transformation on ${\bf C}^n$, associated to the transpose matrix, and
with respect to the dual norm $\|\cdot \|_{p'}$.  One can check that
$M_s$ is a continuous function of $1/s$, $1/s \in [0,1]$, using the
inequalities (\ref{max_{1 le j le n} a_j le (sum_{j=1}^n
a_j^p)^{1/p}}), (\ref{(sum_{j=1}^n a_j^q)^{1/q} le (sum_{j=1}^n
a_j^p)^{1/p}}), (\ref{(sum_{j=1}^n a_j^p)^{1/p} le n^{1/p} max_{1 le j
le n} a_j}), and (\ref{(sum_{j=1}^n a_j^p)^{1/p} le n^{(1/p) -
(1/q)}(sum_{j=1}^n a_j^q)^{1/q}}).

	As in Lemma \ref{lemma about convexity conditions}, we would
like to show that for each $p, q \in [1,\infty]$ there is a $t \in
(0,1)$ such that (\ref{M_r le M_p^t M_q^{1-t}}) holds.  Fix a real
number $r$, $1 < r < \infty$, and let $r'$ be its conjugate exponent.
There exist $x^0, y^0 \in {\bf C}^n$ at which the supremum in the
definition (\ref{def of M_p, 1 < p < infty}) of $M_r$ is attained,
i.e., which satisfy
\begin{equation}
\label{maximizing A(x,y)}
	|A(x^0, y^0)| = M_r
\end{equation}
and the normalizations
\begin{equation}
\label{x^0 has norm 1}
	\Big(\sum_{k=1}^n |x^0_k|^r \Big)^{1/r} = 1
\end{equation}
and
\begin{equation}
\label{y^0 has norm 1}
	\Big(\sum_{j=1}^n |y^0_j|^{r'} \Big)^{1/r'} = 1.
\end{equation}
This follows from standard considerations of continuity and
compactness.

	Observe that
\begin{eqnarray}
\label{|A(x^0, y^0)| = ..., 1}
	|A(x^0, y^0)| 
  & = & \biggl| \sum_{j=1}^n \sum_{k=1}^n y^0_j \, a_{j, k} \, x^0_k \biggr|
									\\
  & = & \Big(\sum_{i=1}^n |y^0_i|^{r'} \Big)^{1/r'}
\Big(\sum_{j=1}^n \biggl| \sum_{k=1}^n a_{j, k} \, x^0_k \biggr|^r \Big)^{1/r}.
								\nonumber
\end{eqnarray}
The first step uses only the definition of $A$.  If the second
equality were replaced with $\le$, then it would be a consequence of
H\"older's inequality.  If equality did not hold, then we could
replace $y^0$ with an element of ${\bf C}^n$ which satisfies (\ref{y^0
has norm 1}) and for which equality does hold, increasing the value of
$|A(x^0, y^0)|$.  Similarly,
\begin{eqnarray}
\label{|A(x^0, y^0)| = ..., 2}
	|A(x^0, y^0)| 
  & = & \biggl| \sum_{j=1}^n \sum_{k=1}^n y^0_j \, a_{j, k} \, x^0_k \biggr|
									\\
  & = &   \Big(\sum_{k=1}^n \biggl| \sum_{j=1}^n y^0_j \, a_{j, k} \biggr|^{r'}
								  \Big)^{1/r'}
	\Big(\sum_{l=1}^n |x^0_l|^r \Big)^{1/r}.
								\nonumber
\end{eqnarray}

	Because of the second equality in (\ref{|A(x^0, y^0)| = ...,
1}), there is a $\mu \ge 0$ such that
\begin{equation}
\label{|sum_{k=1}^n a_{j, k} x^0_k| = mu |y^0_j|^{r' - 1}}
	\biggl|\sum_{k=1}^n a_{j, k} \, x^0_k \biggr| = \mu \, |y^0_j|^{r' - 1}
\end{equation}
for $j = 1, 2, \ldots, n$.  This can be derived from the proof of
H\"older's inequality, by analyzing the conditions in which equality
holds.  Similarly, there is a $\nu \ge 0$ such that
\begin{equation}
\label{|sum_{j=1}^n y^0_j a_{j, k}| = nu |x^0_k|^{r-1}}
	\biggl|\sum_{j=1}^n y^0_j \, a_{j, k} \biggr| = \nu \, |x^0_k|^{r-1}
\end{equation}
for $k = 1, 2, \ldots, n$.  From (\ref{|A(x^0, y^0)| = ..., 1}) and
(\ref{|A(x^0, y^0)| = ..., 2}), we have that
\begin{equation}
	M_r = 
\Big(\sum_{j=1}^n \biggl| \sum_{k=1}^n a_{j, k} \, x^0_k \biggr|^r \Big)^{1/r}
	    =
\Big(\sum_{k=1}^n \biggl| \sum_{j=1}^n y^0_j \, a_{j, k} \biggr|^{r'}
								 \Big)^{1/r'},
\end{equation}
using also (\ref{maximizing A(x,y)}), (\ref{x^0 has norm 1}), and
(\ref{y^0 has norm 1}).  Substituting (\ref{|sum_{k=1}^n a_{j, k}
x^0_k| = mu |y^0_j|^{r' - 1}}) in the first equality, we get that
\begin{equation}
	M_r = \mu \, \Big(\sum_{j=1}^n |y^0_j|^{r (r'-1)}\Big)^{1/r}.
\end{equation}
Because $r (r' - 1) = r'$, since $1/r + 1/r' = 1$, we can apply
(\ref{y^0 has norm 1}) to get that $M_r = \mu$.  For the same reasons,
$M_r = \nu$.

	Now suppose that $p$, $q$, and $t$ are real numbers such that
$1 \le p < q \le \infty$, $0 < t < 1$, and $1/r = t/p + (1-t)/q$.
Thus $p, q' < \infty$, where $q'$ is the conjugate exponent of $q$.
Observe that
\begin{equation}
\label{(sum_{j=1}^n |sum_{k=1}^n a_{j, k} x^0_k|^p)^{1/p} le ..}
\Big(\sum_{j=1}^n \biggl| \sum_{k=1}^n a_{j, k} \, x^0_k \biggr|^p \Big)^{1/p}
	\le M_p \, \Big(\sum_{k=1}^n |x^0_k|^p \Big)^{1/p}
\end{equation}
and
\begin{equation}
\label{(sum_{k=1}^n |sum_{j=1}^n y^0_j a_{j, k}|^{q'})^{1/q'} le ..}
\Big(\sum_{k=1}^n \biggl| \sum_{j=1}^n y^0_j \, a_{j, k} \biggr|^{q'}
								 \Big)^{1/q'}
	\le M_q \, \Big(\sum_{j=1}^n |y^0_j|^{q'} \Big)^{1/q'}.
\end{equation}
Applying the previous computations, we get that
\begin{equation}
\label{M_r -- M_p inequality}
	M_r \, \Big(\sum_{j=1}^n |y^0_j|^{p (r'-1)}\Big)^{1/p}
	   \le M_p \, \Big(\sum_{k=1}^n |x^0_k|^p \Big)^{1/p}
\end{equation}
and
\begin{equation}
\label{M_r -- M_q inequality}
	M_r \, \Big(\sum_{k=1}^n |x^0_k|^{q' (r-1)} \Big)^{1/q'}
	   \le M_q \, \Big(\sum_{j=1}^n |y^0_j|^{q'} \Big)^{1/q'}.
\end{equation}

	We are going to need some identities with indices.  Let us
first check that
\begin{equation}
\label{t (frac{1}{r} - frac{1}{p}) = (1-t) (frac{1}{r'} - frac{1}{q'})}
	t \Big(\frac{1}{r} - \frac{1}{p}\Big)
		= (1-t) \Big(\frac{1}{r'} - \frac{1}{q'}\Big).
\end{equation}
Because $1/r = t/p + (1-t)/q$, we have that
\begin{equation}
	t \Big(\frac{1}{r} - \frac{1}{p}\Big)
		= t \, (1-t) \, \Big(\frac{-1}{p} + \frac{1}{q}\Big).
\end{equation}
Similarly, $1/r' = t/p' + (1-t)/q'$, and 
\begin{equation}
	(1-t) \Big(\frac{1}{r'} - \frac{1}{q'}\Big)
		= (1-t) \, t \, \Big(\frac{1}{p'} - \frac{1}{q'}\Big)
		= (1-t) \, t \, \Big(\frac{-1}{p} + \frac{1}{q}\Big).
\end{equation}
This proves (\ref{t (frac{1}{r} - frac{1}{p}) = (1-t) (frac{1}{r'} -
frac{1}{q'})}).

	Suppose that
\begin{equation}
\label{frac{t}{r} = frac{1-t}{r'}}
	\frac{t}{r} = \frac{1-t}{r'}.
\end{equation}
This implies that
\begin{equation}
	\frac{t}{p} = \frac{1-t}{q'},
\end{equation}
by (\ref{t (frac{1}{r} - frac{1}{p}) = (1-t) (frac{1}{r'} -
frac{1}{q'})}).  Hence
\begin{equation}
	r' - 1 = \frac{1-t}{t} \quad\hbox{and}\quad r - 1 = \frac{t}{1-t},
\end{equation}
because $r \, (r' - 1) = r'$ and $r' \, (r - 1) = r$, since $1/r +
1/r' = 1$.  Therefore
\begin{equation}
\label{p (r'-1) = q', q' (r-1) = p}
	p \, (r'-1) = q' \quad\hbox{and}\quad q' \, (r-1) = p.
\end{equation}

	We can take the $t$th and $(1-t)$th powers of (\ref{M_r -- M_p
inequality}) and (\ref{M_r -- M_q inequality}), respectively, and then
multiply to get
\begin{eqnarray}
\label{M_r -- M_p, M_q inequality}
\lefteqn{M_r \, \Big(\sum_{j=1}^n |y^0_j|^{p (r'-1)}\Big)^{t/p}
	\Big(\sum_{k=1}^n |x^0_k|^{q' (r-1)} \Big)^{(1-t)/q'}} 	\\
	& & \le M_p^t \, M_q^{1-t} \, \Big(\sum_{k=1}^n |x^0_k|^p \Big)^{t/p}
		\Big(\sum_{j=1}^n |y^0_j|^{q'} \Big)^{(1-t)/q'}.
								\nonumber
\end{eqnarray}
Assuming (\ref{frac{t}{r} = frac{1-t}{r'}}), this reduces to
\begin{equation}
\label{M_r le M_p^t M_q^{1-t}, assuming frac{t}{r} = frac{1-t}{r'}}
	M_r \le M_p^t \, M_q^{1-t}.
\end{equation}
because the factors involving $x^0$ and $y^0$ on the left and right
sides of (\ref{M_r -- M_p, M_q inequality}) exactly match up under
these conditions, by the computations in the preceding paragraph.  To
summarize, for each $p$, $q$ with $1 \le p < q \le \infty$, there is a
$t \in (0,1)$ such that (\ref{frac{t}{r} = frac{1-t}{r'}}) holds when
$r$ is given by $1/r = t/p + (1-t)/q$.  For this choice of $t$, we get
the inequality (\ref{M_r le M_p^t M_q^{1-t}, assuming frac{t}{r} =
frac{1-t}{r'}}).  Theorem \ref{M. Riesz' convexity theorem} now
follows from Lemma \ref{lemma about convexity conditions}, with the
small adaptation to functions on closed intervals.

\chapter{Some dyadic analysis}
\label{dyadic analysis}

\section{Dyadic intervals}
\label{dyadic intervals}

	Normally, a reference to ``the unit interval'' in the real
line might suggest the closed interval $[0,1]$, but here it will be
convenient to use $[0,1)$ instead, for minor technical reasons.

\begindefinition
\label{Dyadic intervals in [0,1)}\index{dyadic intervals}
The \emph{dyadic subintervals} of $[0,1)$ are the intervals of the
form $[j \, 2^{-k}, (j+1) \, 2^{-k})$, where $j$ and $k$ are
nonnegative integers, and $j+1 \le 2^k$.  In particular, the length of
a dyadic interval in $[0, 1)$ is of the form $2^{-k}$, where $k$ is a
nonnegative integer.
\end{definition}

	The dyadic intervals in ${\bf R}$ can be defined in the same
way, with arbitrary integers $j$ and $k$.  The half-open, half-closed
condition leads to nice properties in terms of disjointness, as in the
next two lemmas, whose simple proofs are left as exercises.

\beginlemma
\label{Partitions of [0,1)}
For each nonnegative integer $k$, $[0,1)$ is the union of the dyadic
subintervals of length $2^{-k}$, and these subintervals are pairwise
disjoint.
\end{lemma}

\beginlemma
\label{Comparing two intervals}
If $J_1$ and $J_2$ are two dyadic subintervals of $[0,1)$, then either
$J_1 \subseteq J_2$, or $J_2 \subseteq J_1$, or $J_1 \cap J_2 =
\emptyset$.
\end{lemma}

	More precisely, if $J_1$, $J_2$ are dyadic subintervals of
$[0,1)$ such that the length of $J_2$ is less than or equal to the
length of $J_1$, then either $J_2 \subseteq J_1$ or $J_1 \cap J_2 =
\emptyset$.

\beginlemma
\label{Partitions of dyadic intervals}
If $J$ is a dyadic subinterval of $[0,1)$ of length $2^{-k}$, and if
$n$ is an integer greater than $k$, then $J$ is the union of the
dyadic subintervals of $J$ of length $2^{-n}$, and these subintervals
are pairwise disjoint.  Every dyadic subinterval of $[0,1)$ of length
$2^{-n}$ is contained in a unique dyadic subinterval of $[0,1)$ of
length $2^{-k}$ when $n \ge k$.
\end{lemma}

	This is easy to see.

\beginlemma
\label{unions of dyadic intervals}
If $\mathcal{F}$ is an arbitrary collection of dyadic subintervals of
$[0,1)$, then there is a subcollection $\mathcal{F}_0$ of
$\mathcal{F}$ such that
\begin{equation}
\label{bigcup_{J in mathcal{F}_0} J = bigcup_{J in mathcal{F}} J}
	\bigcup_{J \in \mathcal{F}_0} J = \bigcup_{J \in \mathcal{F}} J
\end{equation}
and the elements of $\mathcal{F}_0$ are pairwise disjoint.
\end{lemma}

	To prove this, we take $\mathcal{F}_0$ to be the set of
\emph{maximal} elements of $\mathcal{F}$, i.e., the set of $J \in
\mathcal{F}$ such that $J \subseteq J'$ for some $J' \in \mathcal{F}$
only when $J' = J$.  Every interval in $\mathcal{F}$ is contained in a
maximal interval in $\mathcal{F}$, since every dyadic subinterval of
$[0, 1)$ is contained in only finitely many dyadic subintervals of
$[0, 1)$.  Thus every element of $\mathcal{F}$ is contained in an
element of $\mathcal{F}_0$, which implies (\ref{bigcup_{J in
mathcal{F}_0} J = bigcup_{J in mathcal{F}} J}).  Any two maximal
elements of $\mathcal{F}$ which are distinct are disjoint, by Lemma
\ref{Comparing two intervals}, which implies the second property of
$\mathcal{F}_0$ in the lemma.

	Let $f$ be a real or complex-valued function on the unit
interval $[0,1)$ which is sufficiently well-behaved for integrals of
$f$ over subintervals of $[0, 1)$ to be defined.  One is welcome to
restrict one's attention to step functions here, and we shall simplify
this a bit further in a moment. For each nonnegative integer $k$, let
$E_k(f)$ be the function on $[0,1)$ defined by
\begin{equation}
\label{def of E_k(f)}
	E_k(f)(x) = 2^{-k} \int_J f(y) \, dy,
\end{equation}
where $J$ is the dyadic subinterval of $[0, 1)$ with length $2^{-k}$
that contains $x$.  Of course, $E_k(f)$ is linear in $f$.

\beginlemma
\label{Some properties of E_k(f)}
(a) For each $f$, $E_k(f)$ is constant on the dyadic subintervals of
$[0,1)$ of length $2^{-k}$.

(b) If $f$ is constant on the dyadic subintervals of $[0,1)$ of length
$2^{-k}$, then $E_k(f) = f$.  

(c) For any $f$, $E_j(E_k(f)) = E_k(f)$ and $E_k(E_j(f)) = E_k(f)$
when $j \ge k$.

(d) If $g$ is a function on $[0,1)$ which is constant on the dyadic
subintervals of $[0,1)$ of length $2^{-k}$, then $E_k(g \, f) = g \,
E_k(f)$ for each $f$.
\end{lemma}

	This is easy to verify, directly from the definitions.  Note
that the first part of (c) holds simply because $E_k(f)$ is constant
on dyadic subintervals of $[0,1)$ of length $2^{-j}$ when $j \ge k$.
In the second part of (c), one is first averaging $f$ on the smaller
dyadic intervals of length $2^{-j}$ to get $E_j(f)$, and then
averaging the result on the larger dyadic intervals of length $2^{-k}$
to get $E_k(E_j(f))$, and the conclusion is that this is the same as
averaging over the dyadic intervals of length $2^{-k}$ directly.

\begindefinition
\label{Dyadic step functions}\index{dyadic step functions}
A function $f$ on $[0,1)$ is a \emph{dyadic step function} if it is a
finite linear combination of indicator functions of dyadic
subintervals of $[0,1)$.
\end{definition}

\beginlemma 
\label{characterizations of dyadic step functions}
Let $f$ be a function on $[0,1)$.  The following are equivalent:

	(a) $f$ is a dyadic step function;

	(b) There is a nonnegative integer $k$ such that $f$ is
constant on every dyadic subinterval of $[0,1)$ of length $2^{-k}$;

	(c) $E_k(f) = f$ for some nonnegative integer $k$, and hence
for all sufficiently large integers $k$.
\end{lemma}

	One can check this using the previous lemma.  From now on, one
is welcome to restrict one's attention to dyadic step functions in
this chapter.

\beginlemma
\label{integrals and E_j}
For any functions $f$, $g$ on $[0,1)$ and nonnegative integer $j$,
\begin{eqnarray}
	\int_{[0,1)} E_j(f)(x) \, g(x) \, dx
		& = & \int_{[0,1)} f(x) \, E_j(g)(x) \, dx		\\
		& = & \int_{[0,1)} E_j(f)(x) \, E_j(g)(x) \, dx.   \nonumber
\end{eqnarray}
\end{lemma}

\beginlemma 
\label{orthogonality properties}
For any functions $f$, $g$ on $[0,1)$ and positive integers $j$, $k$
with $j \ne k$,
\begin{equation}
	\int_{[0,1)} E_0(f)(x) \, (E_j(g)(x) - E_{j-1}(g)(x)) \, dx = 0
\end{equation}
and
\begin{equation}
	\int_{[0,1)} (E_j(f)(x) - E_{j-1}(f)(x)) 
			\, (E_k(g)(x) - E_{k-1}(g)(x)) \, dx = 0.
\end{equation}
\end{lemma}

	The computations for these two lemmas are straightforward and
left to the reader.

	Let $I$ be a dyadic subinterval of $[0,1)$, and let $I_l$ and
$I_r$ be the two dyadic subintervals of $I$ of half the size of $I$.
The Haar function\index{Haar functions} $h_I(x)$ on $[0,1)$ associated
to the interval $I$ is defined by
\begin{eqnarray}
\label{def of h_I(x)}
	h_I(x) & = & - |I|^{1/2} \qquad\hbox{when } x \in I_l
								\\
	       & = & |I|^{1/2} \qquad\enspace\, \hbox{when } x \in I_r
							\nonumber \\
	       & = & 0 \qquad\qquad\enspace\hbox{when } 
					x \in [0,1) \backslash I.
							\nonumber
\end{eqnarray}
Observe that
\begin{equation}
\label{int_{[0,1)} h_I(x) dx = 0}
	\int_{[0,1)} h_I(x) \, dx = 0
\end{equation}
and
\begin{equation}
\label{int_{[0,1)} h_I(x)^2 dx = 1}
	\int_{[0,1)} h_I(x)^2 \, dx = 1.
\end{equation}
In addition, there is a special Haar function $h_0(x)$ on $[0,1)$
defined by $h_0(x) = 1$ for every $x \in [0,1)$, for which we also
have
\begin{equation}
	\int_{[0,1)} h_0(x)^2 \, dx = 1.
\end{equation}

	If $I$ and $J$ are distinct dyadic subintervals of $[0,1)$,
then $h_I$ and $h_J$ satisfy the orthogonality property
\begin{equation}
\label{int_{[0,1)} h_I(x) h_J(x) dx = 0}
	\int_{[0,1)} h_I(x) \, h_J(x) \, dx = 0.
\end{equation}
For if $I$ and $J$ are disjoint, then $h_I(x) \, h_J(x) = 0$ for every
$x \in [0,1)$, and the integral vanishes trivially.  Otherwise, one of
the intervals $I$ and $J$ is contained in the other, and we may as
well assume that $J \subseteq I$, since the two cases are completely
symmetric.  Because $J \ne I$, $J \subseteq I_l$ or $J \subseteq I_r$,
$h_I$ is constant on $J$, and (\ref{int_{[0,1)} h_I(x) h_J(x) dx = 0})
follows from (\ref{int_{[0,1)} h_I(x) dx = 0}).  If $I$ is any dyadic
subinterval of $[0,1)$, then
\begin{equation}
	\int_{[0,1)} h_0(x) \, h_I(x) \, dx = 0,
\end{equation}
by (\ref{int_{[0,1)} h_I(x) dx = 0}).

	For each function $f$ on $[0, 1)$ and nonnegative integer $k$,
\begin{equation}
	E_k(f) = \langle f, h_0 \rangle \, h_0 
		+ \sum_{|I| \ge 2^{-k+1}} \langle f, h_I \rangle \, h_I.
\end{equation}
Here the sum is taken over all dyadic subintervals $I$ of $[0,1)$ with
$|I| \ge 2^{-k+1}$, and is interpreted as being $0$ when $k = 0$.
Also, $\langle f, h_0 \rangle$, $\langle f, h_I \rangle$ are the
integrals of $f$ times $h_0$, $h_I$, respectively.  In particular,
dyadic step functions are finite linear combinations of Haar
functions.  If $f$ is a dyadic step function on $[0,1)$, then
\begin{equation}
	f = \langle f, h_0 \rangle \, h_0 
		+ \sum_I \langle f, h_I \rangle \, h_I,
\end{equation}
where the sum is taken over all dyadic subintervals $I$ of $[0,1)$.
The sum is actually a finite sum, since $\langle f, h_I \rangle = 0$
for all but finitely many $I$.  This expression for $f$ follows from
the orthonormality conditions for the Haar functions described
earlier.

\section{Maximal functions}
\label{maximal functions}

	As mentioned in the previous section, one is welcome to
restrict one's attention to real or complex-valued functions on $[0,
1)$ that are dyadic step functions in this chapter.  The \emph{dyadic
maximal function}\index{maximal functions} $M(f)$ associated to a
function $f$ on $[0,1)$ is defined by
\begin{equation}
\label{def of M(f)}
	M(f)(x) = \sup_{k \ge 0} |E_k(f)(x)|.
\end{equation}
Equivalently, $M(f)(x)$ is equal to
\begin{equation}
\label{def of M(f), 2}
 \qquad \sup \bigg\{ \biggl| \frac{1}{|J|} \int_J f(y) \, dy \biggr|
	 : \hbox{ $J$ is a dyadic subinterval of $[0,1)$ and $x \in J$}\bigg\}.
\end{equation}
For each nonnegative integer $l$, put
\begin{equation}
\label{def of M_l(f)}
	M_l(f)(x) = \max_{0 \le k \le l} |E_k(f)(x)|,
\end{equation}
which is the same as
\begin{equation}
\label{def of M_l(f), 2}
 \qquad	M_l(f)(x) = 
 \max \bigg\{ \biggl| \frac{1}{|J|} \int_J f(y) \, dy \biggr|
   : \ J \subseteq [0,1), \ x \in J, \hbox{ and } |J| \ge 2^{-l}\bigg\},
\end{equation}
where the maximum is again taken over dyadic subintervals $J$ of $[0,
1)$.  Thus
\begin{equation}
\label{M_l(f) le M_r(f) when r ge l}
	M_l(f) \le M_r(f)  \quad\hbox{when $r \ge l$}
\end{equation}
and
\begin{equation}
\label{M(f)(x) = sup_{l ge 0} M_l(f)(x) for every x in [0,1)}
	M(f)(x) = \sup_{l \ge 0} M_l(f)(x) \quad\hbox{for every $x \in [0,1)$}.
\end{equation}

	For any pair of functions $f_1$, $f_2$ on $[0,1)$,
\begin{equation}
	M(f_1 + f_2) \le M(f_1) + M(f_2)
\end{equation}
and
\begin{equation}
	M_l(f_1 + f_2) \le M_l(f_1) + M_l(f_2)
\end{equation}
for each $l \ge 0$.  Also,
\begin{equation}
	M(c \, f) = |c| \, M(f)
\end{equation}
and
\begin{equation}
	M_l(c \, f) = |c| \, M_l(f)
\end{equation}
for any function $f$ and constant $c$.  Thus $M(f)$, $M_l(f)$ are
sublinear in $f$.

\beginlemma
\label{M(f), M_l(f), E_l(f), etc.}
If $f$ is constant on the dyadic subintervals of $[0,1)$ of length
$2^{-l}$, then $M(f)$ is constant on the dyadic subintervals of
$[0,1)$ of length $2^{-l}$, and $M(f) = M_l(f)$.  For any function
$f$,
\begin{equation}
\label{M_l(f) = M(E_l(f))}
	M_l(f) = M(E_l(f)),
\end{equation}
and $M_l(f)$ is constant on dyadic intervals of length $2^{-l}$.
\end{lemma}

	Exercise.

\begincorollary
\label{M(f) for dyadic step functions}
If $f$ is a dyadic step function, then $M(f)$ is too, and $M(f) =
M_j(f)$ for sufficiently large $j$.
\end{corollary}

\beginlemma [Supremum bound for M(f)]
\label{supremum bound for M(f)}
If $|f(x)| \le A$ for some $A \ge 0$ and every $x \in [0,1)$, then
$M(f)(x) \le A$ for every $x \in [0,1)$.
\end{lemma}

	This is an easy consequence of the definitions.  Lemma
\ref{supremum bound for M(f)} also works if $|f(x)| \le A$ for every
$x \in [0,1)$ except for a small set that does not affect the
integrals.

\beginproposition [Weak-type estimate for M(f)]
\label{weak-type estimate for M(f)}
For every $\lambda > 0$,
\begin{equation}
\label{weak-type inequality for M(f)}
	|\{x \in [0,1) : M(f)(x) > \lambda\}| 
		\le \frac{1}{\lambda} \int_{[0,1)} |f(w)| \, dw.
\end{equation}
\end{proposition}

	The left-hand side of (\ref{weak-type inequality for M(f)})
refers to the \emph{measure} of the set in question, the meaning of
which is clarified by the proof.

	Let $\lambda > 0$ be given, let $\mathcal{F}$ be the
collection of dyadic intervals subintervals $L$ of $[0,1)$ such that
\begin{equation}
\label{|frac{1}{|L|} int_L f(y) dy| > lambda}
	\biggl| \frac{1}{|L|} \int_L f(y) \, dy \biggr| > \lambda,
\end{equation}
and let us check that
\begin{equation}
\label{{x in [0,1) : M(f)(x) > lambda} = bigcup_{L in mathcal{F}} L}
	\{x \in [0,1) : M(f)(x) > \lambda\} = \bigcup_{L \in \mathcal{F}} L.
\end{equation}
If $x \in [0,1)$ and $M(f)(x) > \lambda$, then there is a dyadic
interval $L$ in $[0,1)$ such that $x \in L$ and $L$ satisfies
(\ref{|frac{1}{|L|} int_L f(y) dy| > lambda}), because of (\ref{def of
M(f), 2}), and hence the left side of (\ref{{x in [0,1) : M(f)(x) >
lambda} = bigcup_{L in mathcal{F}} L}) is contained in the right side
of (\ref{{x in [0,1) : M(f)(x) > lambda} = bigcup_{L in mathcal{F}}
L}).  Conversely, if $L \in \mathcal{F}$, then
\begin{equation}
	M(f)(x) \ge \biggl| \frac{1}{|L|} \int_L f(y) \, dy \biggr|
			> \lambda
\end{equation}
for every $x \in L$, and therefore $L$ is contained in the left side
of (\ref{{x in [0,1) : M(f)(x) > lambda} = bigcup_{L in mathcal{F}}
L}).  Thus the right side of (\ref{{x in [0,1) : M(f)(x) > lambda} =
bigcup_{L in mathcal{F}} L}) is contained in the left side, and
(\ref{{x in [0,1) : M(f)(x) > lambda} = bigcup_{L in mathcal{F}} L})
follows.

	As in Lemma \ref{unions of dyadic intervals}, if
$\mathcal{F}_0$ consists of the maximal elements of $\mathcal{F}$,
then
\begin{equation}
\label{bigcup_{L in mathcal{F}_0} L = bigcup_{L in mathcal{F}} L}
	\bigcup_{L \in \mathcal{F}_0} L = \bigcup_{L \in \mathcal{F}} L,
\end{equation}
and the intervals in $\mathcal{F}_0$ are pairwise disjoint.  Thus
\begin{equation}
\label{|{x in [0,1) : M(f)(x) > lambda}| = sum_{L in mathcal{F}_0} |L|}
	|\{x \in [0,1) : M(f)(x) > \lambda\}| = \sum_{L \in \mathcal{F}_0} |L|.
\end{equation}
If $f$ is a dyadic step function which is constant on the dyadic
subintervals of $[0, 1)$ of length $2^{-l}$, then $M(f)$ is constant
on the dyadic intervals of length $2^{-l}$, and the elements of
$\mathcal{F}_0$ have length $\ge 2^{-l}$.

	Each interval $L \in \mathcal{F}_0$ satisfies
(\ref{|frac{1}{|L|} int_L f(y) dy| > lambda}), which gives
\begin{equation}
\label{|L| < frac{1}{lambda} | int_L f(y) dy|}
	|L| < \frac{1}{\lambda} \biggl| \int_L f(y) \, dy \biggr|
		\le \frac{1}{\lambda} \int_L |f(y)| \, dy,
\end{equation}
and hence
\begin{equation}
\label{sum_{L in mathcal{F}_0} |L| < ...}
	\sum_{L \in \mathcal{F}_0} |L| 
 < \sum_{L \in \mathcal{F}_0} \frac{1}{\lambda}  \int_L |f(y)| \, dy
 = \frac{1}{\lambda} \int_{\bigcup_{L \in \mathcal{F}_0} L} |f(y)| \, dy,
\end{equation}
using the disjointness of the intervals $L \in \mathcal{F}_0$.
Therefore 
\begin{equation}
\label{|{x in [0,1) : M(f)(x) > lambda}| le ...}
	|\{x \in [0,1) : M(f)(x) > \lambda\}| \le 
  \frac{1}{\lambda} \int_{\{y \in [0,1) : M(f)(y) > \lambda\}} |f(y)| \, dy,
\end{equation}
which implies (\ref{weak-type inequality for M(f)}).

\beginlemma
\label{modified weak-type estimate for M(f)}
For each $\lambda > 0$,
\begin{equation}
\label{modified weak-type inequality}
	|\{x \in [0,1): M(f)(x) > 2 \, \lambda\}|
 \le \frac{1}{\lambda} \int_{\{u \in [0,1) : |f(u)| > \lambda\}}  |f(u)| \, du.
\end{equation}
\end{lemma}

	Let $\lambda > 0$ be given, and put
\begin{equation}
	f_1(x) = f(x) \enspace\hbox{when}\enspace |f(x)| \le \lambda,
	\quad f_1(x) = 0 \enspace\hbox{when}\enspace |f(x)| > \lambda,
\end{equation}
and
\begin{equation}
	f_2(x) = f(x) \enspace\hbox{when}\enspace |f(x)| > \lambda,
	\quad f_2(x) = 0 \enspace\hbox{when}\enspace |f(x)| \le \lambda.
\end{equation}
Thus $f(x) = f_1(x) + f_2(x)$ and $M(f_1)(x) \le \lambda$ for every $x
\in [0, 1)$.  This implies that
\begin{equation}
	M(f)(x) \le \lambda + M(f_2)(x)
\end{equation}
for every $x$, and hence
\begin{equation}
\label{|{x : M(f)(x) > 2 lambda}| le |{x : M(f_2)(x) > lambda}|}
	|\{x \in [0,1) : M(f)(x) > 2 \, \lambda\}|
		\le |\{x \in [0,1) : M(f_2)(x) > \lambda\}|.
\end{equation}
We can apply Proposition \ref{weak-type estimate for M(f)} with $f$
replaced by $f_2$ to get that
\begin{equation}
	|\{x \in [0,1) : M(f_2)(x) > \lambda\}| 
		\le \frac{1}{\lambda} \int_{[0,1)} |f_2(u)| \, du,
\end{equation}
and the lemma follows.
	
\beginlemma
\label{integrals from distribution functions}
If $g(x)$ is a nonnegative real-valued function on $[0,1)$, and $p$ is
a positive real number, then
\begin{equation}
\label{formula for integrals, distribution function}
	\int_{[0,1)} g(x)^p \, dx = 
		\int_0^\infty p \, \lambda^{p-1} \, 
			|\{x \in [0,1): g(x) > \lambda\}| \, d\lambda.
\end{equation}
\end{lemma}

	One can see this by integrating $p \, \lambda^{p-1}$ on the
set
\begin{equation}
	\{(x,\lambda) \in {\bf R}^2 : x \in [0,1), 0 < \lambda < g(x)\}
\end{equation}
first in $\lambda$, and then in $x$, and first in $x$, and then in
$\lambda$.

\beginproposition
\label{p-integral inequalities for M(f)}
For each real number $p > 1$,
\begin{equation}
\label{the p-integral inequality}
	\int_{[0,1)} M(f)(x)^p \, dx 
		\le \frac{2^p \, p}{p-1} \int_{[0,1)} |f(y)|^p \, dy.
\end{equation}
\end{proposition}

	To prove this, we apply Lemma \ref{integrals from distribution
functions} with $g = M(f)$ to get that
\begin{equation}
	\int_{[0,1)} M(f)(x)^p \, dx 
	     =  \int_0^\infty p \, \lambda^{p-1} 
			|\{x \in [0,1): M(f)(x) > \lambda\}| \, d\lambda.
\end{equation}
By (\ref{modified weak-type inequality}) with $\lambda$ replaced by
$\lambda/2$,
\begin{eqnarray}
 \qquad \int_{[0,1)} M(f)(x)^p \, dx
	     & \le & \int_0^\infty p \, \lambda^{p-1} 
   \Big(\frac{2}{\lambda} \int_{\{u \in [0,1) : |f(u)| > \lambda/2\}}  
			|f(u)| \, du \Big) \, d\lambda	   	\\
	     & = & \int_0^\infty \int_{\{u \in [0,1) : |f(u)| > \lambda/2\}} 
				2p \, \lambda^{p-2} |f(u)| \, du \, d\lambda.
							\nonumber
\end{eqnarray}
Interchanging the order of integration leads to
\begin{equation}
	\int_{[0,1)} M(f)(x)^p \, dx
  \le \int_{[0,1)} \int_0^{2 \, |f(u)|} 
			2p \, \lambda^{p-2} |f(u)| \, d\lambda \, du.
\end{equation}
Because $p > 1$,
\begin{eqnarray}
	\int_{[0,1)} M(f)(x)^p \, dx					
  & \le & \int_{[0,1)} 2p \, (p-1)^{-1} (2 \, |f(u)|)^{p-1} |f(u)|  \, du
								         \\
  &  =  & \frac{2^p \, p}{p-1} \int_{[0,1)} |f(u)|^p \, du,	\nonumber
\end{eqnarray}
as desired.

\section{Square functions}
\label{square functions}

	The \emph{dyadic square function} \index{square functions}
$S(f)$ associated to a functon $f$ on $[0,1)$ is defined by
\begin{equation}
\label{def of S(f)}
	S(f)(x) = \Big(|E_0(f)(x)|^2 +
                   \sum_{j=1}^\infty |E_j(f)(x) - E_{j-1}(f)(x)|^2 \Big)^{1/2}.
\end{equation}
For each nonnegative integer $l$, put
\begin{equation}
\label{def of S_l(f)}
	S_l(f)(x) = 
	   \Big(|E_0(f)(x)|^2 + 
		\sum_{j=1}^l |E_j(f)(x) - E_{j-1}(f)(x)|^2\Big)^{1/2},
\end{equation}
where the sum on the right side is interpreted as being $0$ when $l =
0$.  Thus
\begin{equation}
	S_l(f)(x) \le S_p(f)(x) \quad\hbox{when $l \le p$}
\end{equation}
and $S(f)(x) = \sup_{l \ge 0} \, S_l(f)(x)$.  It is easy to see that
$S(f)$, $S_l(f)$ are sublinear in $f$, in the sense that
\begin{equation}
	S(f_1 + f_2) \le S(f_1) + S(f_2)
\end{equation}
and $S(c \, f) = |c| \, S(f)$ for all functions $f_1$, $f_2$, and $f$
on $[0,1)$ and all constants $c$, and similarly for $S_l(f)$.

\beginlemma
\label{S(f), S_l(f), E_l(f), etc.}
If $f$ is constant on the dyadic subintervals of $[0,1)$ of length
$2^{-l}$, then $S(f)$ is constant on the dyadic subintervals of
$[0,1)$ of length $2^{-l}$, and $S(f) = S_l(f)$.  For any $f$,
\begin{equation}
	S_l(f) = S(E_l(f)),
\end{equation}
and $S_l(f)$ is constant on dyadic intervals of length $2^{-l}$.
\end{lemma}

	Exercise.  

\begincorollary
If $f$ is a dyadic step function on $[0,1)$, then $S(f)$ is too, and
$S(f) = S_j(f)$ for sufficiently large $j$.
\end{corollary}

\beginlemma
\label{2-norm of S(f)}
For any function $f$ on $[0,1)$, 
\begin{equation}
\label{int S_l(f)^2 = int |E_l(f)|^2}
	\int_{[0,1)} S_l(f)(x)^2 \, dx = \int_{[0,1)} |E_l(f)(x)|^2 \, dx
\end{equation}
for every $l \ge 0$, and
\begin{equation}
\label{int S(f)^2 = int |f|^2}
	\int_{[0,1)} S(f)(x)^2 \, dx = \int_{[0,1)} |f(x)|^2 \, dx.
\end{equation}
\end{lemma}

	Of course
\begin{equation}
	E_l(f) = E_0(f) + \sum_{j=1}^l (E_j(f) - E_{j-1}(f)),
\end{equation}
and to prove the lemma one uses the orthogonality conditions in Lemma
\ref{orthogonality properties}.

\section{Estimates, 1}
\label{Estimates, 1}

\beginproposition
\label{p-norm of S(f) bounded by p-norm of M(f), p < 2}
If $0 < p < 2$, then there is a positive real number $C_1(p)$ such
that
\begin{equation}
\label{int S(f)^p le C_1(p) int M(f)^p}
	\int_{[0,1)} S(f)(x)^p \, dx
		\le C_1(p) \int_{[0,1)} M(f)(x)^p \, dx
\end{equation}
for any function $f$ on $[0,1)$.
\end{proposition}

	Let $p < 2$, a function $f$ on $[0, 1)$, and $\lambda > 0$ be
given, and consider
\begin{equation}
	|\{x \in [0,1) : S(f)(x) > \lambda\}|.
\end{equation}
Let $\mathcal{F}$ denote the set of dyadic subintervals $J$
of $[0,1)$ such that
\begin{equation}
	\frac{1}{|J|} \biggl| \int_J f(y) \, dy \biggr| > \lambda.
\end{equation}
If  $\mathcal{F}_0$ is the set of maximal intervals in
$\mathcal{F}$, then
\begin{equation}
\label{union mathcal{F}_0}
	\bigcup_{J \in \mathcal{F}_0} J = \bigcup_{J \in \mathcal{F}} J
\end{equation}
and
\begin{equation}
\label{disjointness mathcal{F}_0}
	J_1 \cap J_2 = \emptyset 
		\quad\hbox{when } J_1, J_2 \in \mathcal{F}_0, \ J_1 \ne J_2,
\end{equation}
as in Lemma \ref{unions of dyadic intervals}.

	Suppose that $[0,1)$ is not an element of $\mathcal{F}_0$, and
let $\mathcal{F}_1$ be the set of dyadic subintervals $L$ of $[0,1)$
for which there is a $J \in \mathcal{F}_0$ such that
\begin{equation}
\label{J subseteq L and |J| = |L|/2}
	J \subseteq L \hbox{ and } |J| = |L|/2.
\end{equation}
Because $\mathcal{F}_0$ consists of maximal intervals in $\mathcal{F}$,
each $L$ in $\mathcal{F}_1$ does not lie in $\mathcal{F}$, and hence
\begin{equation}
\label{frac{1}{|L|} |int_L f(y) dy| le lambda}
	\frac{1}{|L|} \biggl| \int_L f(y) \, dy \biggr| \le \lambda
\end{equation}
for every $L \in \mathcal{F}_1$.

	The elements of $\mathcal{F}_1$ need not be disjoint, and so
we let $\mathcal{F}_{10}$ be the set of maximal elements of
$\mathcal{F}_1$.  As usual,
\begin{equation}
	\bigcup_{L \in \mathcal{F}_{10}} L = \bigcup_{L \in \mathcal{F}_1} L
\end{equation}
and
\begin{equation}
	L_1 \cap L_2 = \emptyset 
		\quad\hbox{when } L_1, L_2 \in \mathcal{F}_{10}, \ L_1 \ne L_2.
\end{equation}
Let $f_\lambda(x)$ be the function on $[0,1)$ defined by
\begin{eqnarray}
\label{def of f_lambda(x)}
	f_\lambda(x) & = & \frac{1}{|L|} \int_L f(y) \, dy
		   \quad\hbox{when } x \in L, \ L \in \mathcal{F}_{10},
									\\
		     & = & f(x)  \quad\qquad\qquad
 		\, \hbox{when } x \in [0,1) \backslash 
			\Big(\bigcup_{I \in \mathcal{F}_{10}} I \Big).
								\nonumber
\end{eqnarray}

\beginlemma 
\label{averages of f and averages of f_lambda}
If $K$ is a dyadic subinterval of $[0,1)$ such that
\begin{equation}
 K \backslash \Big(\bigcup_{I \in \mathcal{F}_{10}} I \Big) \ne \emptyset
\end{equation}
or $L \subseteq K$ for some $L \in \mathcal{F}_{10}$, then
\begin{equation}
	\frac{1}{|K|} \int_K f(u) \, du 
		= \frac{1}{|K|} \int_K f_\lambda(u) \, du.
\end{equation}
\end{lemma}

	Under these conditions, $K$ is the disjoint union of the
intervals $L \in \mathcal{F}_{10}$ such that $L \subseteq K$ and $K
\backslash \big(\bigcup_{I \in \mathcal{F}_{10}} I \big)$.  The
integral of $f$ over $K$ is equal to the sum of the integrals of $f$
over these sets, which is the same as the integral of $f_\lambda$ over $K$.

\begincorollary
\label{S(f)(x) = S(f_lambda)(x) on a certain set}
If $x \in [0,1) \backslash \big(\bigcup_{L \in \mathcal{F}_{10}} L
\big)$, then $S(f)(x) = S(f_\lambda)(x)$.
\end{corollary}

	For these $x$'s, Lemma \ref{averages of f and averages of
f_lambda} implies that $E_j(f)(x) = E_j(f_\lambda)(x)$ for every
nonnegative integer $j$, and hence $S(f)(x) = S(f_\lambda)(x)$.

	Using the corollary, it is easy to see that
\begin{eqnarray}
\label{|{x : S(f)(x) > lambda}| le ....}
\lefteqn{|\{x \in [0,1) : S(f)(x) > \lambda\}|}		\\
	& & \qquad \le \sum_{L \in \mathcal{F}_{10}} |L| 
			+ |\{x \in [0,1) : S(f_\lambda)(x) > \lambda\}|.
								\nonumber
\end{eqnarray}
For each $L \in \mathcal{F}_{10}$, there is a $J \in \mathcal{F}_0$
such that $J \subseteq L$ and $|J| = |L|/2$, and this leads to
\begin{equation}
  \sum_{L \in \mathcal{F}_{10}} |L| \le 2 \sum_{J \in \mathcal{F}_0} |J|.
\end{equation}
By (\ref{union mathcal{F}_0}) and (\ref{disjointness mathcal{F}_0}),
\begin{equation}
	\sum_{L \in \mathcal{F}_{10}} |L|
		 \le 2 \, \biggl| \bigcup_{J \in \mathcal{F}_0} J \biggr|
		  = 2 \, \biggl| \bigcup_{J \in \mathcal{F}} J \biggr|.
\end{equation}
As in (\ref{{x in [0,1) : M(f)(x) > lambda} = bigcup_{L in mathcal{F}}
L}),
\begin{equation}
	\bigcup_{J \in \mathcal{F}} J = \{x \in [0,1) : M(f)(x) > \lambda\}.
\end{equation}
Therefore
\begin{equation}
	\sum_{L \in \mathcal{F}_{10}} |L| 
		\le 2 \, |\{x \in [0,1) : M(f)(x) > \lambda\}|,
\end{equation}
and hence
\begin{eqnarray}
\label{|{x : S(f)(x) > lambda}| le ...,2}
\lefteqn{\quad |\{x \in [0,1) : S(f)(x) > \lambda\}|}		\\
	& & \le 2 \, |\{x \in [0,1) : M(f)(x) > \lambda\}|
			+ |\{x \in [0,1) : S(f_\lambda)(x) > \lambda\}|.
								\nonumber
\end{eqnarray}
By Lemma \ref{2-norm of S(f)},
\begin{equation}
 \qquad \lambda^2 \, |\{x \in [0,1) : S(f_\lambda)(x) > \lambda\}|
		\le \int_{[0,1)} S(f_\lambda)(x)^2 \, dx
		 =   \int_{[0,1)} |f_\lambda(x)|^2 \, dx.
\end{equation}
Thus
\begin{eqnarray}
\label{|{x : S(f)(x) > lambda}| le ...,3}
\lefteqn{|\{x \in [0,1) : S(f)(x) > \lambda\}|}		\\
	& & \le 2 \, |\{x \in [0,1) : M(f)(x) > \lambda\}|
			+ \lambda^{-2} \int_{[0,1)} |f_\lambda(x)|^2 \, dx.
								\nonumber
\end{eqnarray}

\beginlemma
\label{|f_lambda| le min(lambda, M(f))}
$|f_\lambda| \le \min (\lambda, M(f))$.
\end{lemma}

	Indeed, for each dyadic subinterval $I$ of $[0,1)$, we have that
\begin{equation}
	\frac{1}{I} \biggl|\int_I f(y) \, dy \biggr| \le M(f)(x)
\end{equation}
automatically when $x \in I$, and
\begin{equation}
	\frac{1}{I} \biggl|\int_I f(y) \, dy \biggr| \le \lambda
\end{equation}
when $I \in \mathcal{F}_1$ or $x \in I \backslash \big(\bigcup_{L \in
\mathcal{F}_{10}} L \big)$, since $I \not\in \mathcal{F}$ in these two
cases.  With the help of Lemma \ref{averages of f and averages of
f_lambda}, one can actually get the stronger estimate
\begin{equation}
	M(f_\lambda) \le \min (\lambda, M(f)).
\end{equation}

	Because of the lemma, we may replace (\ref{|{x : S(f)(x) >
lambda}| le ...,3}) with
\begin{eqnarray}
\label{|{x : S(f)(x) > lambda}| le ...,4}
\lefteqn{\quad\enspace |\{x \in [0,1) : S(f)(x) > \lambda\}|}		\\
	& & \le 2 \, |\{x \in [0,1) : M(f)(x) > \lambda\}|
		+ \lambda^{-2} \int_{[0,1)} \min(\lambda, M(f)(x))^2 \, dx.
								\nonumber
\end{eqnarray}
At the beginning of this argument, just before (\ref{J subseteq L and
|J| = |L|/2}), we assumed that $[0,1)$ is not an element of
$\mathcal{F}_0$.  If $[0,1)$ is an element of $\mathcal{F}_0 \subseteq
\mathcal{F}$, then
\begin{equation}
	\biggl|\int_{[0,1)} f(y) \, dy \biggr| > \lambda,
\end{equation}
$M(f)(x) > \lambda$ for every $x \in [0, 1)$, and hence
\begin{equation}
	|\{x \in [0,1) : S(f)(x) > \lambda\}|
		\le |\{x \in [0,1) : M(f)(x) > \lambda\}|.
\end{equation}
Thus (\ref{|{x : S(f)(x) > lambda}| le ...,4}) holds in general.  By
Lemma \ref{integrals from distribution functions},
\begin{equation}
\label{formula for int S(f)^p}
	\int_{[0,1)} S(f)(x)^p \, dx 
		= \int_0^\infty p \, \lambda^{p-1} 
			|\{x \in [0,1): S(f)(x) > \lambda\}| \, d\lambda
\end{equation}
and
\begin{equation}
\label{formula for int M(f)^p}
	\int_{[0,1)} M(f)(x)^p \, dx  
		= \int_0^\infty p \, \lambda^{p-1} 
			|\{x \in [0,1): M(f)(x) > \lambda\}| \, d\lambda.
\end{equation}
Therefore
\begin{eqnarray}
\label{int_{[0,1)} S(f)(x)^p dx le ....}
\lefteqn{\quad \int_{[0,1)} S(f)(x)^p \, dx}	\\
	& & \le 2 \int_{[0,1)} M(f)(x)^p \, dx 
	      + \int_0^\infty p \, \lambda^{p-3} 
		     \int_{[0,1)} \min(\lambda, M(f)(x))^2 \, dx \, d\lambda.
								\nonumber
\end{eqnarray}
We can interchange the order of integration and replace the second
term on the right side of (\ref{int_{[0,1)} S(f)(x)^p dx le ....})
with
\begin{equation}
	\int_{[0,1)} \int_0^\infty p \, \lambda^{p-3}
		\min(\lambda, M(f)(x))^2 \, d\lambda \, dx.
\end{equation}
The integral in $\lambda$ can be computed exactly, since
\begin{equation}
	\int_{M(f)(x)}^\infty p \, \lambda^{p-3} \, M(f)(x)^2 \, d\lambda
		= \frac{p}{2-p} M(f)(x)^p
\end{equation}
and
\begin{equation}
	\int_0^{M(f)(x)} p \, \lambda^{p-3} \, \lambda^2 \, d\lambda
		= M(f)(x)^p.
\end{equation}
Proposition \ref{p-norm of S(f) bounded by p-norm of M(f), p < 2} now
follows by using these formulae in (\ref{int_{[0,1)} S(f)(x)^p dx le
....}).

	The coefficient $(2 - p)^{-1}$ in the previous computations is not very
nice, and one can get bounded constants for $p$ near $2$ using interpolation 
arguments.  One can also start with an estimate for $p = 4$ instead of $p = 2$,
as in Section \ref{another argument for p = 4}, and use the same method as here
to get estimates for $0 < p < 4$ that remain bounded for $p$ near $2$.

\beginproposition [Weak-type estimate for $S(f)$]
\label{weak-type estimate for S(f)}
For any function $f$ on $[0,1)$ and $\lambda > 0$,
\begin{equation}
\label{|{ x in [0,1) : S(f) > lambda }| le ...}
	|\{ x \in [0,1) : S(f) > \lambda \}| 
		\le \frac{3}{\lambda} \int_{[0,1)} |f(x)| \, dx.
\end{equation}
\end{proposition}

	This follows from practically the same arguments as above.  By
(\ref{|{x : S(f)(x) > lambda}| le ...,3}) and Lemma \ref{|f_lambda| le
min(lambda, M(f))},
\begin{eqnarray}
\label{|{x : S(f)(x) > lambda}| le ...,3, revised version}
\lefteqn{|\{x \in [0,1) : S(f)(x) > \lambda\}|}		\\
	& & \le 2 \, |\{x \in [0,1) : M(f)(x) > \lambda\}|
			+ \lambda^{-1} \int_{[0,1)} |f_\lambda(x)| \, dx,
								\nonumber
\end{eqnarray}
where $f_\lambda(x)$ is as in (\ref{def of f_lambda(x)}).  To get
(\ref{|{ x in [0,1) : S(f) > lambda }| le ...}), one can use
Proposition \ref{weak-type estimate for M(f)} and the observation that
\begin{equation}
\label{int_{[0,1)} |f_lambda(x)| dx le int_{[0,1)} |f(x)| dx}
	\int_{[0,1)} |f_\lambda(x)| \, dx \le \int_{[0,1)} |f(x)| \, dx.
\end{equation}

\section{Estimates, 2}
\label{Estimates, 2}

\beginproposition
\label{p-norm of M(f) bounded by p-norm of S(f), p < 2}
If $0 < p < 2$, then there is a positive real number $C_2(p)$ such
that
\begin{equation}
\label{int M(f)^p le C_2(p) int S(f)^p}
	\int_{[0,1)} M(f)(x)^p \, dx 
		\le C_2(p) \int_{[0,1)} S(f)(x)^p \, dx
\end{equation}
for any function $f$ on $[0,1)$.
\end{proposition}

	Let $p < 2$ and $f$ be given, and let $\lambda > 0$ be a
positive real number.

\beginlemma
\label{{x in [0,1) : S(f)(x) > lambda} is a union of dyadic intervals}
The set
\begin{equation}
\label{S(f) > lambda}
	\{x \in [0,1) : S(f)(x) > \lambda\}
\end{equation}
is a union of dyadic subintervals of $[0,1)$.
\end{lemma}

	If $w \in [0,1)$ and
\begin{equation}
 \quad 	S(f)(w) = \Big(|E_0(f)(w)|^2 + 
		\sum_{j=1}^\infty |E_j(f)(w) - E_{j-1}(f)(w)|^2 \Big)^{1/2}
			> \lambda,
\end{equation}
then
\begin{equation}
	\Big(|E_0(f)(w)|^2 + 
		\sum_{j=1}^l |E_j(f)(w) - E_{j-1}(f)(w)|^2 \Big)^{1/2}
			> \lambda
\end{equation}
for some $l$.  Let $I$ be the dyadic subinterval of $[0,1)$ such that
$|I| = 2^{-l}$ and $w \in I$.  Because $E_j(f)$ is constant on dyadic
intervals of length $2^{-j}$, $E_j(f)(y) = E_j(f)(w)$ when $j \le l$
and $y \in I$, and hence
\begin{equation}
	\Big(|E_0(f)(y)|^2 + 
		\sum_{j=1}^l |E_j(f)(y) - E_{j-1}(f)(y)|^2 \Big)^{1/2}
			> \lambda
\end{equation}
for every $y \in I$.  Therefore $S(f)(y) > \lambda$ for every $y \in
I$, and the lemma follows easily.

	Let $\mathcal{G}_0$ be the collection of maximal dyadic
subintervals of $[0,1)$ contained in the set (\ref{S(f) > lambda}).
As usual,
\begin{equation}
\label{bigcup_{J in mathcal{G}_0} J = {x in [0,1) : S(f)(x) > lambda}}
	\bigcup_{J \in \mathcal{G}_0} J = \{x \in [0,1) : S(f)(x) > \lambda\},
\end{equation}
and the intervals in $\mathcal{G}_0$ are pairwise disjoint.  In
particular,
\begin{equation}
\label{sum_{J in mathcal{G}_0} |J| = |{x in [0,1) : S(f)(x) > lambda}|}
	\sum_{J \in \mathcal{G}_0} |J| = |\{x \in [0,1) : S(f)(x) > \lambda\}|.
\end{equation}

	Suppose that (\ref{S(f) > lambda}) is not equal to the whole
unit interval $[0,1)$.  Let $\mathcal{G}_1$ be the collection of
dyadic subintervals $L$ of $[0,1)$ for which there is a $J \in
\mathcal{G}_0$ such that
\begin{equation}
\label{J subseteq L and |J| = |L| / 2}
	J \subseteq L \hbox{ and } |J| = |L| / 2.
\end{equation}
Because the elements of $\mathcal{G}_0$ are maximal dyadic intervals
contained in (\ref{S(f) > lambda}), each $L \in \mathcal{G}_1$ is not
a subset of (\ref{S(f) > lambda}).  Thus for each $L \in
\mathcal{G}_1$ there is a point $u \in L$ such that $S(f)(u) \le
\lambda$.  If $\ell(L)$ denotes the nonnegative integer such that
$2^{-\ell(L)} = |L|$, then
\begin{equation}
\label{... le lambda}
	\Big(|E_0(f)(u)|^2 + 
	    \sum_{j=1}^{\ell(L)} |E_j(f)(u) - E_{j-1}(f)(u)|^2 \Big)^{1/2}
			\le \lambda,
\end{equation}
where the sum on the left is interpreted as being $0$ if $\ell(L) =
0$.  More precisely, this inequality holds for at least one $u \in L$,
and hence at every $u \in L$, because $E_j(f)$ is constant on $L$ when
$j \le \ell(L)$.

	The intervals in $\mathcal{G}_1$ need not be pairwise
disjoint, and we can pass to the subcollection $\mathcal{G}_{10}$ of
maximal elements of $\mathcal{G}_1$ to get
\begin{equation}
\label{bigcup_{L in mathcal{G}_{10}} L = bigcup_{L in mathcal{G}_1} L}
	\bigcup_{L \in \mathcal{G}_{10}} L
		= \bigcup_{L \in \mathcal{G}_1} L
\end{equation}
and $L_1 \cap L_2 = \emptyset$ when $L_1, L_2 \in \mathcal{G}_1$ and
$L_1 \ne L_2$.  The definition of $\mathcal{G}_1$ implies that
$\bigcup_{J \in \mathcal{G}_0} J \subseteq \bigcup_{L \in
\mathcal{G}_1} L$, and therefore
\begin{equation}
\label{{x in [0,1) : S(f)(x) > lambda} subseteq cup_{L in mathcal{G}_{10}} L}
	\{x \in [0,1) : S(f)(x) > \lambda\}
		\subseteq \bigcup_{L \in \mathcal{G}_{10}} L .
\end{equation}
Also,
\begin{equation}
\label{sum_{L in mathcal{G}_{10}} |L| le ....}
	\sum_{L \in \mathcal{G}_{10}} |L| 
		\le \sum_{J \in \mathcal{G}_0} 2 \, |J|
		= 2 \, |\{x \in [0,1) : S(f)(x) > \lambda\}|.
\end{equation}

	Let $g_\lambda(x)$ be the function defined on $[0,1)$ by
\begin{eqnarray}
\label{def of g_lambda}
	g_\lambda(x) & = & \frac{1}{|L|} \int_L f(y) \, dy
			\quad\hbox{when } x \in L, \ L \in \mathcal{G}_{10}  \\
		     & = & f(x)
			\qquad\qquad\quad \hbox{ when } 
    x \in [0,1) \backslash \Big(\bigcup_{I \in \mathcal{G}_{10}} I \Big).
								\nonumber
\end{eqnarray}

\beginlemma
\label{frac{1}{|K|} int_K g_lambda = frac{1}{|K|} int_K g when ....}
If $K$ is a dyadic subinterval of $[0,1)$ such that
\begin{equation}
 K \backslash \Big(\bigcup_{I \in \mathcal{G}_{10}} I \Big) \ne \emptyset
\end{equation}
or $L \subseteq K$ for some $L \in \mathcal{G}_{10}$, then
\begin{equation}
	\frac{1}{|K|} \int_K g_\lambda(u) \, du
		= \frac{1}{|K|} \int_K f(u) \, du.
\end{equation}
\end{lemma}

	This uses the fact that $K$ is the disjoint union of the $L
\in \mathcal{G}_{10}$ with $L \subseteq K$ and $K \backslash
\big(\bigcup_{I \in \mathcal{G}_{10}} I \big)$, as in Lemma
\ref{averages of f and averages of f_lambda}.

\begincorollary
\label{M(f) = M(g_lambda) on the proper set}
If $x \in [0,1) \backslash \big(\bigcup_{I \in \mathcal{G}_{10}} I
\big)$, then $M(f)(x) = M(g_\lambda)(x)$.
\end{corollary}

\begincorollary
\label{S(f), S(g_lambda), etc.}
If $x \in [0,1) \backslash \big(\bigcup_{I \in \mathcal{G}_{10}} I
\big)$, then $S(f)(x) = S(g_\lambda)(x)$.  If $L \in
\mathcal{G}_{10}$, $v \in L$, and $2^{-\ell(L)} = |L|$, then
\begin{equation}
\label{S(g_lambda)(v) = ... on L, L in mathcal{G}_{10}}
	S(g_\lambda)(v) = \Big(|E_0(f)(v)|^2 + 
	    \sum_{j=1}^{\ell(L)} |E_j(f)(v) - E_{j-1}(f)(v)|^2 \Big)^{1/2}.
\end{equation}
\end{corollary}

	These two corollaries follow from Lemma \ref{frac{1}{|K|}
int_K g_lambda = frac{1}{|K|} int_K g when ....} and the relevant
definitions.

\begincorollary
\label{S(g_lambda) le min(lambda, S(f))}
$S(g_\lambda) \le \min(\lambda, S(f))$.
\end{corollary}

	Corollary \ref{S(f), S(g_lambda), etc.} implies that
$S(g_\lambda) \le S(f)$, and we get $S(g_\lambda) \le \lambda$ using
also (\ref{... le lambda}) and (\ref{{x in [0,1) : S(f)(x) > lambda}
subseteq cup_{L in mathcal{G}_{10}} L}).

	Because of Corollary \ref{M(f) = M(g_lambda) on the
proper set},
\begin{eqnarray}
\lefteqn{\{x \in [0,1) : M(f)(x) > \lambda\}}	\\
	& & \subseteq \Big(\bigcup_{I \in \mathcal{G}_{01}} I\Big)
		\cup \{x \in [0,1) : M(g_\lambda)(x) > \lambda\}, \nonumber
\end{eqnarray}
and hence
\begin{eqnarray}
\lefteqn{|\{x \in [0,1) : M(f)(x) > \lambda\}|}   \\
	 & & 	\le \Big(\sum_{I \in \mathcal{G}_{01}} |I| \Big)
		   + |\{x \in [0,1) : M(g_\lambda)(x) > \lambda\}|.
								\nonumber
\end{eqnarray}
Therefore
\begin{eqnarray}
\lefteqn{\quad  |\{x \in [0,1) : M(f)(x) > \lambda\}|}	\\
	& & \le 2 \, |\{x \in [0,1) : S(f)(x) > \lambda\}|
		   + |\{x \in [0,1) : M(g_\lambda)(x) > \lambda\}|,
								\nonumber
\end{eqnarray}
by (\ref{sum_{L in mathcal{G}_{10}} |L| le ....}).  Of course,
\begin{equation}
	|\{x \in [0,1) : M(g_\lambda)(x) > \lambda\}|
		\le \lambda^{-2} \int_{[0,1)} M(g_\lambda)(u)^2 \, du,
\end{equation}
and
\begin{equation}
	\int_{[0,1)} M(g_\lambda)(u)^2 \, du
		\le C \, \int_{[0,1)} |g_\lambda(y)|^2 \, dy
\end{equation}
for some $C > 0$, by Proposition \ref{p-integral inequalities for
M(f)}.  Moreover,
\begin{eqnarray}
	\int_{[0,1)} |g_\lambda(y)|^2 \, dy 
		& = & \int_{[0,1)} S(g_\lambda)(w)^2 \, dw		\\
		& \le & \int_{[0,1)} \min(\lambda, S(f)(w))^2 \, dw,
								\nonumber
\end{eqnarray}
by Lemma \ref{2-norm of S(f)} and Corollary \ref{S(g_lambda) le
min(lambda, S(f))}.  It follows that
\begin{equation}
	\quad |\{x \in [0,1) : M(g_\lambda)(x) > \lambda\}|
	    \le C \, \lambda^{-2} \int_{[0,1)} \min(\lambda, S(f)(w))^2 \, dw,
\end{equation}
and consequently
\begin{eqnarray}
\label{|{x in [0,1) : M(f)(x) > lambda}| le ...., 2}
\lefteqn{\qquad |\{x \in [0,1) : M(f)(x) > \lambda\}|}   \\
	 & & 	\le 2 \, |\{x \in [0,1) : S(f)(x) > \lambda\}|
	      + C \, \lambda^{-2} \int_{[0,1)} \min(\lambda, S(f)(w))^2 \, dw.
							\nonumber
\end{eqnarray}

	We assumed near the beginning of the argument that the set
(\ref{S(f) > lambda}) is not all of $[0,1)$.  If it is, then the
preceding inequality holds trivially.  The rest of the proof of
Proposition \ref{p-norm of M(f) bounded by p-norm of S(f), p < 2}
proceeds via computations like those in the previous section.

\section{Duality, 1}
\label{Duality, 1}

	 If $f_1$, $f_2$ are functions on $[0,1)$, and $l$ is a nonnegative
integer, then
\begin{eqnarray}
	& & \quad \int_{[0,1)} E_l(f_1) \, E_l(f_2) \, dx =		\\
	& & \int_{[0,1)} \Big(E_0(f_1) \, E_0(f_2) 
 + \sum_{j=1}^l (E_j(f_1) - E_{j-1}(f_1)) \, (E_j(f_2) - E_{j-1}(f_2))\Big) 
				\, dx,				\nonumber
\end{eqnarray}
where the sum is interpreted as being $0$ when $l = 0$.  This is a
``bilinear'' version of (\ref{int S_l(f)^2 = int |E_l(f)|^2}) in Lemma
\ref{2-norm of S(f)}, which can be verified in essentially the same
way.  By the Cauchy--Schwarz inequality for sums,
\begin{equation}
	\biggl|\int_{[0,1)} E_l(f_1) \, E_l(f_2) \, dx \biggr|
		\le \int_{[0,1)} S_l(f_1) \, S_l(f_2) \, dx,
\end{equation}
and for suitable functions $f_1$ and $f_2$,
\begin{equation}
\label{|int_{[0,1)} f_1(x) f_2(x) dx| le int_{[0,1)} S(f_1)(x) S(f_2)(x) dx}
	\biggl|\int_{[0,1)} f_1(x) \, f_2(x) \, dx \biggr|
		\le \int_{[0,1)} S(f_1)(x) \, S(f_2)(x) \, dx.
\end{equation}

\beginproposition
\label{q-norm of f bounded by q-norm of S(f), q > 2}
For each $q > 2$, there is a $C_3(q) > 0$ such that
\begin{equation}
	\int_{[0,1)} |f(x)|^q \, dx \le C_3(q) \int_{[0,1)} S(f)(x)^q \, dx.
\end{equation}
\end{proposition}

	Let $q > 2$ be given, and let $p$, $1 < p < \infty$, be the
exponent dual to $q$, so that $1/p + 1/q = 1$ and $p <2$.  
By H\"older's inequality,
\begin{equation}
 \quad 	 \biggl|\int_{[0,1)} f_1(x) \, f_2(x) \, dx \biggr|
		\le \Big(\int_{[0,1)} S(f_1)^q \, dy \Big)^{1/q}
			\, \Big(\int_{[0,1)} S(f_2)^p \, dw \Big)^{1/p}.
\end{equation}
Propositions \ref{p-integral inequalities for M(f)} and \ref{p-norm of
S(f) bounded by p-norm of M(f), p < 2} yield
\begin{eqnarray}
\lefteqn{\biggl|\int_{[0,1)} f_1(x) \, f_2(x) \, dx \biggr|} \\
        & \le & C \, \Big(\int_{[0,1)} S(f_1)(y)^q \, dy \Big)^{1/q}
	     \, \Big(\int_{[0,1)} |f_2(w)|^p \, dw \Big)^{1/p} \nonumber
\end{eqnarray}
for some $C > 0$.  In general, if
\begin{equation}
\label{integral of product le A times L^p norm}
	\biggl|\int_{[0,1)} f_1(x) \, f_2(x) \, dx \biggr|
		\le A \, \Big(\int_{[0,1)} |f_2(w)|^p \, dw \Big)^{1/p}
\end{equation}
for some $A \ge 0$ and arbitrary functions $f_2$ on $[0, 1)$, then
\begin{equation}
\label{L^q norm le A}
	\Big(\int_{[0,1)} |f_1(x)|^q \, dx\Big)^{1/q} \le A,
\end{equation}
and the proposition follows.

\section{Duality, 2}
\label{Duality, 2}

\beginproposition
\label{q-norm of S(f) bounded by q-norm of f, q > 2}
For each $q > 2$, there is a $C_4(q) > 0$ such that
\begin{equation}
	\int_{[0,1)} S(f)(x)^q \, dx \le C_4(q) \int_{[0,1)} |f(x)|^q \, dx.
\end{equation}
\end{proposition}

	Let $q > 2$ be given, and let $p$ be the conjugate exponent to
$q$.  It suffices to show that the proposition holds with $S(f)$
replaced with $S_l(f)$ for every $l$, with a constant that does not
depend on $l$.  To do this, it is enough to show that
\begin{equation}
\label{integral with f and alpha's}
 \quad	\biggl|\int_{[0,1)} \Big(\alpha_0(x) \, E_0(f)(x) 
+ \sum_{j=1}^l \alpha_j(x) \, \big(E_j(f)(x) - E_{j-1}(f)(x)\big)\Big) \, dx
							\biggr|
\end{equation}
is less than or equal to a constant times the product of
\begin{equation}
\label{(int_{[0,1)} |f(y)|^q dy)^{1/q}}
	\Big(\int_{[0,1)} |f(y)|^q \, dy\Big)^{1/q}
\end{equation}
and
\begin{equation}
\label{integral of alpha's}
    \Big(\int_{[0,1)} \Big(\sum_{j=0}^l |\alpha_j(w)|^2\Big)^{p/2} 
							\, dw \Big)^{1/p}
\end{equation}
for arbitrary functions $\alpha_0, \ldots, \alpha_l$ on $[0, 1)$.  By
Lemma \ref{integrals and E_j}, (\ref{integral with f and alpha's}) is
equal to
\begin{equation}
	\biggl|\int_{[0,1)} \Big(E_0(\alpha_0)(x) 
 + \sum_{j=1}^l \big(E_j(\alpha_j)(x) - E_{j-1}(\alpha_j)(x)\big)\Big) 
			\, f(x) \, dx \biggr|.
\end{equation}
H\"older's inequality implies that this is less than or equal to the
product of (\ref{(int_{[0,1)} |f(y)|^q dy)^{1/q}}) and
\begin{equation}
\label{integral of E(alpha)'s}
	\Big(\int_{[0,1)} \biggl|E_0(\alpha_0)(x) +
   \sum_{j=1}^l \big(E_j(\alpha_j)(x) - E_{j-1}(\alpha_j)(x)\big)\biggr|^p 
						\, dx\Big)^{1/p}.
\end{equation}
Thus we would like to show that (\ref{integral of E(alpha)'s}) is less
than or equal to a constant times (\ref{integral of alpha's}).
Proposition \ref{p-norm of M(f) bounded by p-norm of S(f), p < 2}
implies that (\ref{integral of E(alpha)'s}) is bounded by a constant
times
\begin{equation}
\label{integral of S of ... (with E's of alpha's)}
	\Big(\int_{[0,1)} S \Big(E_0(\alpha_0) +
		\sum_{j=1}^l (E_j(\alpha_j) - E_{j-1}(\alpha_j))\Big)(x)^p
						\, dx \Big)^{1/p},
\end{equation}
and so we would like to show that (\ref{integral of S of ... (with
E's of alpha's)}) is bounded by a constant times (\ref{integral of
alpha's}).  One can check that
\begin{equation}
	S \Big(E_0(\alpha_0) +
	\sum_{j=1}^l \big(E_j(\alpha_j) - E_{j-1}(\alpha_j)\big)\Big)(x)
\end{equation}
is equal to
\begin{equation}
	\Big(|E_0(\alpha_0)(x)|^2 + 
   \sum_{j=1}^l |E_j(\alpha_j)(x) - E_{j-1}(\alpha_j)(x)|^2 \Big)^{1/2}.
\end{equation}
It therefore remains to show that
\begin{equation}
    \quad
	\Big(\int_{[0,1)} \Big(|E_0(\alpha_0)(x)|^2 + 
   \sum_{j=1}^l |E_j(\alpha_j)(x) - E_{j-1}(\alpha_j)(x)|^2 \Big)^{p/2}
							\, dx \Big)^{1/p}
\end{equation}
is bounded by a constant times (\ref{integral of alpha's}).  This can
be done using the results discussed in the next section.

\section{Auxiliary estimates}
\label{Auxiliary estimates}

	Let $l$ be a nonnegative integer, and let $\beta_0(x),
\beta_1(x), \ldots, \beta_l(x)$ be nonnegative functions on $[0,1)$.
Given $p, r \ge 1$, consider the problem of bounding
\begin{equation}
\label{(int_{[0,1)} (sum_{j=0}^l E_j(beta_j)(x)^r)^{p/r} dx)^{1/p}}
    \Big(\int_{[0,1)} \Big(\sum_{j=0}^l E_j(\beta_j)(x)^r \Big)^{p/r}
						\, dx \Big)^{1/p}
\end{equation}
by a constant times
\begin{equation}
\label{(int_{[0,1)} (sum_{j=0}^l beta_j(x)^r)^{p/r} dx)^{1/p}}
	\Big(\int_{[0,1)} \Big(\sum_{j=0}^l \beta_j(x)^r \Big)^{p/r}
						\, dx \Big)^{1/p},
\end{equation}
where the constant does not depend on $l$ or $\beta_0(x), \beta_1(x),
\ldots, \beta_l(x)$.  If $r = \infty$, then
\begin{equation}
	\Big(\sum_{j=0}^l \beta_j(x)^r \Big)^{1/r},	\quad
		\Big(\sum_{j=0}^l E_j(\beta_j)(x)^r \Big)^{1/r}
\end{equation}
should be replaced with
\begin{equation}
	\max_{0 \le j \le l} \beta_j(x), \quad
		\max_{0 \le j \le l} E_j(\beta_j)(x),
\end{equation}
as usual.  

\beginlemma
\label{E_j has norm 1 w.r.t. L^p norm}
For each $p \ge 1$, nonnegative integer $j$, and nonnegative function
$\beta$ on $[0,1)$,
\begin{equation}
	\int_{[0,1)} E_j(\beta)(x)^p \, dx
		\le \int_{[0,1)} \beta(x)^p \, dx.
\end{equation} 
\end{lemma}

	If $J$ is any interval in $[0,1)$, then
\begin{equation}
	\Big(\frac{1}{|J|} \int_J \beta(y) \, dy \Big)^p
		\le \frac{1}{|J|} \int_J \beta(y)^p \, dy,
\end{equation}
by Jensen's inequality.  Lemma \ref{E_j has norm 1 w.r.t. L^p norm}
follows by summing this over the dyadic intervals $J$ of length $2^{-j}$.

	Using Lemma \ref{E_j has norm 1 w.r.t. L^p norm}, it is easy
to see that (\ref{(int_{[0,1)} (sum_{j=0}^l E_j(beta_j)(x)^r)^{p/r}
dx)^{1/p}}) is less than or equal to (\ref{(int_{[0,1)} (sum_{j=0}^l
beta_j(x)^r)^{p/r} dx)^{1/p}}) when $p = r$.  When $r = \infty$, we
might as well restrict our attention to the case where the $\beta_j$'s
are all the same, and the question reduces to one about maximal
functions.  Lemma \ref{supremum bound for M(f)} and Proposition
\ref{p-integral inequalities for M(f)} yield suitable estimates for $1
< p \le \infty$.

\beginlemma
\label{earlier inequality equivalent to one with h's too}
Suppose that $1 \le r < p < \infty$, and let $s \in (1,\infty)$ be
conjugate to $p/r$, so that $1/s + r/p = 1$.  For each positive real
number $A_0$, (\ref{(int_{[0,1)} (sum_{j=0}^l E_j(beta_j)(x)^r)^{p/r}
dx)^{1/p}}) is less than or equal to $A_0$ times (\ref{(int_{[0,1)}
(sum_{j=0}^l beta_j(x)^r)^{p/r} dx)^{1/p}}) for arbitrary nonnegative
functions $\beta_0, \beta_1, \ldots, \beta_l$ on $[0,1)$ if and only
if
\begin{eqnarray}
\label{inequality with beta's and h's}
   & & \int_{[0,1)} \Big(\sum_{j=0}^l E_j(\beta_j)(x)^r \Big) \, h(x) \, dx
									\\
	& & \le A_0^r \,
\Big(\int_{[0,1)} \Big(\sum_{j=0}^l \beta_j(y)^r \Big)^{p/r} 
							\, dy \Big)^{r/p}
	\, \Big(\int_{[0,1)} h(w)^s \, dw \Big)^{1/s}
								\nonumber
\end{eqnarray} 
for arbitrary nonnegative functions $\beta_0, \beta_1, \ldots,
\beta_l$ and $h$ on $[0,1)$.
\end{lemma}

	This is basically the same observation as in (\ref{integral of
product le A times L^p norm}) and (\ref{L^q norm le A}), applied to
this situation.  Let us continue to assume that $1 \le r < p <
\infty$, and that $s$ is conjugate to $p/r$.  We would like to show
that (\ref{inequality with beta's and h's}) holds for a suitable
choice of $A_0$.  Because $r \ge 1$, $E_j(\beta_j)^r \le
E_j(\beta_j^r)$, by Jensen's inequality.  Hence
\begin{equation}
	\quad
     \int_{[0,1)} \Big(\sum_{j=0}^l E_j(\beta_j)(x)^r \Big) \, h(x) \, dx
  \le \int_{[0,1)} \Big(\sum_{j=0}^l E_j(\beta_j^r)(x) \Big) \, h(x) \, dx.
\end{equation}
This implies that
\begin{equation}
     \int_{[0,1)} \Big(\sum_{j=0}^l E_j(\beta_j)(x)^r \Big) \, h(x) \, dx
  \le \int_{[0,1)} \sum_{j=0}^l \beta_j(x)^r \, E_j(h)(x) \, dx,
\end{equation}
and thus
\begin{equation}
     \int_{[0,1)} \Big(\sum_{j=0}^l E_j(\beta_j)(x)^r \Big) \, h(x) \, dx
  \le \int_{[0,1)} \sum_{j=0}^l \beta_j(x)^r \, M(h)(x) \, dx.
\end{equation}
By H\"older's inequality,
\begin{eqnarray}
   & & \int_{[0,1)} \Big(\sum_{j=0}^l E_j(\beta_j)(x)^r \Big) \, h(x) \, dx
									\\
	& & \le 
\Big(\int_{[0,1)} \Big(\sum_{j=0}^l \beta_j(y)^r \Big)^{p/r} 
							\, dy \Big)^{r/p}
	\, \Big(\int_{[0,1)} M(h)(w)^s \, dw \Big)^{1/s}.
								\nonumber
\end{eqnarray}
It follows from Proposition \ref{p-integral inequalities for M(f)}
that (\ref{inequality with beta's and h's}) holds for some $A_0 > 0$.
This shows that (\ref{(int_{[0,1)} (sum_{j=0}^l
E_j(beta_j)(x)^r)^{p/r} dx)^{1/p}}) is bounded by a constant times
(\ref{(int_{[0,1)} (sum_{j=0}^l beta_j(x)^r)^{p/r} dx)^{1/p}}) when $1
\le r < p < \infty$.

\beginlemma
\label{same estimates for (p,r) and (p', r')}
If $1 < p < \infty$, $1 < r < \infty$, and $p'$, $r'$ are the
exponents conjugate to $p$, $r$, respectively, then for each positive
real number $A_0$, (\ref{(int_{[0,1)} (sum_{j=0}^l
E_j(beta_j)(x)^r)^{p/r} dx)^{1/p}}) is less than or equal to $A_0$
times (\ref{(int_{[0,1)} (sum_{j=0}^l beta_j(x)^r)^{p/r} dx)^{1/p}})
for arbitrary nonnegative functions $\beta_0, \beta_1, \ldots,
\beta_l$ on $[0,1)$ if and only if the same is true with $p$, $r$
replaced by $p'$, $r'$.
\end{lemma}

	This also works for $p, r = 1, \infty$, with minor adjustments
of the usual type.  To prove the lemma, the main step is to observe
that (\ref{(int_{[0,1)} (sum_{j=0}^l E_j(beta_j)(x)^r)^{p/r}
dx)^{1/p}}) is bounded by $A_0$ times (\ref{(int_{[0,1)} (sum_{j=0}^l
beta_j(x)^r)^{p/r} dx)^{1/p}}) for all nonnegative functions
$\beta_0(x), \beta_1(x), \ldots, \beta_l(x)$ on $[0,1)$ if and only if
\begin{eqnarray}
\label{new inequality with beta_j's and gamma_j's}
\lefteqn{\qquad \int_{[0,1)} 
	\Big(\sum_{j=0}^l E_j(\beta_j)(x) \, \gamma_j(x) \Big) \, dx}  \\
    & & \le A_0 \,
	\Big(\int_{[0,1)} \Big(\sum_{j=0}^l \beta_j(x)^r \Big)^{p/r}
						\, dx \Big)^{1/p} \ 
	\Big(\int_{[0,1)} \Big(\sum_{j=0}^l \gamma_j(x)^{r'} \Big)^{p'/r'}
						\, dx \Big)^{1/p'}
								\nonumber
\end{eqnarray}
for all nonnegative functions $\beta_0(x), \beta_1(x), \ldots,
\beta_l(x)$ and $\gamma_0(x), \gamma_1(x), \ldots, \gamma_l(x)$ on
$[0,1)$.  It follows from the lemma and the remarks preceding it that
(\ref{(int_{[0,1)} (sum_{j=0}^l E_j(beta_j)(x)^r)^{p/r} dx)^{1/p}}) is
less than or equal to a constant times (\ref{(int_{[0,1)} (sum_{j=0}^l
beta_j(x)^r)^{p/r} dx)^{1/p}}) when $1 < p < r < \infty$.

\section{Interpolation}
\label{interpolation}

	Let $T$ be a linear operator acting on real or complex-valued
dyadic step functions on $[0,1)$.  Suppose that $1 \le p < q \le
\infty$,
\begin{equation}
\label{(int_{[0,1)} |T(f)(x)|^p dx)^{1/p} le ..}
	\Big(\int_{[0,1)} |T(f)(x)|^p \, dx \Big)^{1/p}
		\le N_p \, \Big(\int_{[0,1)} |f(x)|^p \, dx \Big)^{1/p},
\end{equation}
and
\begin{equation}
\label{(int_{[0,1)} |T(f)(x)|^q dx)^{1/q} le ..}
	\Big(\int_{[0,1)} |T(f)(x)|^q \, dx \Big)^{1/q}
		\le N_q \, \Big(\int_{[0,1)} |f(x)|^q \, dx \Big)^{1/q}
\end{equation}
when $q < \infty$ or
\begin{equation}
\label{sup_{x in [0,1)} |T(f)(x)| le N_infty sup_{[0,1)} |f(x)|}
	\sup_{x \in [0,1)} |T(f)(x)| \le N_\infty \, \sup_{[0,1)} |f(x)|
\end{equation}
when $q = \infty$, for some $N_p, N_q \ge 0$ and each $f$.  If $0 < t
< 1$ and
\begin{equation}
	\frac{1}{r} = \frac{t}{p} + \frac{1-t}{q},
\end{equation}
then
\begin{equation}
\label{estimate for int |T(f)(x)|^r}
	\Big(\int_{[0,1)} |T(f)(x)|^r \, dx \Big)^{1/r}
   \le N_p^t \, N_q^{1-t} \, \Big(\int_{[0,1)} |f(x)|^r \, dx \Big)^{1/r}.
\end{equation}
This can be derived from Theorem \ref{M. Riesz' convexity theorem}, as
follows.  For each positive integer $l$,
\begin{equation}
\label{(int_{[0,1)} |E_l(T(f))(x)|^p dx)^{1/p} le ...}
	\Big(\int_{[0,1)} |E_l(T(f))(x)|^p \, dx \Big)^{1/p}
		\le N_p \, \Big(\int_{[0,1)} |f(x)|^p \, dx \Big)^{1/p},
\end{equation}
and analogously for $q$ instead of $p$.  Theorem \ref{M. Riesz'
convexity theorem} can be applied to get that
\begin{equation}
\label{estimate for int |E_l(T(f))(x)|^r}
	\Big(\int_{[0,1)} |E_l(T(f))(x)|^r \, dx \Big)^{1/r}
   \le N_p^t \, N_q^{1-t} \, \Big(\int_{[0,1)} |f(x)|^r \, dx \Big)^{1/r}
\end{equation}
for all step functions $f$ on $[0,1)$ that are constant on dyadic
intervals of length $2^{-l}$, by thinking of $E_l \circ T$ as a linear
transformation on that space, which can be identified with ${\bf R}^n$
or ${\bf C}^n$, $n = 2^l$, as appropriate.  Once one has
(\ref{estimate for int |E_l(T(f))(x)|^r}) for step functions that are
constant on dyadic intervals of length $2^{-l}$ for every $l$, it is
easy to derive (\ref{estimate for int |T(f)(x)|^r}) for arbitrary
dyadic step functions.  Of course, one can extend this to other
classes of functions too.

	The maximal and square function operators discussed in this
chapter are not linear, but the same interpolation inequalities can be
applied to them.  One can show this by approximating these operators
by linear operators.  Suppose that $T$ is a not-necessarily-linear
operator acting on dyadic step functions on $[0, 1)$ such that for
each dyadic step function $f$ there is a linear operator $A$ on the
same space of functions with the properties that
\begin{equation}
	|A(h)(x)| \le T(h)(x)
\end{equation}
for every $h$, $x$ and
\begin{equation}
	T(f)(x) = |A(f)(x)|.
\end{equation}
If $T$ satisfies (\ref{(int_{[0,1)} |T(f)(x)|^p dx)^{1/p} le ..})  and
(\ref{(int_{[0,1)} |T(f)(x)|^q dx)^{1/q} le ..}) or (\ref{sup_{x in
[0,1)} |T(f)(x)| le N_infty sup_{[0,1)} |f(x)|}), then the analogous
inequalities hold for these approximating linear operators $A$.  By
interpolation, the approximating linear operators $A$ satisfy
(\ref{estimate for int |T(f)(x)|^r}), and therefore $T$ does too.  One
can approximate maximal functions in this way by linear operators of
the form
\begin{equation}
	E_{\alpha(x)}(f)(x),
\end{equation}
where $\alpha(x)$ takes values in nonnegative integers.  One can
approximate square functions by linear operators of the form
\begin{equation}
	\alpha_0(x) \, E_0(f)(x) 
		+ \sum_{i=1}^l \alpha_i(x) (E_i(f)(x) - E_{i-1}(f)(x)),
\end{equation}
where $\big(\sum_{i=0}^l |\alpha_i(x)|^2 \big)^{1/2} \le 1$.

\section{Another argument for $p = 4$}
\label{another argument for p = 4}

        Let $f$ be a function on $[0, 1)$, and consider
\begin{equation}
 S(f)(x)^4 = \Big(|E_0(f)(x)|^2 +
              \sum_{j = 1}^\infty |E_j(f)(x) - E_{j - 1}(f)(x)|^2\Big)^2.
\end{equation}
Put
\begin{equation}
R_j(f)(x) = \Big(\sum_{k = j}^\infty |E_k(f)(x) - E_{k - 1}(f)(x)|^2\Big)^{1/2}
\end{equation}
for each positive integer $j$.  Thus
\begin{eqnarray}
 \quad  S(f)(x)^4 & = & (|E_0(f)(x)|^2 + R_1(f)(x)^2)^2 \\
& = & |E_0(f)(x)|^4 + 2 \, |E_0(f)(x)|^2 \, R_1(f)(x)^2 + R_1(f)(x)^4,\nonumber
\end{eqnarray}
and
\begin{eqnarray}
 R_1(f)(x)^4 & = & \sum_{j = 1}^\infty |E_j(f)(x) - E_{j - 1}(f)(x)|^4 \\
              & &    + 2 \sum_{j = 1}^\infty |E_j(f)(x) - E_{j - 1}(f)(x)|^2
                                           \, R_{j + 1}(f)(x)^2. \nonumber
\end{eqnarray}

        If $I$ is a dyadic interval of length $2^{-j}$, then
\begin{equation}
\label{int_I (|E_j(f)(x)|^2 + R_{j + 1}(f)(x)^2) dx = int_I |f(x)|^2 dx}
 \int_I (|E_j(f)(x)|^2 + R_{j + 1}(f)(x)^2) \, dx = \int_I |f(x)|^2 \, dx.
\end{equation}
This is analogous to Lemma \ref{2-norm of S(f)}, using orthogonality
properties on $I$ analogous to those in Lemma \ref{orthogonality
properties} on $[0, 1)$.  In particular,
\begin{equation}
\label{int_I R_{j + 1}(f)(x)^2 dx le int_I |f(x)|^2 dx le int_I M(|f|^2)(x) dx}
 \int_I R_{j + 1}(f)(x)^2 \, dx \le \int_I |f(x)|^2 \, dx
                                  \le \int_I M(|f|^2)(x) \, dx,
\end{equation}
where $M(|f|^2)$ is the dyadic maximal function associated to $|f|^2$.
It follows that
\begin{eqnarray}
\lefteqn{\int_{[0, 1)} |E_j(f)(x) - E_{j - 1}(f)(x)|^2 \, R_{j + 1}(f)(x)^2
                                                        \, dx} \\
 & \le & \int_{[0, 1)} |E_j(f)(x) - E_{j - 1}(f)(x)|^2 \, M(|f|^2)(x) \, dx
                                                                   \nonumber
\end{eqnarray}
for each $j \ge 1$, by expressing the integral over $[0, 1)$ as a sum
of integrals over dyadic intervals of length $2^{-j}$, and using the
fact that $|E_j(f)(x) - E_{j - 1}(f)(x)|^2$ is constant on dyadic
intervals of length $2^{-j}$.

        Using estimates like these, one can check that
\begin{equation}
\label{int_{[0, 1)} S(f)(x)^4 dx le C int_{[0, 1)} S(f)(x)^2 M(|f|^2)(x) dx}
        \int_{[0, 1)} S(f)(x)^4 \, dx
                       \le C \int_{[0, 1)} S(f)(x)^2 \, M(|f|^2)(x) \, dx
\end{equation}
for some constant $C \ge 0$ that does not depend on $f$.  This also
uses the fact that $M(f)^2 \le M(|f|^2)$ to deal with the diagonal
terms.  This gives another way to estimate the $L^4$ norm of $S(f)$ in
terms of the $L^4$ norm of $f$.  More precisely, one can first apply
the Cauchy--Schwarz inequality to the right side of (\ref{int_{[0, 1)}
S(f)(x)^4 dx le C int_{[0, 1)} S(f)(x)^2 M(|f|^2)(x) dx}).  This
implies that the $L^4$ norm of $S(f)$ is bounded by a constant times
the product of the square root of the $L^4$ norm of $S(f)$ and the
fourth root of the $L^2$ norm of $M(|f|^2)$.  Dividing both sides by
the square root of the $L^4$ norm of $S(f)$ and then squaring, one
gets that the $L^4$ norm of $S(f)$ is bounded by a constant times the
square root of the $L^2$ norm of $M(|f|^2)$.  The latter is bounded by
a constant multiple of the $L^4$ norm of $f$, as desired, because of
the $L^2$ estimates for the maximal function applied to $|f|^2$.

\section{Rademacher functions}
\label{rademacher functions}

	For each positive integer $j$, the $j$th \emph{Rademacher
function}\index{Rademacher fuctions} $r_j$ is the dyadic step
function on $[0, 1)$ which is constant on dyadic intervals of length
$2^{-j}$ and whose values alternate between $1$ and $-1$.  Thus
\begin{equation}
	r_j(t) = 1
\end{equation}
when $k \, 2^{-j} \le t < (k + 1) \, 2^{-j}$ and $k$ is an even
integer, and
\begin{equation}
	r_j(t) = -1
\end{equation}
when $k$ is odd.  In particular,
\begin{equation}
	|r_j(t)| = 1
\end{equation}
for each $j$ and $t$.  If $I$ is a dyadic subinterval of $[0, 1)$ of
length $|I| > 2^{-j}$, then
\begin{equation}
	\int_I r_j(t) \, dt = 0,
\end{equation}
because the values of $r_j$ alternate between $1$ and $-1$ on $I$.  If
$j$ and $l$ are distinct positive integers, then
\begin{equation}
	\int_0^1 r_j(t) \, r_l(t) \, dt = 0.
\end{equation}
Thus the Rademacher functions are orthogonal with respect to the usual
integral inner product
\begin{equation}
	\langle f_1, f_2 \rangle = \int_0^1 f_1(t) \, f_2(t) \, dt
\end{equation}
for real-valued functions on the unit interval.  Since
\begin{equation}
	\int_0^1 r_j(t)^2 \, dt = 1
\end{equation}
for each $j$, the Rademacher functions are orthonormal with respect to
this inner product.

	If $f$ is a real-valued dyadic step function on $[0, 1)$ and
$p$ is a positive real number, then we put
\begin{equation}
	\|f\|_p = \Big(\int_0^1 |f(x)|^p \, dx\Big)^{1/p}.
\end{equation}
This is a norm when $p \ge 1$, and a quasinorm when $0 < p < 1$.
Jensen's inequality implies that
\begin{equation}
	\|f\|_p \le \|f\|_q
\end{equation}
when $p \le q$.  If
\begin{equation}
	f = \sum_{j = 1}^n a_j \, r_j
\end{equation}
is a linear combination of Rademacher functions, then
\begin{equation}
	\|f\|_2 = \Big(\sum_{j = 1}^n a_j^2\Big)^{1/2},
\end{equation}
by orthonormality.  Hence
\begin{equation}
	\|f\|_p \le \Big(\sum_{j = 1}^n a_j^2\Big)^{1/2} \le \|f\|_q
\end{equation}
when $p \le 2 \le q$.  It turns out that for each $q > 2$ there is a
$B(q) > 0$ such that
\begin{equation}
	\|f\|_q \le B(q) \, \Big(\sum_{j = 1}^n a_j^2\Big)^{1/2},
\end{equation}
and for each $p < 2$ there is a $B(p) > 0$ such that
\begin{equation}
	\Big(\sum_{j = 1}^n a_j^2 \Big)^{1/2} \le B(p) \, \|f\|_p.
\end{equation}
These constants do not depend on $n$ or the coefficients $a_1, \ldots,
a_n$.  By contrast,
\begin{equation}
         \max_{0 \le t < 1} |f(t)| = \sum_{j = 1}^n |a_j|,
\end{equation}
and hence the analogous statement for $q = +\infty$ does not work.

	If $q$ is an even integer, then we can expand $|f(t)|^q$ as a
$q$-fold sum of products of Rademacher functions.  The integral of a
product of Rademacher functions is $1$ when the corresponding indices
are equal in pairs, and is $0$ otherwise.  This permits one to
estimate $\|f\|_q^q$ by a multiple of
\begin{equation}
	\Big(\sum_{j = 1}^n a_j^2\Big)^{q/2},
\end{equation}
as desired.  Actually, it suffices to know that each index of a
Rademacher function in a product is equal to at least one other index
when the integral of the product is different from $0$.  If $q > 2$ is
not an even integer, then we can apply the previous assertion to the
smallest even integer $Q > q$ and use the monotonicity of $\|f\|_q$.
One could use H\"older's inequality instead, in the form
\begin{equation}
	\|f\|_q \le \|f\|_{Q}^a \, \|f\|_2^{1 - a}
\end{equation}
where
\begin{equation}
	\frac{1}{q} = \frac{a}{Q} + \frac{1 - a}{2},
\end{equation}
to get a better constant.  For $p < 2$, we can use H\"older's
inequality in the form
\begin{equation}
	\|f\|_2 \le \|f\|_p^b \, \|f\|_4^{1 - b}
\end{equation}
with
\begin{equation}
	\frac{1}{2} = \frac{b}{p} + \frac{1 - b}{4}
\end{equation}
and replace $\|f\|_4$ by a multiple of $\|f\|_2$ to estimate $\|f\|_2$
in terms of $\|f\|_p$.

	One can also see this as a consequence of the analysis of the
previous sections.  If $E_l$ is as defined in (\ref{def of E_k(f)}),
then $E_l(r_j) = 0$ when $j > l$ and $E_l(r_j) = r_j$ when $j \le l$.
Hence
\begin{equation}
	E_l(r_j) - E_{l - 1}(r_j) = 0
\end{equation}
when $j \ne l$, and
\begin{equation}
	E_j(r_j) - E_{j - 1}(r_j) = r_j.
\end{equation}
Therefore
\begin{equation}
	E_j(f) - E_{j - 1}(f) = a_j \, r_j
\end{equation}
for $j = 1, \ldots, n$, and
\begin{equation}
	S(f) = \Big(\sum_{j = 1}^n a_j^2\Big)^{1/2}.
\end{equation}
In particular, the square function $S(f)$ is constant on $[0, 1)$.

\section{Walsh functions}
\label{walsh functions}

	If $A = \{j_1, \ldots, j_n\}$ is a finite set of positive
integers, then the \emph{Walsh function}\index{Walsh functions} $w_A$
is the dyadic step function on the unit interval which is the product
of the Rademacher functions with these indices, i.e.,
\begin{equation}
	w_A(t) = r_{j_1}(t) \cdots r_{j_n}(t).
\end{equation}
This should be interpreted as the constant function equal to $1$ on
$[0, 1)$ when $A = \emptyset$.  Thus
\begin{equation}
	|w_A(t)| = 1
\end{equation}
for each $A$ and $t$, and
\begin{equation}
	\int_0^1 w_A(t)^2 \, dt = 1.
\end{equation}
The Walsh functions are orthonormal with respect to the usual integral
inner product, because the integral of a product of Rademacher
functions on $[0, 1)$ is nonzero if and only if the indices of the
Rademacher functions are equal in pairs, as in the previous section.
One can check that the Walsh functions form an orthonormal basis for
the space of all dyadic step functions on $[0, 1)$.  More precisely,
the Walsh functions associated to subsets $A$ of $\{1, \ldots, n\}$
form an orthonormal basis for the dyadic step functions that are
constant on dyadic intervals of length $2^{-n}$.  Remember that there
are $2^n$ subsets of $\{1, \ldots, n\}$, which is the same as the
number of dyadic subintervals of $[0, 1)$ of length $2^{-n}$.

	If $A \subseteq \{1, \ldots, n\}$, then $w_A$ is constant
on dyadic intervals of length $2^{-n}$, and
\begin{equation}
	E_l(w_A) = w_A
\end{equation}
when $l \ge n$.  If also $n \in A$, then $E_{n - 1}(w_A) = 0$, and therefore
\begin{equation}
	E_l(w_A) = 0
\end{equation}
for each $l < n$.  It follows that
\begin{equation}
	E_l(w_A) - E_{l - 1}(w_A) = 0
\end{equation}
when $l \ne n$, and
\begin{equation}
	E_n(w_A) - E_{n - 1}(w_A) = w_A.
\end{equation}

	The \emph{Walsh group}\index{Walsh group} can be defined as
the set of sequences of $\pm 1$'s, with respect to coordinatewise
multiplication.  Thus the Walsh group is the Cartesian product of a
sequence of copies of the group with two elements, and is a compact
Hausdorff topological space with respect to the product topology.  The
group structure is compactible with the topology, so that the Walsh
group is a commutative topological group.  Walsh functions correspond
exactly to Fourier analysis on this group.

\appendix


\addtocontents{toc}{{\protect\bigskip\protect\hbox{\protect\bf\protect
Appendix}}}

\chapter{Metric spaces}
\label{metric spaces}

\renewcommand{\theequation}{A.\arabic{equation}}

        A \emph{metric space}\index{metric spaces} is a set $M$
together with a nonnegative real-valued function $d(x, y)$ defined for
$x, y \in M$ such that $d(x, y) = 0$ if and only if $x = y$,
\begin{equation}
\label{d(x, y) = d(y, x)}
	d(x, y) = d(y, x)
\end{equation}
for every $x, y \in M$, and
\begin{equation}
\label{the triangle inequality}
	d(x, z) \le d(x, y) + d(y, z)
\end{equation}
for every $x, y, z \in M$.  The function $d(x, y)$ is known as the
\emph{metric} on $M$, and represents the distance between $x$ and $y$
in the metric space.  If $V$ is a real or complex vector space
equipped with a norm $\|v\|$, then it is easy to see that
\begin{equation}
\label{d(v, w) = ||v - w||, appendix}
        d(v, w) = \|v - w\|
\end{equation}
is a metric on $V$.  In particular, the standard Euclidean metric on
${\bf R}^n$ is the metric that corresponds to the standard Euclidean
norm on ${\bf R}^n$ in this way.  If we identify ${\bf C}^n$ with
${\bf R}^{2 n}$ in the usual way, then the metric on ${\bf C}^n$
determined by the standard Euclidean norm corresponds exactly to the
standard Euclidean metric on ${\bf R}^n$.

        Let $(M, d(x, y))$ be a metric space.  A sequence $\{x_j\}_{j
= 1}^\infty$ of elements of $M$ is said to
\emph{converge}\index{convergent sequences} to $x \in M$ if for each
$\epsilon > 0$ there is an $L \ge 1$ such that
\begin{equation}
\label{d(x_j, x) < epsilon}
        d(x_j, x) < \epsilon
\end{equation}
for each $j \ge L$.  One can check that the limit $x$ of a convergent
sequence $\{x_j\}_{j = 1}^\infty$ is unique when it exists, in which
case we put
\begin{equation}
\label{lim_{j to infty} x_j = x}
        \lim_{j \to \infty} x_j = x.
\end{equation}
A sequence of elements $\{x_j\}_{j = 1}^\infty$ of $M$ is said to be a
\emph{Cauchy sequence}\index{Cauchy sequences} if for each $\epsilon >
0$ there is an $L \ge 1$ such that
\begin{equation}
        d(x_j, x_l) < \epsilon
\end{equation}
for every $j, l \ge L$.  One can also check that every convergent
sequence in $M$ is a Cauchy sequence.

        Conversely, if every Cauchy sequence in $M$ converges to an
element of $M$, then $M$ is said to be \emph{complete}.\index{complete
metric spaces} As in Section \ref{real, complex numbers}, ${\bf R}$
and ${\bf C}$ are complete as metric spaces with respect to their
standard metrics.  This implies that ${\bf R}^n$ and ${\bf C}^n$ are
complete with respect to their standard metrics for each positive
integer $n$, because a sequence of elements of ${\bf R}^n$ or ${\bf
C}^n$ is a Cauchy sequence or a convergent sequence if and only if
their corresponding $n$ sequences of coordinates have the same property.

        Let $E$ be a subset of a metric space $M$.  A point $p \in M$
is said to be in the \emph{closure}\index{closure of a set}
$\overline{E}$ of $E$ in $M$ if for each $\epsilon > 0$ there is a
point $q \in E$ such that
\begin{equation}
        d(p, q) < \epsilon.
\end{equation}
If $p \in E$, then one can simply take $q = p$, so that every element
of $E$ is automatically an element of $\overline{E}$.  If
\begin{equation}
        \overline{E} = E,
\end{equation}
then we say that $E$ is a \emph{closed set}\index{closed sets} in $M$.

        One can check that the closure of any set in $M$ is closed.
If $\{x_j\}_{j = 1}^\infty$ is a sequence of elements of a subset $E$
of $M$ that converges to an element $x$ of $M$, then it is easy to see
that $x \in \overline{E}$.  Conversely, every element of $\overline{E}$
is the limit of a sequence of elements of $E$ that converges in $M$.

        If $x$ is an element of a metric space $M$ and $r \ge 0$, then
the closed ball in $M$ with center $x$ and radius $r$ is defined by
\begin{equation}
        \overline{B}(x, r) = \{y \in M : d(x, y \le r\}.
\end{equation}
One can check that this is always a closed set in $M$, using the
triangle inequality.

        A subset $E$ of a metric space $M$ is said to be
\emph{dense}\index{dense sets} in $M$ if $\overline{E} = M$.  This is
equivalent to saying that every element of $M$ is the limit of a
convergent sequence of elements of $E$.  The set ${\bf Q}$ of rational
numbers is dense in the real line with the standard metric, for instance.

        A subset $E$ of a metric space $M$ is said to be
\emph{bounded}\index{bounded sets} if it is contained in a ball,
which is to say that
\begin{equation}
        E \subseteq \overline{B}(p, r)
\end{equation}
for some $p \in M$ and $r \ge 0$.  In this case, we also have that
\begin{equation}
        E \subseteq \overline{B}(q, r + d(p, q))
\end{equation}
for every $q \in M$, by the triangle inequality.  If $E \subseteq M$
is bounded and nonempty, then the \emph{diameter}\index{diameter of a
set} is defined by
\begin{equation}
\label{diam E = sup {d(x, y) : x, y in E}}
        \diam E = \sup \{d(x, y) : x, y \in E\}.
\end{equation}
One can check that the closure $\overline{E}$ of a bounded set $E
\subseteq M$ is bounded as well, and that the diameter of
$\overline{E}$ is the same as the diameter of $E$.

        If $(M, d(x, y))$ is a metric space and $X$ is a subset of
$M$, then the restriction of the metric $d(x, y)$ on $M$ to $x, y \in
X$ satisfies the requirements of a metric on $X$, so that $X$ becomes
a metric space too.  If $M$ is complete as a metric space and $X$ is a
closed subset of $M$, then $X$ is also complete as a metric space.  To
see this, observe that any Cauchy sequence $\{x_j\}_{j = 1}^\infty$ in
$X$ is a Cauchy sequence in $M$ as well.  If $M$ is complete, then
$\{x_j\}_{j = 1}^\infty$ converges to an element $x$ of $M$, and $x
\in X$ when $X$ is a closed set in $M$.

        Suppose now that $(M, d(x, y))$ and $(N, \rho(u, v))$ are both
metric spaces.  Let $f$ be a function on $M$ with values in $N$, which
is the same as a mapping from $M$ into $N$, and which may be expressed
symbolically by $f : M \to N$.  As usual, $f$ is said to be
\emph{continuous}\index{continuous mappings} at a point $x \in M$
if for every $\epsilon > 0$ there is a $\delta > 0$ such that
\begin{equation}
        \rho(f(x), f(y)) < \epsilon
\end{equation}
for every $y \in M$ such that $d(x, y) < \delta$.  If $f$ is
continuous at $x$, and if $\{x_j\}_{j = 1}^\infty$ is a sequence of
elements of $M$ that converges to $x$, then it is easy to see that
$\{f(x_j)\}_{j = 1}^\infty$ converges to $f(x)$ in $N$.  Conversely,
if $f$ is not continuous at $x$, then one can check that there is an
$\epsilon > 0$ and a sequence $\{x_j\}_{j = 1}^\infty$ of elements of
$M$ that converges to $x$ such that
\begin{equation}
        \rho(f(x), f(x_j)) \ge \epsilon
\end{equation}
for each $j$, so that $\{f(x_j)\}_{j = 1}^\infty$ does not converge to
$f(x)$ in $N$.

        A mapping $f : M \to N$ is is said to be
\emph{continuous}\index{continuous mappings} if it is continuous at
every point in $M$.  Suppose that $(M_1, d_1)$, $(M_2, d_2)$, and
$(M_3, d_3)$ are metric spaces, and that $f_1 : M_1 \to M_2$ and $f_2
: M_2 \to M_3$ are continuous mappings between them.  The composition
$f_2 \circ f_1$ is the mapping from $M_1$ into $M_3$ defined by
\begin{equation}
\label{(f_2 circ f_1)(x) = f_2(f_1(x))}
        (f_2 \circ f_1)(x) = f_2(f_1(x))
\end{equation}
for every $x \in M_1$.  One can check that $f_2 \circ f_1$ is also a
continuous mapping from $M_1$ into $M_3$ under these conditions, using
either the definition of continuity in terms of $\epsilon$'s and
$\delta$'s, or the characterization of continuity in terms of
convergent sequences.

        If $f$ and $g$ be continuous real or complex-valued functions
on a metric space $M$, then their sum $f + g$ and product $f \, g$ are
also continuous functions on $M$.  More precisely, when we say that a
real or complex-valued functions on $M$ is continuous, we mean that it
is continuous as a mapping into ${\bf R}$ or ${\bf C}$ with its
standard metric.  To show that $f + g$ and $f \, g$ are continuous,
one can use the characterization of continuous functions in terms of
convergent sequences to reduce to the analogous statements for sums
and products of convergent sequences of real or complex numbers.  Of
course, one can also prove this more directly, using very similar
arguments.

        Let $M$ be a set, and let $(N, \rho(u, v))$ be a metric space.
A sequence $\{f_j\}_{j = 1}^\infty$ of mappings from $M$ into $N$ is
said to converge \emph{pointwise}\index{pointwise convergence} to a
mapping $f : M \to N$ if $\{f_j(x)\}_{j = 1}^\infty$ converges as a
sequence of elements of $N$ to $f(x)$ for every $x \in M$.  Similarly,
$\{f_j\}_{j = 1}^\infty$ is said to converge to $f$
\emph{uniformly}\index{uniform convergence} on $M$ if for each
$\epsilon > 0$ there is an $L \ge 1$ such that
\begin{equation}
\label{rho(f_j(x), f(x)) < epsilon}
        \rho(f_j(x), f(x)) < \epsilon
\end{equation}
for every $j \ge L$ and $x \in M$.  The difference between uniform and
pointwise convergence is that $L$ depends only on $\epsilon$ and not
on $x$ in the definition of uniform convergence, while $L$ is allowed
to depend on both $\epsilon$ and $x$ in the analogous formulation of
pointwise convergence.  If $(M, d(x, y))$ is also a metric space, and
$\{f_j\}_{j = 1}^\infty$ is a sequence of continuous mappings from $M$
into $N$ that converges uniformly to a mapping $f : M \to N$, then
a well-known theorem states that $f$ is also continuous.

        A function $f$ on a set $M$ with values in a metric space $N$
is said to be \emph{bounded}\index{bounded functions} if
\begin{equation}
\label{f(M) = {f(x) : x in M}}
        f(M) = \{f(x) : x \in M\}
\end{equation}
is a bounded subset of $N$.  Note that the sum and product of bounded
real or complex-valued functions on $M$ are also bounded functions on
$M$.  If a sequence $\{f_j\}_{j = 1}^\infty$ of bounded mappings from
$M$ into a metric space $N$ converges uniformly to a mapping $f : M
\to N$, then it is easy to see that $f$ is also bounded.

        Let $(M, d(x, y))$, $(N, \rho(u, v))$ be metric spaces again,
and let $C_b(M, N)$ be the collection of bounded continuous mappings
from $M$ into $N$.  One can check that
\begin{equation}
\label{theta (f, g) = sup {rho(f(x), g(x)) : x in M}}
        \theta (f, g) = \sup \{\rho(f(x), g(x)) : x \in M\}
\end{equation}
defines a metric on $C_b(M, N)$, known as the \emph{supremum
metric}.\index{supremum metric} Note that a sequence $\{f_j\}_{j =
1}^\infty$ of bounded continuous mappings from $M$ into $N$ converges
to a bounded continuous mapping $f : M \to N$ with respect to the
supremum metric if and only if $\{f_j\}_{j = 1}^\infty$ converges to
$f$ uniformly on $M$.

        If $N$ is complete as a metric space, then $C_b(M, N)$ is also
complete with respect to the supremum metric.  More precisely, if
$\{f_j\}_{j = 1}^\infty$ is a sequence of bounded continuous mappings
from $M$ into $N$ that is a Cauchy sequence with respect to
(\ref{theta (f, g) = sup {rho(f(x), g(x)) : x in M}}), then it is easy
to see that $\{f_j(x)\}_{j = 1}^\infty$ is a Cauchy sequence in $N$
for each $x \in M$.  If $N$ is complete, then it follows that
$\{f_j(x)\}_{j = 1}^\infty$ converges to an element $f(x)$ of $N$ for
every $x \in M$.  Using the Cauchy condition with respect to
(\ref{theta (f, g) = sup {rho(f(x), g(x)) : x in M}}), one can check
that $\{f_j\}_{j = 1}^\infty$ converges to $f$ uniformly on $M$, and
hence that $f$ is bounded and continuous on $M$.

        Let $(M, d(x, y))$ and $(N, \rho(u, v))$ be metric spaces.  A
mapping $f : M \to N$ is said to be \emph{uniformly
continuous}\index{uniform continuity} if for every $\epsilon > 0$
there is a $\delta > 0$ such that
\begin{equation}
\label{rho(f(x), f(y)) < epsilon for every x, y in M with d(x, y) < delta}
	\rho(f(x), f(y)) < \epsilon \quad\hbox{for every } x, y \in M
					\hbox{ with } d(x, y) < \delta.
\end{equation}
As before, the composition of two uniformly continuous mappings is
uniformly continuous, as is the limit of a uniformly convergent
sequence of uniformly continuous mappings.  Similarly, the sum of two
uniformly continuous real or complex-valued functions is also
uniformly continuous, as is the product of a uniformly continuous
function and a constant.  The product of two bounded uniformly
continuous real or complex-valued functions is uniformly continuous as
well, but this does not always work without the additional hypothesis
of boundedness.

        Suppose that $E$ is a dense subset of $M$, and that $f$ is a
uniformly continuous mapping from $E$ into $N$.  If $N$ is complete,
then there is a unique extension of $f$ to a uniformly continuous
mapping from $M$ into $N$.  To see this, let $x$ be any element of
$M$, and let $\{x_j\}_{j = 1}^\infty$ be a sequence of elements of $E$
that converges to $x$, which exists because $E$ is dense in $M$.  In
particular, $\{x_j\}_{j = 1}^\infty$ is a Cauchy sequence, and one can
use the uniform continuity of $f : E \to N$ to show that
$\{f(x_j)\}_{j = 1}^\infty$ is a Cauchy sequence in $N$.  If $N$ is
complete, then it follows that $\{f(x_j)\}_{j = 1}^\infty$ converges
in $N$.  One can also check that the limit of this sequence does not
depend on the specific choice of the sequence $\{x_j\}_{j = 1}^\infty$
of elements of $E$ converging to $x$, so that this defines a mapping
from $M$ into $N$.  This new mapping clearly agrees with the original
one on $E$, and one can use the uniform continuity of the original
mapping on $E$ to show that the new mapping is uniformly continuous on
all of $M$.  The uniqueness of the extension follows from the fact
that two continuous mappings from $M$ into $N$ are the same when they
agree on a dense set.

        A subset $A$ of $M$ is said to be \emph{totally
bounded}\index{totally bounded sets} if for each $\epsilon > 0$, $A$
can be covered by finitely many balls of radius $\epsilon$ in $M$.  It
is easy to see that totally bounded sets are bounded, and that bounded
subsets of ${\bf R}^n$ with the standard metric are totally bounded.
If $f : M \to N$ is uniformly continuous and $A \subseteq M$ is
totally bounded, then $f(A)$ is totally bounded in $N$.

\backmatter

\newpage


\printindex

\end{document}